\documentclass[12pt,reqno]{article}

\usepackage{mathtools} %
\DeclarePairedDelimiter{\floor}{\lfloor}{\rfloor} 

\usepackage{amsmath,amsfonts,amssymb,amsthm} 
\usepackage{stmaryrd} %
\newcommand{\iv}[1]{\llbracket #1 \rrbracket} %
\newcommand{\ivl}[1]{\big\llbracket #1 \big\rrbracket} %

\usepackage{ushort} %
\usepackage{enumitem} %

\usepackage[utf8]{inputenc}

\usepackage{textcomp} %

\newcommand{\seqnum}[1]{\href{https://oeis.org/#1}{\rm \underline{#1}}}

\theoremstyle{plain}%
\newtheorem{theorem}{Theorem}[section]
\theoremstyle{definition}
\newtheorem{definition}[theorem]{Definition}

\usepackage{textcomp} %
\usepackage{thmtools}

\declaretheoremstyle[
  headfont= \bf, %
  numbered=unless unique,
  bodyfont=\normalfont,
  spaceabove=1em plus 0.75em minus 0.25em,
  prefoothook=,%
  spacebelow=1em plus 0.75em minus 0.25em,
  qed={{\tiny$\square$}},
]{exmpstyle}

\declaretheorem[
  style=exmpstyle,
  title=Remark,
  refname={remark,remarks},
  Refname={Remark,Remarks}
]{remark}

\declaretheorem[
  style=exmpstyle,
  title=Necessary Condition,
  refname={condition,conditions},
  Refname={Condition,Conditions}
]{condition}

\makeatletter
\let\save@mathaccent\mathaccent
\newcommand*\if@single[3]{
  \setbox0\hbox{${\mathaccent"0362{#1}}^H$}
  \setbox2\hbox{${\mathaccent"0362{\kern0pt#1}}^H$}
  \ifdim\ht0=\ht2 #3\else #2\fi
  }
\newcommand*\rel@kern[1]{\kern#1\dimexpr\macc@kerna}
\newcommand*\widebar[1]{\@ifnextchar^{{\wide@bar{#1}{0}}}{\wide@bar{#1}{1}}}
\newcommand*\wide@bar[2]{\if@single{#1}{\wide@bar@{#1}{#2}{1}}{\wide@bar@{#1}{#2}{2}}}
\newcommand*\wide@bar@[3]{
  \begingroup
  \def\mathaccent##1##2{
    \let\mathaccent\save@mathaccent
    \if#32 \let\macc@nucleus\first@char \fi
    \setbox\z@\hbox{$\macc@style{\macc@nucleus}_{}$}
    \setbox\tw@\hbox{$\macc@style{\macc@nucleus}{}_{}$}
    \dimen@\wd\tw@
    \advance\dimen@-\wd\z@
    \divide\dimen@ 3
    \@tempdima\wd\tw@
    \advance\@tempdima-\scriptspace
    \divide\@tempdima 10
    \advance\dimen@-\@tempdima
    \ifdim\dimen@>\z@ \dimen@0pt\fi
    \rel@kern{0.6}\kern-\dimen@
    \if#31
      \overline{\rel@kern{-0.6}\kern\dimen@\macc@nucleus\rel@kern{0.4}\kern\dimen@}
      \advance\dimen@0.4\dimexpr\macc@kerna
      \let\final@kern#2
      \ifdim\dimen@<\z@ \let\final@kern1\fi
      \if\final@kern1 \kern-\dimen@\fi
    \else
      \overline{\rel@kern{-0.6}\kern\dimen@#1}
    \fi
  }
  \macc@depth\@ne
  \let\math@bgroup\@empty \let\math@egroup\macc@set@skewchar
  \mathsurround\z@ \frozen@everymath{\mathgroup\macc@group\relax}
  \macc@set@skewchar\relax
  \let\mathaccentV\macc@nested@a
  \if#31
    \macc@nested@a\relax111{#1}
  \else
    \def\gobble@till@marker##1\endmarker{}
    \futurelet\first@char\gobble@till@marker#1\endmarker
    \ifcat\noexpand\first@char A\else
      \def\first@char{}
    \fi
    \macc@nested@a\relax111{\first@char}
  \fi
  \endgroup
}
\makeatother

\usepackage{accents}
\DeclareMathSymbol{\widetildesym}{\mathord}{largesymbols}{"65}
\newcommand\lowerwidetildesym{%
  \text{\smash{\raisebox{-1.3ex}{%
        $\widetildesym$}}}}
\newcommand\fixwidetilde[1]{%
  \mathchoice
  {\accentset{\displaystyle\lowerwidetildesym}{#1}}
  {\accentset{\textstyle\lowerwidetildesym}{#1}}
  {\accentset{\scriptstyle\lowerwidetildesym}{#1}}
  {\accentset{\scriptscriptstyle\lowerwidetildesym}{#1}}
}

\newlength\mylen
\usepackage[usenames]{color}
\definecolor{webgreen}{rgb}{0,.5,0}
\definecolor{webbrown}{rgb}{.6,0,0}

\usepackage[colorlinks=true,
linkcolor=webgreen,
filecolor=webbrown,
citecolor=webgreen,
breaklinks=true]{hyperref}

\allowdisplaybreaks

\begin{document}
\pdfbookmark[1]{Title}{Title}

\begin{center}
  \vskip 1cm
  {\LARGE\bf Systematic Enumeration of Fundamental
    Quantities Involving Runs\\\vskip .1in in Binary Strings}
  \vskip 1cm
  \large
  Félix Balado and Guénolé C.M. Silvestre\\
  School of Computer Science\\
  University College Dublin\\
  Belfield Campus\\
  Dublin 4,  Ireland\\
  \href{mailto: felix@ucd.ie}{\tt felix@ucd.ie}\\
  \href{mailto: guenole.silvestre@ucd.ie}{\tt guenole.silvestre@ucd.ie}\\
\end{center}

\vskip .2 in

\pdfbookmark[1]{Abstract}{Abstract}

\begin{abstract}
  We give recurrences, generating functions and explicit exact
  expressions for the enumeration of fundamental quantities involving
  runs in binary strings. We first focus on enumerations concerning
  runs of ones, and we then analyse the same enumerations when runs of
  ones and runs of zeros are jointly considered. We give the
  connections between these two types of run enumeration, and with the
  problem of compositions. We also analyse the same enumerations with
  a Hamming weight constraint. We discuss which of the many number
  sequences that emerge from these problems are already known and
  listed in the OEIS. Additionally, we extend our main enumerative
  results to the probabilistic scenario in which binary strings are
  outcomes of independent and identically distributed Bernoulli
  variables.
\end{abstract}

\section{Introduction}\label{sec:introduction}
Runs in binary strings are essentially uninterrupted sequences of the
same bit.  In this paper we give recurrences, generating functions and
explicit expressions for the enumeration of fundamental quantities
involving runs in binary strings. A lot of prior work on this topic is
probabilistic in nature, in which case runs of ones are typically
called \textit{success} runs~\cite{balakrishnan01:_runs_scans}
---whereas runs of zeros are called \textit{failure} runs. The very
first historical problem on runs, which was studied by de
Moivre~\cite{moivre38:_doctrine}, was in fact a success runs problem.
Whereas some authors focus on problems concerning success
runs~\cite{balakrishnan01:_runs_scans,philippou86:_successes,godbole90:_specific},
others jointly consider success and failure
runs~\cite{wald40:_test,feller68:_introd_probab,schilling90:_longest,bloom96:_probab,balakrishnan01:_runs_scans}.

We can see enumerative results for runs in binary strings as
special cases of the above corpus of research. However, in spite of
the long history of the topic, a uniform and consistent approach to
enumerating runs-related quantities is currently lacking. Some authors
have provided generating functions for specific
settings~\cite{goulden83:_combin_enumer,apostol88:_binary,sedgewick96:_introd,grimaldi05:_without_odd_runs,sinha09,madden17:_gener_funct},
while others have derived explicit (closed-form) expressions for some
others~\cite{sedgewick96:_introd,magliveras96:_enumer,grimaldi05:_without_odd_runs,makri12:_count_runs_ones_ones_runs,nej19:_binary}
---either through direct combinatorial analysis or through
generating functions. Other authors have given recurrence relations
for different problems involving
runs~\cite{austin78:_binar,apostol88:_binary,bloom96:_probab,schilling90:_longest,grimaldi05:_without_odd_runs,nyblom12:_enumerating,nej19:_binary}. But
many relevant enumeration aspects have apparently not been considered
yet.%

Our goal here is to try to fill this gap by systematically studying
the most fundamental enumerations of quantities involving runs in
binary strings within a common framework. Our strategy always
involves: 1) finding one or more recurrences for each of the problems
at hand; and~2) if possible, solving the recurrences to obtain exact
explicit expressions by means of generating functions. This is hardly
a novel approach, but one that has been frequently sidelined in favour
of standard combinatorial analysis in the study of runs ---which makes
a number of these problems harder. The use of generating functions to
study runs harks back to the earliest problem of this kind studied by
de Moivre~\cite{moivre38:_doctrine}. As far as we are aware,
Laplace~\cite[p.\ 421]{hald03:_history} was the first author who used
recurrences ---or finite difference equations, as they were commonly
known back then--- to derive probability generating functions for
problems involving runs.  Later authors, prominently Feller~\cite[Ch.\
XIII]{feller68:_introd_probab}, also used this approach to deal with
runs, although sometimes with a definition of ``run'' different than
the one adopted here ---see Remark~\ref{rmk:intro} below.

In this work, we recover many established results ---often, though not
always, via more streamlined derivations--- while more frequently we
offer novel contributions in the form of recurrences, generating
functions, and explicit expressions.  Although we are aware of the
limitations of such explicit expressions, we do not pursue here asymptotic
approximations. As we proceed, we highlight which of the enumerations
that we address correspond to known number sequences ---meaning that
they have been documented in Sloane's On-Line Encyclopedia of Integer
Sequences (OEIS)~\cite{oeis}.

Importantly, we also give the most basic connections between the
enumerations of runs of ones, runs of ones and/or zeros, and
compositions. Some of these connections have been partly identified by
previous
authors~\cite{bateman48:_power,schilling90:_longest,grimaldi05:_without_odd_runs},
but, as far as we are aware, they have not been fully
investigated. Last but not least, even though our focus is on
enumerations and not on probabilistic results, we show that our
enumerative results (recurrences, generating functions, and explicit
expressions) can be straightforwardly extended to deal with
probabilistic runs, when assuming bitstrings generated by independent
and identically distributed Bernoulli random variables.

\subsection{Structure of the Paper}
This paper is structured as follows. In
Section~\ref{sec:number-n-strings} we study the number of binary
strings of length $n$ that contain prescribed quantities of runs of
ones under different constraints. In Section~\ref{sec:probability} we
consider the probabilistic extensions of the most relevant among the
Section~\ref{sec:number-n-strings} results.  We then consider in
Section~\ref{sec:number-n-strings-hamming} the same problems as in
Section~\ref{sec:number-n-strings}, but when the binary strings are
constrained to having fixed Hamming
weight. Sections~\ref{sec:oz},~\ref{sec:note-prob-extens}
and~\ref{sec:oz_hamming} address the basic same problems as
Sections~\ref{sec:number-n-strings},~\ref{sec:probability}
and~\ref{sec:number-n-strings-hamming}, respectively, when runs of
ones and runs of zeros are jointly considered ---rather than only runs
of ones.

In Section~\ref{sec:number-runs} we study the total number of runs
(under different constraints) that are found over all binary strings, or
over certain subsets of these. In Section~\ref{sec:total-number-ones}
we essentially address the same problems as in
Section~\ref{sec:number-runs}, but when the goal is instead to count
the number of ones contained in runs of ones, rather than the runs
themselves. Finally, we draw the conclusions of this work in
Section~\ref{sec:conclusions}.

\subsection{Definitions}\label{sec:definitions}
We start by precisely stating the definition of ``run'' that we use,
to avoid any possible confusion with previous works.
\begin{definition}[\textit{Run, Length of a
    Run}]\label{def:mood}
  A \textit{run} in a binary string is an uninterrupted sequence of
  bits of the same kind, flanked on each side either by the opposite
  bit or by the start/end of the string (Mood's counting
  criterion~\cite{mood40:_distr_theor_runs}). The number of bits of
  the same kind in a run is its \textit{length}.
\end{definition}

\begin{definition}[\textit{$n$-String}]
  To minimise excessive wording, a binary string of length $n$ is
  referred to in the following as an \textit{$n$-string}.
\end{definition}

\begin{definition}[\textit{$k$-Run}]
  A \textit{$k$-run} means a run of length~$k$, where $k\ge 0$.
\end{definition}

\begin{definition}[\textit{Null Run, Nonnull Run}]\label{def:nonnull}
  A \textit{null run} is defined to be a $0$-run, i.e., a zero-length
  run, whereas a \textit{nonnull run} is a $k$-run with $k\ge 1$.
\end{definition}

\begin{remark}\label{rmk:nonnull}
  Observe that with the definition above we are just explicitly
  identifying the existence of zero-length runs. Although perhaps
  initially puzzling, null runs ($0$-runs) follow directly from Mood's
  criterion in Definition~\ref{def:mood}.  Letting
  $\tilde{b}=\text{mod}(b+1,2)$ with $b\in\{0,1\}$, a null run of
  $b\,$s in an $n$-string occurs whenever a $\tilde{b}$ is immediately
  followed by another $\tilde{b}$, or whenever the $n$-string
  starts/ends with~$\tilde{b}$. For example, $^{'}\!0^{'}\!010^{'}$
  contains three $0$-runs of ones, $^{'}\!0101$ contains one $0$-run
  of ones, and $^{'}\!0^{'}\!0^{'}\!0^{'}\!0^{'}$ contains five
  $0$-runs of ones, where the apostrophes mark the positions of the
  null runs of ones. Even if stating the obvious, observe that a null
  run of $b\,$s contains no~$b\,$s. The apparent lack of applications
  for null runs is probably the reason why they have not been
  considered by previous authors. However, we will see that
  disregarding null runs is not an option in a complete theory of
  runs, as they play a fundamental role in key theoretical and
  practical results.
\end{remark}

\begin{definition}[\textit{($\ushort{k}\le\widebar{k}$)-Run,  ($\ge\! k$)-Run, ($\le\!k$)-Run}]
  A \textit{($\ushort{k}\le\widebar{k}$)-run} is defined to be a
  $k$-run with $0\le\ushort{k}\le k\le \widebar{k}$. %
  One may verbalise a $(\ushort{k}\le\widebar{k})$-run by saying a
  ``lower $k$ to higher~$k$ run''.  A \textit{$(\ge\! k)$-run} is
  defined to be a $(k\le \infty)$-run, which may include null runs,
  whereas a \textit{$(\le\! k)$-run} is defined to be a $(1\le k)$-run
  if null runs are excluded and a $(0\le k)$-run otherwise. These
  additional definitions can be naturally verbalised as a ``greater or
  equal than $k$ run'' and a ``smaller or equal than $k$ run,''
  respectively. Lastly, observe that we can also express a $k$-run as
  a $(k\le k)$-run.
\end{definition}

\begin{definition}[\textit{Odd and Even Runs, $p$-Parity Runs}] An
  \textit{odd (even) run} is defined to be a run whose length is
  odd (even). More generally, we use the expression \textit{odd (even)
    runs} when describing ensembles of runs whose lengths, which may
  be different, are all odd (even). For the sake of brevity in
  mathematical expressions, we also use the term \textit{$p$-parity
    run}, which refers to an odd run when $p=1$, and to an even run
  when~$p=0$. Throughout the paper we only consider nonnull even
  runs ---odd runs are always nonnull.
\end{definition}

\begin{remark}\label{rmk:intro}
  Mood's criterion is perhaps the most common way to define a
  run.  However, among the authors who adhere to this criterion, some use
  the term ``run'' to denote runs formed by one single kind of bit
  ---typically, runs of ones--- whereas others jointly consider runs
  of ones and zeros.  With the first convention, `$11000111100$'
  contains two runs of ones: a $2$-run and a $4$-run (or two runs of
  zeros: a $3$-run and a $2$-run), whereas with the second one it
  contains four runs (of ones and zeros): two $2$-runs, one $3$-run
  and one $4$-run.  In the first group of authors we would find most
  of those who deal with success runs, such as
  Apostol~\cite{apostol88:_binary}, Balakrishnan and
  Koutras~\cite{balakrishnan01:_runs_scans}, Makri et
  al.~\cite{makri07:_shortest_longest} and many others. Among authors
  simultaneously considering runs of ones and zeros we have Wishart
  and Hirschfeld~\cite{wishart36:_theorem},
  Stevens~\cite{stevens39:_distr}, Wald and
  Wolfowitz~\cite{wald40:_test} ---in their classic runs-based
  test---, Bloom~\cite{bloom96:_probab}, again Balakrishnan and
  Koutras~\cite{balakrishnan01:_runs_scans}, and others. The
  combination of Mood's criterion with two kinds of runs in binary
  strings is summarised in the following statement by Feller~\cite[p.\
  42]{feller68:_introd_probab}: \textit{``In any ordered sequence of
    elements of two kinds, each maximal subsequence of elements of
    like kind is called a run''.}  It should be noted, though, that
  Feller also considers runs of ones (success runs) separately from
  runs of zeros, and that he often handles other definitions of
  ``run''---see below.

  As just mentioned, the reader must be aware that the literature
  contains definitions of ``run'' divergent from Mood's criterion. A
  relevant case for this paper is the definition used by
  Flajolet and
  Sedgewick~\cite{flajolet09:_analytic,sedgewick13:_introduct} or by
  Nyblom~\cite{nyblom12:_enumerating}, for whom a run of ones of
  length~$k$ in a binary string still refers to the appearance of~$k$
  consecutive ones, but not necessarily terminated by a zero or by the
  end of the string (or started by a zero or by the start of the
  string).  Nevertheless, when discussing $n$-strings devoid of such
  runs, as these authors do, this definition of a run of ones of
  length $k$ is completely equivalent to a $(\ge k)$-run of ones
  according to our definition above. 
  
  Other authors consider not only definitions of ``run'' which do not
  conform to Mood's criterion, but also enumeration schemes that count
  overlapping runs. As a relevant example, Balakrishnan and Koutras
  consider Type I, II and III run enumeration schemes in their
  treatise on runs~\cite{balakrishnan01:_runs_scans}, of which only
  type II is relevant to the criterion adopted in this paper. To
  illustrate these schemes, consider the string
  `110111001111101'. With these authors' Type~I scheme ---originally
  proposed by Feller~\cite[p.\ 305]{feller68:_introd_probab}--- a
  $k$-run does not have to end with the opposite bit or with the end
  of the string. Thus, according to this scheme, we may say that the
  aforementioned string contains four runs of ones of length 2:
  `\underbar{11}0\underbar{11}100\underbar{11}\kern0.1em\underbar{11}101',
  which is not compatible with Mood's criterion. However, Feller's and
  Mood's criteria are equivalent when considering the longest
  run. Feller's criterion, which has often been used in runs theory,
  is somewhat unnatural. In fact, it was only introduced by this
  author to be able to study runs through renewal theory, as he
  considered that the classical theory of runs ---based on Mood's
  criterion--- was ``messy''.  Next, using the Type II scheme
  in~\cite{balakrishnan01:_runs_scans} we may say that the
  aforementioned string contains three runs of ones of length at least
  2: `\underbar{11}0\underbar{111}00\underbar{11111}01', which is
  compatible with Mood's criterion, and equivalent to us saying that
  it contains three $(\ge 2)$-runs of ones, or, equivalently, three
  $(2\le 15)$-runs of ones. On the other hand, overlapping of runs is
  allowed with the Type III scheme handled by Balakrishnan and Koutras
  (Ling's criterion~\cite{ling89}), according to which we may say, for
  example, that the string in question contains seven runs of ones of
  length~2:
  `$\underline{11}0\underline{1\overline{1}}\overline{1}00\underline{1\overline{1}}\overline{\hspace{.075cm}\underline{\!1}}\underline{\hspace{.075cm}\overline{\!1}}\overline{1}01$'.
  Of course, this is not compatible with our criterion
  either. 

  Finally, different authors have used names other than ``run'' to
  describe runs that ---usually either \textit{avant la lettre} or
  unwittingly--- comply with Mood's criterion. Some examples are:
  \textit{sequence}~\cite{gold29:_note_frequen},
  \textit{join}~\cite{wishart36:_theorem},
  \textit{group}~\cite{stevens39:_distr},
  \textit{block}~\cite{austin78:_binar,goulden83:_combin_enumer},
  \textit{subset of consecutive ones}~\cite{apostol88:_binary},
  \textit{isolated tuple}~\cite{apostol88:_binary}, or
  \textit{clump}~\cite{bloom96:_probab}. Runs with particular lengths
  have also been called \textit{singles/singletons/isolated
    letters/isolated
    ones}~\cite{bloom98:_singles,apostol88:_binary,gold29:_note_frequen,austin78:_binar}
  (i.e., $1$-runs), \textit{isolated pairs}~\cite{apostol88:_binary}
  (i.e, $2$-runs), \textit{isolated triples}~\cite{apostol88:_binary}
  (i.e., $3$-runs), etc. The expression \textit{maximal
    block/run}~\cite{goulden83:_combin_enumer,madden17:_gener_funct}
  has also been used to denote a run that follows Mood's criterion,
  but we believe that this term is best avoided in order to prevent
  any misunderstandings with maximum length runs. Like other previous
  authors~\cite{fu03:_runs,sinha09}, we are crediting Mood in
  Definition~\ref{def:mood}, but the fact is that the use of the term
  ``run'' with the exact same meaning precedes him (see for
  example~\cite{cochran36:_plants}), and the first explicit statement
  of ``Mood's criterion'' seems to have been subsequently done by
  Mosteller~\cite[p.\ 229]{mosteller41:_note_applic}. Looser uses of
  the term ``run'' in probability are even older~\cite{pearson97}.
\end{remark}

\subsection{Notation}
\label{sec:notation}
The binomial coefficient is defined for any two $a,b\in\mathbb{Z}$ as
${a\choose b}=a(a-1)\cdots(a-b+1)/b!$ if $a\ge b\ge 0$, and as
${a\choose b}=0$ if $b<0$ or $a<b$. Throughout the paper we make extensive use
of Iverson's bracket notation~\cite[p.\ 24]{graham94:_concrete}: a
true-or-false statement enclosed in blackboard bold straight brackets
---i.e., $\iv{\cdot}$--- takes the value~$1$ if the statement is true,
and $0$ if it is false. Finally, we follow the empty summation convention:
$\sum_{i=a}^b=0$ when $b<a$. %

\begin{remark}
  For the sake of brevity, in the remainder of this paper the
  unqualified term ``run'' always implicitly refers to a run of ones
  unless explicitly specified otherwise.  Most of the
  exceptions to the aforementioned naming convention occur in
  Sections~\ref{sec:oz},~\ref{sec:note-prob-extens},~\ref{sec:oz_hamming},
  and~\ref{sec:number-runs-oz-m}, where we jointly consider runs of
  ones and runs of zeros, and where we explicitly indicate the nature
  of the runs as required. Likewise, we assume that runs may be null
  unless explicitly indicated otherwise.
\end{remark}

\section{Number of  $n$-Strings that Contain Prescribed Quantities of
  Runs Under Different Constraints}
\label{sec:number-n-strings}

In this section we address several fundamental enumeration problems
that involve counting the number of $n$-strings that contain
prescribed quantities of runs subject to different constraints.

\subsection[Number of n-Strings that Contain Exactly m (ḵ≤k)-Runs]{Number of $n$-Strings that Contain Exactly $m$ ($\ushort{k}\le\widebar{k}$)-Runs}
\label{sec:w_n_m_klek}
We denote by~$w_{\ushort{k}\le\widebar{k}}(n,m)$ the number of
$n$-strings that contain exactly $m$
$(\ushort{k}\le\widebar{k})$-runs, i.e., $m$ runs whose lengths are
within $\ushort{k}$ and $\widebar{k}$. Notice that such $n$-strings
may also contain other runs longer than $\widebar{k}$ or shorter than
$\ushort{k}$. Because of the general nature of the definition of a
$(\ushort{k}\le\widebar{k})$-run, the results established in this
section form the foundation that allows us to address many other
runs-related enumerations in a simple ---and sometimes even trivial---
way.

We begin by formulating a  necessary condition that underpins
our approach.
\begin{condition}{(Existence of $n$-strings
  containing $m$~$(\ushort{k}\le\widebar{k})$-runs)}
  \begin{equation}
    \label{eq:m_necessary_condition}
    w_{\ushort{k}\le\widebar{k}}(n,m)> 0\quad \Longrightarrow\quad 0\le m \le \floor[\bigg]{\frac{n+1}{\ushort{k}+1}}.
\end{equation}
\end{condition}
The reason for the nonnegativity of $m$ is clear, and perhaps only the
case $m=0$ merits some explanation: observe that we can have
$n$-strings devoid of $(\ushort{k}\le\widebar{k})$-runs. As for the
upper bound in~\eqref{eq:m_necessary_condition}, if we have $m$ runs
with lengths $k_1,\dots,k_m$ such that
$\ushort{k}\le k_j\le \widebar{k}$ for $j=1,\dots,m$, then, because
the number of zeros in an $n$-string that contains $m$ runs must be at
least $m-1$, a necessary condition for these runs to fit in an
$n$-string is
\begin{equation}
  \label{eq:main_bound}
   m\,\ushort{k}+m-1 \le n,
\end{equation}
which is equivalent to the upper bound
in~\eqref{eq:m_necessary_condition}. Condition~\eqref{eq:m_necessary_condition}
is not sufficient, as the validity of its right-hand side does not
guarantee that we can accommodate $m$
$(\ushort{k}\le\widebar{k})$-runs within an $n$-string. However, if
the right-hand side is not true, then such accommodation is certainly
impossible.

In order to streamline some formulas, it is convenient to name the
difference between both sides of inequality~\eqref{eq:main_bound}:
\begin{equation}\label{eq:excess}
  e=n-(m\, \ushort{k}+m-1).
\end{equation}
Of course, $e\ge 0$
is another way to express the upper bound
in~\eqref{eq:m_necessary_condition}.

\begin{remark}
  Necessary condition~\eqref{eq:m_necessary_condition} does not hold
  if $n<-1$. However, consider the two cases
  \begin{equation}
    n=-1\text{ and }n=0,
  \end{equation}
  for which~\eqref{eq:m_necessary_condition} does not negate the
  possibility of existence of $n$-strings containing $m$
  $(\ushort{k}\le\widebar{k})$-runs, as long as $m=0$. Even if the two
  cases above may seem absurd ---especially $n=-1$, since $n=0$ is
  commonly accepted to be the degenerate case of an empty string---
  they play the crucial role of initialising all the recurrences that
  we give in this section ---and also in
  Sections~\ref{sec:probability}
  and~\ref{sec:number-n-strings-hamming}.
\end{remark}

Finally, we also discuss a relevant enumeration that stems from the
main one studied in this section: the number of $n$-strings that
contain at least $m$ $(\ushort{k}\le\widebar{k})$-runs, which we
denote by $w_{\ushort{k}\le\widebar{k}}^{\ge m}(n)$, and which we may
obtain from~$w_{\ushort{k}\le\widebar{k}}(n,m)$ using
\begin{equation}\label{eq:wn_klek_atleast}
  w_{\ushort{k}\le\widebar{k}}^{\ge m}(n)=\sum_{t=m}^{\floor*{\frac{n+1}{\ushort{k}+1}}}w_{\ushort{k}\le\widebar{k}}(n,t).
\end{equation}
The case $m=1$ can also be simply put as
\begin{equation}
  \label{eq:w_n_klek_atleast1_from w_n_0_klek}
  w_{\ushort{k}\le\widebar{k}}^{\ge 1}=2^n-w_{\ushort{k}\le\widebar{k}}(n,0).
\end{equation}
The number of $n$-strings that contains at most $m$
$(\ushort{k}\le\widebar{k})$-runs is
$w_{\ushort{k}\le\widebar{k}}^{\le m}(n)=2^n-w_{\ushort{k}\le\widebar{k}}^{\ge (m+1)}(n)$, so we do not need
to discuss this complementary enumeration.

\subsubsection{Recurrences}
\label{sec:recurrences_w_n_m_klek}

Let us first obtain a recurrence relation for
$w_{\ushort{k}\le\widebar{k}}(n,m)$.  Consider the contribution to
$w_{\ushort{k}\le\widebar{k}}(n,m)$ from the ensemble of $n$-strings
that begin with an $i$-run. If $\ushort{k}\le i\le \widebar{k}$ then they
contribute $w_{\ushort{k}\le\widebar{k}}(n-(i+1),m-1)$ to
$w_{\ushort{k}\le\widebar{k}}(n,m)$. On the other hand, if
$i>\widebar{k}$ or $i<\ushort{k}$ then they contribute
$w_{\ushort{k}\le\widebar{k}}(n-(i+1),m)$ to
$w_{\ushort{k}\le\widebar{k}}(n,m)$. Observe that the ``$+1$'' in
$(i+1)$ is there to guarantee that the first bit immediately after the
$i$-run is a zero ---i.e., to terminate the run.  Thus, considering
all possible lengths of a starting run we obtain the following
bivariate recurrence:
\begin{align}
  \label{eq:recurrence_w_n_m_klek_basic}
  w_{\ushort{k}\le\widebar{k}}(n,m)=&~\sum_{i=\ushort{k}}^{\widebar{k}}
    w_{\ushort{k}\le\widebar{k}}\big(n-(i+1),m-1\big)+\sum_{i=0}^{\ushort{k}-1}
    w_{\ushort{k}\le\widebar{k}}\big(n-(i+1),m\big)\nonumber\\  
  &+\sum_{i=\widebar{k}+1}^{n}
    w_{\ushort{k}\le\widebar{k}}\big(n-(i+1),m\big).  
\end{align}
The natural question is: what if $i\ge n$? As we see later, as long as
we take~\eqref{eq:m_necessary_condition} into account, we do not have
to worry about such questions with the correct initialisation of the
recurrence.

An alternative way to find a recurrence for
$w_{\ushort{k}\le\widebar{k}}(n,m)$ is based on the following
observation: if we remove the first bit of all $n$-strings, then we
can see that $w_{\ushort{k}\le\widebar{k}}(n,m)$ is approximately
equal to $2\,w_{\ushort{k}\le\widebar{k}}(n-1,m)$. The overcounting or
undercounting in this estimate with respect to the true enumeration
only depends on the cases in which an $n$-string starts with a
$\ushort{k}$-run or with a $(\widebar{k}+1)$-run, as these are the
only cases in which the number of $(\ushort{k}\le\widebar{k})$-runs
can be increased or decreased with respect to the $(n-1)$-string
formed by removing the first bit. Thus, in order to make the
aforementioned estimate exact we just need to make the following
adjustments to it:
\begin{enumerate}[label=\alph*)]
\item For every $n$-string that starts with a $\ushort{k}$-run, we
  must add one whenever the remaining $(n-(\ushort{k}+1))$-string
  contains $(m-1)$ $(\ushort{k}\le\widebar{k})$-runs, but we must
  subtract one whenever the remaining $(n-(\ushort{k}+1))$-string
  contains $m$ $(\ushort{k}\le\widebar{k})$-runs already.

\item For every $n$-string that starts with a $(\widebar{k}+1)$-run, we
  must subtract one from the estimate whenever the remaining
  $(n-(\widebar{k}+2))$-string contains $(m-1)$
  $(\ushort{k}\le\widebar{k})$-runs, but we must add one whenever the
  remaining $(n-(\widebar{k}+2))$-string contains $m$
  $(\ushort{k}\le\widebar{k})$-runs .
\end{enumerate}
Collecting the contributions from a) and b) we get the following
alternative recurrence:
\begin{align}
  \label{eq:recurrence_w_n_m_klek_simpler}
  w_{\ushort{k}\le\widebar{k}}(n,m)=2\,w_{\ushort{k}\le\widebar{k}}(n-1,m) & + w_{\ushort{k}\le\widebar{k}}\big(n-(\ushort{k}+1),m-1\big)-w_{\ushort{k}\le\widebar{k}}\big(n-(\ushort{k}+1),m\big)\nonumber\\
            & - w_{\ushort{k}\le\widebar{k}}\big(n-(\widebar{k}+2),m-1\big)+w_{\ushort{k}\le\widebar{k}}\big(n-(\widebar{k}+2),m\big).
\end{align}
Observe that recurrence~\eqref{eq:recurrence_w_n_m_klek_simpler} can
also be derived from recurrence~\eqref{eq:recurrence_w_n_m_klek_basic}
using
$w_{\ushort{k}\le\widebar{k}}(n,m)-w_{\ushort{k}\le\widebar{k}}(n-1,m)$. Recurrence
relation~\eqref{eq:recurrence_w_n_m_klek_simpler} is equivalent to but
simpler than~\eqref{eq:recurrence_w_n_m_klek_basic}: the number of
recursive calls in~\eqref{eq:recurrence_w_n_m_klek_basic}
is~$n$~(full-history recurrence), whereas the number of recursive
calls in~\eqref{eq:recurrence_w_n_m_klek_simpler} is five,
independently of $n$. Also,
unlike~\eqref{eq:recurrence_w_n_m_klek_basic} which contains an
$n$-dependent summation,
recurrence~\eqref{eq:recurrence_w_n_m_klek_simpler} is directly
amenable to the computation of the generating function associated to
$w_{\ushort{k}\le\widebar{k}}(n,m)$, as we see in
Section~\ref{sec:generating-function_klek}. In the remainder we work
with~\eqref{eq:recurrence_w_n_m_klek_simpler}.

To initialise~\eqref{eq:recurrence_w_n_m_klek_simpler} we can use the
case $n=1$, in which we have by inspection that
\begin{align}
  \label{eq:trivial_n1}
  w_{\ushort{k}\le\widebar{k}}(1,m)=&~\big(\iv{\ushort{k}=0}\iv{\widebar{k}=0}+2\,\iv{\ushort{k}>1}\big)\iv{m=0}+\big(\iv{\ushort{k}=1}+\iv{\ushort{k}=0}\iv{\widebar{k}\ge
  1}\big)\iv{m=1}\nonumber\\ &+\iv{\ushort{k}=0}\iv{m=2}.
\end{align}
On the other hand, setting $n=1$ in
recurrence~\eqref{eq:recurrence_w_n_m_klek_simpler}  we get
\begin{align}
  \label{eq:init_helper_w_n_m_klek}
  w_{\ushort{k}\le\widebar{k}}(1,m)=2\,w_{\ushort{k}\le\widebar{k}}(0,m) &+ w_{\ushort{k}\le\widebar{k}}(-\ushort{k},m-1)-w_{\ushort{k}\le\widebar{k}}(-\ushort{k},m)\nonumber\\
&-w_{\ushort{k}\le\widebar{k}}(-(\widebar{k}+1),m-1)+w_{\ushort{k}\le\widebar{k}}(-(\widebar{k}+1),m).
\end{align}
We wish~\eqref{eq:init_helper_w_n_m_klek} to
equal~\eqref{eq:trivial_n1}. Taking~\eqref{eq:m_necessary_condition} into account, we see
that the desired equality is fulfilled for all~$m$ and
$0\le \ushort{k}\le\widebar{k}$ by choosing
\begin{align}\label{eq:w_n_m_klek_init}
  w_{\ushort{k}\le\widebar{k}}(-1,m)&=\iv{m=0},\\
  w_{\ushort{k}\le\widebar{k}}(0,m)&=\ivl{m=\iv{\ushort{k}=0}},\label{eq:w_n_m_klek_init2}
\end{align}
which thus constitute the initialisation of
recurrence~\eqref{eq:recurrence_w_n_m_klek_simpler}. Through the same
procedure as above, the reader may verify
that~\eqref{eq:w_n_m_klek_init} and~\eqref{eq:w_n_m_klek_init2} also
initialise~\eqref{eq:recurrence_w_n_m_klek_basic}. Even if an
uncommon sight, the nested Iversonian brackets in the expression above
are not a typo, and indeed we will meet again this type of expression
in subsequent sections.  

Finally, we can obtain a recurrence for
$w_{\ushort{k}\le\widebar{k}}^{\ge m}(n)$ by adding
recurrences~\eqref{eq:recurrence_w_n_m_klek_basic}
or~\eqref{eq:recurrence_w_n_m_klek_simpler} over the range of $m$
in~\eqref{eq:wn_klek_atleast}. Let us do so with the simpler
recurrence~\eqref{eq:recurrence_w_n_m_klek_simpler}. Taking into
account necessary condition~\eqref{eq:m_necessary_condition}, we get
\begin{equation}
  \label{eq:wn_klek_atleast_rec}
  w_{\ushort{k}\le\widebar{k}}^{\ge m}(n)=2\,w_{\ushort{k}\le\widebar{k}}^{\ge m}(n-1)+w_{\ushort{k}\le\widebar{k}}\big(n-(\ushort{k}+1),m-1\big)-w_{\ushort{k}\le\widebar{k}}\big(n-(\widebar{k}+2),m-1\big).
\end{equation}
Although the  recurrence above depends on
$w_{\ushort{k}\le\widebar{k}}(n,m)$, we show in the next section that
it suffices to obtain the generating function of
$w_{\ushort{k}\le\widebar{k}}^{\ge m}(n)$.

\subsubsection{Generating Functions}
\label{sec:generating-function_klek}
Let us next obtain the bivariate ordinary generating function (ogf)
associated to $w_{\ushort{k}\le\widebar{k}}(n,m)$, i.e.,
\begin{equation}
  W_{\ushort{k}\le\widebar{k}}(x,y)=\sum_n\sum_{m} w_{\ushort{k}\le\widebar{k}}(n,m)\, x^n y^m.\label{eq:bi_gf_klek}
\end{equation}
Notice that we have not set summation limits in~\eqref{eq:bi_gf_klek},
which means that we are adding over all integers $n$ and $m$. This is
because, in order to streamline our task, we follow Graham et al.'s
procedure in~\cite[Sec.\ 7.3]{graham94:_concrete} and
apply~\eqref{eq:bi_gf_klek} to a version of
recurrence~\eqref{eq:recurrence_w_n_m_klek_simpler} valid for all
values of $n$ and $m$.  First of all, taking necessary
condition~\eqref{eq:m_necessary_condition} into account,
recurrence~\eqref{eq:recurrence_w_n_m_klek_simpler} is valid not only
for $n\ge 1$ and $m\ge 0$, but also for $n<-1$ and/or $m<0$ ---in
which cases
$w_{\ushort{k}\le\widebar{k}}(n,m)=0$. So let us see what happens in
the remaining two cases $n=-1$ and $n=0$, by computing
$w_{\ushort{k}\le\widebar{k}}(-1,m)$ and
$w_{\ushort{k}\le\widebar{k}}(0,m)$
using~\eqref{eq:recurrence_w_n_m_klek_simpler} and then comparing the
results with the correct initialisation values
in~\eqref{eq:w_n_m_klek_init} and~\eqref{eq:w_n_m_klek_init2}.

In the first case, from~\eqref{eq:recurrence_w_n_m_klek_simpler}
and~\eqref{eq:m_necessary_condition} we mistakenly have that
$w_{\ushort{k}\le\widebar{k}}(-1,m)=0$ instead of the correct value
$w_{\ushort{k}\le\widebar{k}}(-1,m)=\iv{m=0}$, but we can ``fix'' this
by adding $\iv{n=-1} \iv{m=0}$
to~\eqref{eq:recurrence_w_n_m_klek_simpler}. In the second case,
according to~\eqref{eq:recurrence_w_n_m_klek_simpler} and
using~\eqref{eq:w_n_m_klek_init} and~\eqref{eq:m_necessary_condition}
we wrongly have that $w_{\ushort{k}\le\widebar{k}}(0,m)=2\,\iv{m=0}$
when $\ushort{k}> 0$ and
$w_{\ushort{k}\le\widebar{k}}(0,m)=\iv{m=0}+\iv{m=1}$ when
$\ushort{k}=0$, rather than the correct
value~\eqref{eq:w_n_m_klek_init2}. Once again, we can ``fix'' these cases
just by subtracting $\iv{n=0} \iv{m=0}$
from~\eqref{eq:recurrence_w_n_m_klek_simpler}. Therefore an extended
version of recurrence~\eqref{eq:recurrence_w_n_m_klek_simpler} valid
for all values of $n$ and $m$ when
taking~\eqref{eq:m_necessary_condition} into account is
\begin{align}
  \label{eq:recurrence-gf-klek}
  w_{\ushort{k}\le\widebar{k}}(n,m)=2\,w_{\ushort{k}\le\widebar{k}}(n-1,m)
            & + w_{\ushort{k}\le\widebar{k}}(n-(\ushort{k}+1),m-1)-w_{\ushort{k}\le\widebar{k}}(n-(\ushort{k}+1),m)\nonumber\\
            & - w_{\ushort{k}\le\widebar{k}}(n-(\widebar{k}+2),m-1)+w_{\ushort{k}\le\widebar{k}}(n-(\widebar{k}+2),m)\nonumber\\
&  +\iv{m=0}\big(\iv{n=-1}-\iv{n=0}\big).
\end{align}
This same recurrence is obtained if we first
make~\eqref{eq:recurrence_w_n_m_klek_basic} valid for all~$n$ and~$m$
and then obtain
$w_{\ushort{k}\le\widebar{k}}(n,m)-w_{\ushort{k}\le\widebar{k}}(n-1,m)$
using that extended recurrence.
Considering~\eqref{eq:bi_gf_klek}, we can now get
$W_{\ushort{k}\le\widebar{k}}(x,y)$ just by
multiplying~\eqref{eq:recurrence-gf-klek} on both sides by $x^n y^m$
and then adding over $n$ and $m$. This yields
\begin{equation*}
  \label{eq:gf_eq_klek}
  W_{\ushort{k}\le\widebar{k}}(x,y)=2x\,W_{\ushort{k}\le\widebar{k}}(x,y)+x^{\ushort{k}+1}(y-1)\,W_{\ushort{k}\le\widebar{k}}(x,y) -x^{\widebar{k}+2}(y-1)\,W_{\ushort{k}\le\widebar{k}}(x,y)+x^{-1}-1,
\end{equation*}
and thus the ogf sought is
\begin{equation}\label{eq:ogf_w_n_m_klek}
   W_{\ushort{k}\le\widebar{k}}(x,y)=\frac{1-x}{x\big(1-2x+(1-y)(x^{\ushort{k}+1}-x^{\widebar{k}+2})\big)}.
 \end{equation}
 It is not difficult to identify the coefficient of $W_{\ushort{k}\le\widebar{k}}(x,y)$ that corresponds
 to~$y^m$, i.e., $[y^m]W_{\ushort{k}\le\widebar{k}}(x,y)$. To this end we
 rewrite~\eqref{eq:ogf_w_n_m_klek} as
\begin{equation*}
  \label{eq:ogf-2_klek}
  W_{\ushort{k}\le\widebar{k}}(x,y)=\frac{1-x}{x\big(1-2x+x^{\ushort{k}+1}-x^{\widebar{k}+2}\big)}\cdot\frac{1}{(1-c y)}
\end{equation*}
with
\begin{equation*}
  \label{eq:c_klek}
  c=\frac{x^{\ushort{k}+1}-x^{\widebar{k}+2}}{1-2x+x^{\ushort{k}+1}-x^{\widebar{k}+2}}.
\end{equation*}
Then, because from the negative binomial theorem we have that $(1-c y)^{-1}=\sum_{m\ge 0} (c y)^m$ it follows that
\begin{equation}
  \label{eq:ogf_m_klek}
  [y^m]W_{\ushort{k}\le\widebar{k}}(x,y)=\frac{(1-x)\big(x^{\ushort{k}+1}-x^{\widebar{k}+2}\big)^m}{x\big(1-2x+x^{\ushort{k}+1}-x^{\widebar{k}+2}\big)^{m+1}},
\end{equation}
which is the ogf enumerating the  binary strings that contain
exactly $m$ $(\ushort{k}\le\widebar{k})$-runs.

Finally, we can obtain the ogf
$W_{\ushort{k}\le\widebar{k}}^{\ge m}(x)=\sum_x
w_{\ushort{k}\le\widebar{k}}^{\ge m}(n)\, x^n$ by
multiplying~\eqref{eq:wn_klek_atleast_rec} by $x^n$ on both sides and
then adding on $n$. This gives
\begin{align}
  \label{eq:w_n_klek_atleast_ogf}
  W_{\ushort{k}\le\widebar{k}}^{\ge m}(x)&=\frac{x^{\ushort{k}+1}-x^{\widebar{k}+2}}{1-2x}\,[y^{m-1}]W_{\ushort{k}\le\widebar{k}}(x,y)\nonumber\\
  &=\frac{1-x}{(1-2x)\,x }\Bigg(\frac{x^{\ushort{k}+1}-x^{\widebar{k}+2}}{1-2x+x^{\ushort{k}+1}-x^{\widebar{k}+2}}\Bigg)^m,
\end{align}
where we have used~\eqref{eq:ogf_m_klek} to get the last expression.
This is the ogf enumerating the binary strings that contain at least
$m$
$(\ushort{k}\le\widebar{k})$-runs---see~\eqref{eq:wn_klek_atleast}.

\subsubsection{Explicit Expressions}
\label{sec:w_n_m_klek_explicit}
It is possible, in principle, to get explicit expressions for
$w_{\ushort{k}\le\widebar{k}}(n,m)$ from~\eqref{eq:ogf_m_klek} ---and
for $w^m_{\ushort{k}\le\widebar{k}}(n)$
from~\eqref{eq:w_n_klek_atleast_ogf}--- using the same standard method
that we adopt in many subsequent sections. This method is,
essentially, the repeated application of the (negative) binomial
theorem, plus, in some cases, the solution of one or more simple
Diophantine equations. However, the reader may verify that, in this
general case, this strategy leads to overcomplicated expressions,
which depend on whether $\gcd(\widebar{k}-\ushort{k}+1,\ushort{k})$ is
equal to or greater than one. Consequently, we only obtain explicit
expressions for the special cases of
$w_{\ushort{k}\le\widebar{k}}(n,m)$ considered in
Sections~\ref{sec:w_n_m_k}--\ref{sec:number-wl_n_k_gtk}.

\subsection{Number of $n$-Strings that Contain Exactly $m$ $k$-Runs}
\label{sec:w_n_m_k}

We denote the number of $n$-strings that contain exactly $m$ $k$-runs
by~$w_{k}(n,m)$. Notice that such $n$-strings may also
contain other runs of lengths different than~$k$.

The earliest reference that we know of for this enumeration is the
work of Apostol~\cite{apostol88:_binary}, who gave several recurrences
and a generating function for $w_k(n,m)$ motivated by an electrical
engineering problem. This relevant piece of work appears to have faded
into obscurity, perhaps due to its idiosyncratic naming
conventions. Koutras and Papastavridis~\cite{koutras93:_number}
recovered some of Apostol's results, but most later authors were
unaware of them. They also were unaware of the first closed-form
formula for $w_k(n,m)$ ---a triple-summation expression obtained as a special
case of a more general computation--- given by Magliveras and
Wei~\cite[Thm.\ 2.3]{magliveras96:_enumer}.  Sinha and
Sinha~\cite{sinha09} produced a triple-summation explicit expression for
$w_k(n,m)$, drawing on a generating function for the interstitial gaps
between the $m$ $k$-runs. Soon afterwards, Makri and
Psillakis~\cite{makri11:_bernoul} gave a simpler double-summation explicit
formula for~$w_k(n,m)$. Their derivation exploited an existing
combinatorial result of their own for the distribution of balls in
urns~\cite{makri07:_succes}. More recently,
Madden~\cite{madden17:_gener_funct} has given a generating function
for the equivalent problem of enumerating the $n$-strings \textit{that
  begin with zero} and contain a prescribed number of
runs of a given length. %

As far as we are concerned, this problem is a special case of the
enumeration in Section~\ref{sec:w_n_m_klek} with
$\ushort{k}=\widebar{k}=k$, and thus
\begin{equation}
  \label{eq:w_n_m_k_from_klek}
  w_k(n,m)=w_{k\le k}(n,m).
\end{equation}
We also study the number of $n$-strings that contain at least $m$
$k$-runs, which from~\eqref{eq:w_n_m_k_from_klek}
and~\eqref{eq:wn_klek_atleast} is
\begin{equation}\label{eq:wn_k_atleast}
  w_{k}^{\ge m}(n)=\sum_{t=m}^{\floor*{\frac{n+1}{k+1}}}w_{k}(n,t).
\end{equation}

\subsubsection{Recurrences}
\label{sec:recurrences_w_n_m_k}
From~\eqref{eq:recurrence_w_n_m_klek_basic},~\eqref{eq:recurrence_w_n_m_klek_simpler}
and~\eqref{eq:w_n_m_k_from_klek}, two recurrence relations
for $w_{k}(n,m)$ are
\begin{equation}
  \label{eq:recurrence_w_n_m_k_basic}
  w_{k}(n,m)=w_{k}\big(n-(k+1),m-1\big)+\sum_{i=0\atop  i\neq k}^n w_{k}\big(n-(i+1),m\big),
\end{equation}
and
\begin{align}
  \label{eq:recurrence_w_n_m_k_simpler}
  w_{k}(n,m)=2\,w_{k}(n-1,m)%
            & + w_{k}\big(n-(k+1),m-1\big)-w_{k}\big(n-(k+1),m\big)\nonumber\\
            & - w_{k}\big(n-(k+2),m-1\big)+w_{k}\big(n-(k+2),m\big),
\end{align}
which, from~\eqref{eq:w_n_m_klek_init} and~\eqref{eq:w_n_m_klek_init2},
are both initialised by %
\begin{align}\label{eq:w_n_m_k_init}
  w_{k}(-1,m)&=\iv{m=0},\\
  w_{k}(0,m)&=\ivl{m=\iv{k=0}}.\label{eq:w_n_m_k_init2}
\end{align}

\begin{remark}
  Apostol gives a host of recurrences for this problem~\cite[Thms.\
  2--11]{apostol88:_binary}. In Apostol's notation,
  $A_{m}^{(k)}=w_k(n,m)$, but he also uses $S_m(n)=w_1(n,m)$,
  $P_m(n)=w_2(n,m)$ and $T_m(n)=w_3(n,m)$.  This author gives several
  recurrences for enumerating the $n$-strings with prescribed numbers
  of \textit{isolated singletons} (i.e., $1$-runs): in particular, he
  gives recurrences for $w_1(n,0)$ with $n\ge 3$, for $w_1(n,1)$ with
  $n\ge 3$, for $w_1(n,2)$ with $n\ge 3$, and for $w_1(n,m)$ with
  $n\ge 3$ and $m\ge 2$.  He then gives recurrences for enumerating
  $n$-strings with prescribed numbers of \textit{isolated pairs}
  (i.e., $2$-runs), in particular, recurrences for $w_2(n,0)$ with
  $n\ge 3$, for $w_2(n,1)$ with $n\ge 3$, for $w_2(n,2)$ with
  $n\ge 3$, and for $w_2(n,m)$ with $n\ge 3$ and $m\ge 2$. Finally, he
  extends these recurrences to \textit{isolated $k$-tuples} (i.e.,
  $k$-runs). He first gives a recurrence for $w_k(n,0)$ with
  $n\ge k+2$, then a recurrence for $w_k(n,1)$ with $n\ge k+2$, lastly
  a recurrence for $w_k(n,m)$ with $m\ge 2$, $k\ge 1$ and
  $n\ge m k + m-1$ [in our notation, $e\ge 0$, see~\eqref{eq:excess}].

  Any enumeration that might be obtained through Apostol's various
  recurrences can also be obtained through
  either~\eqref{eq:recurrence_w_n_m_k_basic}
  or~\eqref{eq:recurrence_w_n_m_k_simpler} in a simpler and more
  general way---our recurrences above are valid for all $n\ge 1$,
  $m\ge 0$ and $k\ge 0$ when taking~\eqref{eq:m_necessary_condition}
  into account. All this comes down to the necessity of considering
  the cases $n=-1$ and $n=0$ in the initialisation ---see related
  comments at the very end of Section~\ref{sec:conclusions}. The
  importance of the initial values~\eqref{eq:w_n_m_k_init}
  and~\eqref{eq:w_n_m_k_init2} cannot be understated: with them, all
  of Apostol's recurrences would become a special case
  of~\eqref{eq:recurrence_w_n_m_k_simpler}.
  
  To conclude this remark, in the specific case of $m=0$
  recurrence~\eqref{eq:recurrence_w_n_m_k_simpler} becomes
  \begin{equation}
    \label{eq:recurrence_m0}
    w_{k}(n,0)=2\,w_{k}(n-1,0)  -w_{k}\big(n-(k+1),0\big)+w_{k}\big(n-(k+2),0\big),
  \end{equation}
  i.e., we recover the same recurrence given by Madden in~\cite[Case
  4]{madden17:_gener_funct}---note that $w_k(n,0)=N(n+1,k+1)$ in this
  author's notation. See further comments about Madden's work in
  Remark~\ref{rmk:apostol_madden}. Of course,~\eqref{eq:recurrence_m0}
  is also one of Apostol's recurrences~\cite[Eq.\ (15)]{apostol88:_binary}.
\end{remark}

\subsubsection{Generating Functions}
\label{sec:generating-function}
From~\eqref{eq:ogf_w_n_m_klek}, the ogf  $W_{k}(x,y)=\sum_{n,m} w_k(n,m)\, x^n y^m$ is
\begin{equation}\label{eq:ogf_w_n_m_k}
  W_k(x,y)=\frac{1-x}{x\big(1-2x+(1-y)(x^{k+1}-x^{k+2})\big)},
\end{equation}
and from~\eqref{eq:ogf_m_klek} we get
 \begin{equation}
  \label{eq:ogf_m}
  [y^m]W_k(x,y)=x^{m(k+1)-1}\bigg(\frac{1-x}{1-2x+x^{k+1}-x^{k+2}}\bigg)^{m+1},
\end{equation}
which is the ogf enumerating the binary strings that contain
exactly $m$ $k$-runs.

\begin{remark}\label{rmk:apostol_madden}
  The ogf~\eqref{eq:ogf_m} was first given by Apostol~\cite[Eq.\
  (25)]{apostol88:_binary}, who determined it through his
  aforementioned recurrences. Koutras and Papastavridis~\cite[Sec.\
  5.c]{koutras93:_number} also found~\eqref{eq:ogf_m} through the
  particularisation of a general result based on occupancy models
  (i.e., distributions of balls into cells). Finally, this ogf is also
  Madden's main theorem in~\cite{madden17:_gener_funct}.  Bearing in
  mind that ---except for~\cite{sinha09}--- all prior work on this
  problem was unbeknownst to him, Madden observes
  regarding~\eqref{eq:ogf_m}: ``\textit{This result is quite
    elementary, but we have not been able to find it in any other
    source}''. Rediscoveries and reworkings like this are not uncommon
  in the long history of the theory of runs, and we will find more as
  the paper progresses. In any case, they show that fundamental
  problems concerning runs still merit a closer look. 
\end{remark}  

Finally, from~\eqref{eq:w_n_klek_atleast_ogf}, the ogf $W^{\ge
  m}_k(x)=\sum_{n}w^{\ge m}_k(n)\,x^n$ is
\begin{equation}
  \label{eq:w_n_k_atleast_ogf}
  W^{\ge m}_{k}(x)=\frac{1-x}{(1-2x)\,x }\Bigg(\frac{x^{k+1}-x^{k+2}}{1-2x+x^{k+1}-x^{k+2}}\Bigg)^m,
\end{equation}
which enumerates the  binary strings that contain at
least $m$ $k$-runs ---see~\eqref{eq:wn_k_atleast}.

\subsubsection{Explicit Expression}
\label{sec:explicit-expression}
We find next an explicit expression for $w_k(n,m)$. We proceed by applying successive power series expansions in order to
then identify the coefficient of $x^n$ in~\eqref{eq:ogf_m}. Through
elementary algebraic manipulations, and using $e$ defined
in~\eqref{eq:excess}, we can rewrite~\eqref{eq:ogf_m} as
$[y^m]W_k(x,y)=x^{n-e}\, \big(1-x\,((1-x)^{-1}-x^{k})\big)^{-(m+1)}$.
We now expand this expression using in succession
$(1-x)^{-v}=\sum_{r\ge 0} {r+v-1\choose r} x^r$ (negative binomial
theorem), the binomial theorem, and again the negative binomial theorem, to get
\begin{align}
 [y^m]W_k(x,y) 
  &=\sum_{r\ge 0}{r+m\choose r} \sum_{s\ge 0} {r\choose s}
    (-1)^{r-s}\sum_{t\ge 0}{t+s-1\choose t}\, x^{n-e+r(k+1)-sk+t}.  \label{eq:expansion}
\end{align}
Next, we have to identify the coefficient of $x^n$ in
$[y^m]W_k(x,y)$, which, from~\eqref{eq:expansion}, requires finding
the ranges of indices~$r$, $s$ and~$t$ in that expression for which
\begin{align}\label{eq:mrst}
  r(k+1)-s k+t=e.
\end{align}
Equation~\eqref{eq:mrst} determines a single value of $t$ for any
given $(r,s)$ pair. So, as long as we guarantee that $t\ge 0$, we only
need to focus our attention on~$r$ and~$s$. Because $s\le r$ (as
otherwise ${r\choose s}=0$), we have
from \eqref{eq:mrst} that $t\le e-r$. 
So, the highest value of $r$ for which~$t$ can have a nonnegative
value is~$e$. All this considered, the coefficient of $x^n$
in~\eqref{eq:expansion}, or, equivalently, $[x^ny^m]W_k(x,y)$, is
\begin{equation}
  \label{eq:w_gen}
  w_k(n,m)=\sum_{r=0}^e{r+m\choose r} \sum_{s=0}^r(-1)^{r-s} {r\choose s}
  {e-1-r(k+1)+s(k+1)\choose e-r(k+1)+sk~~}.
\end{equation}
Expression~\eqref{eq:w_gen} is similar to the one previously derived
by Makri and Psillakis~\cite[Eq.\ (8)]{makri11:_bernoul}, which they
obtained as a corollary to an earlier probabilistic result of
theirs~\cite{makri07:_succes}. Like their expression,~\eqref{eq:w_gen}
is simpler than the expressions given by Magliveras and
Wei~\cite[Thm.\ 2.3]{magliveras96:_enumer} ---the earliest explicit
result--- or by Sinha and Sinha~\cite[(3)]{sinha09}, because it
involves two summations rather than three. 

\begin{remark}
  Concurring with Madden's comments in~\cite[p.\ 2]
  {madden17:_gener_funct}, we also found that Sinha and Sinha's
  expression for $w_k(n,m)$~\cite[Eq.\ (3)]{sinha09} (in their
  notation,~$N_n^{m,k}$) is not always in agreement with the true
  counts.  Their expression for~$N_n^{m,k}$ is essentially correct,
  but it fails when $e=0$ [in our notation,
  see~\eqref{eq:excess}]. One can fix this minor oversight by starting
  the summations in~\cite[Eqs.\ (2)\ and\ (3)]{sinha09} at $r=0$ rather
  than at~$r=1$.
\end{remark}

\subsubsection{OEIS}
\label{sec:oeis-wknm}
Madden~\cite{madden17:_gener_funct} has pointed out that most of the
number sequences that emanate from $w_k(n,m)$ are not in the
OEIS. Below we give the few ones that we have found to be documented
already.
  
\begin{itemize}

\item $m=0$ (as indicated by Madden):

  $w_1(n-2,0)$ is \seqnum{A005251} for $n\ge 2$. %

  $w_2(n-3,0)$ is \seqnum{A049856} for $n\ge 3$. %

  $w_3(n-1,0)$ is \seqnum{A108758} for $n\ge 1$. %

\item $m=1$

  $w_1(n,1)$ is \seqnum{A079662} (number of occurrences of $1$ in all
  compositions of $n$ without $2$'s). %

\item $w_1(n-1,m)$ is the $m$th column of \seqnum{A105114} for $n\ge 1$.

  $w_2(n-1,m)$ is the $m$th column of \seqnum{A218796} for $n\ge 1$.
  
\item $k=0$ (null runs)
  
  $w_0(n,0)$ is \seqnum{A000045} (Fibonacci numbers). %

   $w_0(n,1)$ is \seqnum{A006367} (number of binary vectors of length
   $n+1$ beginning with $0$ and containing just $1$ singleton). %

   $w_{0}(n+1,2)$ is \seqnum{A105423} (number of compositions of $n+2$
   having exactly two parts equal to $1$). %

 \item Only sequences from $w^{\ge m}_k(n)$ with $m=1$ are in the
   OEIS; in this case, $w^{\ge 1}_k(n)=2^n-w_k(n,0)$:

  $w^{\ge 1}_{1}(n)$ is \seqnum{A384153} ---see also Section~\ref{sec:oeis_wl_n_m_k}.

  $w^{\ge 1}_{2}(n)$ is \seqnum{A177795} (number of length $n$ binary words that have at least one maximal run of 1's having length two).

\end{itemize}

\subsection[Number of n-Strings that Contain Exactly m
  (≥k)-Runs]{Number of $n$-Strings that Contain Exactly $m$  ($\ge\!k$)-Runs}
\label{sec:w_n_m_gtk}

Let us next address the enumeration of the $n$-strings that contain
exactly $m$ $(\ge\! k)$-runs (i.e., runs of length $k$ or longer),
which we denote by $w_{\ge k}(n,m)$. Notice that such $n$-strings may
also contain other runs shorter than $k$.

Explicit solutions to this problem were previously given by
Muselli~\cite{muselli96:_simple} ---in the context of probabilistic
success runs--- and by Makri and Psillakis~\cite{makri11:_bernoul}. In
the special case $m=0$, which is of particular relevance in
Section~\ref{sec:number-ws_n_k_ltk},
a generating function and an explicit expression were given by
Sedgewick and Flajolet~\cite{sedgewick96:_introd}, and a recurrence by
Nyblom~\cite{nyblom12:_enumerating}. A probabilistic recurrence and a
probability generating function were also given by Balakrishnan and
Koutras~\cite{balakrishnan01:_runs_scans} ---see comments about these
results in Remark~\ref{rmk:balakrishnan} (Section~\ref{sec:pi_n_m_gtk}).

If $k>n$ then $w_{\ge k}(n,m)=0$, and so we assume $k\le n$. Therefore,
the problem is a special case of the enumeration in
Section~\ref{sec:w_n_m_klek} with $\ushort{k}=k$ and $\widebar{k}=n$,
and thus
\begin{equation}
  \label{eq:w_n_m_gtk_from_klek}
  w_{\ge k}(n,m)=w_{k\le n}(n,m).
\end{equation}
We also study the number of $n$-strings that contain at least $m$
$(\ge\! k)$-runs, which from~\eqref{eq:w_n_m_gtk_from_klek}
and~\eqref{eq:wn_klek_atleast}, is
\begin{equation}\label{eq:wn_gtk_atleast}
  w_{\ge k}^{\ge m}(n)=\sum_{t=m}^{\floor*{\frac{n+1}{k+1}}}w_{\ge k}(n,t).
\end{equation}

\subsubsection{Recurrences}
\label{sec:recurrence_gtk}
From~\eqref{eq:recurrence_w_n_m_klek_basic},~\eqref{eq:recurrence_w_n_m_klek_simpler}
and~\eqref{eq:w_n_m_gtk_from_klek}, two recurrence relations for
$w_{\ge k}(n,m)$ are
\begin{equation}  \label{eq:recurrence_gtk_basic}
  w_{\ge k}(n,m)=\sum_{i=0}^{k-1}w_{\ge k}\big(n-(i+1),m\big)+\sum_{i=k}^{n}w_{\ge k}\big(n-(i+1),m-1\big),
\end{equation}
and
\begin{equation}
  \label{eq:recurrence_gtk_alt}
  w_{\ge k}(n,m)= 2\,w_{\ge k}(n-1,m)
                   + w_{\ge k}\big(n-(k+1),m-1\big)-w_{\ge k}\big(n-(k+1),m\big).
\end{equation}
which from~\eqref{eq:w_n_m_klek_init} and~\eqref{eq:w_n_m_klek_init2}  are both initialised by
\begin{align}\label{eq:init_w_n_m_gtk}
  w_{\ge k}(-1,m)=& \iv{m=0}\\
  w_{\ge k}(0,m)=& \ivl{m=\iv{k=0}}.\label{eq:init_w_n_m_gtk2}
\end{align}

\subsubsection{Generating Functions}
\label{sec:generating-function-gtk}
We next obtain the ogf
$W_{\ge k}(x,y)=\sum_{n,m} w_{\ge k}(n,m)\, x^n y^m$. Although a valid
ogf is directly obtained by setting $\ushort{k}=k$ and $\widebar{k}=n$
in~\eqref{eq:ogf_w_n_m_klek}, the resulting expression is
$n$-dependent. We can derive an alternative $n$-independent ogf, valid
for all $k\ge 0$, from~\eqref{eq:recurrence_gtk_alt} using the same
procedure as in Section~\ref{sec:generating-function_klek}. This
yields
\begin{equation}  \label{eq:ogf_gtk}
  W_{\ge k}(x,y)=\frac{1-x}{x\big(1-2x +(1-y)x^{k+1}\big)}.
\end{equation}
Through the same method used
to get~\eqref{eq:ogf_m_klek} from~\eqref{eq:ogf_w_n_m_klek}, we now
obtain from~\eqref{eq:ogf_gtk}
\begin{equation}\label{eq:ogf_m_gtk}
  [y^m]W_{\ge k}(x,y)=x^{m(k+1)-1}\frac{1-x}{\big(1-2x+x^{k+1}\big)^{m+1}},
\end{equation}
which is the ogf enumerating the  binary strings that
contain exactly $m$ $(\ge\!k)$-runs.

\begin{remark}\label{rmk:ogf-gtk-sedgewick}
  We make here some comments for the case $m=0$. In this
  case, recurrence~\eqref{eq:recurrence_gtk_basic} becomes
  \begin{equation}
    \label{eq:recurrence_gtk_basic_m0}
    w_{\ge k}(n,0)=\sum_{i=0}^{k-1}w_{\ge k}\big(n-(i+1),0\big),
  \end{equation}
  which was given by Nyblom~\cite[Thm.\ 2.1]{nyblom12:_enumerating}
  (with validity $k\ge 2$, and with different initialisation). This
  author was unaware of the fact that the same recurrence had been
  previously given by Schilling
  ---see~\eqref{eq:ws_n_ltk_w_n_m_gtk_connection}
  and~\eqref{eq:rec_ws_n_lek_schilling} in
  Section~\ref{sec:number-ws_n_k_ltk}.  It is also interesting to
  compare the ogf~\eqref{eq:ogf_m_gtk} in this special case, i.e.,
  $[y^{0}]W_{\ge k}(x,y)=x^{-1}(1-x)/(1-2x+x^{k+1})$, with the
  equivalent ogf given by Sedgewick and Flajolet~\cite[p.\
  368]{sedgewick96:_introd} %
  for enumerating the binary strings devoid of runs of $k$ consecutive
  zeros: $S_k(x)=(1-x^k)/(1-2x+x^{k+1})$. The comparison is possible
  because, when $m=0$, the definition of a run of $k$ consecutive
  equal bits by Sedgewick and Flajolet is equivalent to our definition
  of a $(\ge\!k)$-run ---see Remark~\ref{rmk:intro}.  Even though
  these two generating functions are different, they have the same
  coefficients in their power series expansion for $n\ge 1$.
\end{remark}

Finally, we may obtain
$W^{\ge m}_{\ge k}(x)=\sum_{n}w^{\ge m}_{\ge k}(n)\,x^n$ directly
from~\eqref{eq:w_n_klek_atleast_ogf}, but this leads to an
$n$-dependent expression. We can derive a more general ogf, valid for
all $k\ge 0$, by working from~\eqref{eq:wn_klek_atleast_rec}
specialised to this case. Doing so yields
\begin{align}
  \label{eq:w_n_gtk_atleast_ogf}
  W_{\ge k}^{\ge m}(x)&=\frac{x^{k+1}}{1-2x}\,[y^{m-1}]W_{\ge k}(x,y)\nonumber\\
  &=\frac{1-x}{(1-2x)\,x }\Bigg(\frac{x^{k+1}}{1-2x+x^{k+1}}\Bigg)^m,
\end{align}
where we have used~\eqref{eq:ogf_m_gtk} in the second step above.
This is the ogf enumerating the binary strings that contain at least
$m$  $(\ge\! k)$-runs ---see~\eqref{eq:wn_gtk_atleast}.

\subsubsection{Explicit Expression}
\label{sec:explicit-expression-1}
To obtain an explicit expression for $w_{\ge k}(n,m)$ we first apply
in succession the negative binomial theorem and the binomial theorem
to the denominator of~\eqref{eq:ogf_m_gtk}. With this expansion we can
express~\eqref{eq:ogf_m_gtk} as
\begin{equation}\label{eq:ogf_m_gtk_2}
  [y^m]W_{\ge k}(x,y)=(1-x)\sum_{r\ge 0}{m+r\choose r}\sum_{s\ge 0}{r\choose s}(-1)^{(r-s)}2^sx^{n-e+(k+1)r-ks},
\end{equation}
where we have used~\eqref{eq:excess} in the exponent of $x$.  In order
to determine the coefficient of $x^n$ in the above expression we have
to find the values of indices $r$ and~$s$ that solve each of the
following two Diophantine equations:
\begin{align}
  \label{eq:diophantine1}
  (k+1)r-ks&=e,\\
  \label{eq:diophantine2}
  (k+1)r-ks&=e-1.
\end{align}
Since
$\gcd(k+1,-k)=1$, we can solve both equations. By inspection, a particular solution to~\eqref{eq:diophantine1} is
$r=s=e$. Thus the general solution is of the form $r=e-k\, t$ and
$s=e-(k+1)\,t$ for integer $t$, which fulfils $s\le r$ for
nonnegative~$t$. As we also need to guarantee $r\ge 0$ and $s\ge 0$,
the valid range of $t$ is $0 \le t \le \floor{e/(k+1)}$. Similarly, a particular solution to~\eqref{eq:diophantine2} is
$r=s=e-1$, and thus the general solution is of the form $r=e-1-k\, t$
and $s=e-1-(k+1)\,t$ for integer $t$, and the range of $t$ is
$0 \le t \le \floor{(e-1)/(k+1)}$.

Combining these solutions we can see that the coefficient of $x^n$
in~\eqref{eq:ogf_m_gtk_2}, or, equivalently, $[x^ny^m]W_{\ge k}(x,y)$,
is
\begin{align}
  \label{eq:w_gtk_gen_diophant}
  w_{\ge k}(n,m)=\sum_{t=0}^{\floor{\frac{e}{k+1}}}  (-1)^{t}\,2^{e-(k+1)t}\Bigg(&{m+e-k t\choose e-kt}{e-kt\choose e-(k+1)t}\nonumber\\
                  &-\frac{1}{2}{m+e-1-kt\choose e-1-kt}{e-1-kt\choose e-1-(k+1)t}\Bigg).
\end{align}
Observe that we use just one summation on $t$, rather than the
difference of two summations on $t$ with different ranges, by taking
advantage of the fact that the last binomial coefficient becomes zero
for $t>\lfloor(e-1)/(k+1)\rfloor$. This single-summation formula is clearly
different from the previous double-summation formulas given by
Muselli~\cite[Thm.\ 1 with $p=q=\frac{1}{2}$ gives
$w_{\ge k}(n,m)/2^n$]{muselli96:_simple} and Makri and
Psillakis~\cite[Eq.\ (23)]{makri11:_bernoul} ---which closely resemble
each other, and which were derived through combinatorial
analysis. However Muselli also simplified his double-summation expression into
a single-summation formula
comparable
to~\eqref{eq:w_gtk_gen_diophant}~\cite[Thm.\ 3]{muselli96:_simple}.

The case $m=0$ is especially relevant, as we see in more detail in
Section~\ref{sec:number-ws_n_k_ltk} ---see also
Remark~\ref{rmk:waiting} below. %
With this argument choice, expression~\eqref{eq:w_gtk_gen_diophant}
becomes
\begin{equation}
  \label{eq:w_gtk_gen_diophant_m0}
  w_{\ge k}(n,0)=\sum_{t=0}^{\floor{\frac{n+1}{k+1}}}  (-1)^{t}\,2^{n+1-(k+1)t}\Bigg({n+1-kt\choose n+1-(k+1)t}-\frac{1}{2}{n-kt\choose n-(k+1)t}\Bigg),
\end{equation}
which is very similar to the explicit expression given by Sedgewick
and Flajolet~\cite[p.\
370]{sedgewick96:_introd}%
---in their notation, $w_{\ge k}(n,0)=[x^n]S_k(x)$.

\begin{remark}\label{rmk:fibs}
  As pointed out by Sedgewick and
  Flajolet~\cite[p.\ 369]{sedgewick96:_introd}, 
  \begin{equation}
    \label{eq:w_n_0_gtk_fibonacci}
    w_{\ge k}(n,0)=F^{(k)}_{n+k},
  \end{equation}
  where $F^{(j)}_i$ are the $j$th order Fibonacci numbers, also called
  Fibonacci $j$-step numbers, $j$-nacci numbers, which are defined
  by the recurrence
  \begin{equation}
    \label{eq:nth_order_fibonacci_standard}
    F^{(j)}_{i}=\sum_{r=1}^j F^{(j)}_{i-r}.
  \end{equation}
  for $j\ge 1$ and $i\ge j$, initialised with $F^{(j)}_{j-1}=1$ and
  $F^{(j)}_i=0$ for $0\le i<j-1$. The case $j=2$ gives the standard
  Fibonacci numbers, but observe that we also allow the degenerate
  case $j=1$, in which case $F_i^{(1)}=1$ for all $i\ge 1$. The
  connection described by~\eqref{eq:w_n_0_gtk_fibonacci} is clear when
  observing recurrence~\eqref{eq:recurrence_gtk_basic_m0} ---see
  also~\cite[p.\ 368]{sedgewick96:_introd}--- or the ogf in
  Remark~\ref{rmk:ogf-gtk-sedgewick}.
  Considering~\eqref{eq:w_n_klek_atleast1_from w_n_0_klek}, a
  consequence of~\eqref{eq:w_n_0_gtk_fibonacci} is
  \begin{equation}
    \label{eq:altleast_fib}
    w^{\ge 1}_{\ge k}(n)=2^n-F^{(k)}_{n+k}.
  \end{equation}
\end{remark}

\begin{remark}\label{rmk:waiting}
  The special case $w_{\ge k}(n,0)$ with $k\ge 1$ can also be used to
  enumerate the $n$-strings that feature their first consecutive
  appearance of $k$ ones at index $i$, which we denote
  by~$\widehat{w}_k(n,i)$. Notice that we are not only referring to
  the first appearance of a $k$-run: for the avoidance of doubt, index
  $i$ marks the position of the $k$th one in the first uninterrupted
  sequence of~$k$~ones, which may be followed by a zero, a one, or the
  end of the string. To obtain $\widehat{w}_k(n,i)$ we just have to
  make the following observation: if the first $k$ consecutive ones
  occur at index $i$, then the initial $(i-(k+1))$-substring must be
  devoid of $(\ge\!  k)$-runs. For every such start of the $n$-string
  we have $2^{n-i}$ possible endings. Combining these two facts, we
  thus have that
  \begin{equation}
    \label{eq:what_n_i}
    \widehat{w}_k(n,i)=w_{\ge k}(i-(k+1),0)\,\,2^{n-i}\,\iv{i\le n}.
  \end{equation}
  Note that this formula works even if $i< k$, in which case
  $w_{\ge k}(i-(k+1),m)=0$ through necessary
  condition~\eqref{eq:m_necessary_condition}. The
  enumeration~$\widehat{w}_k(n,i)$ is relevant in waiting-time type
  problems. Bearing in mind~\eqref{eq:w_n_0_gtk_fibonacci}, the reader
  is referred to~\cite[Sec.\ 2.8]{balakrishnan01:_runs_scans} for a
  literature review  on the role of $j$th order Fibonacci numbers
  in the context of waiting-time problems with
  runs. %
\end{remark}

\subsubsection{OEIS}
\label{sec:oeis-wgtknm}
We report below the number sequences stemming from $w_{\ge k}(n,m)$
that are listed the OEIS.
 
\begin{itemize}
  
\item $m=0$ (see Remark~\ref{rmk:fibs})

  $w_{\ge 2}(n-2,0)$ is \seqnum{A000045} (Fibonacci numbers) for $n\ge 2$.
  
  $w_{\ge 3}(n-3,0)$ is \seqnum{A000073} (tribonacci numbers) for $n\ge 3$.

  $w_{\ge 4}(n-4,0)$ is \seqnum{A000078} (tetranacci numbers) for $n\ge 4$.

  $w_{\ge 5}(n-5,0)$ is \seqnum{A001591} (pentanacci numbers) for $n\ge 5$.

  etc

\item $m=1$

  $w_{\ge 2}(n-1,1)$ is \seqnum{A006478} for $n\ge 1$. %

\item $w_{\ge 1}(n-1,m)$ is the $m$th column of \seqnum{A034839} for
  $n\ge 1$.

  $w_{\ge 2}(n,m)$ is the $m$th column of \seqnum{A334658}.

\item $k=1$: apart from~\seqnum{A034839}, many sequences with this
  parameter are individually documented in the OEIS for different
  values of $m$ ---see~\eqref{eq:w_n_m} and~\eqref{eq:w_n_m_explicit}
  in Section~\ref{sec:anylength}.

  $w_{\ge 1}(n,1)$ is \seqnum{A000217} (triangular numbers, or
  ${n+1\choose 2}$). %
  
  $w_{\ge 1}(n-1,2)$ is \seqnum{A000332} (binomial
  coefficient ${n\choose 4}$) for $n\ge 1$. %
  
  $w_{\ge 1}(n-1,3)$ is \seqnum{A000579} (binomial
  coefficient ${n\choose 6}$) for $n\ge 1$. %

    $w_{\ge 1}(n-1,4)$ is \seqnum{A000581} (binomial
    coefficient ${n\choose 8}$) for $n\ge 1$. %
    
  $w_{\ge 1}(n-1,5)$ is \seqnum{A001287} (binomial
  coefficient ${n\choose 10}$) for $n\ge 1$. %

  For $6\le m\le 24$, $w_{\ge 1}(n-1,m)$ is \seqnum{A010965}$+2(m-6)$ (binomial
  coefficient ${n\choose 2m}$) for $n\ge 1$. %
  
\item $k=0$ (number of $n$-strings that contain exactly $m$ runs,
  including null runs): we can see from the ogf~\eqref{eq:ogf_m_gtk}
  that $w_{\ge 0}(n,m)={n\choose m-1}$, and therefore this sequence
  features in many OEIS entries. See also Section~\ref{sec:anylength}.
  
\item Only $w^{\ge m}_{\ge k}(n)$ sequences having $m=1$ are in the
  OEIS---see~\eqref{eq:wn_gtk_atleast}. Of course,
  from~\eqref{eq:w_n_klek_atleast1_from w_n_0_klek} we have $w^{\ge
    1}_{\ge k}(n)=2^n-w_{\ge k}(n,0)$.

  $w^{\ge 1}_{\ge 1}(n)=2^n-1$ is \seqnum{A000225} (sometimes called
  Mersenne numbers).
  
  $w^{\ge 1}_{\ge 2}(n)$ is \seqnum{A008466}.

  $w^{\ge 1}_{\ge 3}(n)$ is \seqnum{A050231}.

  $w^{\ge 1}_{\ge 4}(n)$ is \seqnum{A050232}. %

  $w^{\ge 1}_{\ge 5}(n)$ is \seqnum{A050233}. 

  $w^{\ge 1}_{\ge 6}(n)$ is \seqnum{A143662}.

  $w^{\ge 1}_{\ge 7}(n)$ is \seqnum{A151975}.

  $w^{\ge 1}_{\ge k}(n)$ is the $k$th column  of \seqnum{A050227}.

  Weisstein discusses $w^{\ge 1}_{\ge k}(n)$ in~\cite{weisstein}  in connection
  with Feller's work on probabilistic success
  runs~\cite{feller68:_introd_probab}. Nevertheless, we should note
  that $w^{\ge 1}_{\ge k}(n)$ was actually first studied by de
  Moivre~\cite{moivre38:_doctrine}, and the probability generating
  function attributed in~\cite{weisstein} to Feller was first found by
  Laplace ---see more details in Remarks~\ref{rmk:demoivre_laplace}
  and~\ref{rmk:wait_prob}  towards the end of Section~\ref{sec:pi_n_m_gtk}.

  Finally, because of~\eqref{eq:altleast_fib}, each of the OEIS
  sequences above can be paired with a corresponding~$k$th order Fibonacci
  sequence (mentioned at the start of this section).
\end{itemize}

\subsection[Number of n-Strings that Contain Exactly m Nonnull
(≤k)-Runs]{Number of $n$-Strings that Contain Exactly $m$ Nonnull
  ($\le\! k$)-Runs}
\label{sec:w_n_m_ltk}

Let us next address what is essentially the counterpart of the problem
in the previous section: the enumeration of the $n$-strings containing
exactly $m$ nonnull $(\le\! k)$-runs (i.e., runs of length~$k$
or shorter, but strictly greater than zero), which we denote
by~$w_{\le k} (n, m)$.  Hence, in this section we assume $k\ge
1$. Notice that the $n$-strings that we enumerate may also contain any
number of null runs, and/or of runs with lengths longer than $k$.

Austin and Guy~\cite{austin78:_binar} gave a recurrence and a
semi-explicit expression for ``\textit{the number of binary sequences
  of length $n$ in which the ones occur only in blocks of length at
  least $k$}'' (where $k\ge 2$), that is to say, for the special case
$w_{\le (k-1)}(n,0)$ of the enumeration that we consider in this
section. Other than in~\cite{austin78:_binar}, $w_{\le k}(n,m)$ seems
not to have received any attention in the literature ---which is
somewhat surprising, given that this enumeration looks like a mere
variation of $w_{\ge k}(n,m)$. A plausible reason for this is that it
may be harder to address this problem through direct combinatorial
analysis. This hypothesis is suggested by the explicit
expression~\eqref{eq:w_n_m_ltk_explicit} for $w_{\le k}(n,m)$ that we
derive later ---a triple-summation expression, as opposed to the single sum
expression~\eqref{eq:w_gtk_gen_diophant} for $w_{\ge k}(n,m)$. Another
possible explanation might be the perception that the problem in this
section has fewer applications.  Whatever the cause,
at least one good reason for studying $w_{\le k}(n,m)$ is that the
special case $w_{\le k}(n,0)$ plays a key role in
Section~\ref{sec:number-wl_n_k_gtk}, just like we have mentioned that
$w_{\ge k}(n,0)$ plays a key role in
Section~\ref{sec:number-ws_n_k_ltk}.

Again, the problem at hand is a special case of the enumeration in
Section~\ref{sec:w_n_m_klek}, here using $\ushort{k}=1$ and
$\widebar{k}=k$, and thus
\begin{equation}
  \label{eq:w_n_m_ltk_from_klek}
  w_{\le k}(n,m)=w_{1\le k}(n,m).
\end{equation}
We also study the number of $n$-strings that contain at least $m$
$(\le\! k)$-runs, which, from~\eqref{eq:w_n_m_ltk_from_klek}
and~\eqref{eq:wn_klek_atleast}, is
\begin{equation}\label{eq:wn_ltk_atleast}
  w_{\le k}^{\ge m}(n)=\sum_{t=m}^{\floor*{\frac{n+1}{2}}}w_{\le k}(n,t).
\end{equation}

\subsubsection{Recurrences}
\label{sec:recurrences_w_n_m_ltk}
From~\eqref{eq:recurrence_w_n_m_klek_basic},~\eqref{eq:recurrence_w_n_m_klek_simpler}
and~\eqref{eq:w_n_m_ltk_from_klek}, two recurrence relations
for $w_{\ge k}(n,m)$ are
\begin{align}  \label{eq:recurrence_w_n_m_ltk_basic}
  w_{\le k}(n,m)=~&w_{\le k}(n-1,m)+\sum_{i=k+1}^{n}w_{\le
                  k}\big(n-(i+1),m\big)%
                +\sum_{i=1}^{k}w_{\le k}\big(n-(i+1),m-1\big),
\end{align}
and
\begin{align}
  \label{eq:recurrence_w_n_m_ltk_simpler}
  w_{\le k}(n,m)=2\,w_{\le k}(n-1,m)%
                  & + w_{\le k}(n-2,m-1)-w_{\le k}(n-2,m)\nonumber\\ 
                  & - w_{\le k}\big(n-(k+2),m-1\big)+w_{\le k}\big(n-(k+2),m\big),
\end{align}
both of which, from~\eqref{eq:w_n_m_ltk_from_klek},~\eqref{eq:w_n_m_klek_init}
and~\eqref{eq:w_n_m_klek_init2}, are initialised by
\begin{equation}
  \label{eq:init_w_n_m_lek}
  w_{\le k}(-1,m)=w_{\le k}(0,m)=\iv{m=0}.
\end{equation}

\begin{remark}
  In the special case $m=0$,
  recurrence~\eqref{eq:recurrence_w_n_m_ltk_simpler} becomes
  \begin{equation}
    \label{eq:austin_guy}
    w_{\le k}(n,0)=2\,w_{\le k}(n-1,0)-w_{\le k}(n-2,0)+w_{\le k}\big(n-(k+2),0\big),    
  \end{equation}
  which was given by Austin and Guy for
  $a_n^{(k)}=w_{\le(k-1)}(n,0)$~\cite[Eq.\
  (1\textquotesingle)]{austin78:_binar}. Importantly, these authors
  used the very same initialisation strategy as above, i.e.,~\eqref{eq:init_w_n_m_lek}.
\end{remark}
\subsubsection{Generating Functions}
\label{sec:generating-function-w_n_m_ltk}
We next obtain $W_{\le k}(x,y)=\sum_{n,m}w_{\le
  k}(n,m)\,x^ny^m$. Setting $\ushort{k}=1$ and $\widebar{k}=k$
in~\eqref{eq:ogf_w_n_m_klek} we have
\begin{equation}
  \label{eq:bivariate_ogf_w_n_m_ltk}
  W_{\le k}(x,y)=\frac{1-x}{x\big(1-2x+(1-y)(x^2-x^{k+2})\big)},
\end{equation}
whereas~\eqref{eq:ogf_m_klek} yields in turn
\begin{equation}
  \label{eq:ym_ogf_w_n_m_ltk}
  [y^m]W_{\le k}(x,y)=\frac{x^{2m-1} (1-x)}{1-x^k}\bigg(\frac{1-x^k}{1-2x+x^2-x^{k+2}}\bigg)^{m+1}.
\end{equation}
This is the ogf enumerating the binary strings that contain
exactly~$m$ nonnull $(\le\! k)$-runs.

Finally, from~\eqref{eq:w_n_klek_atleast_ogf},
$W^{\ge m}_{\le k}(x)=\sum_{n}w^{\ge m}_{\le k}(n)\,x^n$ is
\begin{equation}
  \label{eq:w_n_ltk_atleast_ogf}
  W^{\ge m}_{\le k}(x)=\frac{1-x}{(1-2x)\,x }\Bigg(\frac{x^{2}-x^{k+2}}{1-2x+x^{2}-x^{k+2}}\Bigg)^m,
\end{equation}
which is the ogf enumerating the  binary strings that contain at
least $m$ nonnull $(\le\! k)$-runs ---see~\eqref{eq:wn_ltk_atleast}.

\subsubsection{Explicit Expression}
\label{sec:explicit-expression-w_n_m_lek}
We now derive an explicit expression for $w_{\le k}(n,m)$ by finding
$[x^ny^n]W_{\le k}(x,y)$, or, equivalently, the coefficient of $x^n$ in
the power series expansion of~\eqref{eq:ym_ogf_w_n_m_ltk}. To do so,
we rewrite this expression as
$[y^m]W_{\le k}(x,y)=x^{2m-1}(1-x^k)^m (1-x)^{-(2m+1)}
\big(1-x^{k+2}/(1-x)^2\big)^{-(m+1)}$ and we develop each of the
last three factors in it into power series. Using the negative
binomial series in the first of these factors, the binomial theorem in
the second one, and again the negative binomial theorem (twice) in the
third one, we have that
\begin{align*}
  \big(1-x^k\big)^m&=\sum_{r\ge 0}{m\choose r}(-1)^r x^{k\,r},\\
  (1-x)^{-(2m+1)}&=\sum_{s\ge 0}{s+2m\choose s} x^s,\\
  \bigg(1-\frac{x^{k+2}}{(1-x)^2}\bigg)^{-(m+1)}&=\sum_{t\ge 0}{t+m\choose t}
    x^{(k+2)\,t}\sum_{p\ge 0}{p+2t-1\choose p}x^p.
\end{align*}
Collecting the exponents of $x$, we see that we have to determine next
the ranges of the summation indices $r, s, t$ and $p$ that fulfil
$2m-1+k\,r+s+(k+2)\,t+p=n$, or, using~\eqref{eq:excess},
$k\,r+s+(k+2)\,t+p=e$. Let us define at this point the auxiliary variable
\begin{equation*}
  \label{eq:d}
  d=e-k\,r.
\end{equation*}
As all indices are nonnegative, $s+(k+2)\,t+p\ge 0$. So we must
guarantee $d\ge 0$ in order to have solutions to our problem, which
implies~$r\le \floor{e/k}$. We next have to solve $d=s+(k+2)\,t+p$ for
integers $s, t$ and $p$, where $s,t,p\ge
0$. %
The maximum of $t$ happens when $s=p=0$, and thus
$t\le \floor{d/(k+2)}$.  Given $t$, the maximum of $s$ happens when
$p=0$, and thus $s\le d-(k+2)\,t$. Finally, given $t$ and $s$, we have
that $p=d-(k+2)\,t-s$.

Collecting all these solutions, we can see that the coefficient of $x^n$
in~\eqref{eq:ym_ogf_w_n_m_ltk} is
\begin{align}
  \label{eq:w_n_m_ltk_explicit}
  w_{\le k}(n,m)=&\!\!\sum_{r=0}^{\floor{e/k}}(-1)^r{m\choose
              r}\sum_{t=0}^{\floor{d/(k+2)}}{t+m\choose t}\!\! \sum_{s=0}^{d-(k+2)t}
  \!\!{s+2m \choose s}{d-kt-s-1\choose d-(k+2)t-s}.
\end{align}
In the  case~$m=0$, which is especially relevant in
Section~\ref{sec:number-wl_n_k_gtk}, the triple-summation expression above
simplifies considerably:
\begin{align}
  \label{eq:w_n_m_ltk_explicit_m0}
  w_{\le k}(n,0)%
&=1+\sum_{t=1}^{\floor*{\frac{n+1}{k}}}\frac{n+1-kt}{2t}\, {n-kt
  \choose n+1-(k+2)t}.
\end{align}
This expression also constitutes an alternative to the semi-explicit
expression $w_{\le(k-1)}(n,0)=\text{round}(c_k\gamma_k^n)$ given by
Austin and Guy~\cite[p.\ 85]{austin78:_binar}, which is based on
asymptotic considerations and requires determining a special constant
$c_k$ and finding the real root $\gamma_k$ of a $k$th order
polynomial.

\subsubsection{OEIS}
\label{sec:oeis_w_n_m_ltk}
As in previous sections, we report the few sequences
emanating from $w_{\le k}(n,m)$ that we have been able to find in the
OEIS.

\begin{itemize}
\item $m=0$ 
  
$w_{\le 1}(n-2,0)=w_1(n-2,0)$  is \seqnum{A005251} for $n\ge 2$.

$w_{\le 2}(n-1,0)$ is \seqnum{A005252} for $n\ge 1$. %

$w_{\le 3}(n,0)$ is \seqnum{A005253} (number of binary
words of length $n$ in which the ones occur only in blocks of length at
least 4).

$w_{\le 4}(n-2,0)$  is \seqnum{A005689} (number of Twopins positions)
for $n\ge 6$. %

$w_{\le 5}(n-1,0)$  is \seqnum{A098574} for $n\ge 1$. %

$w_{\le 6}(n+6,0)$ is \seqnum{A217838} (number of $n$ element $0..1$
arrays with each element the minimum of $7$ adjacent elements of a
random $0..1$ array of $n+6$ elements).  %

  These OEIS sequences actually led us to finding reference~\cite{austin78:_binar}.

\item $m=1$

  $w_{\le 1}(n,1)=w_1(n,1)$ is \seqnum{A079662} (number of
  occurrences of $1$ in all compositions of $n$ without $2$'s). %
  
\item Sequences from $w^{\ge m}_{\le k}(n)$ ---see~\eqref{eq:wn_ltk_atleast}: 

  $w^{\ge 1}_{\le 1}(n)=w^{\ge 1}_1(n)$ is \seqnum{A384153} ---see also Section~\ref{sec:oeis_wl_n_m_k}. 

\end{itemize}

\subsection{Number of $n$-Strings that Contain Exactly $m$ Nonnull Runs}
\label{sec:anylength}
In this section we study the number of $n$-strings that contain
exactly $m$ nonnull runs (of arbitrary lengths, all strictly greater
than zero), which we denote by
$w(n,m)$. As may be expected, this problem is also a special case of
the enumeration in Section~\ref{sec:w_n_m_klek} with $\ushort{k}=1$
and $\widebar{k}=n$, and thus
\begin{equation}
  \label{eq:w_n_m_from_klek}
  w(n,m)=w_{1\le n}(n,m).
\end{equation}
Alternatively, using the enumerations in Sections~\ref{sec:w_n_m_gtk}
and~\ref{sec:w_n_m_ltk} we may also write
\begin{equation}
  \label{eq:w_n_m}
  w(n,m)=w_{\ge 1}(n,m)=w_{\le n}(n,m).
\end{equation}
We also study the number of $n$-strings that contain at least $m$
runs, which from~\eqref{eq:w_n_m_from_klek}
and~\eqref{eq:wn_klek_atleast} is
\begin{equation}\label{eq:wn_atleast}
  w^{\ge m}(n)=\sum_{t=m}^{\floor*{\frac{n+1}{2}}}w(n,t).
\end{equation}

\subsubsection{Recurrences}
\label{sec:recurrences-w_n_m}
From~\eqref{eq:w_n_m_from_klek},~\eqref{eq:recurrence_w_n_m_klek_basic}
and~\eqref{eq:recurrence_w_n_m_klek_simpler} two recurrences for
$w(n,m)$ are:
\begin{equation*}
  w(n,m)=w(n-1,m)+\sum_{i=1}^{n}w(n-(i+1),m-1),
\end{equation*}
and
\begin{equation*}
  w(n,m)=2\, w(n-1,m)+w(n-2,m-1)-w(n-2,m),
\end{equation*}
both of which,
from~\eqref{eq:w_n_m_from_klek},~\eqref{eq:w_n_m_klek_init}
and~\eqref{eq:w_n_m_klek_init2}, have initial values
\begin{equation*}
  w(-1,m)=w(0,m)=\iv{m=0}.
\end{equation*}
\subsubsection{Generating Functions}
\label{sec:generating-function-w_n_m}
From~\eqref{eq:w_n_m_from_klek} and~\eqref{eq:w_n_m}, we may get the
ogf $W(x,y)=\sum_{n,m}w(n,m)\,x^ny^m$ in two different ways: by
setting $k=n$ in~\eqref{eq:bivariate_ogf_w_n_m_ltk} ---which is the
same as setting $\ushort{k}=1$ and $\widebar{k}=n$
in~\eqref{eq:ogf_w_n_m_klek}--- or else by setting $k=1$
in~\eqref{eq:ogf_gtk}. The first option gives an $n$-dependent
ogf, but the second option gives a simpler more
general version:
\begin{equation}
  \label{eq:ogf_w_n_m_bivariate}
  W(x,y)=\frac{1-x}{x\big(1-2x +(1-y)x^2\big)}.
\end{equation}
On the other hand, setting $k=1$ in~\eqref{eq:ogf_m_gtk}, yields
\begin{equation}
  \label{eq:ogf_w_n_m}
  [y^m]W(x,y)=\frac{x^{2m-1}}{(1-x)^{2m+1}},
\end{equation}
which is the ogf enumerating the binary strings that contain exactly
$m$ nonnull runs.

Finally, setting $k=1$ in~\eqref{eq:w_n_gtk_atleast_ogf} we have that
$W^{\ge m}(x)=\sum_{n}w^{\ge m}(n)\,x^n$ is
\begin{equation}
  \label{eq:w_n_ltk_atleast_ogf}
  W^{\ge m}(x)=\frac{1}{1-2x}\Bigg(\frac{x}{1-x}\Bigg)^{2m-1},
\end{equation}
which is the ogf enumerating the  binary strings that contain at
least $m$ nonnull runs ---see~\eqref{eq:wn_atleast}.

\subsubsection{Explicit Expressions}
Using~\eqref{eq:w_n_m}, we can get explicit expressions for $w(n,m)$
by setting $k=1$ in~\eqref{eq:w_gtk_gen_diophant} or by setting~$k=n$
in~\eqref{eq:w_n_m_ltk_explicit}, but the resulting formulas are not
obviously simplifiable ---especially the second one.  We may obtain a
much simpler expression by relying on~\eqref{eq:ogf_w_n_m}, which, by
applying the negative binomial theorem, can be expanded as
\begin{equation}
  \label{eq:ogf_w_n_m_expanded}
  [y^m]W(x,y)=\sum_{t\ge 0}{t+2m\choose 2m}x^{2m+t-1}.
\end{equation}
To determine $w(n,m)=[x^ny^m]W(x,y)$ we have to find the value of $t$ for which
$n=2m+t-1$, from which we get
\begin{equation}
  \label{eq:w_n_m_explicit}
  w(n,m)={n+1 \choose 2m}.
\end{equation}

\begin{remark}
  Goulden and Jackson~\cite[Sec.\ 2.4.4]{goulden83:_combin_enumer} found~\eqref{eq:w_n_m_explicit} using the
  following alternative ogf ---cf.~\eqref{eq:ogf_w_n_m_bivariate}:
  \begin{equation*}
    \label{eq:ogf_w_n_m_bivariate_goulden_jackson}
    W(x,y)=\frac{1+(y-1)x}{1-2x +(1-y)x^2}.
  \end{equation*}
  At any rate, expression~\eqref{eq:w_n_m_explicit} suggests
  that~$w(n,m)$ can also be found through basic combinatorial
  reasoning. One such combinatorial explanation is as follows:
  assume that the number of ones in an $n$-string that contains
  exactly $m$ nonnull runs is $t=m,\dots,n-(m-1)$ ---as the minimum is
  $m$ and there must be at least $m-1$ zeros. For a given $t$, there
  are ${t-1\choose t-m}$ ways to choose the lengths of the $m$ runs,
  and ${n-t+1 \choose n-t-m+1}$ ways to place the $n-t$ zeros around
  them. Using  ${a\choose b}={a\choose a-b}$, we can thus write
  \begin{equation}
    \label{eq:w_n_m_comb}
    w(n,m)=\sum_{t=m}^{n-m+1}{t-1\choose m-1}{n-t+1 \choose m}.
  \end{equation}
  We can also put this summation as
  $w(n,m)=\sum_{t=1}^{n+1}{t-1\choose m-1}{n-t+1 \choose m}$, because
  the first binomial coefficient in~\eqref{eq:w_n_m_comb} is zero for
  $t<m$ and the second one is zero for $t>n-m+1$.  Hence,
  $w(n,m)=\sum_{t=0}^n {t\choose m-1}{n-t\choose m}$, which is a
  variation of Vandermonde's convolution that adds up
  to~\eqref{eq:w_n_m_explicit} ---see~\cite[Eq.\
  (5.26)]{graham94:_concrete}. 
\end{remark}

\subsubsection{OEIS}
\label{sec:oeis_w_n_m}
We report here the OEIS sequences relevant to the enumerations in this section.
\begin{itemize}
\item   Many $w(n,m)$ sequences are in the OEIS for different values of $m$
  ---see sequences with $k=1$ in
  Section~\ref{sec:oeis-wgtknm}. %

\item As for $w^{\ge m}(n)$ ---see~\eqref{eq:wn_atleast}---, it yields
  the following OEIS sequences:

$w^{\ge 1}(n)=2^n-1$ is \seqnum{A000225}. 

$w^{\ge 2}(n)$ is \seqnum{A002662}. %

$w^{\ge 3}(n)$ is \seqnum{A002664}. %
 
$w^{\ge 4}(n)$ is \seqnum{A035039}. %

$w^{\ge 5}(n)$ is \seqnum{A035041}. %

These connections are due to the fact that, using
${n+1\choose 2t}={n\choose 2t}+ {n\choose 2t-1}$,
from~\eqref{eq:wn_atleast} and~\eqref{eq:w_n_m_explicit} we can write
\begin{align*}
  \label{eq:1}
  w^{\ge m}(n)& =\sum_{t=m}^{\floor*{(n+1)/2}}\bigg({n\choose 2t}+ {n\choose 2t-1}\bigg)=\sum_{t=2m \atop t\text{ even}}^{n}{n\choose t}+ \sum_{t=2m-1 \atop t\text{ odd}}^{n}{n\choose t}%
    =\sum_{t=2m-1}^{n}{n\choose t}\nonumber\\&=2^n-\sum_{t=0}^{2(m-1)}{n\choose  t},
\end{align*}
where in the second equality we have taken advantage of the fact that
${n\choose t}=0$ for $t>n$ to set $n$ as the upper limit in both summations.
The last expression is the common definition to all five OEIS
sequences above.
\end{itemize}

\subsection[Number of n-Strings Whose Longest Run Is a k-Run or a (≤k)-Run]{Number of $n$-Strings Whose Longest Run Is a $k$-Run or a $(\le\! k$)-Run}
\label{sec:number-ws_n_k_ltk}
We denote the number of $n$-strings whose longest run is a $k$-run by
$\widebar{w}_{k}(n)$, whereas~$\widebar{w}_{\le k}(n)$ represents the
number of $n$-strings whose longest run is a $(\le\! k)$-run. A
recurrence for~$\widebar{w}_{\le k}(n)$ was previously given by
Schilling~\cite{schilling90:_longest}, and, in the context of
probabilistic success runs, explicit expressions were given by
Godbole~\cite{godbole90:_specific} and
Muselli~\cite{muselli96:_simple}. Finally, Flajolet and
Sedgewick~\cite{flajolet09:_analytic} gave a generating function for
$\widebar{w}_{\le k}(n)$ constructed through the symbolic method, and
Prodinger gave its asymptotic behaviour~\cite{prodinger15}. As
for previous work on $\widebar{w}_k(n)$, probabilistic versions of
this enumeration were given by Muselli~\cite{muselli96:_simple} and by
Makri et al.~\cite{makri07:_shortest_longest}, and Schilling provided
asymptotic estimates of its distribution~\cite{schilling12:_surpr}.

For all $k\ge 0$, the two enumerations that we are interested in are
related as follows
\begin{equation}
  \label{eq:connection_ws_n_lek_ws_n_k}
  \widebar{w}_{ \le k}(n)=\sum_{j=0}^k \widebar{w}_{j}(n).
\end{equation}
Conversely, for $k\ge 1$ 
\begin{equation}
  \label{eq:connection_ws_n_lek_ws_n_k_converse}
  \widebar{w}_k(n)=  \widebar{w}_{ \le k}(n)-  \widebar{w}_{ \le (k-1)}(n),
\end{equation}
whereas $\widebar{w}_0(n)=\widebar{w}_{\le 0}(n)$.
Both enumerations follow directly from the results in
Section~\ref{sec:w_n_m_gtk} because of the following fact:
\begin{equation}
  \label{eq:ws_n_ltk_w_n_m_gtk_connection}
  \widebar{w}_{\le k}(n)=w_{\ge (k+1)}(n,0).
\end{equation}
Thus, $\widebar{w}_0(n)=w_{\ge 1}(n,0)=1$ for $n\ge 1$, as the longest
run is a null run only in the all-zeros $n$-string.

Of course, considering~\eqref{eq:ws_n_ltk_w_n_m_gtk_connection}
and~\eqref{eq:w_n_m_gtk_from_klek}, $\widebar{w}_{\le k}(n)$ can also
be expressed as a special case of $w_{\ushort{k}\le\widebar{k}}(n,m)$
in Section~\ref{sec:w_n_m_klek}, but
using~\eqref{eq:ws_n_ltk_w_n_m_gtk_connection} we are able to simplify
our presentation and directly obtain explicit expressions.

\subsubsection{Recurrences}
\label{sec:recurrences_ws}
We start by finding recurrences for $\widebar{w}_{\le k}(n)$.
Schilling's recurrence~\cite[Eq.\ (1)]{schilling90:_longest} is
recovered by directly applying
\eqref{eq:ws_n_ltk_w_n_m_gtk_connection}
to~\eqref{eq:recurrence_gtk_basic}:
\begin{equation} \label{eq:rec_ws_n_lek_schilling}
  \widebar{w}_{\le  k}(n)=\sum_{i=0}^k \widebar{w}_{\le k}(n-(i+1)).
\end{equation}
On the other hand, an alternative recurrence is obtained by directly
applying~\eqref{eq:ws_n_ltk_w_n_m_gtk_connection}
to~\eqref{eq:recurrence_gtk_alt}, which gives
\begin{equation} \label{eq:rec_ws_n_lek_simpler}
  \widebar{w}_{\le k}(n)=2\,\widebar{w}_{\le k}(n-1)-\widebar{w}_{\le k}(n-(k+2)).
\end{equation}
From~\eqref{eq:ws_n_ltk_w_n_m_gtk_connection},~\eqref{eq:init_w_n_m_gtk}
and~\eqref{eq:init_w_n_m_gtk2}, we can see that the initialisation
of both~\eqref{eq:rec_ws_n_lek_schilling}
and~\eqref{eq:rec_ws_n_lek_simpler} is %
\begin{equation}
  \label{eq:init_wlekn}
  \widebar{w}_{\le k}(-1)=\widebar{w}_{\le k}(0)=1.
\end{equation}

\begin{remark}
  Schilling initialised \eqref{eq:rec_ws_n_lek_schilling}
  in~\cite[Eq.\ (1)]{schilling90:_longest} using 
  \begin{equation}\label{eq:schilling_init}
    \widebar{w}_{\le k}(n)=2^n\quad\text{if }n\le k,
  \end{equation}
  but there is no need to resort to~\eqref{eq:schilling_init} when
  initialising~\eqref{eq:rec_ws_n_lek_schilling}
  with~\eqref{eq:init_wlekn}. Also, Bloom gave
  recurrence~\eqref{eq:rec_ws_n_lek_simpler} for $k=4$~\cite[Eq.\
  (9)]{bloom98:_singles}.
\end{remark}

We may also obtain recurrences for $\widebar{w}_k(n)$ for $k\ge 1$ by
relying on~\eqref{eq:rec_ws_n_lek_schilling}
and~\eqref{eq:rec_ws_n_lek_simpler}. Using~\eqref{eq:rec_ws_n_lek_schilling}
in~\eqref{eq:connection_ws_n_lek_ws_n_k_converse}, after some
elementary algebraic manipulations we get
\begin{equation}
  \label{eq:ws_n_k_rec_a}
  \widebar{w}_k(n)=\sum_{i=0}^{k-1}\widebar{w}_k(n-(i+1))+\sum_{j=0}^k \widebar{w}_j(n-(k+1)),
\end{equation}
whereas if we instead use~\eqref{eq:rec_ws_n_lek_simpler}
in~\eqref{eq:connection_ws_n_lek_ws_n_k_converse} we get the
following alternative recurrence:
\begin{equation}
  \label{eq:ws_n_k_rec_b}
  \widebar{w}_{k}(n)=2\,\widebar{w}_{k}(n-1)+\sum_{j=0}^{k-1}\widebar{w}_{j}(n-(k+1))
  -\sum_{j=0}^{k}\widebar{w}_{j}(n-(k+2)).
\end{equation}
From~\eqref{eq:connection_ws_n_lek_ws_n_k_converse}
and~\eqref{eq:init_wlekn}, it follows that
$\widebar{w}_{k}(-1)= \widebar{w}_{k}(0)=0$ for $k\ge 1$. As seen
from~\eqref{eq:ws_n_k_rec_a} and~\eqref{eq:ws_n_k_rec_b}, we also need
initialisation for $k=0$. As $\widebar{w}_0(n)=w_{\ge 1}(n,0)$, we
have from~\eqref{eq:init_w_n_m_gtk} and~\eqref{eq:init_w_n_m_gtk2}
that $\widebar{w}_0(-1)=\widebar{w}_0(0)=1$.  Thus, the initialisation
of recurrences \eqref{eq:ws_n_k_rec_a} and~\eqref{eq:ws_n_k_rec_b} is
\begin{equation*}
  \widebar{w}_{k}(-1)=  \widebar{w}_{k}(0)=\iv{k=0}.
\end{equation*}
The ease with which these recurrences are initialised shows that null
runs are inherent to this problem ---although the stronger reason for
taking null runs into account will be seen in Section~\ref{sec:oz}.

\subsubsection{Generating Functions}
\label{sec:generating-functions_ws}
We may obtain generating functions for $\widebar{w}_{\le k}(n)$ and
$\widebar{w}_{k}(n)$ by using the recurrences in the previous section,
although the reader may verify that this is easier for
$\widebar{w}_{\le k}(n)$ than for $\widebar{w}_k(n)$. Let us follow
instead the path of least resistance by
exploiting~\eqref{eq:ws_n_ltk_w_n_m_gtk_connection}
and~\eqref{eq:connection_ws_n_lek_ws_n_k_converse}.  The ogf
$\widebar{W}_{\!\!\le k}(x)=\sum_n\widebar{w}_{\le k}(x)\,x^n$ follows
directly from~\eqref{eq:ws_n_ltk_w_n_m_gtk_connection}
and~\eqref{eq:ogf_m_gtk}:
\begin{equation}
  \label{eq:ogf_ws_n_ltk2}
  \widebar{W}_{\!\!\le k}(x)=\frac{1-x}{x\big(1-2x+x^{k+2}\big)}.
\end{equation}
Likewise, for $k\ge 1$ the ogf
$\widebar{W}_{\!k}(x)=\sum_n\widebar{w}_{k}(x)\,x^n$ is directly
obtained using~\eqref{eq:ogf_ws_n_ltk2}
in~\eqref{eq:connection_ws_n_lek_ws_n_k_converse}, which yields
\begin{align}\label{eq:ogf_wsk}
  \widebar{W}_{\!k}(x)&=\widebar{W}_{\!\!\le k}(x)-\widebar{W}_{\!\!\le (k-1)}(x)\nonumber\\
  &= \frac{x^k(1-x)^{2}}{(1-2x+x^{k+1})(1-2x+x^{k+2})}.
\end{align}

\begin{remark}
  If we use~\eqref{eq:ogf_wsk}
  and~\eqref{eq:connection_ws_n_lek_ws_n_k} to recover
  $\widebar{W}_{\!\!\le k}(x)$ from $\widebar{W}_{\!k}(x)$, instead
  of~\eqref{eq:ogf_ws_n_ltk2} we get
  \begin{equation}
    \label{eq:ogf_ws_lek_2}
    \widebar{W}_{\!\!\le k}(x)=\frac{1-x^{k+1}}{1-2x+x^{k+2}},
  \end{equation}
  which was given by Flajolet and Sedgewick in~\cite[p.\
  51]{flajolet09:_analytic}. Although \eqref{eq:ogf_ws_lek_2} is
  different from~\eqref{eq:ogf_ws_n_ltk2}, it has the same
  coefficients in its power series expansion for~$n\ge 1$. Prodinger
  gives an asymptotic analysis of $\widebar{w}_{\le k}(n)$ using this
  ogf~\cite{prodinger15}.
\end{remark}

\subsubsection{Explicit Expressions}
\label{sec:explicit-solution}
Explicit expressions for $\widebar{w}_{k}(n)$ and
$\widebar{w}_{\le k}(n)$ are available
through~\eqref{eq:connection_ws_n_lek_ws_n_k_converse},
\eqref{eq:ws_n_ltk_w_n_m_gtk_connection}
and~\eqref{eq:w_gtk_gen_diophant_m0}. The resulting formula for
$\widebar{w}_{\le k}(n)$ is comparable to the simplified single-summation
combinatorial expression given by Muselli~\cite[Cor.\ 1 using
$p=q=1/2$ gives
$\widebar{w}_{\le(k-1)}(n)/2^n$]{muselli96:_simple}. In contrast, our
single-summation expression for $\widebar{w}_k(n)$ is much simpler than the
triple-summation expression given by Makri et al.\ \cite[Thm.\ 2.1.4 using
$p=q=1/2$ gives $\widebar{w}_k(n)/2^n$]{makri07:_shortest_longest}.

Let us comment on some further aspects of our results in this
section. Through\ \eqref{eq:ws_n_ltk_w_n_m_gtk_connection}
and~\eqref{eq:w_n_0_gtk_fibonacci} we have that
\begin{equation}\label{eq:wl_gtk_fibonacci}
  \widebar{w}_{\le k}(n)=F^{(k+1)}_{n+k+1}, 
\end{equation}
where $F^{(j)}_i$ are the $j$th order Fibonacci numbers
---see~\eqref{eq:nth_order_fibonacci_standard}. This is also clear
from
recurrence~\eqref{eq:rec_ws_n_lek_schilling}. Regarding~\eqref{eq:wl_gtk_fibonacci},
Schilling observed already the special case
$\widebar{w}_{\le 1}(n)=F^{(2)}_{n+2}$, when he mentioned that
``(\ldots) \textit{the number} (\ldots) \textit{of sequences of length~$n$ that contain no two
  consecutive heads is the $(n + 2)$nd Fibonacci number}''~\cite{schilling90:_longest}. See as
well Nyblom's comments regarding this case~\cite[Cor.\ 2.1]{nyblom12:_enumerating}. %

Additionally, because
from~\eqref{eq:connection_ws_n_lek_ws_n_k_converse} it follows that
$\widebar{w}_{1}(n)=\widebar{w}_{\le 1}(n)-1$, then we also have that
\begin{equation}
  \label{eq:ws_n_1_fibonacci}
  \widebar{w}_1(n)=F^{(2)}_{n+2}-1.
\end{equation}

\begin{remark}
  For $n$-strings drawn uniformly at random, the expected length of the
  longest run is $\bar{\ell}_n=\big(\sum_{k=0}^n k\, \widebar{w}_k(n)\big)/2^n$.
  Thus, by writing $\widebar{w}_k(n)=[x^n]\widebar{W}_{\!k}(x)$,  we have from
  \eqref{eq:ogf_wsk} that
  \begin{equation}
    \label{eq:avg_length}
    \bar{\ell}_n=\frac{1}{2^n}[x^n]\sum_{k=1}^n k\, \frac{x^k(1-x)^{2}}{(1-2x+x^{k+1})(1-2x+x^{k+2})}.
  \end{equation}
  This expression is equivalent to the one given by Sedgewick and
  Flajolet for the same quantity~\cite[p.\
  426]{sedgewick13:_introduct}, which in fact contains a small
  oversight: their summation must start at~$k=1$ rather than at~$k=0$.
  Of course, we can also get $\bar{\ell}_n$ as a double-summation
  explicit expression through~\eqref{eq:connection_ws_n_lek_ws_n_k_converse},
  \eqref{eq:ws_n_ltk_w_n_m_gtk_connection}
  and~\eqref{eq:w_gtk_gen_diophant_m0}.

  Using~\eqref{eq:avg_length} and a computer algebra system, we may
  also check the elegant asymptotic estimate
  $\hat{\bar{\ell}}_{n}=\gamma \ln 2 + \log_{2}(n/2)$ given by
  Schilling in~\cite{schilling12:_surpr} against the exact value
  $\bar{\ell}_n$ for some moderate values of $n$. For instance,
  $\hat{\bar{\ell}}_{1000}=9.3658$ and
  $\hat{\bar{\ell}}_{3000}=10.9508$, whereas
  $\bar{\ell}_{1000}=9.3000$ and $\bar{\ell}_{3000}=10.8839$.
\end{remark}

\subsubsection{OEIS}
\label{sec:oeis-wskn}
We report here the enumerations in this section that yield OEIS sequences.

\begin{itemize}
\item Sequences emanating from $\widebar{w}_k(n)$:
  
$\widebar{w}_1(n-2)$ is \seqnum{A000071} (Fibonacci numbers minus one)
for $n\ge 2$ ---cf.\ \eqref{eq:ws_n_1_fibonacci}. %

$\widebar{w}_j(n-1)$ for $j=2,3$ and $4$ are \seqnum{A000100},
\seqnum{A000102} and \seqnum{A006979}, respectively (number of
compositions of $n$ in which the maximal part is $j+1$) for $n\ge 1$. %

$\widebar{w}_5(n-1)$ is \seqnum{A006980} for $n\ge 6$. %

$\widebar{w}_k(n)$ is the $k$th column of~\seqnum{A048004}.

Additionally, $2^n\bar{\ell}_n$ is \seqnum{A119706}
---see~\eqref{eq:avg_length}. %

\item Sequences emanating from $\widebar{w}_{\le k}(n)$: see $m=0$ in
  Section~\ref{sec:oeis-wgtknm}
  ---recall~\eqref{eq:ws_n_ltk_w_n_m_gtk_connection}.

\end{itemize}

\subsection[Number of n-Strings Whose Shortest Nonnull Run Is a k-Run or a
(≥k)-Run]{Number of $n$-Strings Whose Shortest Nonnull Run Is a $k$-Run or a $(\ge\! k$)-Run}
\label{sec:number-wl_n_k_gtk}
Let us call $\ushort{w}_{k}(n)$ the number of $n$-strings whose
shortest nonnull run is a $k$-run, and let~$\ushort{w}_{\ge k}(n)$ be
the number of $n$-strings whose shortest run is a $(\ge\! k)$-run. We
therefore assume $k\ge 1$ in this section. These problems have been
studied by Makri et al.\ \cite{makri07:_shortest_longest} in the
context of probabilistic success runs.

Both quantities are related as follows:
\begin{equation*}
  \label{eq:wl_n_gtk_wl_n_k_connection}
  \ushort{w}_{\ge k}(n)=\sum_{j=k}^n \ushort{w}_{j}(n),
\end{equation*}
and so we also have
\begin{equation}
  \label{eq:wl_n_k_from_wl_n_gtk}
  \ushort{w}_{k}(n)=\ushort{w}_{\ge k}(n)-\ushort{w}_{\ge (k+1)}(n).
\end{equation}
The enumeration of $\ushort{w}_{\ge k}(n)$ and that of
$\ushort{w}_{k}(n)$ follow directly from the results in
Section~\ref{sec:w_n_m_ltk}, because for $k>1$ we have the following
relation:
\begin{equation}
  \label{eq:wl_n_gtk_w_n_m_ltk_connection}
  \ushort{w}_{\ge k}(n)=w_{\le (k-1)}(n,0)-1,
\end{equation}
whereas
\begin{equation}\label{eq:wl_n_gt1}
  \ushort{w}_{\ge 1}(n)=2^n-1.
\end{equation}
We subtract one both in \eqref{eq:wl_n_gtk_w_n_m_ltk_connection}
and in~\eqref{eq:wl_n_gt1} to discount the all-zeros $n$-string. Of
course, considering~\eqref{eq:w_n_m_ltk_from_klek}, we have that the
connection in~\eqref{eq:wl_n_gtk_w_n_m_ltk_connection} can also be
expressed using the general enumeration
$w_{\ushort{k}\le\widebar{k}}(n,m)$ studied in
Section~\ref{sec:w_n_m_klek},
but~\eqref{eq:wl_n_gtk_w_n_m_ltk_connection} allows us to simplify our
presentation and to obtain explicit expressions.

\subsubsection{Recurrences}
\label{sec:recurrences_wl}
From~\eqref{eq:wl_n_gtk_w_n_m_ltk_connection}
and~\eqref{eq:recurrence_w_n_m_ltk_basic}, a recurrence for
$\ushort{w}_{\ge k}(n)$ when $k>1$ is
\begin{equation}
  \label{eq:wl_n_gtk_rec}
  \ushort{w}_{\ge k}(n)=\ushort{w}_{\ge k}(n-1)+\sum_{i=k}^n\ushort{w}_{\ge k}(n-(i+1))+\max(n-k+1,0),
\end{equation}
and from~\eqref{eq:wl_n_gtk_w_n_m_ltk_connection}
and~\eqref{eq:recurrence_w_n_m_ltk_simpler} an alternative recurrence
for $k>1$ is
\begin{equation}
  \label{eq:wl_n_gtk_rec_alt}
  \ushort{w}_{\ge k}(n)=\,2\,\ushort{w}_{\ge k}(n-1)-\ushort{w}_{\ge
  k}(n-2)+\ushort{w}_{\ge k}(n-(k+1))+\iv{k\le n}.
\end{equation}
Considering~\eqref{eq:wl_n_gtk_w_n_m_ltk_connection} and~\eqref{eq:init_w_n_m_lek}, the initialisation of these two recurrences is
\begin{equation}
  \label{eq_wl_n_gtk_init}
  \ushort{w}_{\ge k}(-1)=\ushort{w}_{\ge k}(0)=0.
\end{equation}
From~\eqref{eq:wl_n_gtk_rec} and~\eqref{eq:wl_n_gtk_rec_alt} we may
obtain recurrences for $\ushort{w}_{k}(n)$ when $k>1$.
Using~\eqref{eq:wl_n_gtk_rec} in~\eqref{eq:wl_n_k_from_wl_n_gtk} the
first recurrence is
\begin{equation}
  \label{eq:wl_n_k_2_rec}
  \ushort{w}_{k}(n) =  \ushort{w}_{k}(n-1)+\sum_{i=k+1}^n\ushort{w}_{k}(n-(i+1))+\sum_{j=k}^n\ushort{w}_{j}(n-(k+1))+\iv{k\le n},
\end{equation}
where we have used $\max(n-k+1,0)-\max(n-k,0)=\iv{k\le n}$.
Using~\eqref{eq:wl_n_gtk_rec_alt} in~\eqref{eq:wl_n_k_from_wl_n_gtk}
we get the alternative recurrence
\begin{align}
  \label{eq:wl_n_k_2_rec_alt}
  \ushort{w}_{k}(n)
  =~&2\,\ushort{w}_{k}(n-1)-\ushort{w}_{k}(n-2)+\sum_{j=k}^n\ushort{w}_{j}(n-(k+1))
    -\sum_{j=k+1}^n\ushort{w}_{j}(n-(k+2))+\iv{k=n},
\end{align}
where we have used $\iv{k\le n}-\iv{k+1\le n}=\iv{k=n}$.
Inputting~\eqref{eq_wl_n_gtk_init} in~\eqref{eq:wl_n_k_from_wl_n_gtk},
we see that both~\eqref{eq:wl_n_k_2_rec}
and~\eqref{eq:wl_n_k_2_rec_alt} are initialised by
\begin{equation*}
  \ushort{w}_k(-1)=\ushort{w}_k(0)=0.
\end{equation*}
Finally, we address the $k=1$ case.
From~\eqref{eq:wl_n_k_from_wl_n_gtk} and~\eqref{eq:wl_n_gt1} we obtain
$\ushort{w}_1(n)=2^n-1-\ushort{w}_{\ge 2}(n)$, which can be evaluated
using recurrences~\eqref{eq:wl_n_gtk_rec}
or~\eqref{eq:wl_n_gtk_rec_alt}.

\begin{remark}
  The same strategy in~\eqref{eq:wl_n_k_from_wl_n_gtk} was used by
  Makri et al.\ \cite[Thm.\ 2.1.2]{makri07:_shortest_longest} to get,
  in their notation, $P(M_n=k)$ from an explicit computation for
  $P(M_n\ge k)$. These two probabilities correspond to  the
  enumerative quantities $\ushort{w}_k(n)$ and
  $\ushort{w}_{\ge k}(n)$. Also, recurrence~\eqref{eq:wl_n_gtk_rec} is
  the exact enumerative counterpart of the probabilistic recurrence
  given by Makri et al.\ for $P(M_n\ge k)$~\cite[Thm.\
  2.1.3]{makri07:_shortest_longest}. These authors derived their
  recurrence from their explicit expression, rather than the other way
  around, and they did not provide initialisation values.
\end{remark}

\subsubsection{Generating Functions}
\label{sec:generating-functions_wl}
We may derive generating functions from the recurrences in the
previous section, but the simplest way to obtain them is by
exploiting~\eqref{eq:wl_n_gtk_w_n_m_ltk_connection} and
\eqref{eq:wl_n_k_from_wl_n_gtk}.
From~\eqref{eq:wl_n_gtk_w_n_m_ltk_connection}
and~\eqref{eq:ym_ogf_w_n_m_ltk}, the ogf
$\ushort{W}_{\ge k}(x)=\sum_n \ushort{w}_{\ge k}(n)\, x^n$ for $k>1$
is
\begin{equation}
  \label{eq:wl_n_gtk_ogf}
  \ushort{W}_{\ge k}(x)=\frac{1-x}{x (1-2x+x^2-x^{k+1})}-\frac{1}{1-x},
\end{equation} %
whereas from~\eqref{eq:wl_n_gt1} we have
\begin{equation}
  \label{eq:wl_n_gtk_ogf_k1}
  \ushort{W}_{\ge 1}(x)=\frac{1}{1-2x}-\frac{1}{1-x}.
\end{equation} %
Finally, from~\eqref{eq:wl_n_k_from_wl_n_gtk}
and~\eqref{eq:wl_n_gtk_ogf}, for $k>1$ the ogf
$\ushort{W}_k(x)=\sum_n \ushort{w}_k(n)\, x^n$ is
\begin{align}
  \label{eq:wl_n_k_ogf}
  \ushort{W}_k(x)&=\ushort{W}_{\ge k}(x)-\ushort{W}_{\ge(k+1)}(x)\nonumber\\
  &=\frac{x^k (1-x)^2}{((1-x)^2-x^{k+1})((1-x)^2-x^{k+2})},
\end{align} %
whereas
from~\eqref{eq:wl_n_k_from_wl_n_gtk},~\eqref{eq:wl_n_gtk_ogf_k1}
and~\eqref{eq:wl_n_gtk_ogf} we have that
\begin{equation}
  \label{eq:wl_n_k_ogf_k1}
  \ushort{W}_1(x)=\frac{1}{1-2x}-\frac{1-x}{x(1-2x+x^2-x^3)}.
\end{equation} %
To conclude, we observe that both~\eqref{eq:wl_n_gtk_ogf}
and~\eqref{eq:wl_n_k_ogf} happen to be valid not only for $k>1$ but
also for $k=1$. That is to say, for~$k=1$ and $n\ge 1$
\eqref{eq:wl_n_gtk_ogf} and~\eqref{eq:wl_n_k_ogf} have the same
coefficients of $x^n$ as~\eqref{eq:wl_n_gtk_ogf_k1}
and~\eqref{eq:wl_n_k_ogf_k1}, respectively.

\subsubsection{Explicit Expressions}
Explicit expressions for $\ushort{w}_{\ge k}(n)$ and
$\ushort{w}_{k}(n)$ are available
through~\eqref{eq:wl_n_gtk_w_n_m_ltk_connection},~\eqref{eq:wl_n_k_from_wl_n_gtk}
and~\eqref{eq:w_n_m_ltk_explicit_m0}. These single-summation expressions are
simpler than the double-summation combinatorial expressions found by Makri et
al~\cite[Thms.\ 2.1.1 and 2.2.2 with $p=q=1/2$ and multiplied by
$2^n$]{makri07:_shortest_longest}.

\subsubsection{OEIS}
\label{sec:oeis_wl_n_m_k}
We report next any sequences related to $\ushort{w}_k(n)$ and
$\ushort{w}_{\ge k}(n)$ found in the OEIS.

\begin{itemize}
  
\item Sequences emanating from $\ushort{w}_k(n)$: none of the
  $\ushort{w}_k(n)$ sequences, for $k\ge 1$, were in the OEIS previous
  to this work; below are the currently listed sequences:

  $\ushort{w}_1(n)=w^{\ge 1}_1(n)=w^{\ge 1}_{\le 1}(n)$ is \seqnum{A384153}.

  $\ushort{w}_2(n)$ is \seqnum{A384154}.
  
  $\ushort{w}_3(n)$ is \seqnum{A384155}.

  $\ushort{w}_k(n)$ is the $k$th column of \seqnum{A388718}.
  
\item The sequences emanating from $\ushort{w}_{\ge k}(n)$ in the OEIS are:

  $\ushort{w}_{\ge 1}(n)=2^n-1$ is \seqnum{A000225}.

  $\ushort{w}_{\ge 2}(n+2)$ is \seqnum{A077855}.
  
  $\ushort{w}_{\ge 3}(n)$ is \seqnum{A130578} (Number of different
  possible rows ---or columns--- in an $n\times n$ crossword puzzle).
  
  $\ushort{w}_{\ge 4}(n)$ is \seqnum{A209231} (Number of binary words of length $n$ such that there is at least one $0$ and every run of consecutive $0$'s is of length $\ge 4$).

  $\ushort{w}_{\ge k}(n)$ is the $k$th column of \seqnum{A388547}.

\end{itemize}

\subsection{Number of $n$-Strings that Contain Exactly $m$ Nonnull
  $p$-Parity Runs}
\label{sec:wp_n_m}

Let $w_{[p]}(n,m)$ represent the number of $n$-strings that contain
exactly $m$ nonnull $p$-parity runs, i.e., exactly $m$ runs whose
lengths $k_1,\dots,k_m$ are strictly greater than zero and have parity
$\text{mod}(k_i,2)=p$ for $i=1,\dots,m$. The $n$-strings that we
enumerate may have more than~$m$ nonnull runs, as long as the lengths
of all additional runs have parity opposite to $p$, and they may also
have any number of null runs. For the first time in this paper we
address an enumeration unrelated to the results in
Section~\ref{sec:w_n_m_klek}. To the best of our knowledge, the only
authors that have addressed a similar problem are Grimaldi and
Heubach~\cite{grimaldi05:_without_odd_runs}, who studied the number of
$n$-strings devoid of odd runs ---a special case of $w_{[p]}(n,m)$
with $p=1$ and~$m=0$.

First of all, we state a necessary condition similar
to~\eqref{eq:m_necessary_condition}.
\begin{condition}{(Existence of $n$-strings containing $m$ nonnull $p$-parity runs)}
  \begin{equation}
    \label{eq:m_necessary_condition_p}
    w_{[p]}(n,m)> 0\quad \Longrightarrow\quad 0\le m \le \floor[\bigg]{\frac{n+1}{2+\iv{p=0}}}.
\end{equation}
\end{condition}
This is really the same condition as~\eqref{eq:m_necessary_condition}
in Section~\ref{sec:w_n_m_klek}, just noting that the lengths of all
nonnull $p$-parity runs are lower bounded by $\ushort{k}=1+\iv{p=0}$.

We also study the number of $n$-strings that contain at least $m$
nonnull $p$-parity runs, which we denote by $w^m_{[p]}(n)$. This enumeration
can be obtained from $w_{[p]}(n,m)$ as
\begin{equation}\label{eq:wp_n_atleast}
  w^{\ge m}_{[p]}(n)=\sum_{t=m}^{\floor*{\frac{n+1}{2+\iv{p=0}}}}w_{[p]}(n,t).
\end{equation}

\subsubsection{Recurrences}
\label{sec:recurrences_wp_n_m}
We can produce a recurrence to enumerate $w_{[p]}(n,m)$ with a similar
strategy as in Section~\ref{sec:recurrences_w_n_m_klek}. Since we are
not counting null runs, we can split the quantity $w_{[p]}(n,m)$ into
two contributions:
\begin{enumerate}[label=\alph*)]
\item The $n$-strings that begin with $0$  contribute
  $w_{[p]}(n-1,m)$ to $w_{[p]}(n,m)$.
\item As for the $n$-strings that begin with $1$, those that start
  with an odd $i$-run contribute $w_{[p]}(n-(i+1),m-\iv{p=1})$
  to~$w_{[p]}(n,m)$, whereas those that start with an even nonnull
  $i$-run contribute $w_{[p]}(n-(i+1),m-\iv{p=0})$. Equivalently, the
  $n$-strings that start with a nonnull $i$-run contribute
  $w_{[p]}(n-(i+1),m-\iv{p=\text{mod}(i,2)})$ to~$w_{[p]}(n,m)$.
\end{enumerate}
Collecting these two contributions we get the following recurrence:
\begin{equation}
  \label{eq:wp_n_m_rec2}
  w_{[p]}(n,m)=w_{[p]}(n-1,m)+\sum_{i=1}^n w_{[p]}\big(n-(i+1),m-\iv{p=\text{mod}(i,2)}\big).
\end{equation}
To find initial values for recurrence~\eqref{eq:wp_n_m_rec2} we consider
the case $n=1$, in which we know by inspection that
\begin{equation}
  \label{eq:trivial_n1_wp}
  w_{[p]}(1,m)=2\,\iv{p=0}\iv{m=0}+\iv{p=1}\iv{m=1}.
\end{equation}
On the other hand, setting~$n=1$ in~\eqref{eq:wp_n_m_rec2} and taking
into account~\eqref{eq:m_necessary_condition_p} yields
\begin{equation}
  \label{eq:wp_n1}
  w_{[p]}(1,m)=w_{[p]}(0,m)+w_{[p]}(-1,m-\iv{p=1}).
\end{equation}
We wish~\eqref{eq:wp_n1} to equal~\eqref{eq:trivial_n1_wp}. 
Taking into account~\eqref{eq:m_necessary_condition_p}, equality
between the two expressions is 
 achieved by choosing
\begin{equation}
  \label{eq:wp_n_m_init}
  w_{[p]}(-1,m)=w_{[p]}(0,m)=\iv{m=0},
\end{equation}
which are therefore the initial values of~\eqref{eq:wp_n_m_rec2}. %

We produce next a recurrence for $w^{\ge m}_{[p]}(n)$
by applying~\eqref{eq:wp_n_atleast} to~\eqref{eq:wp_n_m_rec2}, which yields
\begin{equation}
  \label{eq:wp_n_atleast_rec}
  w^{\ge m}_{[p]}(n)=w^{\ge m}_{[p]}(n-1)+\sum_{i=1}^n \bigg(w^{\ge m}_{[p]}(n-(i+1))+\iv{p=\text{mod}(i,2)}\,w_{[p]}(n-(i+1),m-1)\bigg).
\end{equation}
Even if this recurrence depends on $w_{[p]}(n,m)$,  we see in the next section
that it suffices to obtain the ogf of $w^{\ge m}_{[p]}(n)$.

\subsubsection{Generating Functions}
\label{sec:generating-function-wp}
Let us next obtain the ogf
$W_{[p]}(x,y)=\sum_{n,m} w_{[p]}(n,m)\, x^ny^m$. First of all, we
obtain a recurrence valid for all $n$ and $m$.  Setting $n=-1$
in~\eqref{eq:wp_n_m_rec2} we get $w_{[p]}(-1,m)=0$ instead of the
correct value $w_{[p]}(-1,m)=\iv{m=0}$ given
by~\eqref{eq:wp_n_m_init}. We can can ``fix'' this by adding
$\iv{n=-1}\iv{m=0}$ to~\eqref{eq:wp_n_m_rec2}. On the other hand,
setting $n=0$ in~\eqref{eq:wp_n_m_rec2} after this change yields
$w_{[p]}(0,m)=\iv{m=0}$, which is correct. We thus have an extended
recurrence valid for all $n$ and $m$, which we use subsequently.

Next, we need a recurrence without a full-history summation to be able
to deduce the ogf. We can obtain one by calculating $w_{[p]}(n,m)-w_{[p]}(n-2,m)$
using the extended recurrence. This yields
\begin{align}
  \label{eq:wp_n_m_rec3_diff}
  w_{[p]}(n,m)-w_{[p]}(n-2,m)=~&w_{[p]}(n-1,m)-w_{[p]}(n-3,m)\nonumber\\
                               &+w_{[p]}\big(n-2,m-\iv{p=1}\big)
                               +w_{[p]}\big(n-3,m-\iv{p=0}\big)\nonumber\\
  &+\big(\iv{n=-1}-\iv{n=1}\big)\iv{m=0}.
\end{align}
\begin{remark}
  Unlike in Sections~\ref{sec:recurrences_w_n_m_klek}
  and~\ref{sec:generating-function_klek}, obtaining the above
  difference before making~\eqref{eq:wp_n_m_rec2} valid for all~$n$
  and~$m$ does not render a valid recurrence. In
  Section~\ref{sec:recurrences_w_n_m_klek} we were able to first
  separately reason each recurrence and then verify that the second
  one could be obtained from the first, but this is not the case
  here.  In general, we can freely operate with a recurrence to
  obtain a new equivalent recursive relation as long as the original
  recurrence is valid for all integer values of its arguments. If this
  is not true, then the manipulation is not guaranteed to render a
  valid relation.
\end{remark}
By multiplying next~\eqref{eq:wp_n_m_rec3_diff} on both sides by
$x^ny^m$ and then adding on $n$ and~$m$ we find
\begin{align*}
  W_{[p]}(x,y)-x^2 W_{[p]}(x,y)=~&x\, W_{[p]}(x,y)-x^3 W_{[p]}(x,y)+ x^2y^{\iv{p=1}} W_{[p]}(x,y)\nonumber\\
                               &+x^3y^{\iv{p=0}} W_{[p]}(x,y)+x^{-1}+x,
\end{align*}
which yields the ogf
\begin{equation}
  \label{eq:wp_bivariate_ogf}
  W_{[p]}(x,y)=\frac{1-x^2}{x\big(1-x-x^2\big(1+y^{\iv{p=1}}\big)+x^3\big(1-y^{\iv{p=0}}\big)\big)}.
\end{equation}
Just like in Section~\ref{sec:generating-function_klek}, obtaining
$[y^m]W_{[p]}(x,y)$ is straightforward by putting~\eqref{eq:wp_bivariate_ogf} as a function of $(1-c y)^{-1}$. Doing so we get
\begin{equation}
  \label{eq:ogf_wp_m}
  [y^m]W_{[p]}(x,y)=\frac{(1-x^2)\,x^{(2+\iv{p=0})m}}{x\,\big(1-x-(1+\iv{p=0})\,x^2+\iv{p=0}\,x^3\big)^{m+1}},
\end{equation}
which is the ogf enumerating the  binary strings that
contain exactly~$m$ nonnull $p$-parity runs.

\begin{remark}\label{rmk:grimaldi}
  In the special case $p=1$ and $m=0$ (i.e., number of $n$-strings
  devoid of odd runs) the ogf~\eqref{eq:ogf_wp_m} is
  $[y^0]W_{[1]}(x,y)=(1-x^2)/\big(x(1-x-x^2)\big)$, which, as indicated
  by Grimaldi and Heubach~\cite[Thm.\
  2]{grimaldi05:_without_odd_runs}, corresponds to
  \begin{equation}
    \label{eq:w_n_m0_p1_fib}
    w_{[1]}(n,0)=F_{n+1}^{(2)},
  \end{equation}
  i.e., the Fibonacci sequence shifted one position.
\end{remark}

To conclude this section, we derive the ogf $W^{\ge m}_{[p]}(x)=\sum_n
w^{\ge m}_{[p]}(n)\, x^n$. We first obtain a recurrence free from $n$-dependent
summations by subtracting $w^{\ge m}_{[p]}(n-2)$ from $w^{\ge m}_{[p]}(n)$
using~\eqref{eq:wp_n_atleast_rec}, which gives
\begin{align}
  \label{eq:wp_n_atleast_rec2}
  w^{\ge m}_{[p]}(n)-w^{\ge m}_{[p]}(n-2)=&~w^{\ge m}_{[p]}(n-1)-w^{\ge m}_{[p]}(n-3)\nonumber\\&+\iv{p=1}w_{[p]}(n-2,m-1)+\iv{p=0}w_{[p]}(n-3,m-1).
\end{align}
As usual, by  multiplying~\eqref{eq:wp_n_atleast_rec2} on both
sides by $x^n$ and then adding over $n$ we get 
\begin{align}
  \label{eq:ogf_wn_p_atleast}
  W^{\ge m}_{[p]}(x)&=\frac{\iv{p=1}\,x^2+\iv{p=0}\,x^3}{1-x-2x^2}\,[y^{m-1}]W_{[p]}(x,y)\nonumber\\
              &=\frac{(1-x^2)\big(\iv{p=1}\,x+\iv{p=0}\,x^2\big)\big(x^{2+\iv{p=0}}\big)^{m-1}}{(1-x-2x^2)\big(1-x-(1+\iv{p=0})\,x^2+\iv{p=0}\,x^3\big)^m},
\end{align}
where we have used~\eqref{eq:ogf_wp_m} in the last
step. This is the ogf enumerating the binary strings that contain at
least $m$ nonnull $p$-parity runs---see~\eqref{eq:wp_n_atleast}.

\subsubsection{Explicit Expressions}\label{sec:w_n_m_p_explicit-expressions}
In this section we find explicit expressions for the two cases of
$w_{[p]}(n,m)$.  We first deal with the case
$p=1$. Rewriting~\eqref{eq:ogf_wp_m} as
$[y^m]W_{[1]}(x,y)=x^{2m-1}(1-x^2)^{-m} (1-x/(1-x^2))^{-m-1}$ and
applying the negative binomial theorem twice, we can expand this ogf as
\begin{equation}
  \label{eq:w_n_m_p_expansion}
  [y^m]W_{[1]}(x,y)=\sum_{s\ge 0} \sum_{t\ge 0} {s+m\choose
    m}{t+s+m-1\choose t} x^{s+2t+2m-1}.
\end{equation}
We now find the coefficient of $x^n$ in this expression by finding the
nonnegative indices $s$ and $t$ such that $s+2t+2m-1=n$. The minimum of
$t$ happens for $s=0$, and therefore $t\le \floor{(n-2m+1)/2}$. Given
a value
of $t$, we have $s=n-2m-t+1$. Considering these solutions
and~\eqref{eq:w_n_m_p_expansion} we thus have that 
\begin{equation}
  \label{eq:w_n_m_p1_explicit}
  w_{[1]}(n,m)=\sum_{t=0}^{\floor{\frac{n-2m+1}{2}}}{n-m-2t+1\choose m}{n-m-t\choose t}.
\end{equation}

As for the case $p=0$, we now rewrite~\eqref{eq:ogf_wp_m} as
$[y^m]W_{[0]}(x,y)=x^{3m-1}(1-x^2)^{-m} (1-x-x^2/(1-x^2))^{-m-1}$
before applying in succession the negative binomial theorem, the
binomial theorem, and again the negative binomial theorem. This yields
the expansion
\begin{equation}
  \label{eq:w_n_m_p_expansion0}
  [y^m]W_{[0]}(x,y)=\sum_{s\ge 0} \sum_{t\ge 0} \sum_{r\ge 0} {s+m\choose
    m}{s\choose t}{r+t+m-1\choose r} x^{s+t+2r+3m-1}.
\end{equation}
To extract the coefficient of $x^n$ we have to find the indices that
fulfil $s+t+2r+3m-1=n$. By setting $s=t=0$ we have the upper bound
$r\le \floor{(n-3m+1)/2}$. Likewise, by setting $t=0$ for a fixed value
of $r$ we have $t\le n-2r-3m+1$. Finally, $t$ is determined by the
equation above given $s$ and $t$. All this considered, we have
from~\eqref{eq:w_n_m_p_expansion0} that
\begin{equation}
  \label{eq:w_n_m_p1_explicit}
  w_{[0]}(n,m)=\sum_{r=0}^{\floor{\frac{n-3m+1}{2}}}\sum_{s=0}^{n-2r-3m+1}{s+m\choose
    s}{s\choose n-2r-s-3m+1}{n-r-s-2m \choose r}.
\end{equation}

\subsubsection{OEIS}
\label{sec:oeis-wp}
The sequences in the OEIS  emanating from the enumeration studied in
this section are:
\begin{itemize}
\item $m=0$

  $w_{[1]}(n-1,0)$ is \seqnum{A000045} (Fibonacci numbers) ---see
  Remark~\ref{rmk:grimaldi}.
  
  $w_{[0]}(n-1,0)$ is \seqnum{A028495}. %

\item $m=1$

  $w_{[1]}(n,1)$ is \seqnum{A029907}. %

  $w_{[0]}(n,1)$ is \seqnum{A384497}. %

\item $w_{[1]}(n,m)$ is the $m$th column of \seqnum{A119473}.

  $w_{[0]}(n,m)$ is the $m$th column of \seqnum{A391669}. %
  
\item Sequences from $w^{\ge m}_{[p]}(n)$ ---see \eqref{eq:wp_n_atleast}.

  $w^{\ge 1}_{[1]}(n)$ is \seqnum{A027934}. %

  $w^{\ge 1}_{[0]}(n)$ is \seqnum{A387332}. %
\end{itemize}

\section{Extensions to Probabilistic Runs}
\label{sec:probability}
At several points in Section~\ref{sec:number-n-strings} we have
specialised probabilistic results about success runs from other
authors (e.g.,~\cite{muselli96:_simple,makri07:_shortest_longest}) in
order to compare them to our enumerations. Thus, one might conclude
that the enumerative results that we have given so far are less
general than their probabilistic counterparts.  But this is a two-way
street: as we see in this section, every enumerative recurrence in
Section~\ref{sec:number-n-strings} also has a direct probabilistic
translation into the case in which the $n$-strings are outcomes
from~$n$ independent and identically distributed (iid) Bernoulli
random variables with parameter $q$, i.e., when each bit in the
$n$-string is independently drawn with probability $0<q<1$ of getting
a $1$.

Since probabilistic results are not our main goal in this paper, we
mainly examine how to extend our most relevant enumerative results in
Sections~\ref{sec:w_n_m_klek} and~\ref{sec:wp_n_m} to the
probabilistic scenario described above. Observe that this also
includes the probabilistic extensions of our results in
Sections~\ref{sec:w_n_m_k}--\ref{sec:number-wl_n_k_gtk}, which are,
essentially, special cases or consequences
of~Section~\ref{sec:w_n_m_klek}.  In fact, we also look into the
probabilistic extension of the enumerative results in
Sections~\ref{sec:w_n_m_k},~\ref{sec:w_n_m_gtk}
and~\ref{sec:number-ws_n_k_ltk}, since these are connected with
several prior findings in the literature.

In order not to overload notation, it is understood that all
recurrences, probability generating functions, explicit expressions,
and moments in this section
implicitly depend on~$q$.

\subsection[Probability that an n-String Contains Exactly m
(ḵ≤k)-Runs]{Probability that an $n$-String Contains Exactly $m$ ($\ushort{k}\le\widebar{k}$)-Runs}
\label{sec:probability-w_m_klek}
We call $\pi_{\ushort{k}\le\widebar{k}}(n,m)$ the probability that an
$n$-string contains exactly $m$ ($\ushort{k}\le\widebar{k}$)-runs.
The obvious case is $q=1/2$, in which
$\pi_{\ushort{k}\le\widebar{k}}(n,m)=w_{\ushort{k}\le\widebar{k}}(n,m)/2^n$.
Likewise, the probability that an $n$-string contains at least $m$
$(\ushort{k}\le\widebar{k})$-runs is
$\pi_{\ushort{k}\le\widebar{k}}^{\ge
  m}(n)=w_{\ushort{k}\le\widebar{k}}^{\ge m}(n)/2^n$ in this case. We
examine next how to get these quantities for arbitrary $q$.

\subsubsection{Recurrences}
\label{sec:recurrences}
We may directly write two equivalent probability recurrences for
$\pi_{\ushort{k}\le\widebar{k}}(n,m)$ by using the two enumerative
recurrences in Section~\ref{sec:recurrences_w_n_m_klek} in conjunction
with the law of total
probabilities. From~\eqref{eq:recurrence_w_n_m_klek_basic} we have the
recurrence
\begin{align}
  \label{eq:recurrence_p_n_m_klek_basic}
  \pi_{\ushort{k}\le\widebar{k}}(n,m)= (1-q)\bigg(&\sum_{i=\ushort{k}}^{\widebar{k}}
    q^i\pi_{\ushort{k}\le\widebar{k}}\big(n-(i+1),m-1\big)+\sum_{i=0}^{\ushort{k}-1}
    q^i\pi_{\ushort{k}\le\widebar{k}}\big(n-(i+1),m\big)\nonumber\\  
  &+\sum_{i=\widebar{k}+1}^{n}
   q^i \pi_{\ushort{k}\le\widebar{k}}\big(n-(i+1),m\big)\bigg).  
\end{align}
whereas from~\eqref{eq:recurrence_w_n_m_klek_simpler} we have the
alternative recurrence
\begin{align}
  \label{eq:recurrence_p_n_m_klek_simpler}
  \pi_{\ushort{k}\le\widebar{k}}(n,m)=\pi_{\ushort{k}\le\widebar{k}}(n-1,m)
             &+ q^{\ushort{k}}(1-q)\Big(\pi_{\ushort{k}\le\widebar{k}}\big(n-(\ushort{k}+1),m-1\big)-\pi_{\ushort{k}\le\widebar{k}}\big(n-(\ushort{k}+1),m\big)\Big)\nonumber\\
            & - q^{\widebar{k}+1}(1-q)\Big(\pi_{\ushort{k}\le\widebar{k}}\big(n-(\widebar{k}+2),m-1\big)-\pi_{\ushort{k}\le\widebar{k}}\big(n-(\widebar{k}+2),m\big)\Big).
\end{align}
Notice that the $2$ factor in the first term
of~\eqref{eq:recurrence_w_n_m_klek_simpler} becomes $1$ in the
equivalent term of~\eqref{eq:recurrence_p_n_m_klek_simpler}. This is
because our probabilistic estimate for the value of
$\pi_{\ushort{k}\le\widebar{k}}(n,m)$ from
$\pi_{\ushort{k}\le\widebar{k}}(n-1,m)$ is
$q\,\pi_{\ushort{k}\le\widebar{k}}(n-1,m)+(1-q)\,\pi_{\ushort{k}\le\widebar{k}}(n-1,m)$. Of
course, we can also get~\eqref{eq:recurrence_p_n_m_klek_simpler}
from~\eqref{eq:recurrence_p_n_m_klek_basic}, using
$\pi_{\ushort{k}\le\widebar{k}}(n,m)-q\,\pi_{\ushort{k}\le\widebar{k}}(n-1,m)$.
In the following we work with the simpler
recurrence~\eqref{eq:recurrence_p_n_m_klek_simpler}, as it is also
directly amenable to obtaining the probability generating function (pgf)
of $\pi_{\ushort{k}\le\widebar{k}}(n,m)$. As we have done before, to
find initial values we use the probability for $n=1$, which, by
inspection, is
\begin{align}
  \label{eq:trivial_n1_p}
  \pi_{\ushort{k}\le\widebar{k}}(1,m)=&~\big(q\,\iv{\ushort{k}=0}\iv{\widebar{k}=0}+(1-q)\iv{\ushort{k}=1}+\iv{\ushort{k}>1}\big)\iv{m=0}\nonumber\\
 &+q\,\big(\iv{\ushort{k}=0}\iv{\widebar{k}\ge 1}+\iv{\ushort{k}=1}\big)\iv{m=1}
  +(1-q)\,\iv{\ushort{k}=0}\iv{m=2}.
\end{align}
Now, setting $n=1$ in~\eqref{eq:recurrence_p_n_m_klek_simpler} yields
\begin{align}\label{eq:init_helper_p_n_m_klek}
  \pi_{\ushort{k}\le\widebar{k}}(1,m)=\pi_{\ushort{k}\le\widebar{k}}(0,m)&+q^{\ushort{k}}(1-q)\Big(\pi_{\ushort{k}\le\widebar{k}}(-\ushort{k},m-1)-\pi_{\ushort{k}\le\widebar{k}}(-\ushort{k},m)\Big)\nonumber\\
            & - q^{\widebar{k}+1}(1-q)\Big(\pi_{\ushort{k}\le\widebar{k}}(-(\widebar{k}+1),m-1)-\pi_{\ushort{k}\le\widebar{k}}(-(\widebar{k}+1),m)\Big),
\end{align}
We want~\eqref{eq:init_helper_p_n_m_klek} to
equal~\eqref{eq:trivial_n1_p}. Taking necessary condition~\eqref{eq:m_necessary_condition} into
account, we can achieve this equality for all $m$ and
$0\le\ushort{k}\le\widebar{k}$ by letting
\begin{align}
  \pi_{\ushort{k}\le\widebar{k}}(-1,m)&=\frac{1}{1-q}\iv{m=0},\label{eq:weird}\\
  \pi_{\ushort{k}\le\widebar{k}}(0,m)&=\ivl{m=\iv{\ushort{k}=0}},\label{eq:not_weird}
\end{align}
which therefore constitute the initialisation
of~\eqref{eq:recurrence_p_n_m_klek_simpler}.  These values
initialise~\eqref{eq:recurrence_p_n_m_klek_basic} as
well.

\begin{remark}\label{rmk:weird}
  As if a binary string of length $-1$ were not weird enough, take a
  moment to ponder that $1/(1-q)$ is greater than one in the
  initialisation~\eqref{eq:weird} of the probability recurrence just
  given.
\end{remark}

As for $\pi_{\ushort{k}\le\widebar{k}}^{\ge m}(n)$, we have
from~\eqref{eq:wn_klek_atleast} that
$\pi_{\ushort{k}\le\widebar{k}}^{\ge
  m}(n)=\sum_{t=m}^{\floor*{(n+1)/(\ushort{k}+1)}}\pi_{\ushort{k}\le\widebar{k}}(n,t)$. Applying
this expression to~\eqref{eq:recurrence_p_n_m_klek_simpler} we get the
recurrence
\begin{equation}
  \label{eq:recurrence_p_n_atleast_klek}
  \pi^{\ge m}_{\ushort{k}\le\widebar{k}}(n)=\pi^{\ge m}_{\ushort{k}\le\widebar{k}}(n-1)
  +
  q^{\ushort{k}}(1-q)\,\pi_{\ushort{k}\le\widebar{k}}(n-(\ushort{k}+1),m-1)-q^{\widebar{k}+1}(1-q)\,\pi_{\ushort{k}\le\widebar{k}}(n-(\widebar{k}+2),m-1),
\end{equation}
which, as we see in the next section, suffices to get the pgf of $\pi^{\ge m}_{\ushort{k}\le\widebar{k}}(n)$.

\subsubsection{Probability Generating Functions}
\label{sec:pgf_pi_n_m_klek}
We can now obtain the bivariate pgf
$\Pi_{\ushort{k}\le\widebar{k}}(x,y)=\sum_{n,m}
\pi_{\ushort{k}\le\widebar{k}}(n,m) \,x^ny^m$ following the same steps as
in Section~\ref{sec:generating-function_klek}. In this case we need to
add $\big(\iv{n=-1}-q\iv{n=0}\big)\iv{m=0}/(1-q)$
to~\eqref{eq:recurrence_p_n_m_klek_simpler} to render it valid for all
$n$ and $m$. Using this extended recurrence, we directly get
\begin{equation}\label{eq:pgf_klek}
  \Pi_{\ushort{k}\le\widebar{k}}(x,y)=\frac{1-q\, x}{(1-q)\,x\,\Big(1-x-(y-1)(1-q)(q^{\ushort{k}}x^{\ushort{k}+1}-q^{\widebar{k}+1}x^{\widebar{k}+2})\Big)}.
\end{equation}
We can extract the coefficient of $y^m$ in this pgf following the same
steps as in Section~\ref{sec:generating-function_klek}, which in this
case lead us to obtain
\begin{equation}\label{eq:pgf_ym}
  [y^m]\Pi_{\ushort{k}\le\widebar{k}}(x,y)=\frac{(1-q)^{m-1}(1-q\,x)\big(q^{\ushort{k}}x^{\ushort{k}+1}-q^{\widebar{k}+1}x^{\widebar{k}+2}\big)^m}{x\Big(1-x+(1-q)\big(q^{\ushort{k}}x^{\ushort{k}+1}-q^{\widebar{k}+1}x^{\widebar{k}+2}\big)\Big)^{m+1}}.
\end{equation}
This is the pgf giving the probability that a binary string contains
exactly $m$ $(\ushort{k}\le\widebar{k})$-runs.

The pgf
$\Pi^{\ge m}_{\ushort{k}\le\widebar{k}}(x)=\sum_n \pi^{\ge
  m}_{\ushort{k}\le\widebar{k}}(n)\, x^n$ is obtained by multiplying
\eqref{eq:recurrence_p_n_atleast_klek} on both sides by $x^n$ and then
adding over~$n$, which yields
\begin{align}
  \label{eq:pgf_p_n_atleast_klek}
  \Pi^{\ge m}_{\ushort{k}\le\widebar{k}}(x)&=\frac{1-q}{1-x} \big(q^{\ushort{k}}x^{\ushort{k}+1}-q^{\widebar{k}+1}x^{\widebar{k}+2}\big) [y^{m-1}]\Pi_{\ushort{k}\le\widebar{k}}(x,y)\nonumber\\
&=\frac{(1-q\,x)}{(1-q)\,x\,(1-x)}\Bigg(\frac{(1-q)\big(q^{\ushort{k}}x^{\ushort{k}+1}-q^{\widebar{k}+1}x^{\widebar{k}+2}\big)}{1-x+(1-q)\big(q^{\ushort{k}}x^{\ushort{k}+1}-q^{\widebar{k}+1}x^{\widebar{k}+2}\big)}\Bigg)^m,
\end{align}
where we have
used~\eqref{eq:pgf_ym}. This is
the pgf giving the probability that a binary string contains at least $m$
$(\ushort{k}\le\widebar{k})$-runs.

\begin{remark}\label{rmk:more_recs_and_pgfs}
  As a rule, pgf~\eqref{eq:pgf_ym} allows us to pursue explicit
  expressions for $\pi_{\ushort{k}\le\widebar{k}}(n,m)$ in particular
  cases of $\ushort{k}$ and $\widebar{k}$, in the same way that we
  have sought explicit expressions for
  $w_{\ushort{k}\le\widebar{k}}(n,m)$ in particular cases of
  $\ushort{k}$ and $\widebar{k}$ in
  Sections~\ref{sec:w_n_m_k}--\ref{sec:anylength}. If we are
  interested instead in the computation of moments for a given~$n$, we
  can do so through~\eqref{eq:pgf_klek} but sometimes the alternative
  route through $[x^n]\Pi_{\ushort{k}\le\widebar{k}}(x,y)$ is cleaner.
  Although getting the coefficient of~$x^n$ is not as straightforward
  as getting the coefficient of~$y^m$ (it cannot usually be done in
  closed form), in special cases one can still get manageable
  expressions for~$[x^n]\Pi_{\ushort{k}\le\widebar{k}}(x,y)$ through
  the binomial theorem. The extra work involved sometimes actually
  pays off, as neater explicit expressions can sometimes be determined
  through $[x^n]\Pi_{\ushort{k}\le\widebar{k}}(x,y)$.
\end{remark}

\subsection[Probability that an n-String Contains Exactly m
k-Runs]{Probability that an $n$-String Contains Exactly $m$
  $k$-Runs}\label{sec:pi_n_m_k}
We denote by $\pi_{k}(n,m)$ the probability that an $n$-string
contains exactly $m$  $k$-runs. This is a special case of the
probability in Section~\ref{sec:probability-w_m_klek} with
$\ushort{k}=\widebar{k}=k$, and thus 
\begin{equation}
  \label{eq:pi_n_m_k_from_klek}
  \pi_{k}(n,m)=\pi_{k\le k}(n,m).
\end{equation}
We discuss this special case in some detail because Makri and
Psillakis~\cite{makri11:_bernoul} previously studied this
probability. Marbe~\cite{marbe16:_mathematische,marbe34:_grundfragen}
and Cochran~\cite{cochran36:_plants} also gave an expectation
connected to it.

\subsubsection{Recurrences}
\label{sec:pi_n_m_k_rec}
Using~\eqref{eq:pi_n_m_k_from_klek},
recurrence~\eqref{eq:recurrence_p_n_m_klek_basic} becomes 
\begin{equation}
  \label{eq:recurrence_p_n_m_k_basic}
  \pi_k(n,m)= (1-q)\bigg(
    q^k\pi_k\big(n-(k+1),m-1\big)+\sum_{i=0\atop i\neq k}^{n}
    q^i\pi_k\big(n-(i+1),m\big)\bigg)
\end{equation}
whereas from~\eqref{eq:recurrence_p_n_m_klek_simpler} we have the
alternative recurrence
\begin{align}
  \label{eq:recurrence_p_n_m_k_simpler}
  \pi_k(n,m)=\pi_k(n-1,m)
             &+ q^{k}(1-q)\Big(\pi_k\big(n-(k+1),m-1\big)-\pi_k\big(n-(k+1),m\big)\Big)\nonumber\\
            & - q^{k+1}(1-q)\Big(\pi_k\big(n-(k+2),m-1\big)-\pi_k\big(n-(k+2),m\big)\Big),
\end{align}
which from~\eqref{eq:weird} and~\eqref{eq:not_weird} are both
initialised by
\begin{align*}
  \pi_{k}(-1,m)&=\frac{1}{1-q}\iv{m=0},\\
  \pi_{k}(0,m)&=\ivl{m=\iv{k=0}}.
\end{align*}

\subsubsection{Probability Generating Functions}
\label{sec:pgf_k}
By setting $\ushort{k}=\widebar{k}=k$
in~\eqref{eq:pgf_klek} we get the pgf $\Pi_k(x,y)=\sum_{n,m}\pi_k(n,m)\,x^ny^m$:
\begin{equation}\label{eq:pgf_k}
  \Pi_k(x,y)=\frac{1-q\, x}{(1-q)\,x\,\Big(1-x-(y-1)(1-q)(q^{k}x^{k+1}-q^{k+1}x^{k+2})\Big)}.
\end{equation}
Similarly, from ~\eqref{eq:pgf_ym} we have
\begin{equation}\label{eq:pgf_k_ym}
  [y^m]\Pi_k(x,y)=\frac{(1-q)^{m-1}(1-q\,x)\big(q^{k}x^{k+1}-q^{k+1}x^{k+2}\big)^m}{x\Big(1-x+(1-q)\big(q^{k}x^{k+1}-q^{k+1}x^{k+2}\big)\Big)^{m+1}},
\end{equation}
which is the pgf giving the probability that a binary string contains
exactly $m$ $k$-runs, whereas from~\eqref{eq:pgf_p_n_atleast_klek} we get
\begin{align}
  \label{eq:pgf_p_n_atleast_k}
  \Pi^{\ge m}_k(x)&=\frac{(1-q\,x)}{(1-q)\,x\,(1-x)}\Bigg(\frac{(1-q)\big(q^{k}x^{k+1}-q^{k+1}x^{k+2}\big)}{1-x+(1-q)\big(q^{k}x^{k+1}-q^{k+1}x^{k+2}\big)}\Bigg)^m.
\end{align}
This is the pgf giving the probability that a binary string contains
at least $m$ $k$-runs.

\subsubsection{Explicit Expression}
\label{sec:explicit-expressions}
In this section we obtain an explicit expression for $\pi_k(n,m)$ by
extracting the coefficient of $x^n$ in~\eqref{eq:pgf_k_ym}. Applying
the negative binomial theorem and the binomial theorem twice, we can
see that~\eqref{eq:pgf_k_ym} can be expanded as
\begin{equation}
  \label{eq:pgf_k_ym_expanded}
  [y^m]\Pi_k(x,y)=\sum_{s\ge 0}\sum_{t\ge 0}\sum_{r\ge 0} {s+m\choose
    m}{s\choose t}{m+t+1\choose r}(-1)^{t+r} (1-q)^{t+m-1}
  q^{k(t+m)+r} x^{(k+1)m-1+s+kt+r}.
\end{equation}
The coefficient of $x^n$ above corresponds to the nonnegative
solutions of indices $s, t$ and $r$ in $(k+1)m-1+s+kt+r=n$, or,
equivalently, in $s+kt+r=e$ using~\eqref{eq:excess} with
$\ushort{k}=k$. For $k>1$, the maximum of $t$ happens for $s=r=0$, and
thus $t\le \floor{e/k}$. Similarly, $r=0$ let us see that $s\le e-kt$,
whereas $r$ is determined by the previous equation for any two values
of $t$ and $s$. If $k=0$ then $t\le e$, as the second binomial is zero
otherwise. Collecting these solutions, we have
from~\eqref{eq:pgf_k_ym_expanded} that
\begin{equation}
  \label{eq:pi_k_n_m_th}
  \pi_k(n,m)=\sum_{t=0}^{\min(\floor{\frac{e}{k}},e)}\sum_{s=0}^{e-kt}{s+m\choose
    m}{s\choose t}{m+t+1\choose e-kt-s}(-1)^{e-(k-1)t-s} (1-q)^{t+m-1}
  q^{k m+e-s}.
\end{equation}
This is equivalent to the explicit expression for $\pi_k(n,m)$ given by
Makri and Psillakis~\cite[Thm.\ 2.1]{makri11:_bernoul} ---in their
notation $P(E_{n,k}=m)=\pi_k(n,m)$.

\subsubsection{Moments}
\label{sec:moments_w_k}
The pgf~\eqref{eq:pgf_k} also allows us to determine the factorial
moments of the random variable (rv) modelling the number of $k$-runs in
an $n$-string drawn at random, which we denote by~$M_{k,n}$. With standard
probabilistic notation we have that
\begin{equation}
  \label{eq:pik_pr}
  \Pr\big(M_{k,n}=m\big)=\pi_k(n,m).
\end{equation}
The first factorial moment (expectation) can be obtained as
$\text{E}(M_{k,n})=[x^n](\partial/\partial y)\Pi_k(x,y)|_{y=1}$.
Thus, we have determine next the coefficient of $x^n$ in
\begin{equation}
  \label{eq:ddy_pik_y1}
  \left.\frac{\partial \Pi_k(x,y)}{\partial y}\right|_{y=1}=\frac{(1-qx)^2q^kx^k}{(1-x)^2}.
\end{equation}
By applying the negative binomial theorem and the binomial theorem we
have that
\begin{equation}
  \left.\frac{\partial \Pi_k(x,y)}{\partial y}\right|_{y=1}=\sum_{t\ge 0}\sum_{s\ge 0}
  {t+1\choose t} {2\choose s} (-1)^sq^{s+k} x^{k+s+t}.
\end{equation}
From $n=k+s+t$ we have that $t\le n-k$ and $s=n-k-t$, and therefore
\begin{equation}
  \label{eq:exp_mkn_prev}
  \text{E}(M_{k,n})=\sum_{t= 0}^{n-k}
  (t+1) {2\choose n-k-t} (-1)^{n-k-t}q^{n-t}.
\end{equation}
There are three cases in this expression, depending on $n-k=0,1$ or
$n-k\ge 2$, but the last two cases have
the same solution. Evaluating them one sees that
\begin{equation}
  \label{eq:exp_mkn_all}
  \text{E}(M_{k,n})=\Big(\big((n-k-1)(1-q)+2\big)(1-q)\,\iv{k<n}+\iv{k=n}\Big)\,q^k.
\end{equation}
The second factorial moment  $\text{E}(M_{k,n}(M_{k,n}-1))$ is the coefficient of $x^n$ in
\begin{equation}
  \label{eq:d2dy2_pik_y1}
  \left.\frac{\partial^2 \Pi_k(x,y)}{\partial y^2}\right|_{y=1}=\frac{2(1-q)(1-qx)^3q^{2k}x^{2k+1}}{(1-x)^3}.
\end{equation}
As before, this can be expanded as
\begin{equation}
  \label{eq:d2dy2_pik_y1}
  \left.\frac{\partial^2 \Pi_k(x,y)}{\partial y^2}\right|_{y=1}=2(1-q)\sum_{t\ge
    0}\sum_{s\ge 0}{t+2\choose t}{3\choose s} (-1)^s q^{2k+s} x^{2k+1+t+s},
\end{equation}
and we now need to solve $2k+1+t+s=n$. Observing that $t\le n-2k-1$
and $s=n-2k-1-t$, we have that 
\begin{equation}
  \label{eq:2ndfact_mkn_prev}
  \text{E}\big(M_{k,n}(M_{k,n}-1)\big)=(1-q)\sum_{t=0}^{n-2k-1}(t+2)(t+1){3\choose
    n-2k-1-t}(-1)^{n-2k-1-t} q^{n-t-1}.
\end{equation}
There are four nonzero cases in~\eqref{eq:2ndfact_mkn_prev}, which,
letting $\xi=n-2k-1$, correspond to $\xi=0, 1, 2$ and $\xi\ge 3$.  The
last three cases have the same solution, and after some algebra we can see that
\begin{align}
  \label{eq:2ndfact_mkn_cases}
  \text{E}\big(M_{k,n}(M_{k,n}-1)\big)=(1 -q)
  q^{n-1}\Big(&(1-q)q^{-\xi}\big(2\, (1 + q + q^2) + 3\,
                (1-q^2) \xi+ (1 - q)^2
                                         \xi^2\big)\iv{\xi> 0}\nonumber\\
  &+2\,\iv{\xi=0}\Big)
\end{align}
With~\eqref{eq:exp_mkn_all} and~\eqref{eq:2ndfact_mkn_cases} we can
obtain the variance of $M_{k,n}$ using
$\text{Var}(M_{k,n})=\text{E}(M_{k,n}(M_{k,n}-1))+\text{E}(M_{k,n})-\text{E}^2(M_{k,n})$. For
example, a trivial case is $\text{E}(M_{n,n})=q^n$ and
$\text{Var}(M_{n,n})=q^n(1-q^n)$, which can also be seen by observing
that $M_{n,n}$ is a Bernoulli rv with parameter~$q^n$.

The expectation~\eqref{eq:exp_mkn_all} was originally given by
Marbe~\cite{marbe16:_mathematische,marbe34:_grundfragen} and then put
on a sounder footing by Cochran~\cite[Eq.\ (5)]{cochran36:_plants}
---see more details in Section~\ref{sec:number-runs},
Remark~\ref{rmk:rk-r-alt}. The expectation and variance of $M_{0,n}$
were implicitly given by Bloom~\cite{bloom96:_probab} ---see discussion
towards the end of Remark~\ref{rmk:isolated_singles} in
Section~\ref{sec:oz}.

\subsection[Probability that an n-String Contains Exactly m
(≥k)-Runs]{Probability that an $n$-String Contains Exactly $m$
  ($\ge k$)-Runs}\label{sec:pi_n_m_gtk}

We denote by $\pi_{\ge k}(n,m)$ the probability that an $n$-string
contains exactly $m$  ($\ge k$)-runs. This is a special case of the
probability in Section~\ref{sec:probability-w_m_klek} with
 $\ushort{k}=k$ and $\widebar{k}=n$, and thus 
 \begin{equation}
   \label{eq:pi_n_m_gtk_from_klek}
   \pi_{\ge k}(n,m)=\pi_{k\le n}(n,m).
 \end{equation}
 We discuss this special case in some detail because several authors
 have previously studied this
 probability~\cite{moivre38:_doctrine,laplace20:_theorie,uspensky32:_probl_runs,
   muselli96:_simple,balakrishnan01:_runs_scans}.

\subsubsection{Recurrences}
\label{sec:pi_n_m_gtk_rec}
Using~\eqref{eq:pi_n_m_gtk_from_klek},
recurrence~\eqref{eq:recurrence_p_n_m_klek_basic} becomes 
  \begin{align}
    \label{eq:recurrence_p_n_m_gtk_basic}
    \pi_{\ge k}(n,m)= (1-q)\bigg(&%
                                   \sum_{i=k}^{n}    q^i\pi_{\ge k}(n-(i+1),m-1)
                                   +\sum_{i=0}^{k-1} q^i\pi_{\ge k}(n-(i+1),m)\bigg),
  \end{align}
whereas the second recurrence~\eqref{eq:recurrence_p_n_m_klek_simpler}
is now
\begin{align}
  \label{eq:recurrence_p_n_m_gtk_simpler}
  \pi_{\ge k}(n,m)=~&\pi_{\ge k}(n-1,m)%
                      + q^{k}(1-q)\Big(\pi_{\ge k}(n-(k+1),m-1)-\pi_{\ge k}(n-(k+1),m)\Big).
\end{align}
From~\eqref{eq:weird} and~\eqref{eq:not_weird} both recurrences are
initialised by
\begin{align}
  \label{eq:p_n_m_gtk_init}
  \pi_{\ge k}(-1,m)&=\frac{1}{1-q}\iv{m=0},\\
  \pi_{\ge k}(0,m)&=\ivl{m=\iv{k=0}}.  \label{eq:p_n_m_gtk_init2}
\end{align}
Lastly,~\eqref{eq:recurrence_p_n_atleast_klek} becomes
\begin{equation}
  \label{eq:recurrence_p_n_m_gtk_demoivre}
  \pi^{\ge m}_{\ge k}(n)=\pi^{\ge m}_{\ge k}(n-1) + q^{k}(1-q)\,\pi_{\ge k}(n-(k+1),m-1).
\end{equation}

\subsubsection{Probability Generating Functions}
\label{sec:pgf_gtk}
Although we can get valid pgfs by
specialising~\eqref{eq:pgf_klek},~\eqref{eq:pgf_ym}
and~\eqref{eq:pgf_p_n_atleast_klek} with $\ushort{k}=k$
and~$\widebar{k}=n$, the resulting expressions depend on $n$ and
require $k\le n$. We obtain next simpler pgfs valid for all $k\ge 0$
by working directly with~\eqref{eq:recurrence_p_n_m_gtk_simpler}
and~\eqref{eq:recurrence_p_n_m_gtk_demoivre}. The pgf
$\Pi_{\ge k}(x,y)=\sum_n \pi_{\ge k}(n,m)\, x^ny^m$ can be obtained by
making recurrence~\eqref{eq:recurrence_p_n_m_gtk_simpler} valid for
all $n$ and $m$, multiplying it on both sides by $x^ny^m$, and then
adding over $n$ and $m$. This yields
\begin{equation}
  \label{eq:pgf_p_n_m_gtk_bivariate}
  \Pi_{\ge k}(x,y)=\frac{1-q\, x}{(1-q)\,x\,\big(1-x-(y-1)(1-q)q^kx^{k+1}\big)}.
\end{equation}
Through the same method as before, we can extract the coefficient of
$y^m$ from from~\eqref{eq:pgf_p_n_m_gtk_bivariate}:
\begin{equation}
  \label{eq:pgf_p_n_m_gtk_m}
  [y^m]\Pi_{\ge k}(x,y)=\frac{(1-q\, x)\big((1-q)q^kx^{k+1}\big)^m}{(1-q)\,x\,\big(1-x+(1-q)q^kx^{k+1}\big)^{m+1}},
\end{equation}
which is the pgf giving the probability that a binary string contains
exactly $m$ $(\ge\! k)$-runs.

Also, multiplying \eqref{eq:recurrence_p_n_m_gtk_demoivre} on both
sides by $x^n$ and then adding over~$n$ yields
\begin{align}
  \label{eq:pgf_demoivre_a}
  \Pi^{\ge m}_{\ge k}(x)&=\frac{1-q}{1-x}\, q^k x^{k+1} [y^{m-1}]\Pi_{\ge k}(x,y)\\
&=\frac{(1-q\,x)}{(1-q)\,x\,(1-x)}\bigg(\frac{(1-q)q^kx^{k+1}}{1-x+(1-q)q^kx^{k+1}}\bigg)^m  \label{eq:pgf_demoivre}
\end{align}
where in the second step above we have
used~\eqref{eq:pgf_p_n_m_gtk_m}.  Expression~\eqref{eq:pgf_demoivre}
is the pgf giving the probability that a binary string contains at
least $m$ $(\ge k)$-runs. 

We finally determine the coefficient of $x^n$
in~\eqref{eq:pgf_p_n_m_gtk_bivariate}, as in this case this allows us
to get cleaner expressions later. By applying the negative binomial
theorem and then the binomial theorem, we can
express~\eqref{eq:pgf_p_n_m_gtk_bivariate} as
\begin{equation}
  \label{eq:pgf_p_n_m_gtk_binomial}
  \Pi_{\ge k}(x,y)=(x^{-1}-q)\sum_{s\ge 0}\sum_{t= 0}^s{s
    \choose t}(1-q)^{t-1} q^{kt}(y-1)^t x^{s+kt}.
\end{equation}
So to get $[x^n]\Pi_{\ge k}(x,y)$ we just need to solve $s+kt-1=n$ and
$s+kt=n$ for nonnegative indices $s$ and $t$. If $k>1$, in both cases
the maximum of $t$ happens when $s=0$, and corresponds to
$\floor{(n+1)/k}$ and $\floor{n/k}$, respectively. Thus,
from~\eqref{eq:pgf_p_n_m_gtk_binomial} we can write
\begin{equation}
  \label{eq:pgf_p_n_m_gtk_xn}
  [x^n]\Pi_{\ge k}(x,y)=\sum_{t= 0}^{\min\big(\floor*{\frac{n+1}{k}},n+1\big)}\Bigg({n-kt+1
    \choose t}-q{n-kt\choose t}\Bigg)(1-q)^{t-1}q^{kt}(y-1)^t.
\end{equation}
Observe that the summation only goes up to $n+1$ when $k=0$, because in this
case of the first binomial coefficient becomes zero for higher values
of $t$.

\subsubsection{Explicit Expression}
We now produce an explicit expression for $\pi_{\ge k}(n,m)$ by extracting the coefficient of $y^m$
from~\eqref{eq:pgf_p_n_m_gtk_xn}. By applying the binomial theorem to
$(y-1)^t$ we immediately see that
$\pi_{\ge k}(n,m)=[x^ny^m]\Pi_{\ge k}(x,y)$ is
\begin{equation}
  \label{eq:pi_gtk_n_m_th}
  \pi_{\ge k}(n,m)=\sum_{t=
    0}^{\min\big(\floor*{\frac{n+1}{k}},n+1\big)}(-1)^{t-m}{t\choose m}\Bigg({n-kt+1
    \choose t}-q{n-kt\choose t}\Bigg)(1-q)^{t-1}q^{kt}.
\end{equation}
which is very close to the simplest explicit expression for
$\pi_{\ge k}(n,m)$ given by Muselli~\cite[Thm.\ 3]{muselli96:_simple}
---in his notation, $P(M_n^{(k)}=m)=\pi_{\ge k}(n,m)$.

\subsubsection{Moments}
\label{sec:moments_w_gtk}
The pgf~\eqref{eq:pgf_p_n_m_gtk_xn} can also be used to obtain moments
of the rv modelling the number of $(\ge k)$-runs in
an $n$-string drawn at random, which we denote by $M_{\ge k,n}$. With
standard probabilistic notation,
\begin{equation}
  \label{eq:pgtk_pr}
  \Pr\big(M_{\ge k,n}=m\big)=\pi_{\ge k}(n,m).
\end{equation}
The expectation of this random
variable is
$\text{E}(M_{\ge k,n})=(d/dy)[x^n]\Pi_{\ge
  k}(x,y)|_{y=1}$. After differentiating~\eqref{eq:pgf_p_n_m_gtk_xn}
with respect to $y$, setting $y=1$ lets us see that the only nonzero
term corresponds to $t=1$, as long as $k\le n$. Therefore,
\begin{equation}
  \label{eq:p_n_gtk_expectation}
  \text{E}(M_{\ge k,n})= \big((n-k)(1-q)+1\big)\,q^k\,\iv{k\le n}.
\end{equation}
The second factorial moment is
$\text{E}(M_{\ge k,n}(M_{\ge k,n}-1))=(d^2/dy^2)[x^n]\Pi_{\ge
  k}(x,y)|_{y=1}$. After differentiating~\eqref{eq:pgf_p_n_m_gtk_xn}
twice with respect to $y$ and then setting $y=1$, we see that the only
nonzero term corresponds to $t=2$. Thus, the second factorial moment
is zero if $k>0$ and $\floor{(n+1)/k}=1$. Therefore we have that %
\begin{equation}
  \label{eq:p_n_gtk_second_fact_moment}
  \text{E}\big(M_{\ge k,n}(M_{\ge    k,n}-1)\big)=%
  \big((n-2k)(1-q)+1+q\big)(n-2k)(1-q)q^{2k}\,\big(\iv{k>0}\iv{k<n}+\iv{k=0}\big).
\end{equation}
The variance of $M_{\ge k,n}$ follows from the two moments above
---see Section~\ref{sec:explicit-expressions}.  The
expectation~\eqref{eq:p_n_gtk_expectation} was apparently first given
by Goldstein~\cite[Eq.\ (5)]{goldstein90:_poisson}, without using a
pgf. Aki and Hirano also obtained~\eqref{eq:p_n_gtk_expectation} and
its corresponding variance~\cite[pp.\ 317--318]{aki93:_markov}
through a pgf similar to~\eqref{eq:pgf_p_n_m_gtk_xn}~\cite[Eq.\
(5)]{aki93:_markov}.  

\begin{remark}\label{rmk:balakrishnan}
  Krishnan Nair~\cite[p.\ 84]{nair42:_prob} gave a pgf which should be
  equivalent to~\eqref{eq:pgf_p_n_m_gtk_m}, but which appears to be
  incorrect.  Also, Balakrishnan and Koutras gave a recurrence for
  $\pi_{\ge k}(n,m)$ very similar
  to~\eqref{eq:recurrence_p_n_m_gtk_basic} \cite[p.\
  144]{balakrishnan01:_runs_scans} but with different
  initialisation. These authors also gave a bivariate pgf for this
  case~\cite[Eq.\ (5.11)]{balakrishnan01:_runs_scans}, which, however,
  they did not get from their aforementioned recurrence.  In our
  notation, the pgf given by Balakrishnan and Koutras is
  \begin{equation}
    \label{eq:pgf_p_n_m_gtk_bivariate_balakrishnan}
    \Pi_{\ge k}(x,y)=\frac{1+(y-1)(q\, x)^k}{1-x-(y-1)(1-q)q^kx^{k+1}}.
  \end{equation}
  It is easier to obtain $[y^m]\Pi_{\ge k}(x,y)$ in closed-form
  from~\eqref{eq:pgf_p_n_m_gtk_bivariate} than
  from~\eqref{eq:pgf_p_n_m_gtk_bivariate_balakrishnan}
  ---\eqref{eq:pgf_p_n_m_gtk_m} is witness to this.  A consequence is
  that is simpler obtaining~\eqref{eq:pgf_demoivre}
  using~\eqref{eq:pgf_p_n_m_gtk_bivariate} than
  using~\eqref{eq:pgf_p_n_m_gtk_bivariate_balakrishnan}.
\end{remark}

\begin{remark}\label{rmk:demoivre_laplace}
  The results in this section also allow us to recover the oldest
  historical result on success runs by de
  Moivre~\cite[Prob. LXXXVIII,\ pp.\ 243--248]{moivre38:_doctrine}:
  ``\textit{To find the Probability of throwing a Chance assigned a
    given number of times without intermission, in any given number of
    Trials},'' or, in terms of our definitions, to find
  $\pi^{\ge 1}_{\ge k}(n)=1-\pi_{\ge k}(n,0)$. Using the same
  numerical values originally employed by de Moivre
  and~\eqref{eq:pi_gtk_n_m_th}, for $n=10$, $k=3$ and $q=1/2$, we may
  see that $\pi^{\ge 1}_{\ge 3}(10)=0.5078=65/128$, which we can also
  verify using $w^{\ge 1}_{\ge 3}(10)/2^{10}$
  ---see~\eqref{eq:wn_gtk_atleast}---, whereas if $q=2/3$
  then~$\pi^{\ge 1}_{\ge 3}(10)=0.8121=592/729$. The truncated decimal
  values are computed using~\eqref{eq:pi_gtk_n_m_th}, whereas the
  fractions can be obtained, for instance, by determining
  $[x^{10}]\pi^{\ge 1}_{\ge 3}(x)$ using~\eqref{eq:pgf_demoivre_a} and
  a computer algebra system. Both probabilities are correctly given
  in~\cite{moivre38:_doctrine}, even though there is a known mistake
  in the procedure followed therein~\cite[p.\
  418]{hald03:_history}. De Moivre did not give a recurrence nor did
  he provide a proof for his solution, which was based on adding a
  given number of terms in the power series expansion of a generating
  function. According to Hald's account~\cite[p.\
  420]{hald03:_history}, the first solution to de Moivre's problem
  based on a finite difference equation (i.e., a recurrence) was given
  by Simpson. Simpson also gave an explicit solution as an infinite
  series. However, as also indicated by Hald, Laplace was the first
  author who essentially derived pgf~\eqref{eq:pgf_demoivre} for
  $m=1$~\cite[p.\ 252]{laplace20:_theorie} ---see
  Remark~\ref{rmk:wait_prob} below.  Finally, Uspensky gave a
  semi-explicit closed-form expression for the probability
  $\pi_{\ge k}(n,0)=1-\pi_{\ge k}^{\ge 1}(n)$~\cite[Eq.\
  (3)]{uspensky32:_probl_runs}.
\end{remark}

\begin{remark}\label{rmk:wait_prob}
  We can also use $\pi_{\ge k}(n,0)$ to obtain the probability
  $\widehat{\pi}_k(n,i)$ that we have to ``wait''~$i$~indices in order
  to observe the first $k\ge 1$ consecutive ones in a randomly drawn
  $n$-string. Reasoning like in Remark~\ref{rmk:waiting}, this
  probability is
  \begin{equation}
    \label{eq:wait_prob}
    \widehat{\pi}_k(n,i)=\pi_{\ge
      k}\big(i-(k+1),0\big)\,(1-q)q^k\,\iv{i\le n}.
  \end{equation}
  When $q=1/2$,~\eqref{eq:wait_prob} follows directly
  from~\eqref{eq:what_n_i} by using
  $\widehat{\pi}_k(n,i)=\widehat{w}_k(n,i)/2^n$.

  Laplace studied the case $\widehat{\pi}_k(n,n)$
  in~\cite[Liv.\ II, Ch.\ II,
  N$^{\text{o}}$12]{laplace20:_theorie}, which we denote in the
  following by $\widehat{\pi}_k(n)$. From~\eqref{eq:wait_prob}
  and~\eqref{eq:pgf_p_n_m_gtk_m} we have that the pgf
  $\widehat{\Pi}_{k}(x)=\sum_n \widehat{\pi}_k(n) \,x^n$ is
  \begin{align}
    \widehat{\Pi}_{k}(x)&=(1-q)q^k x^{k+1}[y^0]\Pi_{\ge k}(x,y)\label{eq:lap_a}\\&=\frac{(1-q\, x)q^kx^k}{1-x+(1-q)q^kx^{k+1}},\label{eq:lap_b}
  \end{align}
  which was given by Laplace~\cite[p.\ 252]{laplace20:_theorie}, and
  which, from~\eqref{eq:lap_a} and~\eqref{eq:pgf_demoivre_a}, can also be put as
  \begin{equation*}
    \label{eq:wait_prob_lap_atleast1}
    \widehat{\Pi}_{k}(x)=(1-x)\,\Pi^{\ge 1}_{\ge k}(x).
  \end{equation*}
  Feller gave~\eqref{eq:lap_b} in~\cite[Ch.\ XIII,
  Eq. (7.6)]{feller68:_introd_probab}, but he was unaware of Laplace's
  earlier computation ---even though he references other results
  from~\cite{laplace20:_theorie}.  This has led a number of authors to
  believe that~\eqref{eq:lap_b} was originally given by
  Feller. Observe that we have not used Feller's criterion to
  produce~\eqref{eq:lap_b}, and, in fact, Laplace did not use renewal
  theory either to derive~\eqref{eq:lap_b} ---see again discussion in
  Remark~\ref{rmk:intro}.
\end{remark}

\subsection{Probability~that~the Longest Run of an~$n$-String Is a
  $(\le\! k)$-Run}
\label{sec:pis_n_ltk}
We denote the probability that the longest run of an
$n$-string is a $(\le\! k)$-run by $\widebar{\pi}_{\le k}(n)$. Of
course, this is a special case of the results in the previous section
because 
\begin{equation}
  \label{eq:pis_n_ltk_from_pi_n_gek}
  \widebar{\pi}_{\le k}(n)=\pi_{\ge (k+1)}(n,0).
\end{equation}
We treat this special case separately because it was previously dealt
with by Muselli~\cite{muselli96:_simple} and
Schilling~\cite{schilling90:_longest}.

\subsubsection{Recurrences}
\label{sec:pis_n_ltk_rec}
Using~\eqref{eq:pis_n_ltk_from_pi_n_gek} in recurrences~\eqref{eq:recurrence_p_n_m_gtk_basic}
and~\eqref{eq:recurrence_p_n_m_gtk_simpler} we have, respectively, the
following two recurrences:
\begin{equation}
  \label{eq:ps_n_ltk_rec1}
  \widebar{\pi}_{\le k}(n)= (1-q)  \sum_{i=0}^{k} q^i\,\widebar{\pi}_{\le k}\big(n-(i+1)\big),
\end{equation}
and
\begin{equation}
  \label{eq:ps_n_ltk_rec2}
  \widebar{\pi}_{\le k}(n)=\widebar{\pi}_{\le k}(n-1)-q^{k+1}(1-q)\, \widebar{\pi}_{\le k}\big(n-(k+2)\big).
\end{equation}
From~\eqref{eq:pis_n_ltk_from_pi_n_gek},~\eqref{eq:p_n_m_gtk_init}
and~\eqref{eq:p_n_m_gtk_init2}, they are both initialised using $\widebar{\pi}_{\le k}(-1)=1/(1-q)$
and $\widebar{\pi}_{\le k}(0)=1$.

\begin{remark}
  Although~\eqref{eq:ps_n_ltk_rec1} is the direct probabilistic
  counterpart of~\eqref{eq:rec_ws_n_lek_schilling}, Schilling did
  not provide this recurrence. This author proposed instead to compute
  $\widebar{\pi}_{\le k}(n)$ as follows~\cite[Eq.\
  (3)]{schilling90:_longest}:
  \begin{equation}
    \label{eq:ps_n_ltk_schilling}
    \widebar{\pi}_{\le k}(n)=\sum_{r=0}^n \widebar{w}_{\le k}(n,r)\, q^r(1-q)^{n-r},
  \end{equation}
  where $\widebar{w}_{\le k}(n,r)$ ---studied later in
  Section~\ref{sec:ws_n_m_k_hamming}--- is the number of $n$-strings
  of Hamming weight $r$ whose longest run is a $(\le\! k)$-run, which
  Schilling showed how to obtain recursively
  using~\eqref{eq:rec_ws_n_r}~\cite[Eq.\ (4)]{schilling90:_longest}.
\end{remark}

\subsubsection{Probability Generating Function}
From~\eqref{eq:pis_n_ltk_from_pi_n_gek} and~\eqref{eq:pgf_p_n_m_gtk_m}, the pgf
$\widebar{\Pi}_{\le k}(x)=\sum_n \widebar{\pi}_{\le k}(n)\, x^n$, is
\begin{equation}
  \label{eq:ps_ltk}
  \widebar{\Pi}_{\le k}(x)=\frac{1-q\, x}{(1-q)\,x\,\big(1-x+(1-q)q^{k+1}x^{k+2}\big)}.
\end{equation}

\subsubsection{Explicit Expressions}
\label{sec:pis_n_ltk_pgf_explicit}
As $\widebar{w}_{\le k}(n,r)$ has the single-summation closed
form~\eqref{eq:ws_n_ltk_r_explicit_good}, Schilling would probably be
pleased to see that~\eqref{eq:ps_n_ltk_schilling} is actually a
double-summation explicit expression for~$\widebar{\pi}_{\le k}(n)$
---although his asymptotic results are a more powerful proposition
indeed. In any case, from~\eqref{eq:pis_n_ltk_from_pi_n_gek}
and~\eqref{eq:pi_gtk_n_m_th} we get a single-summation explicit
expression for~$\widebar{\pi}_{\le k}(n)$, similar to the one given by
Muselli~\cite[Eq.\ (16)]{muselli96:_simple}. Equivalent formulas were
previously found in reliability problems by Hwang~\cite[Thm.\
3]{hwang86:_simpl} and by Lambiris and
Papastavridis~\cite[Eq.\ (1)]{lambiris85:_exact}.

\subsection{Probability that an $n$-String Contains Exactly $m$
  Nonnull $p$-Parity Runs}
\label{sec:p_wp_n_m}
Let $\pi_{[p]}(n,m)$ be the probability that an $n$-string contains
exactly $m$ nonnull $p$-parity runs. Like in
Section~\ref{sec:probability-w_m_klek}, if $q=1/2$ then we simply have that
$\pi_{[p]}(n,m)=w_{[p]}(n,m)/2^n$. We examine the case for arbitrary $q$ next.

\subsubsection{Recurrence}
\label{sec:pi_n_m_p_recurrence}
We may write a probability recurrence for
$\pi_{[p]}(n,m)$ by using the enumerative
recurrence~\eqref{eq:wp_n_m_rec2} in conjunction with the law of total
probabilities:
\begin{equation}
  \label{eq:p_wp_n_m_rec}
  \pi_{[p]}(n,m)=(1-q)\bigg(\pi_{[p]}(n-1,m)+\sum_{i=1}^n q^i \pi_{[p]}\big(n-(i+1),m-\iv{p=\text{mod}(i,2)}\big)\bigg).
\end{equation}
As usual, we find initial values for recurrence~\eqref{eq:p_wp_n_m_rec}
by using the case $n=1$, which we know by inspection to be
\begin{equation}
  \label{eq:trivial_n1_p_wp}
  \pi_{[p]}(1,m)=\big(q\,\iv{p=0}+1-q\big)\iv{m=0}+q\,\iv{p=1}\iv{m=1}.
\end{equation}
On the other hand, setting~$n=1$ in~\eqref{eq:p_wp_n_m_rec} yields
\begin{equation}
  \label{eq:p_wp_n1}
  \pi_{[p]}(1,m)=(1-q)\,\bigg(\pi_{[p]}(0,m)+q\,\pi_{[p]}(-1,m-\iv{p=1})\bigg).
\end{equation}
We wish~\eqref{eq:p_wp_n1} to equal~\eqref{eq:trivial_n1_p_wp}. We can achieve this equality by using
\begin{align*}
  \pi_{[p]}(-1,m)&=\frac{1}{1-q}\iv{m=0},\\
  \pi_{[p]}(0,m)&=\iv{m=0},
\end{align*}
which are therefore the initial values
of~\eqref{eq:p_wp_n_m_rec}. %

From~\eqref{eq:wp_n_atleast} we also have have
$\pi^{\ge
  m}_{[p]}(n)=\sum_{t=m}^{\floor*{(n+1)/(2+\iv{p=0})}}\pi_{[p]}(n,t)$. Applying
this expression to recurrence~\eqref{eq:p_wp_n_m_rec} we get
\begin{equation}
  \label{eq:pi_n_atleast_rec}
  \pi^{\ge m}_{[p]}(n)=(1-q)\Big(\pi^{\ge
    m}_{[p]}(n-1)+\sum_{i=1}^nq^i \big(\pi^{\ge m}_{[p]}(n-(i+1))+\iv{p=\text{mod}(i,2)}\,\pi_{[p]}(n-(i+1),m-1)\big)\Big).
\end{equation}

\subsubsection{Probability Generating Functions}
\label{sec:pgf_pi_n_m_p}
In order to obtain the pgf
$\Pi_{[p]}(x,y)=\sum_{n,m}\pi_{[p]}(n,m)\,x^ny^m$ we first make
recurrence~\eqref{eq:p_wp_n_m_rec} valid for all $n$ and $m$ by adding
$\iv{n=-1}\iv{m=0}/(1-q)$ to~\eqref{eq:p_wp_n_m_rec}. We
then get a recurrence free from $n$-dependent summations by obtaining
$\pi_{[p]}(n,m)-q^2\,\pi_{[p]}(n-2,m)$ using the extended recurrence,
which yields
\begin{align}
  \label{eq:p_wp_n_m_rec2}
  \pi_{[p]}(n,m)=&~q^2\,\pi_{[p]}(n-2,m)+(1-q)\bigg(\pi_{[p]}(n-1,m)-q^2\,\pi_{[p]}(n-3,m)\nonumber\\
                 &+q\,
                   \pi_{[p]}(n-2,m-\iv{p=1})+q^2\,\pi_{[p]}(n-3,m-\iv{p=0})\bigg)\nonumber\\
  &+\big(\iv{n=-1}-q^2\iv{n=1}\big)\iv{m=0}\frac{1}{1-q}.
\end{align}
We can now directly obtain $\Pi_{[p]}(x,y)$
from~\eqref{eq:p_wp_n_m_rec2} by the usual method of multiplying
across by $x^ny^m$ and adding over the range of $n$ and $m$:
\begin{align}
  \label{eq:pgf_p_wp}
  \Pi_{[p]}(x,y)=\frac{1-q^2 x^2}{(1-q)\,x\,\big(1-q^2\,x^2-(1-q)\big(x+qx^2y^{\iv{p=1}}-q^2x^3(1-y^{\iv{p=0}})\big)\big)}.
\end{align}
To conclude, we extract the coefficient of $y^m$ from~\eqref{eq:pgf_p_wp} using
the same approach as in Section~\ref{sec:generating-function_klek}, to
get
\begin{align}
  \label{eq:pgf_p_wp_ym}
  [y^m]\Pi_{[p]}(x,y)=\frac{(1-q^2x^2)\big((1-q)q^{1+\iv{p=0}}x^{2+\iv{p=0}}\big)^m}{(1-q)\,x\,\big(1-(1-q)x-q^{1+\iv{p=1}}\,x^2+\iv{p=0}(1-q)q^2x^3\big)^{m+1}},
\end{align}
which is the pgf giving the probability that a binary string contains
exactly $m$ nonnull $p$-parity runs.

The related pgf
$\Pi^{\ge m}_{[p]}(x)=\sum_n \pi^{\ge m}_{[p]}(n)\, x^n$ can be
obtained from~\eqref{eq:pi_n_atleast_rec} in a similar way
as~\eqref{eq:pgf_demoivre} in the previous section. We first obtain a
recurrence free from $n$-dependent summations by subtracting
$q^2\pi^{\ge m}_{[p]}(n-2)$ from $\pi^{\ge m}_{[p]}(n)$
using~\eqref{eq:pi_n_atleast_rec}, which yields the following
recurrence:
\begin{align}
  \label{eq:pi_n_atleast_rec2}
  \pi^{\ge m}_{[p]}(n)=~&q^2\pi^{\ge m}_{[p]}(n-2)+
                  (1-q)\Big(\pi^{\ge m}_{[p]}(n-1)+q\,\pi^{\ge m}_{[p]}(n-2) \nonumber\\&
  +\iv{p=1}q\,w_{[p]}(n-2,m-1)+\iv{p=0}q^2\,w_{[p]}(n-3,m-1)\Big).
\end{align}
As usual, by  multiplying~\eqref{eq:pi_n_atleast_rec2} on both
sides by $x^n$ and then adding over $n$ we get 
\begin{align}
  \label{eq:ogf_pi_n_p_atleast}
  \Pi^{\ge m}_{[p]}(x)&=\frac{(1-q)\big(\iv{p=1}\,q\,x^2+\iv{p=0}\,q^2\,x^3\big)}{1-(1-q)x-q
                  x^2}\,[y^{m-1}]\Pi_{[p]}(x,y)\nonumber\\
&=\frac{qx\big(\iv{p=1}+\iv{p=0}\,qx\big)(1-q^2x^2)\big((1-q)q^{1+\iv{p=0}}x^{2+\iv{p=0}}\big)^{m-1}}{\big(1-(1-q)x-
                  q x^2\big)\big(1-(1-q)x-q^{1+\iv{p=1}}\,x^2+\iv{p=0}(1-q)q^2x^3\big)^m},
\end{align}
where we have used~\eqref{eq:pgf_p_wp_ym}. This pgf gives the
probability that a binary string contains at least $m$ nonnull $p$-parity
runs.

\subsubsection{Explicit Expressions}
We next determine explicit expressions for $\pi_{[p]}(n,m)$ by
extracting the coefficient of $x^n$ in~\eqref{eq:pgf_p_wp_ym}. We
begin with the case $p=1$, where we first
express~\eqref{eq:pgf_p_wp_ym} as
$[y^m]\Pi_{[1]}(x,y)=(1-q)^{m-1}q^mx^{2m-1}/((1-q^2x^2)^m(1-(1-q)x/(1-q^2x^2))^{m+1})$. Applying
the negative binomial theorem twice to this expression we get
\begin{equation}
  \label{eq:pgf_p1_wp_ym_expanded}
  [y^m]\Pi_{[1]}(x,y)=\sum_{t\ge 0}\sum_{s\ge 0}{t+m\choose
    m}{s+t+m-1\choose s}(1-q)^{t+m-1}q^{2s+m}x^{2m-1+2s+t}.
\end{equation}
The coefficient of $x^n$ is found through the nonnegative indices $t$
and $s$ that satisfy $2m-1+2s+t=n$. The maximum of $s$ happens when
$t=0$, and thus $s\le \floor{(n-2m+1)/2}$, whereas~$t$ is determined
by the previous equation for any given value of $s$. From these
considerations and~\eqref{eq:pgf_p1_wp_ym_expanded} we have that
\begin{equation}
  \label{eq:pi_n_m_p1_explicit}
  \pi_{[1]}(n,m)=\sum_{s=0}^{\floor{\frac{n-2m+1}{2}}}{n-m+1-2s\choose
    m}{n-m-s\choose s}(1-q)^{n-m-2s}q^{2s+m}.
\end{equation}
We now deal with the case $p=0$, which is slightly more involved as it
involves a double summation rather than a single one
---cf. Section~\ref{sec:w_n_m_p_explicit-expressions}. In this
case~\eqref{eq:pgf_p_wp_ym} can be expanded using the
negative binomial theorem and the binomial theorem (twice) as follows:
\begin{equation}
  \label{eq:pgf_p0_wp_ym_expanded}
  [y^m]\Pi_{[0]}(x,y)=\sum_{t\ge 0}\sum_{s\ge 0}\sum_{r\ge 0}{t+m\choose
    m}{t\choose s}{t-s+1\choose r}(1-q)^{t-s+m-1}q^{2m+s+2r}(-1)^rx^{3m-1+t+s+2r}.
\end{equation}
To determine the coefficient of $x^n$ we just determine the
nonnegative indices $t,s$ and $r$ that fulfil $n=3m-1+t+s+2r$. From
$t=s=0$ we have $s\le \floor{(n-3m+1)/2}$, and from $s=0$ we have that
$t\le n-3m+1-2r$. Finally $s=n-3m+1-2r-t$. Thus we have
from~\eqref{eq:pgf_p0_wp_ym_expanded} that
\begin{equation}
  \label{eq:pi_n_m_p1_explicit}
  \pi_{[0]}(n,m)=\sum_{r=0}^{\floor{n-3m+1}{2}}\sum_{t=0}^{n-3m+1-2r}{t+m\choose
    m}{t\choose s}{t-s+1\choose r}(1-q)^{t-s+m-1}q^{2m+s+2r}(-1)^r,
\end{equation}
where $s$ depends on $r$ and $t$ as just indicated.

\subsubsection{Moments}
We denote by $M_{[p],n}$ the random variable that models the number of
$p$-parity runs in an $n$-string drawn at random, i.e.,
\begin{equation}
  \label{eq:prmpeqm}
  \Pr(M_{[p],n}=m)=\pi_{[p]}(n,m).
\end{equation}
Using~\eqref{eq:pgf_p_wp}, we may obtain the moments of
$M_{[p],n}$. First, the expectation is
$\text{E}(M_{[p],n})=[x^n](\partial/\partial y) \Pi_{[p]}(x,y)|_{y=1}$. Since we have that
\begin{equation}
  \label{eq:ddy_pip_y1}
  \left.\frac{\partial \Pi_{[p]}(x,y)}{\partial
      y}\right|_{y=1}=\frac{qx (1-qx)\big(qx \iv{p=0}+\iv{p=1}\big)}{(1-x)^2(1+qx)},
\end{equation}
by applying the negative binomial theorem twice we get
\begin{equation}
  \left.\frac{\partial \Pi_{[p]}(x,y)}{\partial y}\right|_{y=1}=qx (1-qx)\big(qx \iv{p=0}+\iv{p=1}\big)\sum_{t\ge 0}\sum_{r\ge 0}
{t+1\choose t} (-1)^r q^{r+1} x^{t+r}.
\end{equation}
Focussing on the $p=1$ case, the coefficient of $x^n$ is found using
the nonnegative indices $t$ that solve $n=t+r+1$ and $t+r+2$. The
first equation implies that $t\le n-1$, whereas the second one implies
$t\le n-2$. We thus have from~\eqref{eq:ddy_pip_y1} that
\begin{align}
  \label{eq:exp_m1n}
  \text{E}(M_{[1],n})&=\sum_{t=0}^{n-1}(t+1))(-1)^{n-1-t}q^{n-t}-q\sum_{t=0}^{n-2}(t+1))(-1)^{n-2-t}q^{n-t-1}\nonumber\\
  &=\frac{2q^2\big(1-n(1+q)-(-1)^nq^n\big)}{(1+q)^2}+nq.
\end{align}
For large $n$, $\text{E}(M_{[1],n})\approx (1-2q/(1+q))\, n q$. The case
$p=0$ can be determined using~\eqref{eq:p_n_gtk_expectation} and~\eqref{eq:exp_m1n}:
\begin{equation}
  \label{eq:exp_m0n}
  \text{E}(M_{[0],n})=\text{E}(M_{\ge 1,n})-\text{E}(M_{[1],n}).
\end{equation}
The second factorial moment may be computed in a similar fashion.

\subsection{Discussion}
It should be clear from our exposition in this section that it is
relatively straightforward to obtain probabilistic counterparts of
every aspect of the enumeration problems that we have studied in
Section~\ref{sec:number-n-strings}, and that this procedure is
generally simpler and more versatile than addressing success runs from
a combinatorial analysis viewpoint. In particular, pgfs are easy to
obtain through recurrences, and they are generally preferable to
complicated explicit expressions when it comes to obtaining
moments. If desired, explicit expressions can be also obtained from
pgfs, and, importantly, they are also a conduit to asymptotic
probabilistic results.

Last but not least, we see later in Section~\ref{sec:prob-conn} that
all the probabilistic results in this section can also be approximated
through deterministic enumerations via the law of large numbers.

\section{Number of $n$-Strings of Hamming Weight $r$ that Contain
  Prescribed Quantities of Runs Under Different Constraints}
\label{sec:number-n-strings-hamming}
In this section we address the same basic problems as in
Section~\ref{sec:number-n-strings} but when the $n$-strings are
restricted to having Hamming weight $r$ ---i.e., to containing exactly
$r$ ones. In the terminology of Goulden and Jackson~\cite[Sec.\
2.4.7]{goulden83:_combin_enumer} these are enumerations \textit{with
  type restriction}. Many early works in the theory of runs assume
this constraint, particularly those that use runs for hypothesis
tests. The results in this section are complemented by other results
given later in Section~\ref{sec:oz_hamming}, where runs of ones and
zeros are jointly considered.

We  cover the counterparts of the problems
in Sections~\ref{sec:w_n_m_klek} and~\ref{sec:wp_n_m} in
Sections~\ref{sec:w_n_m_klek_hamming} and~\ref{sec:wp_n_m_r},
respectively. The former enumeration also solves the same problems as
in Sections~\ref{sec:w_n_m_k}--\ref{sec:number-wl_n_k_gtk} when the
fixed Hamming weight constraint is observed.  We actually discuss two
of these special cases in more detail in Sections~\ref{sec:w_n_m_r}
and~\ref{sec:ws_n_m_k_hamming}, since they have been studied by a considerable
number of previous
authors~\cite{stevens39:_distr,mood40:_distr_theor_runs,mosteller41:_note_applic,bateman48:_power,burr61:_longest,philippou86:_successes,schilling90:_longest,schuster94:_exchan}.

Whereas the enumerative recurrences in Sections~\ref{sec:w_n_m_klek}
and~\ref{sec:wp_n_m} are bivariate, solving fixed Hamming weight
versions of these enumerations involves establishing trivariate
recurrence relations. This means that the corresponding ogfs are also
trivariate, which in some scenarios may hinder the determination of
compact explicit expressions. Nevertheless, this is possible in a
number of special cases.

\subsection[Number of n-Strings of Hamming Weight r that Contain
Exactly m (ḵ≤k)-Runs]{Number of $n$-Strings of Hamming Weight $r$ that Contain Exactly $m$ ($\ushort{k}\le\widebar{k}$)-Runs}
\label{sec:w_n_m_klek_hamming}

Let $w_{\ushort{k}\le\widebar{k}}(n,m,r)$ represent the number of
$n$-strings of Hamming weight $r$ that contain exactly $m$
$(\ushort{k}\le\widebar{k})$-runs. As usual, the $n$-strings that we
enumerate may have other runs, as long as they are longer than
$\widebar{k}$ or shorter than $\ushort{k}$. Observe that this is a
conditional version of the general enumeration in
Section~\ref{sec:w_n_m_klek} in which we only consider the
${n\choose r}$ $n$-strings that have Hamming weight~$r$.

Next, we need to tighten necessary
condition~\eqref{eq:m_necessary_condition} in order to take the
Hamming weight constraint into account.
\begin{condition}{(Existence of $n$-strings of Hamming weight $r$
    containing $m$~$(\ushort{k}\le\widebar{k})$-runs)}
  \begin{equation}
    \label{eq:m_necessary_condition_r}
    w_{\ushort{k}\le\widebar{k}}(n,m,r)> 0\quad \Longrightarrow\quad 0\le m \le \min\Bigg(\floor[\bigg]{\frac{n+1}{\ushort{k}+1}},\floor[\bigg]{\frac{r}{\ushort{k}}}\Bigg).
  \end{equation}
  Observe that if $r>0$ and $\ushort{k}=0$ then the upper bound on $m$
  is $n+1$. When $r=0$ and $\ushort{k}=0$ the ratio~$r/\ushort{k}$ is undefined, but the
  upper bound is still $n+1$, as this is the number of null runs in the
  all-zeros $n$-string.
\end{condition}
\subsubsection{Recurrences}
\label{sec:recurrences_w_n_m_k_r}
We can resort to a nearly identical reasoning as in
Section~\ref{sec:recurrences_w_n_m_klek} to get a recurrence for this
enumeration: the quantity $w_{\ushort{k}\le\widebar{k}}(n,m,r)$ can be
broken down into different contributions corresponding to the ensemble
of $n$-strings that start with an $i$-run. If $i>\widebar{k}$ or
$i<\ushort{k}$ then the contribution is
$w_{\ushort{k}\le\widebar{k}}(n-(i+1),m,r-i)$, but if
$\ushort{k}\le i\le \widebar{k}$ then the contribution is
$w_{\ushort{k}\le\widebar{k}}(n-(i+1),m-1,r-i)$. Thus, considering all
possible lengths of the starting run we get the following trivariate
recurrence:
\begin{align}
  \label{eq:rec_w_n_m_klek_hamming}
  w_{\ushort{k}\le\widebar{k}}(n,m,r)=&~\sum_{i=\ushort{k}}^{\widebar{k}}
    w_{\ushort{k}\le\widebar{k}}\big(n-(i+1),m-1,r-i\big)%
  +\sum_{i=0}^{\ushort{k}-1}
    w_{\ushort{k}\le\widebar{k}}\big(n-(i+1),m,r-i\big)\nonumber \\&+\sum_{i=\widebar{k}+1}^{n}
    w_{\ushort{k}\le\widebar{k}}\big(n-(i+1),m,r-i\big).  
\end{align}
Observe that the only difference with respect to
recurrence~\eqref{eq:recurrence_w_n_m_klek_basic} are the updates of
the Hamming weights in the recursive invocations on the right hand
side.  Following the same steps as in
Section~\ref{sec:recurrences_w_n_m_klek}, we can argue an alternative
recurrence equivalent to~\eqref{eq:rec_w_n_m_klek_hamming}:
\begin{align}
  \label{eq:rec_w_n_m_klek_r_simpler}
  w_{\ushort{k}\le\widebar{k}}(n,m,r)=&~w_{\ushort{k}\le\widebar{k}}(n-1,m,r)+w_{\ushort{k}\le\widebar{k}}(n-1,m,r-1)\nonumber\\
            & + w_{\ushort{k}\le\widebar{k}}\big(n-(\ushort{k}+1),m-1,r-\ushort{k}\big)-w_{\ushort{k}\le\widebar{k}}\big(n-(\ushort{k}+1),m,r-\ushort{k}\big)\nonumber\\
            & - w_{\ushort{k}\le\widebar{k}}\big(n-(\widebar{k}+2),m-1,r-(\widebar{k}+1)\big)+w_{\ushort{k}\le\widebar{k}}\big(n-(\widebar{k}+2),m,r-(\widebar{k}+1)\big),
\end{align}
which parallels recurrence~\eqref{eq:recurrence_w_n_m_klek_simpler},
and which can also be obtained by calculating
$w_{\ushort{k}\le\widebar{k}}(n,m,r)-w_{\ushort{k}\le\widebar{k}}(n-1,m,r-1)$
using~\eqref{eq:rec_w_n_m_klek_hamming}.

In the following we work with the simpler
recurrence~\eqref{eq:rec_w_n_m_klek_r_simpler}, as it is free from
$n$-dependent summations and thus directly amenable to determining the
ogf. In order to find initialisation values we rely on the case $n=1$,
in which we know by inspection that
\begin{align}
  \label{eq:w_1_m_k_r_trivial}
  w_{\ushort{k}\le\widebar{k}}(1,m,r)=&~\iv{r=0}\big(\iv{m=0}\iv{\ushort{k}\ge
                                        1}+\iv{m=2}\iv{\ushort{k}=0}\big)\nonumber\\&
  +\iv{r=1}\Big(\iv{m=0}\big(\iv{\ushort{k}>1}+\iv{\ushort{k}=0}\iv{\widebar{k}=0}\big)+\iv{m=1}\iv{\ushort{k}=1}\Big),
\end{align}
Setting next $n=1$ in~\eqref{eq:rec_w_n_m_klek_r_simpler} yields
\begin{align}
  \label{eq:rec_w_1_m_k_r_simpler}
  w_{\ushort{k}\le\widebar{k}}(1,m,r)=&~w_{\ushort{k}\le\widebar{k}}(0,m,r-1)+w_{\ushort{k}\le\widebar{k}}(0,m,r)\nonumber\\
             &+ w_{\ushort{k}\le\widebar{k}}(-\ushort{k},m-1,r-\ushort{k})-w_{\ushort{k}\le\widebar{k}}(-\ushort{k},m,r-\ushort{k})\nonumber\\
            & - w_{\ushort{k}\le\widebar{k}}\big(-(\widebar{k}+1),m-1,r-(\widebar{k}+1)\big)+w_{\ushort{k}\le\widebar{k}}\big(-(\widebar{k}+1),m,r-(\widebar{k}+1)\big),
\end{align}
Taking~\eqref{eq:m_necessary_condition_r} into account, we may verify
that~\eqref{eq:rec_w_1_m_k_r_simpler}
equals~\eqref{eq:w_1_m_k_r_trivial} for all $m,r$ and
$0\le \ushort{k}\le \widebar{k}$ if we choose the values
\begin{align}
  \label{eq:w_n_m_k_r_init}
  w_{\ushort{k}\le\widebar{k}}(-1,m,r)=&\iv{m=0}\iv{r=0},\\
  w_{\ushort{k}\le\widebar{k}}(0,m,r)=&\ivl{m=\iv{\ushort{k}=0}}\iv{r=0}, \label{eq:w_n_m_k_r_init2}
\end{align}
which therefore constitute the initialisation
of~\eqref{eq:rec_w_n_m_klek_r_simpler}.  The reader may check that
these values also initialise the equivalent
recurrence~\eqref{eq:rec_w_n_m_klek_hamming}.

\subsubsection{Generating Function}
\label{sec:ogf_w_n_m_k_r}
We obtain next the trivariate ogf
$W_{\ushort{k}\le\widebar{k}}(x,y,z)=\sum_{n,m,r}
w_{\ushort{k}\le\widebar{k}}(n,m,r)\, x^n y^mz^r$ using
recurrence~\eqref{eq:rec_w_n_m_klek_r_simpler}.  We first need to make
this recurrence valid for all values of $n,m$ and $r$, which, like in
Section~\ref{sec:generating-function_klek}, simply involves
extending~\eqref{eq:rec_w_n_m_klek_r_simpler} so that it works for $n=-1$ and
$n=0$. Through the same procedure as in that section, the reader may
verify that we now have to add
$\iv{n=-1}\iv{m=0}\iv{r=0}-\iv{n=0}\iv{m=0}\iv{r=1}$
to~\eqref{eq:rec_w_n_m_klek_r_simpler}. We then multiply both sides of
the extended recurrence by $x^n y^mz^r$ and sum on $n,m,$ and $r$,
from which we readily get 
\begin{equation}
  \label{eq:ogf_w_n_m_k_r}
  W_{\ushort{k}\le\widebar{k}}(x,y,z)=\frac{1-xz}{x\Big(1-x(z+1)+(1-y)\big(x^{\ushort{k}+1}z^{\ushort{k}}-x^{\widebar{k}+2}z^{\widebar{k}+1}\big)\Big)}.
\end{equation}
Like in Section~\ref{sec:generating-function_klek}, if we write this
ogf as a function of $(1-cy)^{-1}$ we can see that the coefficient of $y^m$ in~\eqref{eq:ogf_w_n_m_k_r} is
\begin{equation}
  \label{eq:ym_ogf_w_n_m_klek_r}
  [y^m]W_{\ushort{k}\le\widebar{k}}(x,y,z)=\frac{(1-xz)\Big(x^{\ushort{k}+1}z^{\ushort{k}}-x^{\widebar{k}+2}z^{\widebar{k}+1}\Big)^m}{x\Big(1-x(z+1)+x^{\ushort{k}+1}z^{\ushort{k}}-x^{\widebar{k}+2}z^{\widebar{k}+1}\Big)^{m+1}}.
\end{equation}

\subsection[Number of n-Strings of Hamming Weight r that Contain
Exactly m Runs]{Number of $n$-Strings of Hamming Weight $r$ that
  Contain Exactly $m$ Nonnull Runs}
\label{sec:w_n_m_r}

We denote by $w(n,m,r)$ the number of $n$-strings of Hamming weight
$r$ that contain exactly~$m$ nonnull runs. This enumeration ---in
fact, the probability of drawing $m$ nonnull runs when $n$-strings
of Hamming weight $r$ are drawn uniformly at random, which we denote
by $\pi(n,m,r)=w(n,m,r)/{n\choose r}$--- was explicitly solved by a
number of
authors~\cite{stevens39:_distr,mood40:_distr_theor_runs,charalambides02:_enumer_combin,schuster94:_exchan}.

This is a particular case of the analysis in
Section~\ref{sec:w_n_m_klek_hamming} with $\ushort{k}=1$ and
$\widebar{k}=n$, and therefore
\begin{equation}
  \label{eq:wnmr}
  w(n,m,r)=w_{1\le n}(n,m,r).
\end{equation}

\subsubsection{Recurrences}
\label{sec:w_n_m_r_rec}
Specialising recurrences~\eqref{eq:rec_w_n_m_klek_hamming}
and~\eqref{eq:rec_w_n_m_klek_r_simpler} to this case yields the
following two recursive relations:
\begin{align}
  \label{eq:rec_w_n_m_r_rec}
  w(n,m,r)=&    w(n-1,m,r)+\sum_{i=1}^{n}    w(n-(i+1),m-1,r-i),
\end{align}
and
\begin{align}
  \label{eq:w_n_m_1len_r_simpler}
    w(n,m,r)=&~w(n-1,m,r-1)+w(n-1,m,r)\nonumber\\
             & + w(n-2,m-1,r-1)-w(n-2,m,r-1),
\end{align}
both of which from~\eqref{eq:w_n_m_k_r_init}
and~\eqref{eq:w_n_m_k_r_init2} are initialised by
\begin{equation}
  \label{eq:w_n_m_r_init}
  w(-1,m,r)=w(0,m,r)=\iv{m=0}\iv{r=0}.
\end{equation}

\subsubsection{Generating Function}
\label{sec:w_n_m_r_ogf}
We obtain next the ogf
$W(x,y,z)=\sum_{n,m,r}w(n,m,r)\,x^ny^mz^r$. Although we may indeed
specialise~\eqref{eq:ym_ogf_w_n_m_klek_r} to this case with
$\ushort{k}=1$ and $\widebar{k}=n$, this gives an $n$-dependent ogf. A
simpler, more general ogf is obtained by proceeding from
recurrence~\eqref{eq:w_n_m_1len_r_simpler} just like in we did to
get~\eqref{eq:ym_ogf_w_n_m_klek_r}
from~\eqref{eq:rec_w_n_m_klek_r_simpler}. This yields
\begin{equation}
  \label{eq:ogf_w_n_m_1len_r}
  W(x,y,z)=\frac{1-xz}{x(1-x-xz+x^2z-x^2yz)},
\end{equation}
from which we have
\begin{equation}
  \label{eq:ym_ogf_w_n_m_1len_r}
  [y^m]W(x,y,z)=\frac{x^{2m-1}z^m}{(1-xz)^{m}(1-x)^{m+1}}.
\end{equation}
\subsubsection{Explicit Expression}
\label{sec:w_n_m_r_explicit}
Next, we get an explicit expression for $w(n,m,r)$.  Applying the
negative binomial theorem twice, we can
express~\eqref{eq:ym_ogf_w_n_m_1len_r} as
\begin{equation}
  \label{eq:ym_ogf_w_n_m_1len_r_expanded_good}
  [y^m]W(x,y,z)=\sum_{s\ge 0}\sum_{t\ge 0}{s+m-1\choose m-1}{t+m\choose m}\, x^{s+t+2m-1}z^{s+m}.
\end{equation}
Thus the coefficient of $x^nz^r$ corresponds to the indices that
fulfil $n=s+t+2m-1$ and $r=s+m$, i.e., $s=r-m$ and
$t=n-s-2m+1=n-r-m+1$. Therefore,
from~\eqref{eq:ym_ogf_w_n_m_1len_r_expanded_good} we have that
$w(n,m,r)=[x^ny^mz^r]W(x,y,z)$ is
\begin{equation}
  \label{eq:w_n_m_1len_r_hypergeo}
  w(n,m,r)={r-1\choose m-1}{n-r+1\choose m}.
\end{equation}
Either in this enumerative form or in its probabilistic version
$\pi(n,m,r)=w(n,m,r)/{n\choose r}$,
expression~\eqref{eq:w_n_m_1len_r_hypergeo} was previously given by
Stevens~\cite[Eq.\ (3.05)]{stevens39:_distr}, Mood~\cite[Eq.\
(2.11)]{mood40:_distr_theor_runs}, Charalambides~\cite[p.\
97]{charalambides02:_enumer_combin}, Gibbons and
Chakraborti~\cite[Cor.\ 3.2.1]{gibbons11:_nonparam}
Schuster~\cite[Cor.\ 3.5]{schuster94:_exchan} and others. Schuster
explicitly notes that $\pi(n,m,r)$ is a hypergeometric
distribution. %

\subsection[Number of n-Strings of Hamming Weight r Whose Longest Run
Is a $k$-Run or a $(\le k)-Run$]{Number of $n$-Strings of Hamming
  Weight $r$ Whose Longest Run Is a  $k$-Run or a $(\le k)$-Run}
\label{sec:ws_n_m_k_hamming}
   
We consider here the enumeration of the $n$-strings of Hamming weight
$r$ whose longest run is a $k$-run or a $(\le k)$-run, which
we respectively denote by~$\widebar{w}_k(n,r)$ and
$\widebar{w}_{\le k}(n,r)$. These are the Hamming weight constrained
versions of the enumerations in Section~\ref{sec:number-ws_n_k_ltk}.

Many authors have studied these problems, both from enumerative and
probabilistic viewpoints. In the latter case we refer to the situation
in which the $n$-strings of Hamming weight~$r$ are drawn uniformly at
random, where we are interested in the probabilities
$\widebar{\pi}_k(n,r)=\widebar{w}_k(n,r)/{n\choose r}$ and
$\widebar{\pi}_{\le k}(n,r)=\widebar{w}_{\le k}(n,r)/{n\choose r}$
---i.e., the cumulative probability.  The earliest study is by
Mosteller, who gave
$\widebar{\pi}_{\ge k}(n,r)=1-\widebar{\pi}_{\le (k-1)}(n,r)$ for the
case $n=2r$~\cite[Eq.\ (4)]{mosteller41:_note_applic}.
Bateman~\cite{bateman48:_power} first considered an arbitrary ratio
$r/n$ in the same problem and gave an explicit formula for
$\widebar{\pi}_{\ge k}(n,r)$. David and Barton~\cite{david62:_chance},
Burr and Cane~\cite{burr61:_longest}, Philippou and
Makri~\cite{philippou86:_successes}, Gibbons and
Chakraborti~\cite{gibbons11:_nonparam}, and
Schuster~\cite{schuster94:_exchan} gave exact formulas for
$\widebar{\pi}_{\le k}(n,r)$. Recurrences
for~$\widebar{w}_{\le k}(n,r)$ and~$\widebar{\pi}_{\le k}(n,r)$ where
given by Schilling~\cite{schilling90:_longest} and
Schuster~\cite{schuster94:_exchan}, respectively. Finally, and Nej and
Reddy~\cite{nej19:_binary} gave a recurrence and a closed-form
explicit expression for~$\widebar{w}_k(n,r)$.

We address this enumeration using the exact same approach as in
Section~\ref{sec:number-ws_n_k_ltk}. Like in that section, we first
state the connection between $\widebar{w}_k(n,r)$ and
$\widebar{w}_{\le k}(n,r)$. For all $k\ge 0$ we have that
\begin{equation}
  \label{eq:connection_ws_n_r_lek_ws_n_r_k}
  \widebar{w}_{\le k}(n,r)=\sum_{j=0}^k \widebar{w}_{j}(n,r),
\end{equation}
and, conversely, for all $k\ge 1$ it holds that
\begin{equation}
  \label{eq:connection_wsr_n_lek_ws_n_k_r_converse}
  \widebar{w}_k(n,r)=\widebar{w}_{\le k}(n,r)-\widebar{w}_{\le (k-1)}(n,r),
\end{equation}
whereas $\widebar{w}_0(n,r)=\widebar{w}_{\le 0}(n,r)$.  The expression
that allows us to get these two quantities through the general
enumeration $w_{\ushort{k}\le\widebar{k}}(n,m,r)$ is
\begin{equation}
  \label{eq:ws_n_r_ltk_w_n_m_r_gtk_connection}
  \widebar{w}_{\le k}(n,r)=w_{\ge (k+1)}(n,0,r),
\end{equation}
where $w_{\ge k}(n,m,r)=w_{k\le n}(n,m,r)$.  Thus
$\widebar{w}_0(n,r)=w_{\ge 1}(n,0,r)=\iv{r=0}$ for $n\ge 1$, as the
longest run is a null run only in the all-zeros $n$-string, which can
only happen if $r=0$.

\subsubsection{Recurrences}
\label{sec:ws_n_m_r_rec}
We first find recurrences for $\widebar{w}_{\le k}(n,r)$. By
applying~\eqref{eq:ws_n_r_ltk_w_n_m_r_gtk_connection} to
~\eqref{eq:rec_w_n_m_klek_hamming} we get
\begin{align}
  \label{eq:rec_ws_n_r}
  \widebar{w}_{\le k}(n,r)&=\sum_{i=0}^{k}    \widebar{w}_{\le k}(n-(i+1),r-i),
\end{align}
which parallels~\eqref{eq:rec_ws_n_lek_schilling}. Also, by
applying~\eqref{eq:ws_n_r_ltk_w_n_m_r_gtk_connection}
to~\eqref{eq:rec_w_n_m_klek_r_simpler} we get the following
alternative recurrence:
\begin{align}
  \label{eq:rec_ws_n_r_simpler}
  \widebar{w}_{\le k}(n,r)&=\widebar{w}_{\le k}(n-1,r)+\widebar{w}_{\le k}(n-1,r-1)
                            -\widebar{w}_{\le k}\big(n-(k+2),r-(k+1)\big),
\end{align}
which
parallels~\eqref{eq:rec_ws_n_lek_simpler}. From~~\eqref{eq:ws_n_r_ltk_w_n_m_r_gtk_connection},~\eqref{eq:w_n_m_k_r_init}
and~\eqref{eq:w_n_m_k_r_init2}, both~\eqref{eq:rec_ws_n_r}
and~\eqref{eq:rec_ws_n_r_simpler} are initialised using
\begin{equation}
  \label{eq:ws_n_ltk_r_init}
  \widebar{w}_{\le k}(-1,r)=\widebar{w}_{\le k}(0,r)=\iv{r=0}.
\end{equation}
We can also get recurrences for $\widebar{w}_{k}(n,r)$ when $k\ge 1$. By
inputting~\eqref{eq:rec_ws_n_r}
in~\eqref{eq:connection_wsr_n_lek_ws_n_k_r_converse} we find the recurrence
\begin{equation}
  \label{eq:ws_n_k_r_rec_a}
  \widebar{w}_k(n,r)=\sum_{i=0}^{k-1}\widebar{w}_k(n-(i+1),r-i)+\sum_{j=0}^k \widebar{w}_j(n-(k+1),r-k),
\end{equation}
whereas if we input~\eqref{eq:rec_ws_n_r_simpler}
in~\eqref{eq:connection_wsr_n_lek_ws_n_k_r_converse} we get the
alternative recurrence
\begin{equation}
  \label{eq:ws_n_k_r_rec_b}
  \widebar{w}_{k}(n,r)=\widebar{w}_{k}(n-1,r)+\widebar{w}_{k}(n-1,r-1)+\sum_{j=0}^{k-1}\widebar{w}_{j}(n-(k+1),r-k)
  -\sum_{j=0}^{k}\widebar{w}_{j}(n-(k+2),r-(k+1)).
\end{equation}
From~\eqref{eq:connection_wsr_n_lek_ws_n_k_r_converse}
and~\eqref{eq:ws_n_ltk_r_init}, the initialisation of both recurrences
is $\widebar{w}_{k}(-1,r)=\widebar{w}_{k}(0,r)=0$ for $k\ge 1$, but we
can see that we also need initialisation for the case $k=0$. Since
$\widebar{w}_{0}(n,r)=w_{\ge 1}(n,0,r)$,
from~\eqref{eq:w_n_m_k_r_init} and~\eqref{eq:w_n_m_k_r_init2} we have
that $\widebar{w}_{0}(-1,r)=\widebar{w}_0(0,r)=\iv{r=0}$. Thus
recurrences~\eqref{eq:ws_n_k_r_rec_a} and~\eqref{eq:ws_n_k_r_rec_b}
are both initialised by
\begin{equation*}
  \widebar{w}_{k}(-1,r)=  \widebar{w}_{k}(0,r)=\iv{k=0}\iv{r=0}.
\end{equation*}
Again,~\eqref{eq:ws_n_k_r_rec_a} and~\eqref{eq:ws_n_k_r_rec_b}
parallel their unconstrained counterparts~\eqref{eq:ws_n_k_rec_a}
and~\eqref{eq:ws_n_k_rec_b}, respectively.

\begin{remark}
  Recurrence~\eqref{eq:rec_ws_n_r} was given by Schilling
  in~\cite[Eq.\ (4)]{schilling90:_longest}. This author examines in
  some detail the case $k=3$ of~\eqref{eq:rec_ws_n_r}, but is
  nevertheless somewhat vague about its general initialisation.
  Schuster also gave a recurrence for
  $\widebar{\pi}_{\le k}(n,r)=\widebar{w}_{\le k}(n,r)/{n\choose r}$
  based on Pascal's triangle of order $k+1$~\cite[Cor.\
  4.3]{schuster94:_exchan}, which is essentially different from
  either~\eqref{eq:rec_ws_n_r} or~\eqref{eq:rec_ws_n_r_simpler}.  In
  any case, both recurrences above are simpler. Also,
  recurrence~\eqref{eq:ws_n_k_r_rec_a} was previously given by Nej and
  Reddy~\cite[Thm.\ 3.2]{nej19:_binary} ---in their notation,
  $\widebar{w}_k(n,r)=F_n(r,k)$. These authors state that their
  recurrence is valid for $n\ge 3$, $1\le r \le n-2$ and
  $\floor{n/(r+1)}\le k\le n$.

  Regarding the piecemeal validity of the recurrences given by
  previous authors mentioned in this remark: notice that all the
  recurrences given here are valid for all $n\ge 1$ and $r, k\ge 0$
  with the initialisations given, just taking necessary
  condition~\eqref{eq:m_necessary_condition_r} into account. As in
  other such fragmentary recurrences we have met before (see
  Section~\ref{sec:w_n_m_gtk}), this shows the necessity of considering
  the cases $n=-1$ and $n=0$ in the initialisation in order to obtain
  the simplest and most general recurrences.
\end{remark}

\subsubsection{Generating Functions}
\label{sec:ws_n_m_r_ogf}
Considering~\eqref{eq:ws_n_r_ltk_w_n_m_r_gtk_connection}, we may
obtain the ogf
$\widebar{W}_{\!\!\le k}(x,z)=\sum_{n,r}\widebar{w}_{\le
  k}(n,r)\,x^nz^r$ by setting $\ushort{k}=k+1$, $\widebar{k}=n$ and
$m=0$ in \eqref{eq:ym_ogf_w_n_m_klek_r}, but this leads to an
$n$-dependent ogf. A simpler, more general ogf is obtained from
recurrence~\eqref{eq:rec_ws_n_r_simpler} as we show next. Adding
$\iv{n=-1}\iv{r=0}-\iv{n=0}\iv{r=1}$ to~\eqref{eq:rec_ws_n_r_simpler}
to make it valid for all values of the parameters, and then
multiplying it on both sides by $x^nz^r$ and summing over $n$ and $r$
we obtain
\begin{equation}
  \label{eq:ogf_ws_n_ltk_r}
  \widebar{W}_{\!\!\le k}(x,z)=\frac{1-xz}{x\big(1-x-xz+x^{k+2}z^{k+1}\big)}.
\end{equation}
It is possible to derive the ogf
$\widebar{W}_{\!k}(x,z)=\sum_{n,r}\widebar{w}_k(n,r)\,x^nz^r$ from any of
the two recurrences~\eqref{eq:ws_n_k_r_rec_a}
or~\eqref{eq:ws_n_k_r_rec_b}, but this is not straightforward due to
the summations on~$j$ leading to recursive relations on the generating
function itself. Nevertheless it is a simple matter to get this ogf
directly from~\eqref{eq:connection_wsr_n_lek_ws_n_k_r_converse}
and~\eqref{eq:ogf_ws_n_ltk_r}:
\begin{equation}
  \label{eq:ym_ogf_ws_n_k_r}
  \widebar{W}_{\!k}(x,z)=\frac{1-xz}{x}\Bigg(\frac{1}{1-x-xz+x^{k+2}z^{k+1}}-\frac{1}{1-x-xz+x^{k+1}z^{k}}\Bigg).
\end{equation}

\subsubsection{Explicit Expressions}
\label{sec:ws_n_m_r_explicit}
We find next an explicit expression for $\widebar{w}_{\le k}(n,r)$,
for which we rewrite~\eqref{eq:ogf_ws_n_ltk_r} first as
$\widebar{W}_{\!\!\le  k}(x,z)=x^{-1}\big(1-(x-x^{k+2}z^{k+1})/(1-xz)\big)^{-1}$. Applying
the negative binomial theorem twice and then the binomial theorem to
this expression, we can expand the ogf as
\begin{equation}
  \label{eq:ogf_ws_n_ltk_r_expanded2}
  \widebar{W}_{\!\!\le  k}(x,z)=\sum_{t\ge 0}\sum_{s\ge 0}\sum_{p\ge 0}
  {s+t-1\choose t-1}{t\choose p} (-1)^p (xz)^{p(k+1)+s} x^{t-1}.
\end{equation}
To extract the coefficient of $x^nz^r$ we have to find nonnegative
indices $t,s$ and $p$ that fulfil
\begin{align*}
  n&=p(k+1)+s+t-1,\\
  r&=p(k+1)+s.
\end{align*}
From the second equation, the maximum of $p$ happens when $s=0$, and
thus $p\le \floor{r/(k+1)}$. One value of $p$ determines $s=r-p(k+1)$ and 
$t=n+1-p(k+1)-s$. We thus have
from~\eqref{eq:ogf_ws_n_ltk_r_expanded2} that $\widebar{w}_{\le
  k}(n,r)=[x^nz^r]\widebar{W}_{\le  k}(x,z)$ is
\begin{equation}
  \label{eq:ws_n_ltk_r_explicit_good}
  \widebar{w}_{\le    k}(n,r)=\sum_{p=0}^{\floor{\frac{r}{k+1}}}
  (-1)^p   {n-p(k+1)\choose n-r}{n-r+1\choose p}.
\end{equation}
Bateman was the first author who gave an expression parallel
to~\eqref{eq:ws_n_ltk_r_explicit_good}~\cite[p.\
101]{bateman48:_power}.  Her formula gives the probability %
$\widebar{\pi}_{\ge k}(n,r)$, and it is based on a generating function
for compositions restricted to a maximum part ---see start of
Section~\ref{sec:oz} for the connection between runs and compositions.
David and Barton~\cite[p.\ 230]{david62:_chance} and then
Schuster~\cite[Cor.\ 4.2]{schuster94:_exchan} gave alternative
derivations of~\eqref{eq:ws_n_ltk_r_explicit_good} in the form
$\widebar{\pi}_{\le k}(n,r)$. The former authors used combinatorial
arguments and characteristic (Bernoulli) random variables, whereas the
latter author used the theory of exchangeability of random variables
and generating functions. Schuster was driven by the computational
impracticality of some earlier explicit expressions
for~$\widebar{\pi}_{\le
  k}(n,r)$~\cite{burr61:_longest,philippou86:_successes}, but he must
have missed the expression by Bateman that we mentioned above.

From~\eqref{eq:ws_n_ltk_r_explicit_good} we also get a
single-summation explicit expression for $\widebar{w}_{k}(n,r)$
through~\eqref{eq:connection_wsr_n_lek_ws_n_k_r_converse}. Regarding
this case, Nej and Reddy found a more compact closed-form expression
(without summations) for $\widebar{w}_k(n,r)$ through combinatorial
analysis, although restricted to the special case $r< n-1$ and
$r<2k$~\cite[Thm.\ 3.1]{nej19:_binary}. This suggests that it might be
possible to evaluate~\eqref{eq:ws_n_ltk_r_explicit_good} in closed
form, perhaps using the ``snake oil''
method~\cite{wilf06:_generatingf}, although we have not succeeded in
doing so.

\subsection{Number of $n$-Strings of Hamming Weight $r$ that Contain
  Exactly $m$ Nonnull $p$-Parity Runs}
\label{sec:wp_n_m_r}
In this section we consider the counterpart of
Section~\ref{sec:wp_n_m} under a Hamming weight constraint. We denote
by $w_{[p]}(n,m,r)$ the number of $n$-strings of Hamming
weight $r$ that contain exactly $m$ nonnull $p$-parity runs.

As usual, we need to state first a necessary condition similar
to~\eqref{eq:m_necessary_condition_r}.
\begin{condition}{(Existence of $n$-strings of Hamming weight $r$
    containing $m$ nonnull $p$-parity runs)}
  \begin{equation}
    \label{eq:m_necessary_condition_p_r}
    w_{[p]}(n,m,r)> 0\quad \Longrightarrow\quad 0\le m \le \min\Bigg(\floor[\bigg]{\frac{n+1}{2+\iv{p=0}}},\floor*{\frac{r}{1+\iv{p=0}}}\Bigg).
\end{equation}
\end{condition}

\subsubsection{Recurrence}
\label{sec:recurrences_w_n_m_p_r}
We obtain a recurrence with a rationale similar to that in
Section~\ref{sec:recurrences_wp_n_m}: $n$-strings that start with $0$
contribute $w_{[p]}(n-1,m,r)$ to $w_{[p]}(n,m,r)$; on the other hand,
$n$-strings that start with a nonnull $i$-run contribute
$w_{[p]}(n-(i+1),m,r-i)$ if $p=\text{mod}(i,2)$ but
$w_{[p]}(n-(i+1),m-1,r-i)$ if $p\neq\text{mod}(i,2)$. These two
contributions lead to the following trivariate recurrence:
\begin{equation}
  \label{eq:rec_w_n_m_p_r}
  w_{[p]}(n,m,r)=w_{[p]}(n-1,m,r)+ \sum_{i=1}^{n} w_{[p]}(n-(i+1),m-\iv{p=\text{mod}(i,2)},r-i).
\end{equation}
Again, the only difference with respect to~\eqref{eq:wp_n_m_rec2} are
the Hamming weight updates. We initialise~\eqref{eq:rec_w_n_m_p_r} by
first determining its correct value for $n=1$.  By inspection, the
value is
\begin{equation}
  \label{eq:wp_n_m_p_r_trivialn1}
  w_{[p]}(1,m,r)=\iv{r=0}\iv{m=0}+\iv{r=1}\Big(\iv{m=0}\iv{p=0}+\iv{m=1}\iv{p=1}\Big).
\end{equation}
Setting $n=1$ in~\eqref{eq:rec_w_n_m_p_r} now yields
\begin{equation}
  \label{eq:rec_w_n_m_p_r_n1}
  w_{[p]}(1,m,r)=w_{[p]}(0,m,r)+w_{[p]}(-1,m-\iv{p=1},r-1).
\end{equation}
It is readily verified that~\eqref{eq:rec_w_n_m_p_r_n1}
equals~\eqref{eq:wp_n_m_p_r_trivialn1} for the following
initialisation values:
\begin{equation}
  \label{eq:rec_w_n_m_p_r_init}
  w_{[p]}(-1,m,r)=w_{[p]}(0,m,r)=\iv{m=0}\iv{r=0}.
\end{equation}

\subsubsection{Generating Function}
\label{sec:w_n_m_p_r_generating-function}
We obtain next the ogf
$W_{[p]}(x,y,z)=\sum_{n,m,r} w_{[p]}(n,m,r)\,x^ny^mz^r$.  To this end,
we first add $\iv{n=-1}\iv{m=0}\iv{r=0}$ to~\eqref{eq:rec_w_n_m_p_r}
to make it valid for all $n$, $m$ and $r$. The next step is obtaining
a recurrence without $n$-dependent summations, for which we compute
$w_{[p]}(n,m,r)-w_{[p]}(n-2,m,r-2)$ using the extended
recurrence. Doing so gives
\begin{align}
  \label{eq:wp_n_m_p_r_rec2}
  w_{[p]}(n,m,r)=&~w_{[p]}(n-2,m,r-2)+w_{[p]}(n-1,m,r)-w_{[p]}(n-3,m,r-2)\nonumber\\
                                    &+w_{[p]}(n-2,m-\iv{p=1},r-1)+w_{[p]}(n-3,m-\iv{p=0},r-2)\nonumber\\
  &+\iv{m=0}\big(\iv{n=-1}\iv{r=0}-\iv{n=1}\iv{r=2}\big).
\end{align}
We can now determine $W_{[p]}(x,y,z)$ by multiplying both sides of~\eqref{eq:wp_n_m_p_r_rec2} by $x^ny^mz^r$
and then adding on $n,m,$ and $r$. This yields
\begin{equation}
  \label{eq:ogf_wp_n_m_p_r}
  W_{[p]}(x,y,z)=\frac{1-x^2z^2}{x\Big(1-x-x^2z\big(z+y^{\iv{p=1}}\big)+x^3z^2\big(1-y^{\iv{p=0}}\big)\Big)}.
\end{equation}
With the usual strategy, we can see that the coefficient of $y^m$
in~\eqref{eq:ogf_wp_n_m_p_r} is
\begin{equation}
  \label{eq:ogf_wp_n_m_p_r_ym}
  [y^m]W_{[p]}(x,y,z)=\frac{(1-x^2z^2)(x^{2+\iv{p=0}} z^{1+\iv{p=0}})^m}{x\Big(1-x-x^2z^2-\iv{p=0}\,x^2z\,(1-xz)\Big)^{m+1}}.
\end{equation}

\subsection{Probabilistic Connections}
\label{sec:prob-conn}
If the probability of drawing a one is $q$, then an $n$-string drawn
at random will roughly contain $nq$ ones for large $n$ with high
probability ---by the law of large numbers. Thus, for large~$n$, all
the probabilistic results in Section~\ref{sec:probability} can be
approximated through their deterministic counterparts in this section, by
assuming that the Hamming weight of an $n$-string drawn at random is
always~$\floor*{nq}$. With this assumption we have that, for large $n$,
\begin{equation}
  \label{eq:lln_connection}
  \pi_{\ushort{k}\le\widebar{k}}(n,m)\approx
  w_{\ushort{k}\le\widebar{k}}(n,m,\floor*{n q})\,{n\choose \floor*{n q}}^{-1},
\end{equation}
and
\begin{equation}
  \label{eq:lln_connection}
  \pi_{[p]}(n,m)\approx  w_{[p]}(n,m,\floor*{n q})\,{n\choose \floor*{n q}}^{-1}.
\end{equation}
Of course, these connections may also be used in reverse, that is to
say, to approximate Hamming-constrained enumerations using
probabilistic expressions.

\section{Number of $n$-Strings that Contain Prescribed Quantities of
  Nonnull Runs of Ones and/or Zeros Under Different Constraints}
\label{sec:oz}

In Section~\ref{sec:number-n-strings} we studied enumerations of
$n$-strings containing prescribed numbers of runs of ones in binary
strings. In this section we see how analogous techniques allow us to
enumerate the $n$-strings that contain prescribed quantities of runs
of ones \textit{and/or zeros} ---i.e., for the sake of clarity, the
case where each run may indistinctly be a run of ones or a run of
zeros--- in several similar scenarios. We focus on extending the main
enumerations in Sections~\ref{sec:w_n_m_klek} and~\ref{sec:wp_n_m} to
the same scenarios for runs of ones and zeros, but we also examine the
extension of the special enumerations in
Sections~\ref{sec:w_n_m_k},~\ref{sec:anylength}
and~\ref{sec:number-ws_n_k_ltk}, as these problems have been addressed
by other
authors~\cite{bloom98:_singles,goulden83:_combin_enumer,schilling90:_longest,suman94:_longest}. The
main change in our approach in this section with respect to previous
sections is that we always establish two mutual recurrences, rather
than a single one ---as we have always done up to this point. These
mutual recurrences correspond, respectively, to enumerations of
$n$-strings that start with a zero or with a one. Importantly, we only
consider nonnull runs of ones or zeros.

As mentioned in the introduction, a number of authors have considered
runs of ones and runs of zeros
jointly~\cite{stevens39:_distr,wald40:_test,feller68:_introd_probab,goulden83:_combin_enumer,schilling90:_longest,bloom96:_probab,bloom98:_singles,balakrishnan01:_runs_scans}.
On first impression, this setting is not so relevant to the study of
success runs, but, in contrast, it can lead to better statistics in
runs-based hypothesis tests (as more information is taken into account
when considering both kinds of runs). Importantly, it also has clear
direct implications for the problem of compositions (ordered
partitions) of $n$, which are the different ways in which we can
partition the sum $1 + 1 + \dots + 1 = n$ into nonzero ordered
parts. The intimate connection between runs of ones and zeros in
$n$-strings and compositions is easy to understand: we can visualise
the parts of a composition of $n$ as the ordered lengths of an
alternating series of runs of ones and zeros in an $n$-string, and in
its ones' complement. For example, consider the following composition
of $n=9$ into five parts: $1+3+2+2+1=9$. We can represent this
composition using the sequence of lengths of five alternating runs of
ones and zeros in two related $9$-strings: `011100110', and its ones'
complement, `100011001' ---see also Figure~\ref{fig:starsbars}, where
the top row can be interpreted as a partition of $1+1+\dots+1$ into
nonzero parts. Two other consequences of this observation are: 1) all
enumerations in this section must be even valued, and 2) the
aforementioned mutual recurrences are always symmetric. To the best of
our knowledge, the earliest author who saw a connection between
compositions and runs was Bateman~\cite{bateman48:_power}.

The relationship between runs \textit{of ones} and compositions is
perhaps less transparent. However, the results in
Section~\ref{sec:number-n-strings} also hint at a close association
between runs of ones and compositions, as shown by the many OEIS
sequences cited that are simultaneously related to both. Grimaldi and
Heubach~\cite{grimaldi05:_without_odd_runs} have in fact explicitly
described one case of this relationship.  The reason for these
connections is the intrinsic link that exists in certain scenarios
between the enumerations of $(n-1)$-strings containing prescribed
quantities of runs of ones and the enumerations of $n$-strings
containing prescribed quantities of runs of ones and/or zeros
---which, as we have discussed in the previous paragraph, are
themselves directly connected to compositions of $n$. We explicitly
give the simplest of such links in Theorems~\ref{thm:chave_das_nozes}
and~\ref{thm:chave_das_nozes_2} in Sections~\ref{sec:s_n_m_klek}
and~\ref{sec:s_n_m_p}, respectively.

Remark~\ref{rmk:compositions} at end of this section overviews the
main consequences of our results for the problem of
compositions.
 
\subsection[Number of n-Strings that Contain Exactly m Nonnull (ḵ≤k)-Runs of
Ones or Zeros]{Number of $n$-Strings that Contain Exactly $m$ Nonnull
  ($\ushort{k}\le\widebar{k}$)-Runs of Ones and/or Zeros}
\label{sec:s_n_m_klek}

We denote by $s_{\ushort{k}\le\widebar{k}}(n,m)$ the number of
$n$-strings that contain exactly $m$ nonnull ($\ushort{k}\le\widebar{k}$)-runs
of ones and/or zeros, which may also contain other runs with lengths
not in the prescribed range.  If we denote the number of compositions
of $n$ that contain exactly $m$ parts between $\ushort{k}$ and
$\widebar{k}$ by $c_{\ushort{k}\le \widebar{k}}(n,m)$, then, from the
observations in the previous section, we also have that
\begin{equation}
  \label{eq:compositions}
  c_{\ushort{k}\le \widebar{k}}(n,m)=\frac{1}{2}\,s_{\ushort{k}\le\widebar{k}}(n,m).
\end{equation}

As in Section~\ref{sec:w_n_m_klek}, we establish first a necessary
condition that guides our enumeration.

\begin{condition}{(Existence of $n$-strings
  containing $m$ nonnull $(\ushort{k}\le\widebar{k})$-runs of ones and/or zeros)}
  \begin{equation}
    \label{eq:m_necessary_condition_s}
    s_{\ushort{k}\le\widebar{k}}(n,m)> 0\quad \Longrightarrow\quad 0\le m \le \floor[\bigg]{\frac{n}{\ushort{k}}}.
  \end{equation}
\end{condition}
The upper bound in~\eqref{eq:m_necessary_condition_s} is due to the
fact that it must always hold that $m\,\ushort{k}\le n$ for an
$n$-string to be potentially able to hold $m$ nonnull
$(\ushort{k}\le\widebar{k})$-runs of ones and/or zeros.  Thus, unlike
in the enumerations of runs of ones, we now necessarily have that the
enumerations must be zero when $n=-1$, which implies that only the
case $n=0$ is involved in the initialisations.

\subsubsection{Recurrence}
\label{sec:s_rec}
We may obtain a recurrence for $s_{\ushort{k}\le\widebar{k}}(n,m)$ by
defining first two mutual recurrences for the number of $n$-strings
that start with bit $b$ and contain exactly $m$
$(\ushort{k}\le\widebar{k})$-runs of ones and/or  zeros, which
we denote by $s^b_{\ushort{k}\le\widebar{k}}(n,m)$ with
$b\in\{0,1\}$. Considering all the $n$-strings that start with an $i$-run
of each kind, we can see through the usual strategy ---see
Section~\ref{sec:recurrences_w_n_m_klek}--- that
\begin{align}
  \label{eq:s0}
  s^b_{\ushort{k}\le\widebar{k}}(n,m)=&~\sum_{i=\ushort{k}}^{\widebar{k}}s^{\tilde{b}}_{\ushort{k}\le\widebar{k}}(n-i,m-1)+\sum_{i=1}^{\ushort{k}-1}s^{\tilde{b}}_{\ushort{k}\le\widebar{k}}(n-i,m)+  \sum_{i=\widebar{k}+1}^{n}s^{\tilde{b}}_{\ushort{k}\le\widebar{k}}(n-i,m),
\end{align}
where $b\in\{0,1\}$ and $\tilde{b}=\text{mod}(b+1,2)$. Of course,
these recurrences closely
parallel~\eqref{eq:recurrence_w_n_m_klek_basic}. To initialise these
two mutual recurrences we can use the case $n=1$, in which we
known by inspection that
\begin{align}
  \label{eq:s0_n1_trivial}
  s^b_{\ushort{k}\le\widebar{k}}(1,m)=\iv{m=0}\iv{\ushort{k}>1}+\iv{m=1}\iv{\ushort{k}=1}.
\end{align}
Setting $n=1$ in~\eqref{eq:s0} we get
\begin{align}
  \label{eq:s0_n1}
  s^b_{\ushort{k}\le\widebar{k}}(1,m)=s^{\tilde{b}}_{\ushort{k}\le\widebar{k}}(0,m-\iv{\ushort{k}=1}).
\end{align}
In can be easily verified that~\eqref{eq:s0_n1_trivial}
and~\eqref{eq:s0_n1} are equal when
\begin{align}
  \label{eq:s0_init}
  s^b_{\ushort{k}\le\widebar{k}}(0,m)=\iv{m=0},
\end{align}
which thus initialises the two recurrences. Since the enumeration we
are interested in can be put as
\begin{equation}
  \label{eq:sum_mutual}
  s_{\ushort{k}\le\widebar{k}}(n,m)=  s^0_{\ushort{k}\le\widebar{k}}(n,m)+  s^1_{\ushort{k}\le\widebar{k}}(n,m),
\end{equation}
from this expression and \eqref{eq:s0} we have the recurrence
\begin{align}
  \label{eq:s_n_m_klek_rec}
  s_{\ushort{k}\le\widebar{k}}(n,m)=&~\sum_{i=\ushort{k}}^{\widebar{k}}s_{\ushort{k}\le\widebar{k}}(n-i,m-1)+\sum_{i=1}^{\ushort{k}-1}s_{\ushort{k}\le\widebar{k}}(n-i,m)+  \sum_{i=\widebar{k}+1}^{n}s_{\ushort{k}\le\widebar{k}}(n-i,m),
\end{align}
which from~\eqref{eq:sum_mutual} and~\eqref{eq:s0_init} is initialised by
\begin{align}
  \label{eq:s_n_m_klek_rec_init}
  s_{\ushort{k}\le\widebar{k}}(0,m)=2\,\iv{m=0}.
\end{align}

\subsubsection{Generating Function}
\label{sec:s_gf}
To get the generating function of $s_{\ushort{k}\le\widebar{k}}(n,m)$
we make first~\eqref{eq:s_n_m_klek_rec} valid for all $n$ and $m$. Since
setting $n=0$
in~\eqref{eq:s_n_m_klek_rec} gives zero, we achieve our goal by adding $2\iv{n=0}\iv{m=0}$ to it. We can now subtract $s_{\ushort{k}\le\widebar{k}}(n-1,m)$
from
$s_{\ushort{k}\le\widebar{k}}(n,m)$ using this extended recurrence to
get a recurrence free from $n$-dependent summations:
\begin{align}
  \label{eq:s_n_m_klek_rec_2}
  s_{\ushort{k}\le\widebar{k}}(n,m)=&~2\,s_{\ushort{k}\le\widebar{k}}(n-1,m)%
                                    +s_{\ushort{k}\le\widebar{k}}(n-\ushort{k},m-1)-s_{\ushort{k}\le\widebar{k}}(n-\ushort{k},m)\nonumber\\
                                    &-s_{\ushort{k}\le\widebar{k}}\big(n-(\widebar{k}+1),m-1\big)+s_{\ushort{k}\le\widebar{k}}\big(n-(\widebar{k}+1),m\big)\nonumber\\
  &+2\iv{m=0}\big(\iv{n=0}-\iv{n=1}\big).
\end{align}
We get next the ogf
$S_{\ushort{k}\le\widebar{k}}(x,y)=\sum_{n,m}s_{\ushort{k}\le\widebar{k}}(n,m)\,x^ny^m$
through the usual procedure of
multiplying~\eqref{eq:s_n_m_klek_rec_2} on both sides by $x^ny^m$ and
then summing over $n$ and $m$, which yields
\begin{equation}
  \label{eq:ogf_s_n_m_klek}
 S_{\ushort{k}\le\widebar{k}}(x,y)= \frac{2(1-x)}{1-2x+(1-y)(x^{\ushort{k}}-x^{\widebar{k}+1})}.
\end{equation}
 We can now state a simple but relevant theorem.
 \begin{theorem}[Fundamental link between the enumeration of $n$-strings
   containing prescribed quantities of
   $(\ushort{k}\le\widebar{k})$-runs of ones and/or zeros and its
   counterpart for runs of ones]\label{thm:chave_das_nozes}
  \begin{equation}
    \label{eq:ogfs_chave_das_nozes}
    S_{\ushort{k}\le\widebar{k}}(x,y)=2x\,W_{(\ushort{k}-1)\le(\widebar{k}-1)}(x,y).
  \end{equation}
  \textit{Proof.} Substitute~\eqref{eq:ogf_w_n_m_klek}
  into~\eqref{eq:ogfs_chave_das_nozes}.\hfill {\tiny$\square$}
\end{theorem}
We can also produce this theorem by comparing
recurrences~\eqref{eq:recurrence_w_n_m_klek_basic}
and~\eqref{eq:s_n_m_klek_rec} together with their initialisations, but
any possible hesitation about the relationship between these two
families of enumerations vanishes after comparing their respective
ogfs.  As far as we know, the earliest sign of
Theorem~\ref{thm:chave_das_nozes} appears in a result by Schilling
---see start of Section~\ref{sec:ss_n_ltk}

Theorem~\ref{thm:chave_das_nozes} is the main reason why we have
restricted our analysis in this section (Section~\ref{sec:oz}) to
nonnull runs of ones and zeros. Since null runs of ones are well
defined (see again Definition~\ref{def:nonnull} and
Remark~\ref{rmk:nonnull} in the introduction), results for
nonnull runs of ones and zeros can always be translated into results for
runs of ones through~\eqref{eq:ogfs_chave_das_nozes} when
$\ushort{k}\ge 1$. 

\begin{remark}\label{rmk:smirnov}
  A \textit{Smirnov word} is a string of symbols from an alphabet in
  which no two adjacent symbols are equal ---equivalently, a string without
  \textit{levels}~\cite[Secs.\
  2.4.13, 2.4.14]{goulden83:_combin_enumer}. The number of binary Smirnov words of length $n$ is
  trivially two for $n\ge 1$. We verify next that this can indeed be
  seen through~\eqref{eq:ogfs_chave_das_nozes}. Using this expression
  and~\eqref{eq:ogf_m_gtk} we can write
  \begin{equation}
    \label{eq:s_w_gk}
    [y^m]S_{\ge k}(x,y)=[y^m]\,2\,x\,W_{\ge (k-1)}(x,y)=2\,x^{mk}\frac{1-x}{\big(1-2x+x^k\big)^{m+1}}.
  \end{equation}
  So $[y^0]S_{\ge 2}(x,y)=2/(1-x)$, and thus $s_{\ge
    2}(n,0)=2$. 
\end{remark}

\subsection[Number of n-Strings that Contain Exactly m Nonnull k-Runs of Ones or
Zeros]{Number of $n$-Strings that Contain Exactly $m$ Nonnull $k$-Runs of Ones and/or Zeros}
\label{sec:s_n_m_k}
We denote by $s_k(n,m)$ the number of $n$-strings that contain exactly
$m$ nonnull $k$-runs of ones and/or zeros.  This enumeration is a
special case of $s_{\ushort{k}\le\widebar{k}}(n,m)$, because
\begin{equation}
  \label{eq:s_n_m_k_special_case}
  s_k(n,m)=s_{k\le k}(n,m).
\end{equation}
We consider this special enumeration separately because the case $k=1$
was previously studied by Bloom~\cite{bloom98:_singles}.

\subsubsection{Recurrences}
\label{sec:s_n_m_k_rec}
From~\eqref{eq:s_n_m_klek_rec} we have the recurrence
\begin{align}
  \label{eq:s_n_m_k_rec}
  s_{k}(n,m)=&~s_{k}(n-k,m-1)+\sum_{i=0\atop  i\neq k}^{n}s_{k}(n-i,m),
\end{align}
which, from~\eqref{eq:s_n_m_klek_rec_init}, is initialised by
$s_k(0,m)=2\,\iv{m=0}$. We can also
specialise~\eqref{eq:s_n_m_klek_rec_2} to get the following
alternative recurrence:
\begin{align}
  \label{eq:s_n_m_k_rec_2}
  s_{k}(n,m)=&~2\,s_{k}(n-1,m) +s_{k}(n-k,m-1)-s_{k}(n-k,m)
               -s_{k}(n-(k+1),m-1)\nonumber\\
  &+s_{k}(n-(k+1),m)+2\iv{m=0}\big(\iv{n=0}-\iv{n=1}\big),
\end{align}
which does not require initialisation as it is valid for all values of
$n$ and $m$ ---when taking~\eqref{eq:m_necessary_condition_s} into account.

\subsubsection{Generating Functions}
\label{sec:s_n_m_k_ogf}
Considering~\eqref{eq:s_n_m_k_special_case}, the coefficient of $y^m$ in
$S_k(x,y)=\sum_{n,m}s_k(n,m)\,x^ny^m$ can be obtained through
Theorem~\ref{thm:chave_das_nozes} and~\eqref{eq:ogf_m}:
\begin{equation}
  \label{eq:S_n_m_k_ogf_ym}
  [y^m]S_k(x,y)=2\,x^{m k}\bigg(\frac{1-x}{1-2x+x^{k}-x^{k+1}}\bigg)^{m+1}.
\end{equation}
This ogf enumerates the $n$-strings that contain exactly $m$ $k$-runs of
ones and/or zeros.

\subsubsection{Explicit Expression}
\label{sec:s_n_m_k_explicit}
From~\eqref{eq:ogfs_chave_das_nozes} we have
$s_k(n,m)=2\,w_{k-1}(n-1,m)$, which gives us a double-summation explicit
expression through~\eqref{eq:w_gen}. This expression does of course
work for $k=1$, where we have
\begin{equation}
  \label{eq:null_in_action}
  s_1(n,m)=2\,w_0(n-1,m).
\end{equation}
Observe that this mapping makes explicit use of an enumeration of null
runs of ones. In any case, we see in the next remark that a simpler
single-summation expression is possible when~$k=1$.

\begin{remark}\label{rmk:isolated_singles}
  Bloom used the term \textit{single} to mean a $1$-run of ones or 
  zeros ---a name reminiscent of Apostol's \textit{isolated singleton}
  for a $1$-run of ones~\cite{apostol88:_binary}--- and obtained a
  recurrence for the number of $n$-strings containing exactly $m$
  $1$-runs of ones and/or zeros~\cite{bloom98:_singles}. Setting $k=1$
  in~\eqref{eq:s_n_m_k_rec_2} we get
  \begin{align}
    \label{eq:rec_isolated singles}
    s_{1}(n,m)=&~s_{1}(n-1,m) +s_{1}(n-1,m-1) +s_{1}(n-2,m)-s_{1}(n-2,m-1)\nonumber\\
              &+2\,\iv{m=0}\big(\iv{n=0}-\iv{n=1}\big).
  \end{align}
  This is essentially Bloom's recurrence~\cite[Eq.\
  (1)]{bloom98:_singles} but in an even more general form,
  since~\eqref{eq:rec_isolated singles} is valid for all $n$ and $m$
  thanks to the inhomogeneous term ---absent in Bloom's
  expression. Due to this, Bloom's
  recursion does not work whenever both arguments of
  $s_1(\cdot,\cdot)$ equal zero in any of its five instances
  in~\eqref{eq:rec_isolated singles}, which the author deals with
  through the initialisation procedure.  Bloom also observes that
  $s_1(n,0)=2\,F_{n-1}^{(2)}$, which, for example, can also be seen by
  taking into account~\eqref{eq:ogfs_chave_das_nozes} and the comment
  about $w_0(n,0)$ in Section~\ref{sec:oeis-wknm}. Lastly, Bloom
  obtains, through counting arguments, a single-summation explicit
  expression for $s_1(n,m)$. We can find the same expression from the
  ogf~\eqref{eq:S_n_m_k_ogf_ym}, which in this case takes the form
  \begin{equation}
    \label{eq:s_n_m_k1_ogf}
    [y^m]S_1(x,y)=2\,x^m\Big(\frac{1-x}{1-x-x^2}\Big)^{m+1}.
  \end{equation}
  Applying the negative binomial theorem twice, we can express~\eqref{eq:s_n_m_k1_ogf} as
  \begin{equation}
    \label{eq:s_n_m_k1_ogf_bin}
    [y^m]S_1(x,y)=2\sum_{t\ge 0}\sum_{l \ge 0} {t+m\choose
      m}{l+t-1\choose t-1} x^{m+2t+l}.
  \end{equation}
  To get the coefficient of $x^n$ in~\eqref{eq:s_n_m_k1_ogf} we find
  the nonnegative indices $t$ and $l$ such that $n=m+2t+l$. From this relation,
  the maximum of $t$ happens when $l=0$, which implies that
  $t\le \lfloor(n-m)/2\rfloor$. For a given value of~$t$, we have that
  $l=n-m-2t$. Therefore, from~\eqref{eq:s_n_m_k1_ogf_bin} we have that
  \begin{equation}
    \label{eq:s_n_m_k1_explicit}
    s_1(n,m)=2\sum_{t=0}^{\floor{\frac{n-m}{2}}}{t+m\choose
      m}{n-m-t-1\choose t-1},
  \end{equation}
  which is the same as~\cite[Eq.\ (7)]{bloom98:_singles}, just
  noting that for $t=0$ the second binomial coefficient is zero
  ---which means that the summation may start at $t=1$.

  To conclude this remark, we define a rv modelling the number of
  $1$-runs of ones and/or zeros in an $n$-string drawn uniformly at
  random, which we call $\fixwidetilde{M}_{1,n}$. In standard
  probability notation,
  $\Pr(\fixwidetilde{M}_{1,n}=m)=s_1(n,m)/2^{n}$.  Bloom gave the
  expectation and variance of $\fixwidetilde{M}_{1,n}$ in~\cite[Eqs.\
  (3) and (4)]{bloom98:_singles}. From~\eqref{eq:null_in_action}, we
  get the first and second factorial moments of
  $\fixwidetilde{M}_{1,n}$ by using $n-1$, $k=0$, and $q=1/2$
  in~\eqref{eq:exp_mkn_all} and~\eqref{eq:2ndfact_mkn_cases},
  respectively ---see Section~\ref{sec:moments_w_k}. The expectation
  is
  \begin{equation}
    \label{eq:exp_tildem_1n}
    \text{E}\big(\fixwidetilde{M}_{1,n}\big)=\bigg(\frac{n+2}{4}\bigg)\,\iv{n>1}+\iv{n=1},
  \end{equation}
  whereas the second factorial moment is
  \begin{align}
    \label{eq:2ndfact_tildem_1n}
    \text{E}\big(\fixwidetilde{M}_{1,n}(\fixwidetilde{M}_{1,n}-1)\big)=\bigg(\frac{n(n+5)}{16}\bigg)\,\iv{n>2}+\frac{1}{2}\,\iv{n=2}.
  \end{align}
  Thus $\text{Var}(\fixwidetilde{M}_{1,n})=(5n+4)/16$ for $n> 2$, as
  Bloom proves by induction. Of course, \eqref{eq:exp_tildem_1n}
  and~\eqref{eq:2ndfact_tildem_1n} can be readily extended to
  $\fixwidetilde{M}_{k,n}$ in the same scenario. However, with unequal
  bit probabilities we must follow a different procedure ---see
  Section~\ref{sec:note-prob-extens}.
\end{remark}

\subsection[Number of n-Strings that Contain Exactly m Nonnull Runs of Ones or
Zeros]{Number of $n$-Strings that Contain Exactly $m$ Nonnull Runs of Ones and/or Zeros}
\label{sec:s_n_m}
We denote by $s(n,m)$ the number of $n$-strings that contain exactly
$m$ nonnull runs (of any lengths, all strictly greater than zero) of
ones and/or zeros. This enumeration is a special case of
$s_{\ushort{k}\le\widebar{k}}(n,m)$, as
\begin{equation}
  \label{eq:s_n_m_special_case}
  s(n,m)=s_{1\le n}(n,m).
\end{equation}
We consider this special enumeration separately because it was previously
studied by Goulden and Jackson~\cite[p.\
76]{goulden83:_combin_enumer}.

\subsubsection{Recurrences}
\label{sec:s_n_m_rec}
From~\eqref{eq:s_n_m_special_case} and~\eqref{eq:s_n_m_klek_rec} we have the recurrence
\begin{align}
  \label{eq:s_n_m_rec}
  s(n,m)=&~\sum_{i=1}^{n}s(n-i,m-1),
\end{align}
which, from~\eqref{eq:s_n_m_klek_rec_init}, is initialised by
$s(0,m)=2\,\iv{m=0}$. An alternative recurrence is found by specialising~\eqref{eq:s_n_m_klek_rec_2}:
\begin{align}
  \label{eq:s_n_m_rec_2}
  s(n,m)=&~s(n-1,m)+s(n-1,m-1)+2\iv{m=0}\big(\iv{n=0}-\iv{n=1}\big).
\end{align}
This recurrence does not require initialisation and it is valid for
all $m$ and $n$.

\subsubsection{Generating Function}
\label{sec:s_n_m_ogf}
The coefficient of $y^m$ in $S(x,y)=\sum_{n,m}s(n,m)\,x^ny^m$ can be
obtained by setting $k=1$ in~\eqref{eq:s_w_gk}:
  \begin{equation}
    \label{eq:s_n_m_ogf_ym}
    [y^m]S(x,y)=\frac{2x^m}{(1-x)^m}.
  \end{equation}

\subsubsection{Explicit Expression}
\label{sec:s_n_m_explicit}
To get $s(n,m)=[x^ny^m]S(x,y)$ we apply the negative binomial
theorem to~\eqref{eq:s_n_m_ogf_ym}, which allows us to express it as
\begin{equation}
  \label{eq:s_n_m_ogf_ym_bin}
  [y^m]S(x,y)=2\sum_{t\ge 0}{t+m-1\choose m-1} x^{t+m}.
\end{equation}
Thus, from $n=t+m$, the coefficient of $x^n$ in~\eqref{eq:s_n_m_ogf_ym_bin} is
\begin{equation}
  \label{eq:s_n_m_explicit}
  s(n,m)=2{n-1\choose m-1},
\end{equation}
which was given by Goulden and Jackson~\cite[p.\
76]{goulden83:_combin_enumer}.

\begin{remark}
  Expression~\eqref{eq:s_n_m_explicit} can alternatively be obtained
  through simple combinatorial reasoning.  Consider the partitioning
  of an $n$-string into $m$ nonempty substrings, where $1\le m\le
  n$. As illustrated in the ``stars and bars'' example in
  Figure~\ref{fig:starsbars}, in which the $n$-string is represented
  by $n$ asterisks and the $m$ partitions by $m+1$ vertical bars, we
  can put this partition into a bijection with $2m$ runs of ones and
  zeros ---corresponding to an $n$-string and to its ones' complement.
  The number of ways in which we can partition an $n$-string into
  $m$~nonempty substrings is ${n-1\choose m-1}$, from
  which~\eqref{eq:s_n_m_explicit} follows.

  \begin{figure}[t!]
  \begin{center}
    \begin{tabular}{@{}c@{}c@{}c@{}c@{}c@{}c@{}c@{}c@{}c@{}c@{}c@{}c@{}c@{}c@{}c@{}c@{}c@{}c@{}c@{}c@{}c@{}c@{}c@{}c@{}c@{}c@{}c@{}}
      $\mid$&$*$&$*$&$*$&$\mid$&$*$&$*$&$*$&$*$&$\mid$& $*$&$*$ &$\mid$& $*$&$*$&$*$&$*$&$*$ &$\mid$&$\cdots$ &$\mid$&$*$&$\mid$&$*$&$*$&$\mid$\\
      \hline
            &1&1&1&  &0&0&0&0 & & 1&1 & & 0&0&0&0&0 &&$\cdots$ && 0&& 1&1&\\
            &0&0&0&  &1&1&1&1 & & 0&0 & & 1&1&1&1&1 && $\cdots$ && 1&& 0&0&\\
    \end{tabular}
  \end{center}
  \caption{Partitioning of an $n$-string into nonempty substrings
    versus runs.}
  \label{fig:starsbars}
\end{figure}

  Finally, note that we can use~\eqref{eq:s_n_m_special_case} and Theorem~\ref{thm:chave_das_nozes} to write
  $s(n,m)=2\,w_{0\le n-1}(n-1,m)$. Thus,
  from~\eqref{eq:s_n_m_explicit}, the number of $n$-strings that
  contain exactly $m$ runs \textit{of ones}, including null runs, is
  ${n\choose m-1}$ ---cf.~\eqref{eq:w_n_m_explicit}, which does not
  count null runs.
\end{remark}

\subsection[Number of n-Strings Whose Longest Nonnull Run,
Is a (≤k)-Run of Ones or Zeros]{Number of $n$-Strings Whose Longest Nonnull
  Run Is a  $(\le k)$-Run of Ones or Zeros}
\label{sec:ss_n_ltk}

We denote by $\widebar{s}_{\le k}(n)$ the number of $n$-strings whose
longest nonnull run is a $(\le k)$-run of ones or zeros. Thus, we
assume $k\ge 1$. We can study this problem as a special case of the
results in Section~\ref{sec:s_n_m_klek}, because
\begin{equation}
  \label{eq:ss_n_ltk_s_n_m0_kp1}
  \widebar{s}_{\le k}(n)=s_{(k+1)\le n}(n,0).
\end{equation}
Again, we look at this special case in some detail because 
previous authors studied this problem. Schilling deduced through
simple reasoning that~\cite[Eq.\ (2)]{schilling90:_longest}
\begin{equation}
  \label{eq:schilling_chave_das_nozes}
  \widebar{s}_{\le k}(n)=2\,\widebar{w}_{\le (k-1)}(n-1),
\end{equation}
which is Theorem~\ref{thm:chave_das_nozes} in action in this special
case. Also, Suman studied the number of $b$-ary strings of length $n$
whose longest run is a $(\le\!k)$-run of any of the $b$-ary
symbols~\cite{suman94:_longest}, giving a recurrence, an ogf, and
three explicit expressions (one of them asymptotic). We only
consider here Suman's nonasymptotic results in the binary case, in
which his $r$ and~$k$ correspond in our setting to $k$ and $2$,
respectively.  In Suman's notation,
$\widebar{s}_{\le k}(n)=c_n$. Finally, Bloom also considered the
asymptotics of $\widebar{s}_{\le k}(n)$ in a particular
case~\cite{bloom98:_singles}.

\subsubsection{Recurrences}
\label{sec:ss_n_ltk_ec}
By applying~\eqref{eq:ss_n_ltk_s_n_m0_kp1} to
recurrence~\eqref{eq:s_n_m_klek_rec} we get
\begin{align}
  \label{eq:ss_n_ltk_rec}
  \widebar{s}_{\le k}(n)=&\sum_{i=1}^{k}\widebar{s}_{\le k}(n-i),
\end{align}
initialised by $\widebar{s}_{\le k}(0)=2$. This recurrence was given
by Suman~\cite[p.\ 121]{suman94:_longest}. Of
course,~\eqref{eq:ss_n_ltk_rec} is also
Schilling's recurrence~\eqref{eq:rec_ws_n_lek_schilling} transformed
by~\eqref{eq:schilling_chave_das_nozes}.

An alternative recurrence is obtained by
applying~\eqref{eq:ss_n_ltk_s_n_m0_kp1}
to~\eqref{eq:s_n_m_klek_rec_2}, which yields
\begin{align}
  \label{eq:ss_n_ltk_rec_2}
  \widebar{s}_{\le k}(n)=&~2\,\widebar{s}_{\le k}(n-1) -\widebar{s}_{\le k}\big(n-(k+2)\big) +2\big(\iv{n=0}-\iv{n=1}\big),
\end{align}
and does not require initialisation ---other than
taking~\eqref{eq:m_necessary_condition_s} into account. Bloom
discussed the asymptotics of $\widebar{s}_{\le 5}(n)$ using the
recurrence for runs of ones equivalent to~\eqref{eq:ss_n_ltk_rec_2}
---i.e., considering~\eqref{eq:schilling_chave_das_nozes},
recurrence~\eqref{eq:rec_ws_n_lek_simpler} with $k=4$, relying on its
characteristic equation~\cite[p.\ 126]{bloom98:_singles}.

\subsubsection{Generating Functions}
\label{sec:ss_n_m_ogf}
Through~\eqref{eq:ogfs_chave_das_nozes} and~\eqref{eq:ogf_ws_n_ltk2} we
have that the ogf $\widebar{S}_{\le
  k}(x)=\sum_n\widebar{s}_{\le k}(n)\,x^n$  is
\begin{equation}
  \label{eq:ogf_ss_ltk2}
  \widebar{S}_{\le k}(x)=\frac{2(1-x)}{1-2x+x^{k+1}}.
\end{equation}
Suman~\cite[p.\ 121]{suman94:_longest} gave the following alternative
equivalent ogf:
\begin{equation}
  \label{eq:ogf_ss_ltk2_suman}
  \widebar{S}_{\le k}(x)=\frac{1-x^{k+1}}{1-2x+x^{k+1}},
\end{equation}
which is, in Suman's notation, $A(x)$ in the binary case with some
straightforward algebraic simplification.

\subsubsection{Explicit Expression}
\label{sec:s_nexplicit}
Suman extracted a double-summation explicit expression for
$\widebar{s}_{\le k}(n)$ from~\eqref{eq:ogf_ss_ltk2_suman}~\cite[Thm.\
1]{suman94:_longest}, and he then simplified it into a single-summation
version~\cite[Thm.\ 2]{suman94:_longest}. Regarding our own results,
we can see that
using~\eqref{eq:schilling_chave_das_nozes},~\eqref{eq:ws_n_ltk_w_n_m_gtk_connection}
and~\eqref{eq:w_gtk_gen_diophant_m0} we can directly obtain a
single-summation explicit expression.

\subsection{Number of $n$-Strings that Contain Exactly $m$ Nonnull $p$-Parity
  Runs of Ones and/or Zeros}
\label{sec:s_n_m_p}

We denote by $s_{[p]}(n,m)$ the number of $n$-strings that contain
exactly $m$ nonnull $p$-parity runs of ones and/or zeros, i.e.,
exactly $m$ runs of ones and/or zeros whose lengths $k_1,\dots,k_m$
have parity $\text{mod}(k_i,2)=p$ for $i=1,\dots,m$ and are strictly
greater than zero. The $n$-strings that we enumerate may have more
than~$m$ nonnull runs of ones and/or zeros, as long as the lengths of
all these additional runs have parity opposite to $p$, and any number
of null runs (whose parity is even). For the same reason as in
Section~\ref{sec:s_n_m_klek}, if we call $c_{[p]}(n,m)$ the number of
compositions of $n$ that contain exactly $m$ ``$p$-parity parts''
---i.e., parts whose parity is $p$--- then we also have that
\begin{equation}
  \label{eq:compositions_p}
  c_{[p]}(n,m)=\frac{1}{2}\,s_{[p]}(n,m).
\end{equation}
As far as we are aware, Goulden and Jackson are the only authors who
have addressed enumeration problems related to the one considered in
this section. These authors have enumerated the $n$-strings in which
an odd run of zeros is never followed by an odd run of
ones~\cite[Sec.\ 2.4.6]{goulden83:_combin_enumer}, and the $n$-strings
in which all runs of ones are even and all runs of zeros are
odd~\cite[Ex.\ 2.4.3]{goulden83:_combin_enumer}. Clearly, the results
in this section do not encompass these two enumerations.

Before continuing, we state a necessary condition similar
to~\eqref{eq:m_necessary_condition_s}.
\begin{condition}{(Existence of $n$-strings containing $m$ nonnull $p$-parity
    runs of ones and/or zeros)}
  \begin{equation}
    \label{eq:m_necessary_condition_p_s}
    s_{[p]}(n,m)> 0\quad \Longrightarrow\quad 0\le m \le \floor[\bigg]{\frac{n}{1+\iv{p=0}}}.
\end{equation}
\end{condition}
This is really the same necessary condition
as~\eqref{eq:m_necessary_condition_s}, just noting that the lengths of
all nonnull $p$-parity runs of ones and/or zeros are lower bounded by
$\ushort{k}=1+\iv{p=0}$.

\subsubsection{Recurrence}
\label{sec:recurrences_sp_n_m}
We can produce a recurrence to enumerate $s_{[p]}(n,m)$ with the same
strategy used in Section~\ref{sec:s_n_m_klek}: we first produce two
mutual recurrences for the number of $n$-strings that start with bit
$b$ and contain exactly $m$ $p$-parity runs of ones and/or zeros,
denoted by $s^b_{[p]}(n,m)$, where $b\in\{0,1\}$. Considering all the
$n$-strings that start with a nonnull $i$-run of each length, we can
readily see through a similar reasoning as in
Section~\ref{sec:recurrences_wp_n_m} that the mutual recurrences
sought are given by
\begin{align}
  \label{eq:s0_p}
  s^b_{[p]}(n,m)%
  &=\sum_{i=1}^{n}s^{\tilde{b}}_{[p]}(n-i,m-\iv{p=\text{mod}(i,2)}),
\end{align}
for $b\in\{0,1\}$, and where $\tilde{b}=\text{mod}(b+1,2)$.  To find
initial values for the mutual recurrences in~\eqref{eq:s0_p} we use the
case $n=1$, in which we know by inspection that
\begin{equation}
  \label{eq:trivial_n1_sp}
  s^b_{[p]}(1,m)=\iv{m=0}\iv{p=0}+\iv{m=1}\iv{p=1}.
\end{equation}
On the other hand, setting~$n=1$ in~\eqref{eq:s0_p} yields
\begin{equation}
  \label{eq:s0_p_n1}
  s^b_{[p]}(1,m)=s^{\tilde{b}}_{[p]}(0,m-\iv{p=1}).
\end{equation}
We may verify that~\eqref{eq:s0_p_n1} equals~\eqref{eq:trivial_n1_sp} when
\begin{equation}
  \label{eq:spb_n_m_init}
  s^b_{[p]}(0,m)=\iv{m=0}.
\end{equation}
which are therefore the initial values of the mutual
recurrences in~\eqref{eq:s0_p}. Since the enumeration we are interested in can be written as
\begin{equation}
  \label{eq:sum_mutual_p}
  s_{[p]}(n,m)=  s^0_{[p]}(n,m)+  s^1_{[p]}(n,m),
\end{equation}
from this expression and \eqref{eq:s0_p} we have the recurrence
\begin{align}
  \label{eq:s_n_m_p_rec}
  s_{[p]}(n,m)=\sum_{i=1}^{n}s_{[p]}(n-i,m-\iv{p=\text{mod}(i,2)}),
\end{align}
which from~\eqref{eq:sum_mutual_p} and~\eqref{eq:spb_n_m_init} is initialised by
\begin{align}
  \label{eq:s_n_m_p_init}
  s_{[p]}(0,m)=2\iv{m=0}.
\end{align}

\subsubsection{Generating Functions}
\label{sec:generating-function-sp}
We now obtain the ogf
$S_{[p]}(x,y)=\sum_{n,m} s_{[p]}(n,m)\, x^ny^m$. First of all, we make
$s_{[p]}(n,m)$ valid for all values of $n$ and $m$. Setting $n=0$
in~\eqref{eq:s_n_m_p_rec} we get $s_{[p]}(0,m)=0$ instead of the correct value
$2\iv{m=0}$. Thus we just need to add $2\iv{m=0}\iv{n=0}$
to~\eqref{eq:s_n_m_p_rec} to obtain an extended recurrence valid for
all $n$ and $m$, which we use subsequently.

Next, we obtain a recurrence without an $n$-dependent summation by
determining $s_{[p]}(n,m)-s_{[p]}(n-2,m)$ using the extended
recurrence. This yields
\begin{align}
  \label{eq:sp_n_m_rec_diff}
  s_{[p]}(n,m)=&~s_{[p]}(n-2,m)+s_{[p]}(n-1,m-\iv{p=1})+s_{[p]}(n-2,m-\iv{p=0})\nonumber\\&+2\iv{m=0}\big(\iv{n=0}-\iv{n=2}\big).
\end{align}
By multiplying next~\eqref{eq:sp_n_m_rec_diff} on both sides by
$x^ny^m$ and then adding on $n$ and~$m$ we find that
\begin{equation}
  \label{eq:sp_bivariate_ogf}
  S_{[p]}(x,y)=\frac{2\,(1-x^2)}{1-x\, y^{\iv{p=1}}-x^2\big(1+y^{\iv{p=0}}\big)}.
\end{equation}
In general, this ogf cannot always be connected to its
counterpart~\eqref{eq:wp_bivariate_ogf} that counts only runs of ones,
unlike $S_{\ushort{k}\le\widebar{k}}(x,y)$ in the previous section
---see Theorem~\ref{thm:chave_das_nozes}. But not all hope is lost, as
the next simple theorem shows.
\begin{theorem}{Fundamental link between enumerations of $n$-strings
    that contain prescribed quantities of nonnull even runs of ones
    and/or zeros and enumerations of $n$-strings containing prescribed
    quantities of odd runs of ones.}\label{thm:chave_das_nozes_2}
  \begin{equation}
    \label{eq:chave_das_nozes_2}
    S_{[0]}(x,y)=2x\, W_{[1]}(x,y).
  \end{equation}
  \textit{Proof.} Substitute~\eqref{eq:wp_bivariate_ogf}
  into~\eqref{eq:chave_das_nozes_2}.\hfill {\tiny$\square$}
\end{theorem}
\begin{remark}
  There is no similar converse mapping between $S_{[1]}(x,y)$ and
  $W_{[0]}(x,y)$, as one might expect. This asymmetry is rather
  unsatisfying, but it has a satisfying explanation: we have only
  considered \textit{nonnull} $p$-parity runs of ones to get
  $W_{[p]}(x,y)$. This assumption is immaterial when $p=1$, as odd
  runs of ones are always nonnull. However, when $p=0$ we are
  disregarding all null runs of ones, which have even parity. Thus,
  the symmetry of the setting with respect to runs of ones and zeros
  is broken. We will not delve further into this issue here, but it is
  not difficult to see that if we modify $W_{[p]}(x,y)$ to include
  null runs of ones, then symmetry is restored and
  \eqref{eq:chave_das_nozes_2} becomes
  $S_{[p]}(x,y)=2x\, W_{[\tilde{p}]}(x,y)$ with $p\in\{0,1\}$ and
  $\tilde{p}=\text{mod}(p+1,2)$. This is yet another example of the
  relevant role played by null runs.
\end{remark}
Extracting the coefficient of $y^m$ from~\eqref{eq:sp_bivariate_ogf}
using the procedure that we have repeatedly used throughout the paper
yields
\begin{equation}
  \label{eq:sp_bivariate_ogf_ym}
  [y^m]S_{[p]}(x,y)=\frac{2\,(1-x^2)\,x^{(1+\iv{p=0})\,m}}{\Big(1-\iv{p=0}x-\big(1+\iv{p=1}\big)x^2\Big)^{m+1}}.
\end{equation}
This is the ogf enumerating the binary strings that contain exactly
$m$ nonnull $p$-parity runs of ones and/or zeros.

\begin{remark}\label{rmk:compositions}
  In this remark we briefly discuss the consequences of our results in
  Section~\ref{sec:oz} for compositions. We start with the results in
  Section~\ref{sec:s_n_m_klek}. Putting together
  Theorem~\ref{thm:chave_das_nozes} and the
  connection~\eqref{eq:compositions} between runs of ones and zeros
  and compositions, we find that the ogf
  $C_{\ushort{k}\le\widebar{k}}(x,y)=\sum_{n,m}c_{\ushort{k}\le\widebar{k}}(n,m)\,x^ny^m$
  is
  \begin{equation}
    \label{eq:comp_runs}
    C_{\ushort{k}\le\widebar{k}}(x,y)=\frac{1}{2}S_{\ushort{k}\le\widebar{k}}(x,y)=x\,W_{\ushort{k}-1\le\widebar{k}-1}(x,y).
  \end{equation}
  Thus, not only all results in
  Sections~\ref{sec:w_n_m_klek}--\ref{sec:number-wl_n_k_gtk} ---which,
  we remind, only concern runs of ones--- can be directly applied to
  the corresponding enumerations of $n$-strings with prescribed
  quantities of nonnull runs of ones and/or zeros, but they can also
  be directly applied to the corresponding enumerations of
  compositions of $n$ with prescribed quantities of parts. We
  give next some examples that recover known results due to other
  authors. We mainly focus on generating functions, but of course
  recurrences and explicit results also follow from the connections
  given.

  \begin{itemize}
  \item From~\eqref{eq:comp_runs}, by dividing~\eqref{eq:s_n_m_ogf_ym}
    and~\eqref{eq:s_n_m_explicit} by $2$ we get the ogf and the
    explicit expression for compositions with exactly $m$ parts (of
    any sizes), i.e., $[y^n]C(x,y)=x^m/(1-x)^m$ and
    $c(n,m)={n-1\choose m-1}$, both of which were given by
    Riordan~\cite[p.\ 124]{riordan58:_comb}.  Goulden and
    Jackson~\cite[p.\ 53]{goulden83:_combin_enumer} gave as well
    $c(n,m)$ and $C(x,y)=(1-x)/(1-(1+y)x)$, which we can also obtain
    using $\ushort{k}=1$ from~\eqref{eq:comp_runs}
    and~\eqref{eq:ogf_gtk}. Observe that Goulden and Jackson also
    found~\eqref{eq:s_n_m_explicit}, but apparently they did not make
    the connection between runs (in their nomenclature, maximal
    blocks) and compositions.

\item Using~\eqref{eq:comp_runs} and~\eqref{eq:ogf_ws_lek_2}, the ogf for
  the number of compositions with no part greater than~$k$ is
  \begin{equation}
    \label{eq:comp_none_gtk}
    \widebar{C}_{\le k}(x)=x\, \widebar{W}_{\!\!\le (k-1)}(x)=\frac{x-x^{k+1}}{1-2x+x^{k+1}},
  \end{equation}   %
  which was also given by Riordan~\cite[Eq.\
  (40)]{riordan58:_comb}. An alternative ogf
  for~\eqref{eq:comp_none_gtk} was also given by Heubach and
  Mansour~\cite[Ex.\ 2.7]{heubach04:_compositions}.

\item From~\eqref{eq:comp_runs} and~\eqref{eq:ogf_wsk}, the number of
  compositions with no part greater than $k$ and at least one part
  equal to $k$ has ogf
  \begin{equation}
    \label{eq:comp_greatest_k}
    \widebar{C}_{k}(x)=x\, \widebar{W}_{(k-1)}(x)=\frac{x^k(1-x)^{2}}{(1-2x+x^{k})(1-2x+x^{k+1})},
  \end{equation}   %
  which was given by Riordan as well~\cite[p.\ 155]{riordan58:_comb}.

  Expression~\eqref{eq:comp_greatest_k} can be read as a bijection
  between the compositions of $n+1$ with largest part equal to $k$ and
  the $n$-strings with longest run (of ones) a $(k-1)$-run.  Nej and
  Reddy~\cite{nej19:_binary} say that there is a bijection between the
  compositions of $n+1$ with largest part $k$ and the $n$-strings with
  Hamming weight $r$ whose longest run is a $(k-1)$-run, but,
  given~\eqref{eq:comp_greatest_k}, this assertion requires
  considering all possible Hamming weights $0\le r \le n$.
  
\item Finally, from~\eqref{eq:comp_runs} and~\eqref{eq:ogf_m_gtk} we
  have that the number of compositions with no part equal to $k$ has
  ogf
  \begin{equation}
    \label{eq:comp_no_k}
    [y^0]C_k(x,y)=[y^0]x\,W_{(k-1)}(x,y)=\frac{1-x}{1-2x+x^k-x^{k+1}}.
  \end{equation}
  This ogf was given by Chinn and Heubach~\cite[Thm.\
  1]{chinn03:compositions}. The recurrences given by these authors for
  this enumeration~\cite[Thm.\ 1]{chinn03:compositions} can also be
  found by using $c_k(n,m)=w_{k-1}(n-1,m)$ with $m=0$
  in~\eqref{eq:recurrence_w_n_m_k_basic} or
  in~\eqref{eq:recurrence_w_n_m_k_simpler}.`
\end{itemize}

We consider subsequently the consequences for compositions of the
results in Section~\ref{sec:s_n_m_p}. Letting
$C_{[p]}(x,y)=\sum_{n,m} c_{[p]}(n,m)\,x^n y^m$,
from~\eqref{eq:compositions_p} we have that
Theorem~\ref{thm:chave_das_nozes_2} carries over to the enumeration of
compositions having prescribed quantities of even parts as follows:
  \begin{equation}
    \label{eq:even_comp_even_runs}
    C_{[0]}(x,y)=\frac{1}{2}S_{[0]}(x,y)=x\,W_{[1]}(x,y).
  \end{equation}
  Thus, from these equalities and~\eqref{eq:wp_bivariate_ogf} we have
  \begin{equation}
    \label{eq:comp_m_even_parts}
    C_{[0]}(x,y)=\frac{1-x^2}{1-x-x^2(1+y)},
  \end{equation}
  which was previously given by Goulden and Jackson~\cite[p.\
  54]{goulden83:_combin_enumer}.

  Grimaldi and Heubach identified the bijection between the
  compositions of $n+1$ with only odd parts and the $n$-strings
  without odd runs of zeros~\cite[Sec.\
  5]{grimaldi05:_without_odd_runs}. Now,
  from~\eqref{eq:even_comp_even_runs} we can write
  \begin{align}
    \label{eq:even_comp_runs_ones_m}
    [y^m]C_{[0]}(x,y)=[y^m]x\,W_{[1]}(x,y),
  \end{align}
  which allows us to make a stronger statement: there is a bijection
  between the compositions of $n+1$ with exactly $m$ even parts and
  the $n$-strings that contain exactly~$m$ odd runs of ones (or,
  indeed, exactly $m$ odd runs of zeros), of which Grimaldi and
  Heubach's observation is the special case $m=0$. Finally, in this
  case we have from~\eqref{eq:even_comp_even_runs}
  and~\eqref{eq:sp_bivariate_ogf_ym} that
  \begin{equation}
    \label{eq:comp_no_even_parts}
    [y^0]C_{[0]}(x,y)=\frac{1-x^2}{1-x-x^2},
  \end{equation}
  i.e., the ogf enumerating compositions with only odd parts, which
  was given by Heubach and Mansour~\cite[Ex.\
  2.9]{heubach04:_compositions}.

  Last but not least, Heubach and Mansour obtained many of their
  results for compositions as a specialisation of a general theorem
  of theirs that gives the ogf of the number of compositions with $m$
  parts in a general subset $A\subseteq \mathbb{N}$~\cite[Thm.\
  2.1]{heubach04:_compositions}. Through
  Theorems~\ref{thm:chave_das_nozes} and~\ref{thm:chave_das_nozes_2},
  most ogfs for runs in the current paper can in principle be recovered
  by relying on Heubach and Mansour's theorem, and also new results
  for runs be obtained.

  Some further comments about compositions and runs are given in
  Remark~\ref{rmk:more_compositions}.
\end{remark}

  \subsubsection{OEIS}
\label{sec:sp_oeis}

The sequences connected to $s_{[p]}(n,m)$ in the OEIS are:
\begin{itemize}
\item $m=0$
  
  $s_{[1]}(n,0)$ is \seqnum{A077957} for $n\ge 1$. %
  
  $s_{[0]}(n-1,0)$ is \seqnum{A006355} for $n\ge 2$. %
  
\item $m=1$
  
  $s_{[1]}(2n-3,1)$ is \seqnum{A057711} for $n\ge 2$. %

  $s_{[0]}(n+1,1)$ is \seqnum{A320947} for $n\ge 2$. %

\end{itemize}

\section{Extensions to Probabilistic Runs of Ones and/or Zeros}
\label{sec:note-prob-extens}

Just as in Section~\ref{sec:probability}, where we discussed the
probabilistic extension of the enumerative results in
Section~\ref{sec:number-n-strings}, we can extend the enumerative
results in Section~\ref{sec:oz} to the same probabilistic scenario:
the case where the $n$-strings are outcomes from $n$ iid Bernoulli
random variables with parameter $q$, where $0<q<1$ is the probability
of drawing a~$1$. As it may be intuited, in this scenario there is no
easy mapping between the results for runs of ones and the results for
runs of ones and zeros, i.e., there are no probabilistic analogues of
Theorems~\ref{thm:chave_das_nozes} and~\ref{thm:chave_das_nozes_2},
and the expressions can get rather involved.

We mainly deal in this section with the extension of the results in
Section~\ref{sec:s_n_m_klek} to the probabilistic scenario. We also
examine the extension of the special case in Section~\ref{sec:s_n_m},
as this was previously analysed in detail by Wishart and
Hirschfeld~\cite{wishart36:_theorem}.  The same approach can obviously
be used to extend the results in Section~\ref{sec:s_n_m_p}, although
we omit it here.

Like in Section~\ref{sec:probability}, it is understood that all
expressions in this section implicitly depend on $q$.

\subsection[Probability that an n-String Contains Exactly m
Nonnull (ḵ≤k)-Runs of Ones and/or Zeros]{Probability that an $n$-String
  Contains Exactly $m$ Nonnull ($\ushort{k}\le\widebar{k}$)-Runs of Ones
  and/or Zeros}
\label{sec:probability-s_m_klek}
We call $\lambda_{\ushort{k}\le\widebar{k}}(n,m)$ the probability that
an $n$-string contains exactly $m$ nonnull
($\ushort{k}\le\widebar{k}$)-runs of ones and/or zeros.  The obvious
case is $q=1/2$, in which
$\lambda_{\ushort{k}\le\widebar{k}}(n,m)=s_{\ushort{k}\le\widebar{k}}(n,m)/2^n$. This
was in fact the approach followed in special cases of this expression
by Bloom~\cite{bloom98:_singles} and by Suman~\cite{suman94:_longest},
who produced $\lambda_1(n,m)=s_1(n,m)/2^n$ and
$\widebar{\lambda}_{\le k}(n)=\widebar{s}_{\le k}(n)/2^n$,
respectively, using explicit expressions ---Suman's results actually
apply to the longest run of any $b$-ary symbol in $b$-ary strings of
length $n$.  We examine next how to get
$\lambda_{\ushort{k}\le\widebar{k}}(n,m)$ for
arbitrary~$q$.

\subsubsection{Recurrence}
\label{sec:lambda_recurrence}
As per our notation conventions, we denote the joint probability that
an $n$-string starts with~$b$ and contains exactly $m$ nonnull
$(\ushort{k}\le\widebar{k})$-runs of ones and/or zeros by
$\lambda^b_{\ushort{k}\le\widebar{k}}(n,m)$. By the law of total
probabilities, the probability that we wish to determine can be
expressed in terms of these joint probabilities as
\begin{equation}
  \label{eq:lambda_klek}
  \lambda_{\ushort{k}\le\widebar{k}}(n,m)=\lambda^0_{\ushort{k}\le\widebar{k}}(n,m)+\lambda^1_{\ushort{k}\le\widebar{k}}(n,m).
\end{equation}
Letting
\begin{equation}
  \label{eq:qb}
  q_b=(1-q)\iv{b=0}+q\iv{b=1}
\end{equation}
and invoking again the law of total probabilities, we have that the
probabilistic version of the mutual recurrences in~\eqref{eq:s0} is
\begin{align}
  \label{eq:lambda_b_klek}
  \lambda^b_{\ushort{k}\le\widebar{k}}(n,m)=&~\sum_{i=\ushort{k}}^{\widebar{k}}q_b^i\,\lambda^{\tilde{b}}_{\ushort{k}\le\widebar{k}}(n-i,m-1)+\sum_{i=1}^{\ushort{k}-1}q_b^i\,\lambda^{\tilde{b}}_{\ushort{k}\le\widebar{k}}(n-i,m)+  \sum_{i=\widebar{k}+1}^{n}q_b^i\,\lambda^{\tilde{b}}_{\ushort{k}\le\widebar{k}}(n-i,m),
\end{align}
where $b\in\{0,1\}$ and $\tilde{b}=\text{mod}(b+1,2)$. By inspection, the value of $\lambda^b_{\ushort{k}\le\widebar{k}}(n,m)$ when $n=1$ is
\begin{align}
  \label{eq:lambda0_n1_trivial}
  \lambda^b_{\ushort{k}\le\widebar{k}}(1,m)=q_b\,\big(\iv{m=0}\iv{\ushort{k}>1}+\iv{m=1}\iv{\ushort{k}=1}\big).
\end{align}
On the other hand, setting
$n=1$ in~\eqref{eq:lambda_b_klek} gives
\begin{align}
  \label{eq:lambda0_n1}
  \lambda^b_{\ushort{k}\le\widebar{k}}(1,m)=q_b\,\lambda^{\tilde{b}}_{\ushort{k}\le\widebar{k}}(0,m-\iv{\ushort{k}=1}),
\end{align}
and we can see that~\eqref{eq:lambda0_n1} equals~\eqref{eq:lambda0_n1_trivial}
for
\begin{align}
  \label{eq:lambda0_init}
  \lambda^b_{\ushort{k}\le\widebar{k}}(0,m)=\iv{m=0}.
\end{align}
This is therefore the initialisation
of~\eqref{eq:lambda_b_klek}. Although we can get
$\lambda_{\ushort{k}\le\widebar{k}}(n,m)$ using~\eqref{eq:lambda_klek}
and~\eqref{eq:lambda_b_klek}, the asymmetry of the mutual recurrences
in~\eqref{eq:lambda_b_klek} prevents us from obtaining a recurrence
for $\lambda_{\ushort{k}\le\widebar{k}}(n,m)$ itself. This theme
---asymmetry between mutual recurrences--- will resurface in all the
enumerations in Section~\ref{sec:oz_hamming}, which can be, in fact,
related asymptotically to the results in this section through the law
of large numbers.

\begin{remark}
  Schilling mentioned that a recurrence for
  $\widebar{\lambda}_{\le k}(n)=\lambda_{(k+1)\le n}(n,0)$ can be
  obtained~\cite[p.\ 200]{schilling90:_longest} in a similar way
  as~\eqref{eq:ps_n_ltk_schilling}. While he did not furnish the
  recurrence, this author indicated that the approximation
  $\widebar{\lambda}_{\le k}(n)\approx \widebar{\pi}_{\le k}(n)$
  ---see Section~\ref{sec:pis_n_ltk}--- works well for very large $n$
  when using~$\check{q}=\max(q,1-q)$ to compute
  $\widebar{\pi}_{\le k}(n)$.  Empirically, for fixed~$n$, the
  accuracy of this approximation increases
  as~$k$~increases. 
\end{remark}

\subsubsection{Probability Generating Function}
\label{sec:lambda_prob-gener-funct}
In spite of not having a recurrence for
$\lambda_{\ushort{k}\le\widebar{k}}(n,m)$, it is still possible to get
the pgf
$\Lambda^b_{\ushort{k}\le\widebar{k}}(x,y)=\sum_{n,m}\lambda^b_{\ushort{k}\le\widebar{k}}(n,m)
\,x^n y^m$, and then, through~\eqref{eq:lambda_klek}, obtain
$\Lambda_{\ushort{k}\le\widebar{k}}(x,y)=\sum_{n,m}\lambda_{\ushort{k}\le\widebar{k}}(n,m)
\,x^n y^m$ using
\begin{equation}
  \label{eq:Lambda_klek}
  \Lambda_{\ushort{k}\le\widebar{k}}(x,y)=\Lambda^0_{\ushort{k}\le\widebar{k}}(x,y)+\Lambda^1_{\ushort{k}\le\widebar{k}}(x,y).
\end{equation}

As usual, we first make~\eqref{eq:lambda_b_klek} valid for all values of
$n$ and $m$. Setting $n=0$ in~\eqref{eq:lambda_b_klek} gives
$\lambda^b_{\ushort{k}\le\widebar{k}}(0,m)=0$ instead
of~\eqref{eq:lambda0_init}, so we just need to add
$\iv{n=0}\iv{m=0}$ to~\eqref{eq:lambda_b_klek} to achieve
our goal. Using this extended recurrence, we compute next
$\lambda^b_{\ushort{k}\le\widebar{k}}(n,m)-q_b\,\lambda^b_{\ushort{k}\le\widebar{k}}(n-1,m)$
to produce a mutual recurrence free from $n$-dependent summations. This yields
\begin{align}
  \label{eq:lambda_b_rec_free}
  \lambda^b_{\ushort{k}\le\widebar{k}}(n,m)-q_b\,\lambda^b_{\ushort{k}\le\widebar{k}}(n-1,m)=&~{q_b}^{\ushort{k}}\,\lambda^{\tilde{b}}_{\ushort{k}\le\widebar{k}}(n-\ushort{k},m-1)-{q_b}^{\widebar{k}+1}\,\lambda^{\tilde{b}}_{\ushort{k}\le\widebar{k}}(n-(\widebar{k}+1),m-1)\nonumber\\
                                                                                       &+q_b\,\lambda^{\tilde{b}}_{\ushort{k}\le\widebar{k}}(n-1,m)-{q_b}^{\ushort{k}}\,\lambda^{\tilde{b}}_{\ushort{k}\le\widebar{k}}(n-\ushort{k}),m)\nonumber\\
                                                                                      &+{q_b}^{\widebar{k}+1}\,\lambda^{\tilde{b}}_{\ushort{k}\le\widebar{k}}(n-(\widebar{k}+1),m)\nonumber\\
  &+\big(\iv{n=0}-q_b\,\iv{n=1}\big)\iv{m=0}.
\end{align}
Multiplying now~\eqref{eq:lambda_b_rec_free} on both sides by $x^ny^m$ and
adding over $n$ and $m$ we obtain
\begin{equation}
  \label{eq:ogf_lambda_b}
  \Lambda^{b}_{\ushort{k}\le\widebar{k}}(x,y)(1-q_b\,x)=\Lambda^{\tilde{b}}_{\ushort{k}\le\widebar{k}}(x,y)\Big((y-1)\big({q_b}^{\ushort{k}}x^{\ushort{k}}-{q_b}^{\widebar{k}+1}x^{\widebar{k}+1}\big)+q_b\,x\Big)+1-q_b\,x.
\end{equation}
This is a system of two equations with two unknowns, i.e., the pgfs
$\Lambda^{0}_{\ushort{k}\le\widebar{k}}(x,y)$ and
$\Lambda^{1}_{\ushort{k}\le\widebar{k}}(x,y)$. In order to streamline
the upcoming expressions we now define
\begin{equation}
  \label{eq:xb}
  x_b=q_b\,x.
\end{equation}
Solving the system in~\eqref{eq:ogf_lambda_b} for
$\Lambda^{1}_{\ushort{k}\le\widebar{k}}(x,y)$ we get
\begin{equation}
  \label{eq:ogf_lambda_1}
  \Lambda^{1}_{\ushort{k}\le\widebar{k}}(x,y)=\frac{(1-x_1)(1-x_0)+\Big((y-1)\big(x_1^{\ushort{k}}-x_1^{\widebar{k}+1}\big)+x_1\Big)(1-x_0)}{(1-x_1)(1-x_0)-\Big((y-1)\big(x_1^{\ushort{k}}-x_1^{\widebar{k}+1}\big)+x_1\Big)\Big((y-1)\big(x_0^{\ushort{k}}-x_0^{\widebar{k}+1}\big)+x_0\Big)}.
\end{equation}
Inputting~\eqref{eq:ogf_lambda_1} in~\eqref{eq:ogf_lambda_b} we get
$\Lambda^{0}_{\ushort{k}\le\widebar{k}}(x,y)$, and 
using then~\eqref{eq:Lambda_klek} we arrive after some algebra
at
\begin{equation}
  \label{eq:ogf_lambda_klek}
  \Lambda_{\ushort{k}\le\widebar{k}}(x,y)=\frac{2-x_1-x_0+(y-1)\Big(\big(x_1^{\ushort{k}}-x_1^{\widebar{k}+1}\big)(1-x_0)+\big(x_0^{\ushort{k}}-x_0^{\widebar{k}+1}\big)(1-x_1)\Big)}{(1-x_1)(1-x_0)-\Big((y-1)\big(x_1^{\ushort{k}}-x_1^{\widebar{k}+1}\big)+x_1\Big)\Big((y-1)\big(x_0^{\ushort{k}}-x_0^{\widebar{k}+1}\big)+x_0\Big)}.
\end{equation}
This pgf allows us to obtain
$\lambda_{\ushort{k}\le\widebar{k}}(n,m)=[x^ny^m]\Lambda_{\ushort{k}\le\widebar{k}}(x,y)$
using a computer algebra system, in a more efficient manner than
through recurrences~\eqref{eq:lambda_b_klek}
or~\eqref{eq:lambda_b_rec_free}. We must surely give up hope of
finding a reasonably simple closed-form general expression here, but
this is not the case in special scenarios as we see next.

\subsection[Probability that an n-String Contains Exactly m
Nonnull Runs of Ones and/or Zeros]{Probability that an $n$-String
  Contains Exactly $m$ Nonnull Runs of Ones
  and/or Zeros}
\label{sec:probability-s_n_m}

Let $\lambda(n,m)$ be the probability that an $n$-string contains
exactly $m$ nonnull runs of ones and/or zeros. This is a special case
of the probability studied in the previous section with $\ushort{k}=1$
and $\widebar{k}=n$, and thus $\lambda(n,m)=\lambda_{1\le
  n}(n,m)$. Wishart and Hirschfeld previously studied this problem in
some detail~\cite{wishart36:_theorem}.
  
\subsubsection{Recurrences}
Letting $\lambda^b(n,m)=\lambda^b_{1\le n}(n,m)$, the mutual
recurrences~\eqref{eq:lambda_b_klek} become
\begin{align}
  \label{eq:lambda_b_n_m}
  \lambda^b(n,m)=&~\sum_{i=1}^{n}q_b^i\,\lambda^{\tilde{b}}(n-i,m-1),
\end{align}
which, from~\eqref{eq:lambda0_init}, are initialised by
$\lambda^b(0,m)=\iv{m=0}$. On the other hand, the alternative mutual
recurrences~\eqref{eq:lambda_b_rec_free} now become
\begin{align}
  \label{eq:lambda_b_rec_free_wis}
  \lambda^b(n,m)-q_b\,\lambda^b(n-1,m)=&~q_b\,\lambda^{\tilde{b}}(n-1,m-1)+\big(\iv{n=0}-q_b\,\iv{n=1}\big)\iv{m=0},
\end{align}
which do not need initialisation thanks to the inhomogeneous term.
The two mutual recurrences in~\eqref{eq:lambda_b_rec_free_wis} were
given by Wishart and Hirschfeld~\cite[Eqs.\ (1) and
(2)]{wishart36:_theorem} ---in their notation,
$\lambda(n,m)=P_n(m-1)$--- but without the inhomogeneous term, which
in fact makes~\eqref{eq:lambda_b_rec_free_wis} more general, as it
renders it valid for all $n\ge 1$ and $m\ge 0$ taking
\eqref{eq:m_necessary_condition_s} into account.

\subsubsection{Probability Generating Function}
Although we may particularise~\eqref{eq:ogf_lambda_klek} using
$\ushort{k}=1$ and $\widebar{k}=n$, as in other similar cases this
approach leads to an $n$-dependent pgf. A simpler, more general pgf is
possible in this case.  Working from~\eqref{eq:lambda_b_rec_free_wis}
and following the same steps as in
Section~\ref{sec:lambda_prob-gener-funct} (i.e., first obtaining a
system of two equations with the two unknowns
$\Lambda^b(x,y)=\sum_{n,m}\lambda^b(x,y)\, x^ny^m$, then solving it,
and finally adding the two solutions) is not difficult to see that the
pgf $\Lambda(x,y)=\sum_{n,m}\lambda(n,m)\,x^ny^m$ is
\begin{equation}
  \label{eq:lambda_pgf}
  \Lambda(x,y)=\frac{2(1-q_0x)(1-q_1x)+\big(q_0(1-q_1x)+q_1(1-q_0x)\big)\,xy}{(1-q_0x)(1-q_1x)-q_0q_1x^2y^2}.
\end{equation}
As usual, we may extract the coefficient of $y^m$
from~\eqref{eq:lambda_pgf} using the negative binomial theorem.  For
$m$ even we have
\begin{equation}
  \label{eq:lambda_pgf_m_even}
  [y^m]\Lambda(x,y)=2\Bigg(\frac{q_0q_1\, x^2}{1-x+q_0q_1 x^2}\Bigg)^{\frac{m}{2}},
\end{equation}
while for $m$ odd the expression is
\begin{equation}
  \label{eq_1:lambda_pgf_m_odd}
  [y^m]\Lambda(x,y)=\bigg(\frac{1}{q_0q_1x}-2\bigg)\Bigg(\frac{q_0q_1\, x^2}{1-x+q_0q_1 x^2}\Bigg)^{\frac{m+1}{2}}.
\end{equation}

\subsubsection{Explicit Expressions}
We next determine the coefficient of $x^n$ in
pgfs~\eqref{eq:lambda_pgf_m_even} and~\eqref{eq_1:lambda_pgf_m_odd} in
order to get single-summation explicit expressions for
$\lambda(n,m)=[x^ny^m]\Lambda(x,y)$. We start by developing the common
term in~\eqref{eq:lambda_pgf_m_even} and~\eqref{eq_1:lambda_pgf_m_odd}
into a power series using the negative binomial theorem followed by
the binomial theorem:
\begin{equation}
  \label{eq:common}
  \Bigg(\frac{q_0q_1\,
    x^2}{1-x+q_0q_1 x^2}\Bigg)^{v}=\sum_{t\ge
    0}\sum_{p\ge 0}{t+v-1\choose t}{t\choose p}(-1)^p (q_0q_1)^{p-v} x^{2v+t+p}.
\end{equation}
The coefficient of $x^n$ in~\eqref{eq:common} is found by solving
$2v+t+p=n$ for the nonnegative indices~$t$ and $p$. The maximum of $t$
happens when $p=0$, and thus $t\le n-2v$. Therefore, letting
$\omega(n,v)=[x^n]\big(q_0q_1\, x^2/(1-x+q_0q_1 x^2)\big)^{v}$ and
using~\eqref{eq:qb} in the following, we have from~\eqref{eq:common}
that
\begin{equation}
  \label{eq:omega}
  \omega(n,v)=\sum_{t=0}^{n-2v}
   {t+v-1\choose t}{t\choose n-2v-t}(-1)^{n-2v-t} ((1-q)q)^{n-v-t}.
\end{equation}
Hence, from~\eqref{eq:lambda_pgf_m_even} and~\eqref{eq:omega} we
have that for $m$ even
\begin{equation}
  \label{eq:lambda_m_even}
  \lambda(n,m)=2\,\omega\Big(n,\frac{m}{2}\Big),
\end{equation}
whereas from~\eqref{eq_1:lambda_pgf_m_odd} and~\eqref{eq:lambda_m_odd}
the $m$ odd case is
\begin{equation}
  \label{eq:lambda_m_odd}
  \lambda(n,m)=\frac{1}{q(1-q)}\,\omega\bigg(n+1,\frac{m+1}{2}\bigg)-2\,\omega\bigg(n,\frac{m+1}{2}\bigg).
\end{equation}
These two single-summation expressions for $\lambda(n,m)$ are similar to the
ones found by Wishart and Hirschfeld~\cite[Eqs.\ (22) and
(23)]{wishart36:_theorem}. However, the analysis by these authors
---which relies on the moment generating function (mgf) centered about
the mean, rather than the pgf--- also allows them to study the
asymptotics of $\lambda(n,m)$ in a more direct manner.

When $q=1/2$ we are back to a symmetric scenario, where we have 
$\lambda(n,m)=s(n,m)/2^n$. Therefore, using~\eqref{eq:s_n_m_explicit}
we can see that
\begin{equation}
  \label{eq:binomial}
  \lambda(n,m)=\frac{{n-1\choose m-1}}{2^{n-1}},
\end{equation}
which is the binomial distribution with parameters $n-1$ and $1/2$, as
observed by Wishart and Hirschfeld~\cite[p.\ 231]{wishart36:_theorem}.

\subsubsection{Moments}
We now let $\fixwidetilde{M}_n$ be the rv modelling the
number of nonnull runs of ones and/or zeros in an $n$-string drawn at random,
whose distribution can be expressed as
$\Pr(\fixwidetilde{M}_n=m)=\lambda(n,m)$.
Using~\eqref{eq:lambda_pgf}, the first moment of $\fixwidetilde{M}_n$
is the coefficient of $x^n$ in
\begin{equation}
  \label{eq:ddy_lambda_y1}
  \left.\frac{\partial \Lambda(x,y)}{\partial
      y}\right|_{y=1}=\frac{x\big(1-x+2(1-q)q x\big)}{(1-x)^2}.
\end{equation}
Expanding this expression using the negative binomial coefficient we
have that
\begin{equation}
  \label{eq:ddy_lambda_y1_exp}
  \left.\frac{\partial \Lambda(x,y)}{\partial
      y}\right|_{y=1}=\frac{x}{1-x}+2(1-q)\,q x^2\sum_{t\ge 0}{t+1\choose t}x^t,
\end{equation}
and, hence, the coefficient of $x^n$ in this expression is
\begin{equation}
  \label{eq:exp_mn}
  \text{E}\big(\fixwidetilde{M}_n\big)=1+2(1-q)q (n-1).
\end{equation}
Another way to obtain this expectation is
through~\eqref{eq:p_n_gtk_expectation}. Explicitly denoting by
$M_{\ge k,n}^{(q)}$ the rv in Section~\ref{sec:moments_w_gtk} with
parameter $q$, we have that
$\text{E}(\fixwidetilde{M}_n)=\text{E}(M_{\ge 1,n}^{(q)}+M_{\ge
  1,n}^{(1-q)})=\text{E}(M_{\ge 1,n}^{(q)})+\text{E}(M_{\ge
  1,n}^{(1-q)})$, even if $M_{\ge 1,n}^{(q)}$ and
$M_{\ge 1,n}^{(1-q)}$ are not independent. This strategy was used by
Cochran~\cite[Eq.\ (6)]{cochran36:_plants} to obtain the expected
number of $k$-runs of ones and/or zeros in an $n$-string drawn at
random from $\text{E}(M_{k,n})$ ---see~\eqref{eq:exp_mkn_all}.

For the second factorial moment we follow the same procedure as above
but with the second derivative of~\eqref{eq:lambda_pgf} evaluated at
$y=1$:
\begin{equation}
  \label{eq:ddy_lambda_y1}
  \left.\frac{\partial^2 \Lambda(x,y)}{\partial
      y^2}\right|_{y=1}=\frac{2(1-q)q x^2\big(2-x-(1-2q)^2 x^2\big)}{(1-x)^3}.
\end{equation}
Applying the negative binomial theorem, the coefficient of $x^n$ in this expression is
\begin{equation}
  \label{eq:2nd_m_mn}
  \text{E}\big(\fixwidetilde{M}_n(\fixwidetilde{M}_n-1)\big)=(1-q)q\big(2n(n-1)-(n-1)(n-2)-(1-2q)^2(n-2)(n-3)\big),
\end{equation}
from which we may obtain $\text{Var}(\fixwidetilde{M}_n)$ in
conjunction with~\eqref{eq:exp_mn}. Wishart and Hirschfeld
obtained~\eqref{eq:exp_mn} directly from their version of
recurrence~\eqref{eq:lambda_b_rec_free_wis} ---recall that in their
notation $\lambda(n,m)=P_n(m-1)$--- and then used the mgf to get the
second, third and fourth \textit{semi-invariants} (i.e., cumulants) of
the distribution~\cite[Eqs. (13) and (14)]{wishart36:_theorem}.

\section{Number of $n$-Strings of Hamming Weight $r$ that Contain
  Prescribed Quantities of Nonnull Runs of Ones and/or Zeros Under
  Different Constraints}
\label{sec:oz_hamming}
In this section we address the enumeration of the $n$-strings of
Hamming weight $r$ that contain prescribed quantities of nonnull runs
of ones and/or zeros under different constraints.  This section is to
Section~\ref{sec:oz} what Section~\ref{sec:number-n-strings-hamming}
is to Section~\ref{sec:number-n-strings}.  One important implication
of the Hamming weight constraint is that, unlike in
Section~\ref{sec:oz}, no enumeration in this section is related to a
counterpart enumeration concerning the compositions of~$n$. This is
because an $n$-string and its ones' complement do not have the same
Hamming weight in general. For the same reason, we do not have
analogues of Theorems~\ref{thm:chave_das_nozes}
and~\ref{thm:chave_das_nozes_2} in this setting: the results in this
section are not related to the results in
Section~\ref{sec:number-n-strings-hamming}. Like in
Sections~\ref{sec:oz} and~\ref{sec:note-prob-extens}, we only consider
nonnull runs of ones or zeros, and we approach the enumerations
through mutual recurrences. Importantly, the Hamming weight constraint
creates a fundamental asymmetry between these mutual recurrences like
the one we observed in Section~\ref{sec:note-prob-extens}, which, in
general, makes it harder to obtain overall recurrences and generating
functions. Due to this difficulty, we only get explicit
expressions in two particular cases.

Regarding previous work on this topic, Stevens~\cite{stevens39:_distr}
and Wald and Wolfowitz~\cite{wald40:_test} first addressed the
enumeration of the $n$-strings that have a prescribed number $m$ of
nonnull runs or ones and/or zeros under the Hamming weight
constraint. Their results were later rederived by  other
authors~\cite{gibbons11:_nonparam,schuster94:_exchan}. Also,
Bateman~\cite{bateman48:_power}, Schuster~\cite{schuster94:_exchan},
Bloom~\cite{bloom96:_probab}, and Jackson (see~\cite{bloom96:_probab})
tackled the enumeration of the $n$-strings of Hamming weight $r$ whose
longest nonnull run is a $(\le k)$-run of ones or zeros.

\subsection[Number of n-Strings of Hamming Weight r that Contain
Exactly m (ḵ≤k)-Runs of Ones and Zeros]{Number of $n$-Strings of
  Hamming Weight $r$ that Contain Exactly $m$ Nonnull
  ($\ushort{k}\le\widebar{k}$)-Runs of Ones and/or Zeros}
\label{sec:s_n_m_klek_hamming}

In this section we enumerate the $n$-strings of Hamming weight $r$
that contain exactly $m$ nonnull $(\ushort{k}\le\widebar{k})$-runs of
ones and/or zeros, which we denote by
$s_{\ushort{k}\le\widebar{k}}(n,m,r)$. Thus, we assume
$\ushort{k}\ge 1$.

As in Section~\ref{sec:w_n_m_klek_hamming}, we first narrow down
necessary condition~\eqref{eq:m_necessary_condition_s} to take into
account the Hamming weight constraint.
\begin{condition}{(Existence of $n$-strings of Hamming weight $r$
    containing $m$ nonnull $(\ushort{k}\le\widebar{k})$-runs of ones and/or zeros)}
  \begin{equation}
    \label{eq:s_m_necessary_condition_r}
    s_{\ushort{k}\le\widebar{k}}(n,m,r)> 0\quad \Longrightarrow\quad 0\le m \le\min\Bigg(\floor[\bigg]{\frac{n}{\ushort{k}}},\floor[\bigg]{\frac{r}{\ushort{k}}}+\floor[\bigg]{\frac{n-r}{\ushort{k}}}\Bigg).
  \end{equation}
\end{condition}

\subsubsection{Recurrence}
\label{sec:recurrence_s_n_m_k_r}
We can enumerate $s_{\ushort{k}\le\widebar{k}}(n,m,r)$ with the same
strategy used in Section~\ref{sec:s_n_m_klek}: we first produce two
mutual recurrences for the number of $n$-strings of Hamming weight $r$
that start with bit~$b$ and contain exactly $m$
${\ushort{k}\le\widebar{k}}$-runs of ones and/or zeros, denoted by
$s^b_{\ushort{k}\le\widebar{k}}(n,m,r)$, where~$b\in\{0,1\}$. The two
trivariate mutual recurrences sought are just like in~\eqref{eq:s0}
but updating the Hamming weight constraints:
\begin{align}
  \label{eq:s0_n_m_r_klek}
  s^b_{\ushort{k}\le\widebar{k}}(n,m,r)=&~\sum_{i=\ushort{k}}^{\widebar{k}}s^{\tilde{b}}_{\ushort{k}\le\widebar{k}}\big(n-i,m-1,r-\iv{b=1}\,i\big)+\sum_{i=1}^{\ushort{k}-1}s^{\tilde{b}}_{\ushort{k}\le\widebar{k}}\big(n-i,m,r-\iv{b=1}\,i\big)\nonumber\\&+  \sum_{i=\widebar{k}+1}^{n}s^{\tilde{b}}_{\ushort{k}\le\widebar{k}}\big(n-i,m,r-\iv{b=1}\,i\big),
\end{align}
where $b\in\{0,1\}$ and $\tilde{b}=\text{mod}(b+1,2)$.  To initialise
these two mutual recurrences we use the case $n=1$, in which we
know by inspection that
\begin{align}
  \label{eq:s0_n_m_r_klek_n1_trivial}
  s^b_{\ushort{k}\le\widebar{k}}(1,m,r)=&~\iv{r=b}
                                          \Big(\iv{m=0}\iv{\ushort{k}>1}+\iv{m=1}\iv{\ushort{k}=1}\Big).
\end{align}
Setting $n=1$ in~\eqref{eq:s0_n_m_r_klek} we have
\begin{equation}
  \label{eq:s0_n_m_r_klek_n1_rec}
  s^b_{\ushort{k}\le\widebar{k}}(1,m,r)=  s^b_{\ushort{k}\le\widebar{k}}\big(0,m-\iv{\ushort{k}=1},r-\iv{b=1}\big).
\end{equation}
We see that~\eqref{eq:s0_n_m_r_klek_n1_rec} 
equals~\eqref{eq:s0_n_m_r_klek_n1_trivial}  if we choose the following initialisation values:
\begin{equation}
  \label{eq:s_r_init}
  s^b_{\ushort{k}\le\widebar{k}}(0,m,r)=\iv{m=0}\iv{r=0}.
\end{equation}
The main enumeration in this section can now be expressed as
\begin{equation}
  \label{eq:s_n_m_klek_r}
  s_{\ushort{k}\le\widebar{k}}(n,m,r)=s^0_{\ushort{k}\le\widebar{k}}(n,m,r)+s^1_{\ushort{k}\le\widebar{k}}(n,m,r).
\end{equation}
Importantly, unlike in Section~\ref{sec:s_rec} the mutual recurrences
in~\eqref{eq:s0_n_m_r_klek} for the $n$-strings that start with $0$ or
with $1$ are not symmetric. This is analogous to what we saw in
Section~\ref{sec:note-prob-extens}. Because of this, we cannot simply
input~\eqref{eq:s0_n_m_r_klek} into~\eqref{eq:s_n_m_klek_r} and then
directly work out a recurrence
for~$s_{\ushort{k}\le\widebar{k}}(n,m,r)$. Therefore we do not attempt
to obtain such a general recurrence, but we see later that this
asymmetry can be wrestled with in special cases.

\subsubsection{Generating Functions}
\label{sec:s_n_m_klek_r_generating-functions}
While we have not been able to produce a general recurrence for
$s_{\ushort{k}\le\widebar{k}}(n,m,r)$, we can still determine its ogf
$S_{\ushort{k}\le\widebar{k}}(x,y,z)=\sum_{n,m,r} s_{\ushort{k}\le\widebar{k}}(n,m,r)\,x^ny^mz^r$
through the mutual recurrences for
$s^b_{\ushort{k}\le\widebar{k}}(n,m,r)$. First of all, setting $n=0$
in~\eqref{eq:s0_n_m_r_klek} we get $s^b_{\ushort{k}\le\widebar{k}}(n,m,r)=0$
instead of~\eqref{eq:s_r_init}. Thus, to get recurrences valid for all
$n, m$, and $r$ we add $\iv{n=0}\iv{m=0}\iv{r=0}$
to~\eqref{eq:s0_n_m_r_klek}. We next get recurrences without $n$-dependent summations
by calculating
$s^b_{\ushort{k}\le\widebar{k}}(n,m,r)-s^b_{\ushort{k}\le\widebar{k}}\big(n-1,m,r-\iv{b=1}\big)$
using the extended recurrences. This yields
\begin{align}
  \label{eq:sb_n_m_klek_r_extended_rec}
  s^b_{\ushort{k}\le\widebar{k}}(n,m,r)=&~s^b_{\ushort{k}\le\widebar{k}}\big(n-1,m,r-\iv{b=1}\big)+s^{\tilde{b}}_{\ushort{k}\le\widebar{k}}\big(n-\ushort{k},m-1,r-\iv{b=1}\,\ushort{k}\big)\nonumber\\&-s^{\tilde{b}}_{\ushort{k}\le\widebar{k}}\big(n-(\widebar{k}+1),m-1,r-\iv{b=1}\,(\widebar{k}+1)\big)+s^{\tilde{b}}_{\ushort{k}\le\widebar{k}}\big(n-1,m,r-\iv{b=1}\big)\nonumber\\&-s^{\tilde{b}}_{\ushort{k}\le\widebar{k}}\big(n-\ushort{k},m,r-\iv{b=1}\,\ushort{k}\big)+s^{\tilde{b}}_{\ushort{k}\le\widebar{k}}\big(n-(\widebar{k}+1),m,r-\iv{b=1}\,(\widebar{k}+1)\big)\nonumber\\
  &+\iv{m=0}\big(\iv{r=0}\iv{n=0}-\iv{r=b}\iv{n=1}\big).
\end{align}
We are now ready to obtain the trivariate ogfs
$S^b_{\ushort{k}\le\widebar{k}}(x,y,z)=\sum_{n,m,r}
s^b_{\ushort{k}\le\widebar{k}}(n,m,r)\,x^ny^mz^r$ for
$b\in\{0,1\}$. Multiplying~\eqref{eq:sb_n_m_klek_r_extended_rec} on
both sides by $x^ny^mz^r$ and adding over $n,m$, and $r$ yields the
following system of equations:
\begin{align}
  \label{eq:S0_S1_system_n_m_klek_r}
  S^0_{\ushort{k}\le\widebar{k}}(x,y,z)(1-x)&=S^1_{\ushort{k}\le\widebar{k}}(x,y,z)\Big(\big(x^{\ushort{k}}-x^{\widebar{k}+1}\big)(y-1)+x\Big)+1-x,\\ 
  S^1_{\ushort{k}\le\widebar{k}}(x,y,z)(1-xz)&=S^0_{\ushort{k}\le\widebar{k}}(x,y,z)\Big(\big(x^{\ushort{k}}z^{\ushort{k}}-x^{\widebar{k}+1}z^{\widebar{k}+1}\big)(y-1)+xz\Big)+1-xz. 
\end{align}
Solving the system we get
\begin{align}
  \label{eq:S1_n_m_klek_r}
  S^1_{\ushort{k}\le\widebar{k}}(x,y,z)&=\frac{(1-x)\Big(\big(x^{\ushort{k}}z^{\ushort{k}}-x^{\widebar{k}+1}z^{\widebar{k}+1}\big)(y-1)+1\Big)}{(1-x)(1-xz)-\Big(\big(x^{\ushort{k}}z^{\ushort{k}}-x^{\widebar{k}+1}z^{\widebar{k}+1}\big)(y-1)+xz\Big)\Big(\big(x^{\ushort{k}}-x^{\widebar{k}+1}\big)(y-1)+x\Big)},
\end{align}
whereas $S^0_{\ushort{k}\le\widebar{k}}(x,y,z)$ is obtained by
substituting~\eqref{eq:S1_n_m_klek_r}
in~\eqref{eq:S0_S1_system_n_m_klek_r}, which yields
\begin{align}
  \label{eq:S0_n_m_klek_r}
  S^0_{\ushort{k}\le\widebar{k}}(x,y,z)&=\frac{\Big(\big(x^{\ushort{k}}z^{\ushort{k}}-x^{\widebar{k}+1}z^{\widebar{k}+1}\big)(y-1)+1\Big)\Big(\big(x^{\ushort{k}}-x^{\widebar{k}+1}\big)(y-1)+x\Big)}{(1-x)(1-xz)-\Big(\big(x^{\ushort{k}}z^{\ushort{k}}-x^{\widebar{k}+1}z^{\widebar{k}+1}\big)(y-1)+xz\Big)\Big(\big(x^{\ushort{k}}-x^{\widebar{k}+1}\big)(y-1)+x\Big)}+1.
\end{align}
Thus, from~\eqref{eq:s_n_m_klek_r}, ~\eqref{eq:S1_n_m_klek_r}
and~\eqref{eq:S0_n_m_klek_r} the desired ogf is
\begin{multline}
  \label{eq:s_n_m_klek_r_ogf}
  S_{\ushort{k}\le\widebar{k}}(x,y,z)=\\\frac{(1-x)(1-xz)+(1-x)\Big(\big(x^{\ushort{k}}z^{\ushort{k}}-x^{\widebar{k}+1}z^{\widebar{k}+1}\big)(y-1)+1\Big)+(1-xz)\Big(\big(x^{\ushort{k}}-x^{\widebar{k}+1}\big)(y-1)+x\Big)}{(1-x)(1-xz)-\Big(\big(x^{\ushort{k}}z^{\ushort{k}}-x^{\widebar{k}+1}z^{\widebar{k}+1}\big)(y-1)+xz\Big)\Big(\big(x^{\ushort{k}}-x^{\widebar{k}+1}\big)(y-1)+x\Big)}.
\end{multline}
Although finding an explicit expression
is impracticable, one can  use a computer algebra system to obtain
$s_{\ushort{k}\le\widebar{k}}(n,m,r)=[x^ny^mz^r]S_{\ushort{k}\le\widebar{k}}(x,y,z)$
from~\eqref{eq:s_n_m_klek_r_ogf}.

\subsection[Number of n-Strings of Hamming Weight r that Contain
Exactly m Nonnull Runs of Ones and Zeros]{Number of $n$-Strings of Hamming
  Weight $r$ that Contain Exactly $m$ Nonnull Runs of Ones and/or Zeros}
\label{sec:s_n_m_hamming}

We call $s(n,m,r)$ the number of $n$-strings of Hamming weight $r$
that contain exactly $m$ nonnull runs of ones and/or zeros (of
arbitrary lengths, all strictly greater than zero). This enumeration
---in fact, the probability $\lambda(n,m,r)=s(n,m,r)/{n\choose r}$---
was first given by Stevens~\cite{stevens39:_distr}. However, its
relevance in the theory of runs is mainly due to Wald and Wolfowitz,
who independently derived the same result one year
later~\cite{wald40:_test} to use it as the basis of their runs-based
hypothesis test for the identity between the distributions from which
two samples are
drawn. 
The probability $\lambda(n,m,r)$ was also rederived later in
alternative ways by several other authors, such as
Gibbons~\cite{gibbons11:_nonparam} and
Schuster~\cite{schuster94:_exchan}. The latter author, who denotes the
probability $\lambda(n,m,r)$ by $R$, comments~\cite[Rem.\
2]{schuster94:_exchan}): ``(\ldots) \textit{much of the difficulty in
  studying the theoretical properties of $R$ is due to the fact that
  $R$ is not symmetrically defined} (\ldots)''. This difficulty, which
was in fact the stumbling block preventing us from producing a general
recurrence in Section~\ref{sec:s_n_m_klek_hamming}, also rears its
ugly head in our approach to this particular enumeration, although in
a less severe manner.

As in previous cases, $s(n,m,r)$ is a special case of our general
analysis in Section~\ref{sec:s_n_m_klek_hamming} with $\ushort{k}=1$
and $\widebar{k}=n$, and therefore
\begin{equation}
  \label{eq:s_n_m}
  s(n,m,r)=s_{1\le n}(n,m,r).
\end{equation}

\subsubsection{Recurrence}
\label{sec:s_n_m_r_rec}
In order to delve deeper into this special enumeration, we start by
specialising recurrence~\eqref{eq:s0_n_m_r_klek} to get
\begin{align}
  \label{eq:sb_n_m_r_rec}
  s^b(n,m,r)=&~\sum_{i=1}^{n}s^{\tilde{b}}\big(n-i,m-1,r-\iv{b=1}\,i\big).
\end{align}
We certainly face the same issue here as in the general case: the two
mutual recurrences defined by~\eqref{eq:sb_n_m_r_rec} are not
symmetric. However, in this particular case we may symmetrise them,
which allows us to produce a recurrence for $s(n,m,r)$. To this end,
we fist make \eqref{eq:sb_n_m_r_rec} valid for all values of the
parameters by adding $\iv{n=0}\iv{m=0}\iv{r=0}$ to it. We then apply
this extended recurrence to itself, which yields
\begin{align}
  s^b(n,m,r)=\sum_{i=1}^{n}\Bigg(&\sum_{j=1}^ns^{b}\big(n-(i+j),m-2,r-\iv{b=1}\,i-\iv{b=0}\,j\big)\nonumber\\&+\iv{n=i}\iv{m=1}\ivl{r=\iv{b=1}\,i}\Bigg)\nonumber\\
               &+\iv{n=0}\iv{m=0}\iv{r=0},  \label{eq:sb_n_m_r_rec_sym}
\end{align}
where we have extended the summation on $j$ from $n-i$ to $n$ to make
the symmetry due to the summations on $i$ and $j$ clear. This does not
alter the recurrence, as from necessary
condition~\eqref{eq:s_m_necessary_condition_r} we have that
$s^b(n,m,r)=0$ for $n<0$.

We can now directly obtain
$s(n,m,r)=s^0(n,m,r)+s^1(n,m,r)$ using~\eqref{eq:sb_n_m_r_rec_sym},
which yields the following recurrence:
\begin{align}
  s(n,m,r)=&~\sum_{i=1}^{n}\sum_{j=1}^{n-i}s\big(n-(i+j),m-2,r-i\big)\nonumber\\%
               &+\iv{n>0}\iv{m=1}\big(\iv{r=n}+\iv{r=0}\big)+2\iv{n=0}\iv{m=0}\iv{r=0}.  \label{eq:s_n_m_r_rec_sym}
\end{align}
Observe that, after the addition of the two mutual recurrences, we
have set the upper limit of the summation on $j$ back to the more efficient $n-i$.

\subsubsection{Generating Function}
\label{sec:s_n_m_r_ogf}
One way to work out the ogf $S(x,y,z)=\sum_{n,m,r}s(n,m,r)\,x^ny^mz^r$
is simply to specialise~\eqref{eq:s_n_m_klek_r_ogf}, but as in other
similar cases this gives an $n$-dependent ogf and is anyway
unwieldy. As we see next, we can get a much simpler and general
expression. While we have the option of finding an ogf through a
version of recurrence~\eqref{eq:s_n_m_r_rec_sym} without $n$-dependent
summations, a gentler approach is to produce instead ogfs from the
recurrences~\eqref{eq:sb_n_m_klek_r_extended_rec} specialised to
this case, i.e., $S^b(x,y,z)=\sum_{n,m,r}s^b(n,m,r)\, x^ny^mz^r$, and
then obtain $S(x,y,z)=S^0(x,y,z)+S^1(x,y,z)$ like we did in
Section~\ref{sec:s_n_m_klek_r_generating-functions}. The
specialisation of~\eqref{eq:sb_n_m_klek_r_extended_rec} yields
\begin{align}
  s^b(n,m,r)=&~s^b\big(n-1,m,r-\iv{b=1}\big)+s^{\tilde{b}}\big(n-1,m-1,r-\iv{b=1}\big)\nonumber\\
             &+\iv{m=0}\big(\iv{r=0}\iv{n=0}-\iv{r=b}\iv{n=1}\big).  \label{eq:sb_n_m_r_extended_rec}
\end{align}
This recurrence is valid already for all values of the
parameters. Thus, multiplying it on both sides by $x^ny^mz^r$ and then
adding over $n,m$ and $r$ we get
\begin{align}
  \label{eq:Sb_n_m_r_ogf}
  S^b(x,y,z) (1-xz^{\iv{b=1}})=&~S^{\tilde{b}}(x,y,z)\, xyz^{\iv{b=1}}+(1-xz^b).
\end{align}
Again, we have a linear system of two equations with two unknowns, i.e.,
$S^0(x,y,z)$ and $S^1(x,y,z)$. Solving it we find that
\begin{equation}
  \label{eq:S_n_m_r_ogf}
 S(x,y,z) =\frac{2 - 2 x + x y - 2 x z + 2 x^2 z + x y z - 2 x^2 y z}{1 - x - x z + x^2 z - x^2 y^2 z}.
\end{equation}

\subsubsection{Explicit Expressions}
We obtain next $s(n,m,r)=[x^ny^mz^r]S(x,y,z)$ 
from~\eqref{eq:S_n_m_r_ogf}. To simplify
the procedure it is convenient to rewrite~\eqref{eq:S_n_m_r_ogf} as follows:
\begin{equation}
  \label{eq:S_alternative}
  S(x,y,z) =\frac{2+\big(-1+(1-x^2z)(1-xz)^{-1}(1-x)^{-1}\big)y}{1-x^2z(1-xz)^{-1}(1-x)^{-1}y^2}.
\end{equation}
By applying the negative binomial theorem, we can
express~\eqref{eq:S_alternative} as
\begin{equation}
  \label{eq:S_alternative_exp}
  S(x,y,z)
  =\bigg(2+\bigg(-1+\frac{1-x^2z}{(1-xz)(1-x)}\bigg)y\bigg)\sum_{p\ge
    0}\bigg(\frac{x^2z}{(1-xz)(1-x)}\bigg)^{p} y^{2p}.
\end{equation}
Thus, defining
\begin{equation}
  \label{eq:sigmasigma}
  \sigma_{\{a,b\}}(x,z)=\frac{(x^2z)^a}{\Big((1-xz)(1-x)\Big)^b},
\end{equation}
we have that the coefficient of $y^m$ for $m$ even  is
\begin{equation}
  \label{eq:S_ym_even}
  [y^m]S(x,y,z)  =2\,\sigma_{\left\{\frac{m}{2},\frac{m}{2}\right\}}(x,z),
\end{equation}
whereas for $m$ odd
\begin{equation}
  \label{eq:S_ym_odd}
  [y^m]S(x,y,z) =-\sigma_{\left\{\frac{m-1}{2},\frac{m-1}{2}\right\}}(x,z)+\sigma_{\left\{\frac{m-1}{2},\frac{m+1}{2}\right\}}(x,z)-\sigma_{\left\{\frac{m+1}{2},\frac{m+1}{2}\right\}}(x,z).
\end{equation}
On the other hand, we can expand~\eqref{eq:sigmasigma} by applying the
negative binomial theorem twice to get
\begin{equation}
  \label{eq:sigma_expanded}
  \sigma_{\{a,b\}}(x,z)=\sum_{s\ge 0}\sum_{t\ge 0}
  {s+b-1\choose b-1}{t+b-1\choose b-1} x^{2a+s+t} z^{a+s}.
\end{equation}
In order to determine $[x^nz^r]\sigma_{\{a,b\}}(x,z)$ we need to find
the indices $s$ and $t$ that fulfil $n=2a+s+t$ and $r=a+s$. From the
second equation we see that $s=r-a$, and thus $t=n-r-a$. Therefore,
from these solutions and~\eqref{eq:sigma_expanded} we have that
\begin{equation}
  \label{eq:xz_sigma}
  [x^nz^r]\sigma_{\{a,b\}}(x,z)={r-a+b-1\choose b-1}{n-r-a+b-1\choose b-1}.
\end{equation}
So, from~\eqref{eq:S_ym_even} and~\eqref{eq:xz_sigma} we obtain the
following expression for $m$ even:
\begin{equation}
  \label{eq:s_n_m_rm_even}
  s(n,m,r)=2\,{r-1\choose \frac{m}{2}-1}{n-r-1\choose \frac{m}{2}-1}.
\end{equation}
From~\eqref{eq:S_ym_odd} and~\eqref{eq:xz_sigma} we have that,
for $m$ odd, the expression is
\begin{equation}
  \label{eq:s_n_m_rm_odd_p}
  s(n,m,r)=-{r-1\choose \frac{m-1}{2}-1}{n-r-1\choose \frac{m-1}{2}-1}+{r\choose \frac{m-1}{2}}{n-r\choose \frac{m-1}{2}}-{r-1\choose \frac{m-1}{2}}{n-r-1\choose \frac{m-1}{2}}.
\end{equation}
Applying ${a\choose b}={a-1\choose b}+{a-1\choose b-1}$ to both
binomial coefficients in the positive term
in~\eqref{eq:s_n_m_rm_odd_p} we see that the negative terms cancel
out, which leads to a simpler expression for $m$ odd:
\begin{equation}
  \label{eq:s_n_m_rm_odd}
  s(n,m,r)={r-1\choose \frac{m-1}{2}-1}{n-r-1\choose \frac{m-1}{2}}+{r-1\choose \frac{m-1}{2}}{n-r-1\choose \frac{m-1}{2}-1}.
\end{equation}
In the form $\lambda(n,m,r)=s(n,m,r)/{n\choose r}$, the explicit
expressions~\eqref{eq:s_n_m_rm_even} and~\eqref{eq:s_n_m_rm_odd} are
well-known since the works of Stevens~\cite[Eqs.\ (3.31) and
(3.32)]{stevens39:_distr} and Wald and Wolfowitz~\cite[Eqs.\ (7) and
(8)]{wald40:_test}. See also the rederivations by Gibbons~\cite[Thm.\
3.2.2]{gibbons11:_nonparam} and Schuster~\cite[Cor.\
3.7]{schuster94:_exchan}, and the alternative moment computations
by Guenther~\cite{guenther78:_remarks}.

\begin{remark}
  Like in Section~\ref{sec:prob-conn}, we can invoke the law of large
  numbers to write
  $\lambda(n,m)\approx \lambda(n,m,\floor{nq})=
  s(n,m,\floor{nq})/{n\choose \floor{nq}}$ for large $n$, where
  $\lambda(n,m)$ is the probability discussed in
  Section~\ref{sec:probability-s_n_m}.  Observe that both in
  $s(n,m,r)$ and in $\lambda(n,m)$ there is a dichotomy between two
  different explicit expressions depending on the parity of $m$.  For
  example, if we let $\fixwidetilde{M}_{n,r}$ be the rv with
  distribution $\Pr(\fixwidetilde{M}_{n,r}=m)=s(n,m,r)/{n\choose r}$,
  then we can approximate its expectation using~\eqref{eq:exp_mn} and
  $q=r/n$ to see that, for large $n$,
  $\text{E}(\fixwidetilde{M}_{n,r})\approx 1+2(n-r)r(n-1)/n^2\approx
  1+2(n-r)r/n$, an expression that is well known~\cite[Eq.\
  (12)]{wald40:_test}.
  
  This connection between $\lambda(n,m)$ and $\lambda(n,m,r)$ implies
  that Wishart and Hirschfeld's work~\cite{wishart36:_theorem} somehow
  predated that of Wald and Wolfowitz~\cite{wald40:_test}: Wald and
  Wolfowitz's test could, in principle, have been based on the earlier
  results in~\cite{wishart36:_theorem} for large $n$.
\end{remark}

\subsection[Number of n-Strings of Hamming Weight r Whose Longest Nonnull Run,
Is a (≤k)-Run of Ones or Zeros]{Number of $n$-Strings of Hamming
  Weight $r$ Whose Longest Nonnull Run Is a  $(\le k)$-Run of Ones or Zeros}
\label{sec:ss_n_m_r}

We denote by $\widebar{s}_{\le k}(n,r)$ the number of $n$-strings of
Hamming weight $r$ whose longest nonnull run is a $(\le k)$-run of
ones or zeros, where $k\ge 1$. Bateman~\cite{bateman48:_power} gave an
explicit expression for the complementary probability of
$\widebar{\lambda}_{\le k}(n,r)=\widebar{s}_{\le k}(n,r)/{n\choose r}$,
and tables for
$\widebar{\lambda}_{k}(n,r)=\widebar{\lambda}_{\le
  k}(n,r)-\widebar{\lambda}_{\le (k-1)}(n,r)$ for small values
of $n, r$ and $k$. Schuster~\cite{schuster94:_exchan},
Bloom~\cite{bloom96:_probab}, and Jackson (see~\cite{bloom96:_probab})
gave recurrences for $\widebar{s}_{\le k}(n,r)$. Schuster
also extended and corrected Bateman's tables.

We can study this problem as a special case of the results in
Section~\ref{sec:s_n_m_klek_hamming}, because
\begin{equation}
  \label{eq:ss_s}
  \widebar{s}_{\le k}(n,r)=s_{(k+1)\le n}(n,0,r).
\end{equation}

\subsubsection{Recurrence}
\label{sec:ss_n_m_r_rec}
Let us see how we can deduce a recurrence for
$\widebar{s}_{\le k}(n,r)$ from the general expressions. Calling
$\widebar{s}^{\,b}_{\le k}(n,r)=s^b_{(k+1)\le n}(n,0,r)$, we can
express the quantity for which we want a recurrence as
\begin{equation}\label{eq:ss_n_ltk_r_main}
  \widebar{s}_{\le k}(n,r)=\widebar{s}^{\,0}_{\le k}(n,r)+\widebar{s}^{\,1}_{\le k}(n,r).
\end{equation}
With the parameters of this special
case~\eqref{eq:s0_n_m_r_klek} becomes
\begin{align}
  \label{eq:b_ss_n_ltk_r}
  \widebar{s}^{\,b}_{\le k}(n,r)=&\sum_{i=1}^{k}\widebar{s}^{\,\tilde{b}}_{\le k}\big(n-i,r-\iv{b=1}\,i\big).
\end{align}
Again, the main obstacle to finding a recurrence for
$\widebar{s}_{\le k}(n,r)$ is the asymmetry between the two mutual
recurrences in~\eqref{eq:b_ss_n_ltk_r}. Nevertheless, like in
Section~\ref{sec:s_n_m_r_rec}, in this case it is also possible to
symmetrise them in a straightforward manner.  We start by
making~\eqref{eq:b_ss_n_ltk_r} valid for all $n$ and $r$ by adding
$\iv{n=0}\iv{r=0}$ to it. We then apply this extended recurrence to
itself, which gives
\begin{align}
  \label{eq:b_ss_n_ltk_r_sym}
  \widebar{s}^{\,b}_{\le k}(n,r)=&~\sum_{i=1}^{k}\Bigg(\sum_{j=1}^{k}\widebar{s}^{\,b}_{\le
                                   k}\big(n-(i+j),r-\iv{b=1}\,i-\iv{b=0}\,j\big)+\iv{n=i}\ivl{r=\iv{b=1}\,i}\Bigg)\nonumber\\
                                 &+\iv{n=0}\iv{r=0}.
\end{align}
We can now input~\eqref{eq:b_ss_n_ltk_r_sym}
into~\eqref{eq:ss_n_ltk_r_main} to obtain the following recursive
relation:
\begin{align}
  \label{eq:s0_n_m_r_klek_our}
  \widebar{s}_{\le k}(n,r)=&~\sum_{i=1}^{k}\Bigg(\sum_{j=1}^{k} \widebar{s}_{\le
                             k}\big(n-(i+j),r-i\big)+\iv{n=i}\big(\iv{r=i}+\iv{r=0}\big)\Bigg)\nonumber\\
                           &+2\iv{n=0}\iv{r=0}.
\end{align}
Taking~\eqref{eq:s_m_necessary_condition_r} into
account,~\eqref{eq:s0_n_m_r_klek_our} is valid for all $n$ and $r$, and so it
does not need initialisation.

\begin{remark}\label{rmk:bloom1}
  At first sight,~\eqref{eq:s0_n_m_r_klek_our} may look unrelated to
  the recurrence for the same enumeration first cleverly surmised and
  then proved by Bloom~\cite[Eq.\ (9)]{bloom96:_probab} ---in his
  notation, $\widebar{s}_{\le k}(n,r)=C_{k+1}(n-r,r)$. However, if we
  obtain $\widebar{s}_{\le k}(n,r)-\widebar{s}_{\le k}(n-1,r)$
  using~\eqref{eq:s0_n_m_r_klek_our}, then we essentially recover
  Bloom's recurrence:
  \begin{align}
    \label{eq:s0_n_m_r_klek_our_bloom}
    \widebar{s}_{\le k}(n,r)=&~\sum_{i=0}^{k}\widebar{s}_{\le k}\big(n-(i+1),r-i\big)-\sum_{i=1}^{k}\widebar{s}_{\le
                               k}\big(n-(i+k+1),r-i\big)\nonumber\\
                             &+\sum_{i=1}^{k}\big(\iv{n=i}-\iv{n=i+1}\big)\big(\iv{r=i}+\iv{r=0}\big)\nonumber\\
                             &+2\big(\iv{n=0}-\iv{n=1}\big)\iv{r=0}.
  \end{align}
  In Bloom's recurrence, index $i$ only affects the second argument of
  $C_t(\cdot,\cdot)$ whereas above it affects both arguments of
  $\widebar{s}_{\le k}(\cdot,\cdot)$. This just a notational artifact:
  in Bloom's notation the addition of the two arguments yields the
  length of the binary string, and thus the homogeneous part of both
  recurrences is exactly the same. The inhomogeneous term in Bloom's
  recurrence does not always equal the inhomogeneous term
  in~\eqref{eq:s0_n_m_r_klek_our_bloom}, which simply indicates a
  different initialisation strategy: 
  both recurrences return the same values.
  
  Finally, if we obtain
  $\widebar{s}_{\le k}(n,r)-\widebar{s}_{\le k}(n-1,r-1)$
  using~\eqref{eq:s0_n_m_r_klek_our_bloom}, we also 
  essentially recover Jackson's recurrence~\cite[Eq.\
  (13)]{bloom96:_probab}:
  \begin{align}
    \label{eq:s0_n_m_r_klek_our_jackson}
    \widebar{s}_{\le k}(n,r)=&~ \widebar{s}_{\le k}\big(n-1,r-1\big)+\widebar{s}_{\le k}\big(n-1,r\big)-\widebar{s}_{\ge
                               k}\big(n-(k+2),r-(k+1)\big)\nonumber\\
                             &-\widebar{s}_{\le k}\big(n-(k+2),r-1\big)+\widebar{s}_{\le k}\big(n-2(k+1),r-(k+1)\big)\nonumber\\
                             &+\sum_{i=1}^{k}\big(\iv{n=i}-\iv{n=i+1}\big)\big(\iv{r=i}+\iv{r=0}\big)\nonumber\\
                             &-\sum_{i=1}^{k}\big(\iv{n=i+1}-\iv{n=i+2}\big)\big(\iv{r=i+1}+\iv{r=1}\big)\nonumber\\
                             &+2\big(\iv{n=0}-\iv{n=1}\big)\iv{r=0}-2\big(\iv{n=1}-\iv{n=2}\big)\iv{r=1}.
  \end{align}
  Again the homogeneous parts of~\eqref{eq:s0_n_m_r_klek_our_jackson}
  and~\cite[Eq.\ (13)]{bloom96:_probab} are identical, but the
  inhomogeneous parts are different. At any rate, both recurrences
  deliver the same numerical values.  The number of recurrent calls
  made by~\eqref{eq:s0_n_m_r_klek_our},
  \eqref{eq:s0_n_m_r_klek_our_bloom}
  and~\eqref{eq:s0_n_m_r_klek_our_jackson} is $k^2$, $2k+1$ and~$5$,
  respectively. Thus, \eqref{eq:s0_n_m_r_klek_our_jackson} is the most
  efficient recurrence when $k> 2$, only improved
  by~\eqref{eq:s0_n_m_r_klek_our} when $k=1$ or
  $k=2$.  

  To conclude, we should mention that Schuster was the first author
  who found a recurrence for~$\widebar{s}_{\le k}(n,r)$ ---in fact,
  for
  $\widebar{\lambda}_{\le k}(n,r)=\widebar{s}_{\le k}(n,r)/{n\choose
    r}$~\cite[Cor.\ 5.4]{schuster94:_exchan}. His recurrence depends
  on an involved ancillary function, and it does not seem to bear a
  close relationship to the recurrences discussed above.
\end{remark}

\subsubsection{Generating Functions}
\label{sec:ss_n_m_r_ogf}
We now obtain
$\widebar{S}_{\le k}(x,z)=\sum_{n,r}\widebar{s}_{\le
  k}(n,r)\,x^nz^r$. We can do so using~\eqref{eq:s0_n_m_r_klek_our},
\eqref{eq:s0_n_m_r_klek_our_bloom}
or~\eqref{eq:s0_n_m_r_klek_our_jackson} without further ado, since
none of these equivalent recurrences contains $n$-dependent summations and
all of them are valid for all arguments. Let us
use~\eqref{eq:s0_n_m_r_klek_our_jackson}, as it involves the least
amount of algebra. Multiplying this recurrence on both sides by
$x^nz^r$ and then adding over $n$ and $r$ we can see that
\begin{equation}
  \label{eq:Ss_n_ltk_r_ogf}
  \widebar{S}_{\le k}(x,z)=\frac{2-x(1+z)+x^{k+2}(z+z^{k+1})-x^{k+1}(1+z^{k+1})}{1-x(1+z)+x^{k+2}(z+z^{k+1})-x^{2(k+1)}z^{k+1}}.
\end{equation}
The reader may verify that the same ogf is obtained from
recurrences~\eqref{eq:s0_n_m_r_klek_our}
or~\eqref{eq:s0_n_m_r_klek_our_bloom}.
Using instead Jackson's original recurrence~\cite[Eq.\
  (13)]{bloom96:_probab} we get a somewhat simpler ogf:
\begin{equation}
  \label{eq:Ss_n_ltk_r_ogf_jackson}
\widebar{S}_{\le k}(x,z)=\frac{1-x^{k+1}(1+z^{k+1})+x^{2(k+1)}z^{k+1}}{1-x(1+z)+x^{k+2}(z+z^{k+1})-x^{2(k+1)}z^{k+1}}.
\end{equation} 
Using any of these two ogfs, we may 
obtain $\widebar{s}_{\le k}(n,r)=[x^nz^r]\widebar{S}_{\le
  k}(x,z)$ using a computer algebra system.

\subsubsection{Explicit Expression}
In order to extract the coefficient of $x^nz^r$ from
$\widebar{S}_{\le k}(x,z)$ we define
\begin{equation}
  \label{eq:axz}
  a^{(u,v,w)}(x,z)= \frac{x^vz^w\big(1-x-xz+x^{k+2}z+x^{k+2}z^{k+1}\big)^{-u}}{1-\big(x^{2(k+1)}z^{k+1}\big)\big(1-x-xz+x^{k+2}z+x^{k+2}z^{k+1}\big)^{-1}}.
\end{equation}
which we use to rewrite~\eqref{eq:Ss_n_ltk_r_ogf} as
\begin{equation}
  \label{eq:SS_ltk_ogf_rewritten}
  \widebar{S}_{\le k}(x,z)=a^{(0,0,0)}(x,z)+a^{(1,0,0)}(x,z)-a^{(1,k+1,0)}(x,z)-a^{(1,k+1,k+1)}(x,z).
\end{equation}
Applying the negative binomial theorem (twice) and then the binomial
theorem (three times), we can expand~\eqref{eq:axz} as
\begin{equation}
  \label{eq:bxz_expand}
  a^{(u,v,w)}(x,z)=\!\!\sum_{t,p,i,j,l\ge 0}\!\!(-1)^i{p+t+u-1\choose p}{p\choose
    i}{p-i\choose j}{i\choose l}
    x^{2(k+1)t+p+(k+1)i+v}\, z^{(k+1)t+i+j+l k+w}.
\end{equation}
Thus, the coefficient of $x^nz^r$ in this expression is obtained by
determining the nonnegative indices $t, p, i, j,$ and $l$ that fulfil
$n=2(k+1)t+p+(k+1)i+v$ and $r=(k+1)t+i+j+lk+w$. This implies that
$t\le \min(\floor{(n-v)/(2(k+1))},\floor{(r-w)/(k+1)})$,
$l\le \floor{(r-w-(k+1)t)/k}$ and $i\le \floor{(n-2(k+1)t-v)/(k+1)}$,
whereas, for any $(t,l,i)$ triple, $j$ is determined by the second
equation given and $p$ by the first one. Hence, we have
\begin{equation}
  \label{eq:bxz_xnzr}
  [x^nz^r]a^{(u,v,w)}(x,z)=\sum_{t,l,i\ge 0} (-1)^i{p+t+u-1\choose p}{p\choose    i}{p-i\choose j}{i\choose l},
\end{equation}
where the upper limits of the summations on $t,l,$ and $i$ are given
by the three aforementioned inequalities, and $p$ and $j$ are
determined using the two equations above. We do not need to check for
nonnegativity of $p$ and $j$, as if that were the case the
corresponding terms in~\eqref{eq:bxz_xnzr} would cancel out.
Using~\eqref{eq:bxz_xnzr} and~\eqref{eq:SS_ltk_ogf_rewritten}, we have
an explicit triple-summation expression
for~$\widebar{s}_{\le k}(n,r)=[x^nz^r]\widebar{S}_{\le k}(x,y)$,
comparable to the triple-summation expression originally given by
Bateman for~$\widebar{\lambda}_{\ge k}(n,r)$~\cite[p.\
100]{bateman48:_power}.

\begin{remark}\label{rmk:bloom2}
  The ogfs~\eqref{eq:Ss_n_ltk_r_ogf}
  or~\eqref{eq:Ss_n_ltk_r_ogf_jackson}, or the explicit expression in
  the previous section, allow us to address the same problem that
  motivated Bloom's work in~\cite{bloom96:_probab}. Bloom's research
  was spurred by his desire to determine the accuracy of the following
  assertion by M. Gardner in~\cite[p.\ 124]{gardner82}:
  ``\textit{(\ldots) a shuffled deck of cards will contain
    coincidences. For instance, almost always there will be a clump of
    six or seven cards of the same color.}''  Gardner, like Bloom,
  uses the word ``clump'' to mean a run.  More specifically, what
  Bloom calls a ``$k$-clump'' is a run of length $k$ or longer of
  either colour of the deck, and he wishes to determine the likelihood
  that a shuffled deck will contain a $6$-clump.
  
  If we just look at the colour of the cards, then there is a
  bijection between the shuffles of a standard deck and the binary
  strings of length $52$ and Hamming weight $26$ (where~$0$ means
  ``red'' and $1$ means ``black'', or vice versa).  The number of
  $n$-strings devoid of $(\ge 6)$-runs of ones or zeros is
  $\widebar{s}_{\le 5}(52,26)=265,692,662,174,100$, and, thus, the
  probability that a shuffled deck of cards drawn uniformly at random
  contains at least one $6$-clump is
  $1-\widebar{s}_{\le 5}(52,26)/{52\choose 26}=0.464241$.

  This certainly coincides with Bloom's analysis, and it confirms this
  author's finding that Gardner's claim is hardly supported by the
  numbers. Bloom also used the expectation~\eqref{eq:exp_bloom} in
  Section~\ref{sec:total-number-runs-klek-hamming} as an indirect way
  to support his findings. However, as seems to have happened so many
  times in the history of the theory of runs, this author missed
  Bateman's results. Had he not missed them, he would simply have
  implemented the explicit expression given by Bateman~\cite[p.\
  100]{bateman48:_power} to answer his question. But had he done so,
  we would have missed on his insightful paper.

  Lastly, the tables for
  $\widebar{\lambda}_{\le k}(n,r)=\widebar{s}_{\le k}(n,r)/{n\choose
    r}$ given by Schuster~\cite[pp. 112--116]{schuster94:_exchan}
  ---in his notation, $\widebar{\lambda}_{\le k}(n,r)=F_M(k,n-r,r)$---
  can also be easily reproduced through Bateman's explicit expression
  (or, of course, through our own explicit expression or through the
  ogfs above).
\end{remark}

\subsubsection{OEIS}
Integer sequence ${52\choose 26}-\widebar{s}_{(k-1)}(52,26)$ is
\seqnum{A086438}.  On the other hand, the numerators and denominators
of the (simplified) fraction sequence
$1-\widebar{s}_{(k-1)}(52,26)/{52\choose 26}$ are, respectively,
\seqnum{A086439} and \seqnum{A086440}.

\subsection[Number of n-Strings of Hamming Weight r that Contain
Exactly m Nonnull p-Parity Runs of Ones and/or Zeros]{Number of
  $n$-Strings of Hamming Weight $r$ that Contain Exactly $m$ Nonnull
  $p$-Parity Runs of Ones and/or Zeros}
\label{sec:s_n_m_p_r}

To conclude this section, we enumerate the $n$-strings of Hamming weight $r$ that
exactly contain $m$ nonnull $p$-parity runs of ones and/or zeros, which we
denote by $s_{[p]}(n,m,r)$. The necessary condition that we require in
this case, which is just the appropriate variation
of~\eqref{eq:s_m_necessary_condition_r}, is given next.

\begin{condition}{(Existence of $n$-strings of Hamming weight $r$
    containing $m$ nonnull $p$-parity runs of ones and/or zeros)}
  \begin{equation}
    \label{eq:s_m_p_necessary_condition_r}
    s_{[p]}(n,m,r)> 0\quad \Longrightarrow\quad 0\le m \le \min\Bigg(\floor[\bigg]{\frac{n}{1+\iv{p=0}}},\floor[\bigg]{\frac{r}{1+\iv{p=0}}}+\floor[\bigg]{\frac{n-r}{1+\iv{p=0}}}\Bigg).
\end{equation}
\end{condition}

\subsubsection{Recurrence}
\label{sec:s_n_m_p_r_recs}
As in all previous cases in Section~\ref{sec:oz_hamming}, we start
by finding two mutual recurrences for this enumeration corresponding
to the $n$-strings that start with bit $b\in\{0,1\}$, which we denote
by $s^b_{[p]}(n,m,r)$. Through the standard procedure, these two
mutual recurrences are given by
\begin{equation}
  \label{eq:sb_n_m_p_r_rec}
  s^b_{[p]}(n,m,r)=\sum_{i=1}^ns^{\tilde{b}}_{[p]}\big(n-i,m-\iv{p=\text{mod}(i,2)},r-\iv{b=1}\,i\big),
\end{equation}
where the only difference with respect to~\eqref{eq:s0_p} is the
Hamming weight update when~$b=1$.

When $n=1$ we know by inspection that
\begin{equation}
  \label{eq:sb_n1_m_p_r_trivial}
  s^b_{[p]}(1,m,r)=\iv{r=b}\big(\iv{m=0}\iv{p=0}+\iv{m=1}\iv{p=1}\big),
\end{equation}
whereas setting $n=1$ in~\eqref{eq:sb_n_m_p_r_rec} yields
\begin{equation}
  \label{eq:sb_n1_m_p_r}
  s^b_{[p]}(1,m,r)=s^{\tilde{b}}_{[p]}\big(0,m-\iv{p=1},r-\iv{b=1}\big).
\end{equation}
We can verify that~\eqref{eq:sb_n1_m_p_r}
equals~\eqref{eq:sb_n1_m_p_r_trivial} for
\begin{equation}
  \label{eq:sb_n_m_p_r_rec_init}
  s^b_{[p]}(0,m,r)= \iv{m=0}\iv{r=0},
\end{equation}
which therefore constitute the initialisation
of~\eqref{eq:sb_n_m_p_r_rec}. Using these recurrences we can now compute
\begin{equation}
  \label{eq:s_n_m_r_from_mutual}
  s_{[p]}(n,m,r)= s^0_{[p]}(n,m,r)+s^1_{[p]}(n,m,r).
\end{equation}

Like in Sections~\ref{sec:s_n_m_r_rec} and~\ref{sec:ss_n_m_r_rec},  it is
also possible to obtain a recurrence on $s_{[p]}(n,m,r)$ by
symmetrising~\eqref{eq:sb_n_m_p_r_rec}. After adding
$\iv{n=0}\iv{m=0}\iv{r=0}$ to~\eqref{eq:sb_n_m_p_r_rec} to make it
valid for all values of the parameters, we may apply this recurrence
to itself to get
\begin{align}
  \label{eq:sb_n_m_p_r_rec_symm}
  s^b_{[p]}(n,m,r)=\sum_{i=1}^n\Bigg(&\sum_{j=1}^n s^{b}_{[p]}\big(n-(i+j),m-\iv{p=\text{mod}(i,2)}-\iv{p=\text{mod}(j,2)},r-\iv{b=1}i-\iv{b=0}j\big)\nonumber\\
&+\iv{n=i}\ivl{m=\iv{p=\text{mod}(i,2)}}\ivl{r=\iv{b=1}\,i}\Bigg)\nonumber\\
                               &+\iv{n=0}\iv{m=0}\iv{r=0},
\end{align}
where we have extended the summation on $j$ from $n-i$ to $n$ to
highlight the symmetries of the expression with respect to $i$ and
$j$. This is possible because of necessary
condition~\eqref{eq:s_m_p_necessary_condition_r}.  Using
next~\eqref{eq:sb_n_m_p_r_rec_symm} in~\eqref{eq:s_n_m_r_from_mutual}
we get the recurrence
\begin{align}
  \label{eq:s_n_m_p_r_rec_symm}
  s_{[p]}(n,m,r)=&\sum_{i=1}^{n}\sum_{j=1}^{n-i} s_{[p]}\big(n-(i+j),m-\iv{p=\text{mod}(i,2)}-\iv{p=\text{mod}(j,2)},r-i\big)\nonumber\\
&+\iv{n>0}\ivl{m=\iv{p=\text{mod}(n,2)}}\big(\iv{r=0}+\iv{r=n}\big)
                                   +2\iv{n=0}\iv{m=0}\iv{r=0}.
\end{align}
Like we did in~\eqref{eq:s_n_m_r_rec_sym}, in the recurrence above we have
again set the upper limit of the summation on~$j$ to $n-i$.

\subsubsection{Generating Function}
\label{sec:s_n_m_p_r_ogf}
The simplest way to obtain the ogf
$S_{[p]}(x,y,z)=\sum_{n,m,r}s_{[p]}(n,m,r)\,x^ny^mz^r$ is to proceed
as in Sections~\ref{sec:s_n_m_klek_r_generating-functions}
and~\ref{sec:s_n_m_r_ogf}. That is to say, we first get the ogfs
$S^b_{[p]}(x,y,z)=\sum_{n,m,r}s^b_{[p]}(n,m,r)$ in order to then
obtain
\begin{equation}
  \label{eq:Sp_ogf_from_mutual}
  S_{[p]}(x,y,z)=S^0_{[p]}(x,y,z)+S^1_{[p]}(x,y,z).
\end{equation}
First, we
 obtain versions of~\eqref{eq:sb_n_m_p_r_rec} without
$n$-dependent summations. We start by adding
$\iv{n=0}\iv{m=0}\iv{r=0}$ to~\eqref{eq:sb_n_m_p_r_rec} to make it
valid for all values of the parameters. We can then obtain
$s^b(n,m,r)-s^b(n-2,m,r-2\iv{b=1})$ using this extended
recurrence, which yields
\begin{align}
  \label{eq:sb_n_m_p_r_rec2}
  s^b_{[p]}(n,m,r)=&~s^b_{[p]}\big(n-2,m,r-2\iv{b=1}\big)+s^{\tilde{b}}_{[p]}\big(n-1,m-\iv{p=1},r-\iv{b=1}\big)\nonumber\\&+s^{\tilde{b}}_{[p]}\big(n-2,m-\iv{p=0},r-2\iv{b=1}\big)\nonumber\\
  &+\iv{n=0}\iv{m=0}\iv{r=0}-\iv{n=2}\iv{m=0}\ivl{r=2\iv{b=1}}.
\end{align}
Thus, multiplying this recurrence on both sides by $x^ny^mz^r$ and
adding over $n,m,$ and $r$ we get the following system of two
equations with two unknowns, i.e., $S^0_{[p]}(x,y,z)$ and $S^1_{[p]}(x,y,z)$:
\begin{align}
  \label{eq:sb_n_m_p_r_rec2}
  S^b_{[p]}(x,y,z)\big(1-x^2z^{2\iv{b=1}}\big)=&~S^{\tilde{b}}_{[p]}(x,y,z)\big(xy^{\iv{p=1}}z^{\iv{b=1}}+x^2y^{\iv{p=0}}z^{2\iv{b=1}}\big)+1-x^2z^{2\iv{b=1}}.
\end{align}
Solving the system and using~\eqref{eq:Sp_ogf_from_mutual}, it can be
seen that
\begin{equation}
  \label{eq:Sn_m_p_r_ogf}
  S_{[p]}(x,y,z)=\frac{2x^4z^2\big(y^{\iv{p=0}}-1\big)+x^3z(z+1)y^{\iv{p=1}}-x^2(z^2+1)\big(y^{\iv{p=0}}-2\big)-x(z+1)y^{\iv{p=1}}-2}{x^4z^2(y^{2\iv{p=0}}-1)+x^3z(z+1)y^{\iv{p=0}-\iv{p=1}}+x^2(zy^{2\iv{p=1}}+z^2+1)-1}.
\end{equation}

\section{Number of Runs Under Different Constraints Over All $n$-Strings and Over Restricted Subsets}
\label{sec:number-runs}

Our main goal in this section is enumerating how many runs under
different constraints are found over all $n$-strings. The most basic
results are the enumerations of the $(\ushort{k}\le\widebar{k})$-runs
over all $n$-strings in Section~\ref{sec:total-number-runs-klek} and
over all $n$-strings with Hamming weight $r$ in
Section~\ref{sec:total-number-runs-klek-hamming}, which we then
specialise to several particular cases.  In
Section~\ref{sec:number-runs-m} we study the number of runs over the
$n$-strings enumerated in Sections~\ref{sec:w_n_m_klek}
and~\ref{sec:wp_n_m}. These same problems are revisited in
Section~\ref{sec:number-runs-oz-m}, but with the goal of determining
the number of runs of ones \textit{and zeros} instead.  Several of the
questions addressed in this section, or variations of them, were
previously studied by Marbe~\cite{marbe16:_mathematische},
Gold~\cite{gold29:_note_frequen}, Cochran~\cite{cochran36:_plants},
Bloom~\cite{bloom96:_probab}, Sinha and
Sinha~\cite{sinha09,sinha12:_energy}, Makri et
al.~\cite{makri11:_bernoul,makri12:_count_runs_ones_ones_runs},
Nyblom~\cite{nyblom12:_enumerating}, and Grimaldi and
Heubach~\cite{grimaldi05:_without_odd_runs} ---and perhaps by other
authors we do not know of, since these are fairly common problems. In
fact, many of the authors just mentioned were not aware of previous
results.

Since runs of zeros appear at some points throughout this section, the
reader is again reminded that, unless explicitly said otherwise, the
term ``run'' without any qualifier refers to a run of ones.

\subsection[Number of (ḵ≤k)-Runs Over All $n$-Strings]{Number of
  $(\ushort{k}\le\widebar{k})$-Runs Over All $n$-Strings}
\label{sec:total-number-runs-klek}

We denote the number of $(\ushort{k}\le\widebar{k})$-runs over
all $n$-strings by
$\rho_{\ushort{k}\le\widebar{k}}(n)$. Considering~\eqref{eq:m_necessary_condition},
it is possible to compute
$\rho_{\ushort{k}\le\widebar{k}}(n)$ by relying on the enumeration
$w_{\ushort{k}\le\widebar{k}}(n,m)$ studied in Section~\ref{sec:w_n_m_klek} as follows:
\begin{equation}\label{eq:r_n_klek_from_w_n_m_klek}
  \rho_{\ushort{k}\le\widebar{k}}(n)=\sum_{m=1}^{\big\lfloor
      \frac{n+1}{\ushort{k}+1}\big\rfloor} m\,  w_{\ushort{k}\le\widebar{k}}(n,m).
\end{equation}
This summation requires an explicit expression for
$w_{\ushort{k}\le\widebar{k}}(n,m)$, which we have not
provided. However, it is also possible to obtain
$\rho_{\ushort{k}\le\widebar{k}}(n)$ directly from the
ogf~\eqref{eq:ogf_w_n_m_klek}:
\begin{equation}
  \label{eq:r_klek_through_ogf}
  \rho_{\ushort{k}\le\widebar{k}}(n)=[x^n]\left.\frac{\partial W_{\ushort{k}\le \widebar{k}}(x,y)}{\partial
  y}\right|_{y=1}.
\end{equation}
See that this is essentially the same as the computation of the first
moment with a pgf
---cf. Section~\ref{sec:moments_w_k}. From~\eqref{eq:ogf_w_n_m_klek}
we have that
\begin{equation}
  \label{eq:ddy_w_n_m_klek_ogf_y1}
  \left.\frac{\partial
    W_{\ushort{k}\le\widebar{k}}(x,y)}{\partial
    y}\right|_{y=1}=\frac{(1-x)\big(x^{\ushort{k}}-x^{\widebar{k}+1}\big)}{(1-2x)^2}.
\end{equation}
By  applying the negative binomial theorem, we can expand~\eqref{eq:ddy_w_n_m_klek_ogf_y1} as follows:
\begin{equation}
  \label{eq:ddy_w_n_m_klek_ogf_y1_expand}
  \left.\frac{\partial
    W_{\ushort{k}\le\widebar{k}}(x,y)}{\partial
    y}\right|_{y=1}=\big(x^{\ushort{k}}-x^{\ushort{k}+1}-x^{\widebar{k}+1}+x^{\widebar{k}+2}\big)\sum_{t\ge
  0}{t+1 \choose t} 2^t x^t.
\end{equation}
To extract the coefficient of $x^n$ we have to solve four equations
for the nonnegative index $t$: a) $n=\ushort{k}+t$, b)
$n=\ushort{k}+1+t$, c) $n=\widebar{k}+1+t$, and d)
$n=\widebar{k}+2+t$. 
When each of these equations has a solution,
from~\eqref{eq:ddy_w_n_m_klek_ogf_y1_expand} their respective contributions to the
coefficient of $x^n$ are:

\begin{enumerate}[label=\alph*)]
  
\item  $(n-\ushort{k}+1)\,2^{n-\ushort{k}}$
  
\item $-(n-\ushort{k})\,2^{n-(\ushort{k}+1)}$
 
\item  $-(n-\widebar{k})\,2^{n-(\widebar{k}+1)}$
  
\item $(n-\widebar{k}-1)\,2^{n-(\widebar{k}+2)}$
  
\end{enumerate}
\begin{itemize}
\item If $\widebar{k}+2\le n$, as
  $\ushort{k}\le\widebar{k}$ there is a solution in each of the four
  cases, so all four contributions must be added.
\item If $\widebar{k}+1=n$ only the first three equations have a
  solution, but all four contributions can still be added because the
  fourth one is zero anyway.

\item If $\widebar{k}\ge n$ then the last two equations do not have
  a solution, and we have the following cases for the first two:
\begin{itemize}[label=$\circ$]
\item If $\ushort{k}>n$ no equation has a solution.
\item If $\ushort{k}< n$  the first two equations have a
  solution, so the first two contributions must be added.
\item If $\ushort{k}=n$ only the first equation has a
  solution, but the first two contributions can still be added
  because the second one is zero anyway.
\end{itemize}
\end{itemize}
Collecting all these contributions we finally can see that
\begin{equation}
  \label{eq:r_n_klek_explicit}%
  \rho_{\ushort{k}\le\widebar{k}}(n)= (n-\ushort{k}+2)\,2^{n-(\ushort{k}+1)}\,\iv{\ushort{k}\le
    n}-(n-\widebar{k}+1)\,2^{n-(\widebar{k}+2)}\,\iv{\widebar{k}<n}.
\end{equation}
Observe that~\eqref{eq:r_n_explicit} is valid for all
$0\le \ushort{k}\le \widebar{k}$ and $n\ge 1$ ---for example, it can
enumerate the null runs over all $n$-strings if desired.

Although we will not show it here, \eqref{eq:r_n_klek_explicit} can
also be obtained by establishing a recurrence for
$\rho_{\ushort{k}\le\widebar{k}}(n)$ from first principles (i.e,
without resorting to any of the results in previous sections), and
then solving it directly by unrolling it. In any case, the approach
that we
have given above is quicker and cleaner.

\subsubsection{Special Cases}
\label{sec:number-runs-klek-special-cases}
In this section we discuss several specialisations of
$\rho_{\ushort{k}\le \widebar{k}}(n)$. We denote the number of $k$-runs over all $n$-strings by $\rho_k(n)$,
and the number of nonnull runs (of arbitrary lengths, all strictly
greater than zero) over all $n$-strings by $\rho(n)$. These
enumerations are special cases of our analysis in the previous
section, as we can write $\rho_k(n)=\rho_{k\le k}(n)$ and
$\rho(n)=\rho_{1\le n}(n)$. Thus, using~\eqref{eq:r_n_klek_explicit} we have that
\begin{equation}
  \label{eq:r_n_k_explicit}
  \rho_k(n)=(n-k+3)\,2^{n-(k+2)}\,\iv{k<n}+\iv{k=n},
\end{equation}
and
\begin{equation}
  \label{eq:r_n_explicit}
  \rho(n)=(n+1)\,2^{n-2}.
\end{equation}
Next, we denote the number of $(\ge\! k)$-runs and of nonnull
$(\le\! k)$-runs over all $n$-strings by $\rho_{\ge k}(n)$ and
by $\rho_{\le k}(n)$, respectively. Once again, these are special cases of
our analysis in the previous section, because we can write
$\rho_{\ge k}(n)=\rho_{k\le n}(n)$ and
$\rho_{\le k}(n) =\rho_{1\le k}(n)$.  Therefore,
from~\eqref{eq:r_n_klek_explicit} we have 
\begin{equation}\label{eq:rgtk_explicit}
  \rho_{\ge k}(n)=(n-k+2)\,2^{n-(k+1)}\,\iv{k\le n}
\end{equation}
and
\begin{equation}
  \label{eq:rltk_explicit}
  \rho_{\le k}(n) =%
  \rho(n)-(n-k+1)\,2^{n-(k+2)}\,\iv{k<n},
\end{equation}
where $\rho(n)$ is given by~\eqref{eq:r_n_explicit}.

\begin{remark} \label{rmk:rk-r-alt} Leaving aside Marbe's formula for
  a moment ---see below--- the first author who
  produced~\eqref{eq:r_n_k_explicit} was
  Gold~\cite{gold29:_note_frequen,cochran38:_golds}, in his analysis
  of the predictability of two-state meteorological series. This
  author gave $\rho_k(n)/2^n$ as the expected number of $k$-runs in an
  $n$-string drawn uniformly at random ---he actually obtained twice
  the value of this expectation, as he addressed runs of ones and
  zeros jointly. Gold arrived at his result by combining
  $\sum_{k=1}^n k\, \rho_k(n)=n 2^{n-1}$ (see
  Remark~\ref{rmk:onesandones}) and $\rho_{k-1}(n-1)=\rho_k(n)$, an
  identity which can easily be argued.  Gold's expectation was then
  extended by Cochran~\cite[Eq.\ (5)]{cochran36:_plants} to the case
  where the bits are not equally likely ---Cochran, who was studying
  the spread of diseases in plants arranged in rows, acknowledges
  Marbe's priority regarding his result~\cite{marbe16:_mathematische},
  but he indicates that the original proof~\cite[p.\
  9]{marbe34:_grundfragen} omits essential steps, and gives his own
  proof by induction. Cochran's expectation in~\cite[Eq.\
  (5)]{cochran36:_plants} is in fact~\eqref{eq:exp_mkn_all} in
  Section~\ref{sec:moments_w_k}, which simplifies
  to~\eqref{eq:r_n_k_explicit} using $q=1/2$ and multiplying by $2^n$.
  
  We briefly discuss other approaches for the derivation of some of
  these enumerations, as we retrace the work of other authors who
  previously dealt with them. These authors were unaware of the
  aforementioned results. As indicated by Sinha and
  Sinha~\cite{sinha09} and by Makri and Psillakis~\cite[Eq.\
  (9)]{makri11:_bernoul}, one way to obtain $\rho_k(n)$ is
  \begin{equation}\label{eq:rnk_wknm}
    \rho_k(n)=\sum_{m=1}^{\floor*{ \frac{n+1}{k+1}}} m\, w_k(n,m),
  \end{equation}
  which is the specialisation
  of~\eqref{eq:r_n_klek_from_w_n_m_klek}. However, this is not an easy
  approach ---consider using~\eqref{eq:w_gen}
  in~\eqref{eq:rnk_wknm}. Tellingly, Makri et al. were only able to
  evaluate an expression similar to~\eqref{eq:rnk_wknm} to
  obtain~\eqref{eq:r_n_k_explicit}~\cite[Eq.\ (7),
  $a=e$]{makri12:_count_runs_ones_ones_runs} by using an earlier
  probabilistic result of theirs, but not by using their own explicit
  expression for $w_k(n,m)$~\cite[Eq. (8)]{makri11:_bernoul}
  in~\eqref{eq:rnk_wknm}. 
  As mentioned at the end of the previous section, it is possible to
  establish a recurrence for $\rho_{\ushort{k}\le\widebar{k}}(n)$
  \textit{ab initio}. In the special case $\ushort{k}=\widebar{k}=k$,
  this recurrence is
  \begin{equation}
    \rho_k(n) =2\,\rho_k(n-1)+2^{n-(k+2)}  \label{eq:recurrence_rk}
  \end{equation}
  for $k+2\le n$, whereas $\rho_{n-1}(n)=2$ and $\rho_n(n)=1$.
  Recurrence~\eqref{eq:recurrence_rk} and explicit
  expression~\eqref{eq:r_n_k_explicit} were worked out from first
  principles by Sinha and Sinha in~\cite{sinha12:_energy} ---according
  to~\cite{sinha09}, an earlier appearance of these authors' results
  regarding $\rho_k(n)$ is in K.~Sinha's PhD
  thesis~\cite{sinha07:_locat}. Sinha and Sinha
  cite~\cite{cochran38:_golds} in passing in~\cite{sinha12:_energy},
  without realising that Cochran's article already
  contains the closed-form solution~\eqref{eq:r_n_k_explicit} they
  find.  Makri et al.\ \cite[Eq. (9),
  $a=e$]{makri12:_count_runs_ones_ones_runs} also noticed
  recurrence~\eqref{eq:recurrence_rk} after finding the explicit
  expression~\eqref{eq:r_n_k_explicit} for $\rho_k(n)$. We would also
  like to mention that we gave two other \textit{ab initio} approaches
  to the problem of finding $\rho_k(n)$ in~\cite{balado2023runs}.
  
  We should mention that, after determining the closed-form
  expression~\eqref{eq:r_n_k_explicit} for $\rho_k(n)$ through any of
  the methods just mentioned, one can then obtain
  $\rho_{\ushort{k}\le
    \widebar{k}}(n)=\sum_{k=\ushort{k}}^{\widebar{k}} \rho_k(n)$
  ---which only involves summations of the form $\sum_k a^k$ and
  $\sum_k k\, a^k$. While this makes the approach in the previous
  section unnecessary, we believe that the method that we have adopted
  is interesting because it connects with the enumerations in
  Section~\ref{sec:number-n-strings}. It also allows for a consistent
  methodology throughout the paper
  ---cf. Sections~\ref{sec:w_n_m_klek},~\ref{sec:probability-w_m_klek},~\ref{sec:w_n_m_klek_hamming},~\ref{sec:s_n_m_klek},~\ref{sec:probability-s_m_klek},~\ref{sec:s_n_m_klek_hamming}
  and~\ref{sec:total-number-ones-klek}.

  As regards $\rho(n)$, this enumeration may be obtained in a number
  of different ways as well. One possibility is to use $\rho_k(n)$ to
  compute $\rho(n)=\sum_{k=1}^n \rho_k(n)$, as pointed out by Sinha
  and Sinha~\cite{sinha12:_energy}. Using~\eqref{eq:r_n_k_explicit} we
  can thus write
  \begin{equation}
    \label{eq:total_1}
    \rho(n)=1+\sum_{k=1}^{n-1} (n-k+3)\,2^{n-(k+2)},
  \end{equation} 
  which, of course, evaluates to~\eqref{eq:r_n_explicit}. But we may also
  get $\rho(n)$ without resorting to $\rho_k(n)$.  One way is by
  specialising~\eqref{eq:r_n_klek_from_w_n_m_klek} with $\ushort{k}=1$
  and $\widebar{k}=n$ to get
  $\rho(n)=\sum_{m=1}^{\floor*{(n+1)/2}} m\, w(n,m)$, which allows us
  to obtain this enumeration by using~\eqref{eq:w_n_m_explicit}:
  \begin{align*}\label{eq:total}
    \rho(n)    &=\sum_{m=1}^{\floor*{ \frac{n+1}{2}}} m\, {n+1 \choose 2m}.
  \end{align*}
  This summation yields~\eqref{eq:r_n_explicit}, a fact that one may
  verify using generating functions.  Similarly, we can get $\rho(n)$
  from the number of $n$-strings with exactly $m$ nonnull runs of ones
  and/or zeros, $s(n,m)$. By symmetry,
  $\rho(n)=(1/2)\sum_{m=1}^n m\,s(n,m)$.
  Inputting~\eqref{eq:s_n_m_explicit} in this expression we get
  \begin{equation}
    \label{eq:total_alt}
    \rho(n)=\sum_{m=1}^{n} m \,{n-1\choose m-1},
  \end{equation}
  and we may verify that this summation also
  gives~\eqref{eq:r_n_explicit}.  Lastly, Nyblom~\cite[Lem.\
  3.1]{nyblom12:_enumerating} also found a recurrence for $\rho(n)$
  which, of course, can also be derived from~\eqref{eq:recurrence_rk},
  and solved it to get~\eqref{eq:r_n_explicit} ---in actual fact, Nyblom's
  results are for $2\,\rho(n)$, as this author counts both runs of
  ones and of zeros.

  Finally, we make a few comments about $\rho_{\ge
    k}(n)$. Expression~\eqref{eq:rgtk_explicit} can also be obtained
  from expectation~\eqref{eq:p_n_gtk_expectation} in
  Section~\ref{sec:moments_w_gtk}, by setting $q=1/2$ and then
  multiplying by $2^n$.  Bloom gave~\eqref{eq:rgtk_explicit} divided
  by $2^{n-1}$ ---i.e., the expected number of $(\ge k)$-runs of ones
  and/or zeros in an $n$-string drawn uniformly at random~\cite[Eq.\
  (12)]{bloom96:_probab}. He most likely also used a probabilistic
  approach, but he gives no details.  Another way to obtain
  $\rho_{\ge k}(n)$ is by
  specialising~\eqref{eq:r_n_klek_from_w_n_m_klek} with $\ushort{k}=k$
  and $\widebar{k}=n$, which yields
  $\rho_{\ge k}(n)=\sum_{m=1}^{\floor*{(n+1)(k+1)}} m\, w_{\ge
    k}(n,m)$.  Although it is again very difficult to handle this
  summation using an explicit expression for $w_{\ge k}(n,m)$, such
  as~\eqref{eq:w_gtk_gen_diophant}, a probabilistic version of this
  expression was used by Makri et al. to
  obtain~\eqref{eq:rgtk_explicit}~\cite[Eq.\ (7),
  $a=g$]{makri12:_count_runs_ones_ones_runs}.  These authors also
  deduced a recurrence working backwards from the explicit
  expression~\cite[Eq.\ (9),
  $a=g$]{makri12:_count_runs_ones_ones_runs}. 
\end{remark}

\begin{remark}\label{rmk:more_compositions}
  There are some immediate connections between the results in this
  section and compositions. For example, the total number of parts in
  all compositions of $n$ is the same as the total number of nonnull
  runs over all $n$-strings, $\rho(n)$, and thus given
  by~\eqref{eq:r_n_explicit}. In similar fashion, the total number of
  parts equal to $k\ge 1$ in all compositions of $n$ equals the total
  number of runs of length $k$ over all $n$-strings, $\rho_k(n)$,
  which is given by~\eqref{eq:r_n_k_explicit}. In the special case
  $k=1$, both recurrence~\eqref{eq:recurrence_rk} and closed-form
  expression~\eqref{eq:r_n_k_explicit} were given by Chin et
  al.~\cite[Eqs.\ (1) and (2)]{chinn92:_cuisenaire}. These authors
  studied a problem concerning Cuisenaire rods (number rods) which is
  equivalent to enumerating all parts equal to $1$ over all
  compositions of $n$.
\end{remark}

\subsubsection{OEIS}
\label{sec:oeis-rgtkn}
We report next sequences connected to the enumerations in this
section that are found in  the OEIS:
\begin{itemize}
\item All sequences emanating from $\rho_k(n)$, $\rho_{\ge k}(n)$, and $\rho(n)$
  are in the OEIS:

  $\rho_k(n+k)$ is \seqnum{A045623} 
  (number of $1$'s in all compositions of $n+1$).

  $\rho_{\ge k}(n+k)$ and $\rho(n+1)$ are \seqnum{A001792}. %

\item Sequences emanating from $\rho_{\le k}(n)$:

  $\rho_{\le 1}(n+1)$ is \seqnum{A045623} %
  ---just like $\rho_k(n+k)$.
  
  $\rho_{\le 2}(n+1)$ is \seqnum{A106472}. %

  $\rho_{\le 3}(n)$ is \seqnum{A386878}. %

\end{itemize}

\subsection[Number of (ḵ≤k)-Runs Over All n-Strings of Hamming
Weight r]{Number of
  $(\ushort{k}\le\widebar{k})$-Runs Over All $n$-Strings of Hamming
  Weight~$r$}
\label{sec:total-number-runs-klek-hamming}
In this section we obtain the total number of
$(\ushort{k}\le\widebar{k})$-runs over all $n$-strings of Hamming
weight $r$, which we denote by  $\rho_{\ushort{k}\le\widebar{k}}(n,r)$. We
can obtain this quantity directly from \eqref{eq:ogf_w_n_m_k_r} as
follows:
\begin{equation}
  \label{eq:r_n_klek_r_from_w_n_m_klek_r_ogf}
  \rho_{\ushort{k}\le\widebar{k}}(n,r)=[x^nz^r]\left.\frac{\partial
      W_{\ushort{k}\le\widebar{k}}(x,y,z)}{\partial y}\right|_{y=1}.
\end{equation}
From~\eqref{eq:ogf_w_n_m_k_r} we have that
\begin{equation}
  \label{eq:ddy_w_n_m_klek_r_ogf_y1}
  \left.\frac{\partial
    W_{\ushort{k}\le\widebar{k}}(x,y,z)}{\partial
    y}\right|_{y=1}=\frac{(1-xz)\big((xz)^{\ushort{k}}-(xz)^{\widebar{k}+1}\big)}{(1-x-xz)^2}.
\end{equation}
By rewriting this expression as
$((xz)^{\ushort{k}}-(xz)^{\widebar{k}+1})/\big((1-xz)(1-x/(1-xz))^2\big)$
and then applying the negative binomial theorem twice, we can
expand it as follows:
\begin{equation}
  \label{eq:ddy_w_n_m_klek_r_ogf_y1_expand}
  \left.\frac{\partial
      W_{\ushort{k}\le\widebar{k}}(x,y,z)}{\partial
      y}\right|_{y=1}=\big((xz)^{\ushort{k}}-(xz)^{\widebar{k}+1}\big)\sum_{t\ge
    0}\sum_{s\ge 0}{t+1\choose t}{s+t\choose t} x^t (xz)^s.
\end{equation}
To extract the coefficient of $x^nz^r$ we just have to find the
nonnegative indices $t$ and $s$ that fulfil two sets of equations: on
the one hand, $n=\ushort{k}+t+s$ and $r=\ushort{k}+s$; on the other
hand, $n=\widebar{k}+1+t+s$ and $r=\widebar{k}+1+s$. The solutions of
the first set are $t=n-r$ and $s=r-\ushort{k}$, whereas the solutions
of the second one are $t=n-r$ and $s=r-\widebar{k}-1$. For a solution
to exist, both sets require $n\ge r$, whereas $r\ge \ushort{k}$ and
$r\ge \widebar{k}+1$ are required in the first and second set,
respectively. From these considerations and
\eqref{eq:ddy_w_n_m_klek_r_ogf_y1_expand} we thus have that
\begin{equation}
  \label{eq:r_n_klek_r}
  \rho_{\ushort{k}\le\widebar{k}}(n,r)=(n-r+1)\bigg({n-\ushort{k}\choose
    r-\ushort{k}}-{n-\widebar{k}-1\choose
    r-\widebar{k}-1}\bigg)\,\iv{n\ge r}.
\end{equation}
Observe that the constraints on $r$ are taken care of by the binomial
coefficients.  Like~\eqref{eq:r_n_klek_explicit}, this expression is
valid for all $0\le \ushort{k}\le \widebar{k}$ and $n\ge 1$, apart
from all $r$.

\subsubsection{Special Cases}
\label{sec:total-number-runs-klek-hamming-special}
We consider next several specialisations of
$\rho_{\ushort{k}\le \widebar{k}}(n,r)$. We denote the number of
$(\ge k)$-runs and the total number of nonnull runs over all
$n$-strings of Hamming weight $r$ by $\rho_{\ge k}(n,r)$ and
$\rho(n,r)$, respectively. These enumerations are special cases
of~\eqref{eq:r_n_klek_r}, as $\rho_{\ge k}(n,r)=\rho_{k\le n}(n,r)$
and $\rho(n,r)=\rho_{1\le n}(n,r)=\rho_{\ge 1}(n,r)$. Thus,
from~\eqref{eq:r_n_klek_r} we see that
\begin{equation}
  \label{eq:r_n_gtk_r}
  \rho_{\ge k}(n,r)=(n-r+1){n-k\choose r-k},
\end{equation}
whereas
\begin{equation}
  \label{eq:r_n_r}
  \rho(n,r)=(n-r+1){n-1\choose r-1}.
\end{equation}
Notice that we do not need Iverson brackets in these two expressions.
Also, denote the number of nonnull $(\le k)$-runs and the number of
$k$-runs over all $n$-strings of Hamming weight $r$ by
$\rho_{\le k}(n,r)$ and $\rho_k(n,r)$, respectively. As above, these
enumerations are special cases of~\eqref{eq:r_n_klek_r}, because
$\rho_{\le k}(n,r)=\rho_{1\le k}(n,r)$ and
$\rho_k(n,r)=\rho_{k\le k}(n,r)=\rho_{\ge 1}(n,r)$. Thus,
from~\eqref{eq:r_n_klek_r} we have that
\begin{equation}
  \label{eq:r_n_ltk_r}
  \rho_{\le k}(n,r)=(n-r+1)\bigg({n-1\choose
    r-1}-{n-k-1\choose r-k-1}\bigg)\,\iv{n\ge r},
\end{equation}
and
\begin{equation}
  \label{eq:r_n_k_r}
  \rho_k(n,r)=(n-r+1){n-k-1\choose r-k}\,\iv{n\ge r}.
\end{equation}

\begin{remark}
  If we denote by $M_{\ge k,n,r}$ the rv that models the number of
  $(\ge k)$-runs in an $n$-string of Hamming weight $r$ drawn
  uniformly at random, then its expectation is
  $\text{E}(M_{\ge k,n,r})=\rho_{\ge k}(n,r)/{n\choose r}$. Similarly,
  if $\fixwidetilde{M}_{\ge k,n,r}$ models the number of
  $(\ge k)$-runs of ones and/or zeros in an $n$-string of Hamming
  weight $r$ drawn uniformly at random then
\begin{align}
  \label{eq:exp_bloom_like}
  \text{E}\big(\fixwidetilde{M}_{\ge k,n,r}\big)&=\frac{1}{{n\choose r}}\big(\rho_{\ge
                           k}(n,r)+\rho_{\ge k}(n,n-r)\big)\nonumber\\
  &=\frac{1}{{n\choose r}}\bigg((n-r+1){n-k\choose n-r}+(r+1){n-k\choose r}\bigg).
\end{align}
Bloom used a probabilistic rationale to give the following expression
for $\text{E}\big(\fixwidetilde{M}_{\ge k,n,r}\big)$~\cite[p.\ 370]{bloom96:_probab}:
\begin{align}
  \label{eq:exp_bloom}
  \text{E}\big(\fixwidetilde{M}_{\ge k,n,r}\big)
  &=\frac{1}{{n\choose k}}\bigg((n-r+1){r\choose k}+(r+1){n-r\choose k}\bigg),
\end{align}
where we have divided by $k!$ both the numerator and the denominator
of Bloom's original expression, for ease of comparison
with~\eqref{eq:exp_bloom_like}.  As
${n-k\choose n-r}{n\choose k}={r\choose k}{n\choose r}$ and
${n-k\choose r}{n\choose k}={n-r\choose k}{n\choose r}$ ---which is
verified simply by developing the factorials in all these binomial
coefficients--- then we can see that~\eqref{eq:exp_bloom_like} is the
same as Bloom's expression~\eqref{eq:exp_bloom}.
\end{remark}

\subsection{Number of Nonnull Runs Over All $n$-Strings that Contain
  Exactly $m$ Nonnull Runs Under Different Constraints}
\label{sec:number-runs-m}
Next, we deal with the problems of enumerating the number of nonnull
runs over all $n$-strings that contain exactly $m$ nonnull
$(\ushort{k}\le\overline{k})$-runs, or exactly $m$ nonnull $p$-parity
runs. In keeping with our notation conventions, we call these
quantities $\rho_{\ushort{k}\le\overline{k}}(n,m)$ and
$\rho_{[p]}(n,m)$, respectively.

These problems are not as elementary as the ones we have dealt with in
the previous two sections. Relatively simple explicit expressions do
not seem possible in general, although they are achievable in
particular cases ---see Section~\ref{sec:oeis-rnm}. However, we can
obtain recurrence relations and generating functions for these
enumerations by relying on their counterparts in
Sections~\ref{sec:w_n_m_klek} and~\ref{sec:wp_n_m}, i.e.,
$w_{\ushort{k}\le\overline{k}}(n,m)$ and $w_{[p]}(n,m)$, respectively.

\subsubsection[Number of Nonnull Runs Over All the n-Strings that
Contain Exactly m Nonnull (ḵ≤k)-Runs]{Number of Nonnull Runs Over All
  $n$-Strings that Contain Exactly $m$ Nonnull
  $(\ushort{k}\le \widebar{k})$-Runs}
\label{sec:r_n_m_klek}

We obtain first a recurrence for $\rho_{\ushort{k}\le\widebar{k}}(n,m)$,
by splitting the contributions due $n$-strings with different initial
bit. As indicated, we assume $\ushort{k}\ge 1$. The $n$-strings that
start with $0$ contribute $\rho_{\ushort{k}\le\widebar{k}}(n-1,m)$, to
$\rho_{\ushort{k}\le\widebar{k}}(n,m)$. On the other hand, the
$n$-strings that start with a nonnull $i$-run with
$\ushort{k}\le i\le \widebar{k}$ contribute
$\rho_{\ushort{k}\le\widebar{k}}(n-(i+1),m-1)$ runs plus the
$w_{\ushort{k}\le\widebar{k}}(n-(i+1),m-1)$ initial runs
themselves. If they start with a nonnull $i$-run with $i<\ushort{k}$
or $i>\widebar{k}$, then they contribute
$\rho_{\ushort{k}\le\widebar{k}}(n-(i+1),m)$ runs plus the
$w_{\ushort{k}\le\widebar{k}}(n-(i+1),m)$ initial runs. This yields
the recurrence
\begin{align}
  \label{eq:r_klek_n_m_rec}
  \rho_{\ushort{k}\le\widebar{k}}(n,m)= \rho_{\ushort{k}\le\widebar{k}}(n-1,m)&+\sum_{i=\ushort{k}}^{\widebar{k}} \Big(
              \rho_{\ushort{k}\le\widebar{k}}(n-(i+1),m-1)+w_{\ushort{k}\le\widebar{k}}(n-(i+1),m-1)\Big)\nonumber\\
            &+ \sum_{i=1}^{\ushort{k}-1} \Big(
              \rho_{\ushort{k}\le\widebar{k}}(n-(i+1),m)+w_{\ushort{k}\le\widebar{k}}(n-(i+1),m)\Big)\nonumber\\
            &+ \sum_{i=\widebar{k}+1}^{n} \Big(
              \rho_{\ushort{k}\le\widebar{k}}(n-(i+1),m)+w_{\ushort{k}\le\widebar{k}}(n-(i+1),m)\Big).
\end{align}
To initialise the recurrence we use the case $n=1$, in which we
can see by inspection that
\begin{align}
  \label{eq:r_klek_init_n1}
  \rho_{\ushort{k}\le\widebar{k}}(1,m)=&~\iv{m=0}\iv{\ushort{k}>1}+\iv{m=1}\iv{\ushort{k}=1}
\end{align}
On the other hand, specialising recurrence~\eqref{eq:r_klek_n_m_rec}
for $n=1$ and considering~\eqref{eq:m_necessary_condition} we get
\begin{align}\label{eq:r_klek_rec_n1}
  \rho_{\ushort{k}\le\widebar{k}}(1,m)=~&\rho_{\ushort{k}\le\widebar{k}}\big(0,m\big)+
 \rho_{\ushort{k}\le\widebar{k}}\big(-1,m-\iv{\ushort{k}=1}\big)+w_{\ushort{k}\le\widebar{k}}\big(-1,m-\iv{\ushort{k}=1}\big).
\end{align}
We wish~\eqref{eq:r_klek_rec_n1} to equal~\eqref{eq:r_klek_init_n1}.
Taking~\eqref{eq:w_n_m_klek_init} and~\eqref{eq:w_n_m_klek_init2} into
account, equality between these two expressions is fulfilled for all
$m$ and $1\le\ushort{k}\le \widebar{k}$ by choosing
\begin{align}
  \label{eq:r_klek_n_m_init}
  \rho_{\ushort{k}\le\widebar{k}}(-1,m)=\rho_{\ushort{k}\le\widebar{k}}(0,m)=0,
\end{align}
which are therefore the initial values of~\eqref{eq:r_klek_n_m_rec}.

Recurrence~\eqref{eq:r_klek_n_m_rec} also allows us to find the ogf
$R_{\ushort{k}\le\widebar{k}}(x,y)=\sum_{n,m}\rho_{\ushort{k}\le\widebar{k}}(n,m)\,x^ny^m$
in terms of the ogf $W_{\ushort{k}\le\widebar{k}}(x,y)$
in~\eqref{eq:ogf_w_n_m_klek}. Recurrence~\eqref{eq:r_klek_n_m_rec} is
valid already for all values of the parameters, so we just obtain the
difference
$\rho_{\ushort{k}\le\widebar{k}}(n,m)-\rho_{\ushort{k}\le\widebar{k}}(n-1,m)$
using~\eqref{eq:r_klek_n_m_rec} which yields
\begin{align*}
  \rho_{\ushort{k}\le\widebar{k}}(n,m)=~&2\, \rho_{\ushort{k}\le\widebar{k}}(n-1,m)+w_{\ushort{k}\le\widebar{k}}(n-2,m)\\
                       &+\rho_{\ushort{k}\le\widebar{k}}(n-(\ushort{k}+1),m-1)+w_{\ushort{k}\le\widebar{k}}(n-(\ushort{k}+1),m-1)\\
                       &-\rho_{\ushort{k}\le\widebar{k}}(n-(\ushort{k}+1),m)-w_{\ushort{k}\le\widebar{k}}(n-(\ushort{k}+1),m)\\
                       &-\rho_{\ushort{k}\le\widebar{k}}(n-(\widebar{k}+2),m-1)-w_{\ushort{k}\le\widebar{k}}(n-(\widebar{k}+2),m-1)\\
                       &+\rho_{\ushort{k}\le\widebar{k}}(n-(\widebar{k}+2),m)+w_{\ushort{k}\le\widebar{k}}(n-(\widebar{k}+2),m).
\end{align*}
By multiplying this expression on both sides by $x^ny^m$ and then
adding on $n$ and $m$, we get
\begin{equation}
  \label{eq:r_n_m_klek_ogf}
  R_{\ushort{k}\le\widebar{k}}(x,y)=\frac{x^2-(1-y)(x^{\ushort{k}+1}-x^{\widebar{k}+2})}{1-2x+(1-y)(x^{\ushort{k}+1}-x^{\widebar{k}+2})} W_{\ushort{k}\le\widebar{k}}(x,y),
\end{equation}
where $W_{\ushort{k}\le\widebar{k}}(x,y)$ is given
by~\eqref{eq:ogf_w_n_m_klek}.

\subsubsection{Number of Nonnull Runs Over All $n$-Strings that
  Contain Exactly $m$ Nonnull $p$-Parity Runs}
\label{sec:r_n_m_p}
We now give a recurrence relation for $\rho_{[p]}(n,m)$. Following the
same strategy as in the previous section, the $n$-strings that start
with $0$ contribute $\rho_{[p]}(n-1,m)$ runs to $\rho_{[p]}(n,m)$. As for
the $n$-strings that start with $1$, if they start with an odd $i$-run
then they contribute $\rho_{[p]}(n-(i+1),m-\iv{p=1})$ runs plus the
$w_{[p]}(n-(i+1),m-\iv{p=1})$ initial runs themselves. If they start
with an even $i$-run then they contribute
$\rho_{[p]}(n-(i+1),m-\iv{p=0})$ runs plus the
$w_{[p]}(n-(i+1),m-\iv{p=0})$ initial runs. This yields the recurrence
\begin{align}
  \label{eq:r_n_m_p_rec}
  \rho_{[p]}(n,m)=~&\rho_{[p]}(n-1,m)+ \sum_{i=1}^{n}
                  \rho_{[p]}(n-(i+1),m-\iv{p=\text{mod}(i,2)})\nonumber\\
  &+ \sum_{i=1}^{n}w_{[p]}(n-(i+1),m-\iv{p=\text{mod}(i,2)}).
\end{align}
To initialise~\eqref{eq:r_n_m_p_rec} we use the value
at $n=1$, which, by inspection, is
\begin{align}
  \label{eq:r_n_m_p_trivial_n1}
  \rho_{[p]}(1,m)=\iv{p=0}\iv{m=0}+\iv{p=1}\iv{m=1}.
\end{align}
On the other hand, setting $n=1$ in~\eqref{eq:r_n_m_p_rec} we have
that
\begin{align}
  \label{eq:r_n_m_p_rec_n1}
  \rho_{[p]}(1,m)=\rho_{[p]}(0,m)+\rho_{[p]}(-1,m-\iv{p=1})+w_{[p]}(-1,m-\iv{p=1}).
\end{align}
We want~\eqref{eq:r_n_m_p_rec_n1} to
equal~\eqref{eq:r_n_m_p_trivial_n1}. 
Taking~\eqref{eq:m_necessary_condition} and~\eqref{eq:wp_n_m_init}
into account, we see that equality between the two expressions is
achieved for all $m$ by choosing
\begin{equation*}
  \label{eq:r_n_m_p_init}
  \rho_{[p]}(-1,m)=\rho_{[p]}(0,m)=0,
\end{equation*}
which are therefore the initial values of~\eqref{eq:r_n_m_p_rec}.

We may next obtain the ogf
$R_{[p]}(x,y)=\sum_{n,m} \rho_{[p]}(n,m)\,x^n y^m$ in terms of the ogf
$W_{[p]}(x,y)$ in~\eqref{eq:wp_bivariate_ogf}. We first obtain the
following difference using~\eqref{eq:r_n_m_p_rec}:
\begin{align}\label{eq:r_n_m_p_rec2}
  \rho_{[p]}(n,m)-\rho_{[p]}(n-2,m)=~&\rho_{[p]}(n-1,m)-\rho_{[p]}(n-3,m)\nonumber\\
                               &+\rho_{[p]}\big(n-2,m-\iv{p=1}\big)%
                               +\rho_{[p]}\big(n-3,m-\iv{p=0}\big)\nonumber\\
                               &+w_{[p]}\big(n-2,m-\iv{p=1}\big)%
                               +w_{[p]}\big(n-3,m-\iv{p=0}\big).
\end{align}
As usual, by multiplying by $x^ny^m$ on both sides of this expression
and then adding on $n$ and $m$ and using~\eqref{eq:wp_bivariate_ogf} we obtain
\begin{equation*}
  \label{eq:r_n_m_p_ogf}
  R_{[p]}(x,y)=\frac{x^2y^{\iv{p=1}}+x^3y^{\iv{p=0}}}{1-x-x^2\big(1+y^{\iv{p=1}}\big)+x^3\big(1-y^{\iv{p=0}}\big)} W_{[p]}(x,y).
\end{equation*}

\begin{remark}\label{rmk:grimaldi_r}
  Grimaldi and Heubach~\cite[Thm.\ 4]{grimaldi05:_without_odd_runs}
  gave a recurrence, an explicit expression and a generating function
  for the total number of runs of ones \textit{and zeros} in the
  $w_{[1]}(n,0)$ $n$-strings devoid of odd runs of zeros---in their
  notation,~$t_n$. Notice that $\rho_{[1]}(n,0)$ solely gives the total
  number of runs of ones in $n$-strings devoid of odd runs of ones. In
  any case, the corresponding generating functions are very similar.
  In the specific case of $m=0$ and $p=1$,~\eqref{eq:r_n_m_p_rec2}
  becomes
  \begin{equation*} %
    \rho_{[1]}(n,0)-\rho_{[1]}(n-2,0)= \rho_{[1]}(n-1,0)+w_{[1]}\big(n-3,0\big).
  \end{equation*}
  From here,  we have that
  \begin{align}
    [y^0]R_{[1]}(x,y)&=\frac{x^3}{1-x-x^2}\, [y^0]W_{[1]}(x,y)\label{eq:r1xy_a}\\
                     &=\frac{x^2-x^4}{(1-x-x^2)^2},\label{eq:r1xy_b}
  \end{align}
  where we have used~\eqref{eq:ogf_wp_m} in the second step.  For
  comparison's sake, the related ogf given
  in~\cite{grimaldi05:_without_odd_runs} is
  $G_{t_n}(x)=(x-x^4)/(1-x-x^2)^2$. See also
  Remark~\ref{rmk:grimaldi_r_2}.
\end{remark}

\subsubsection{OEIS}
\label{sec:oeis-rnm}

The only sequences arising from the enumerations in
Section~\ref{sec:number-runs-m} that we have been able to find in the
OEIS are:

$\rho_{1}(n-1,0)$ is \seqnum{A136444}  $\big(\sum_{k=0}^n k{n-k\choose  2k}\big)$
for $n\ge 1$. %

$\rho_{\ge 1}(n,1)=w_{\ge 1}(n,1)$ is \seqnum{A000217}. %

$\rho_{\ge 1}(n-1,2)$ is  \seqnum{A034827}: $2{n\choose 4}$, for $n\ge
1$. %

$\rho_{\ge 2}(n-1,0)$ is \seqnum{A001629}  (Self-convolution of Fibonacci
numbers), for $n\ge 1$. %

$\rho_{\le 1}(n-1,0)=\rho_1(n-1,0)$ is \seqnum{A136444} for $n\ge 1$.

$\rho_{[1]}(n-1,0)$ is \seqnum{A029907} for $n\ge 1$, or
$\rho_{[1]}(n+1,0)=w_{[1]}(n,1)$ for $n\ge 0$.

$\rho_{[0]}(n-1,0)=w_{[0]}(n,1)$ is \seqnum{A384497} for $n\ge 1$.

\subsection{Number of Nonnull Runs of Ones and/or Zeros Over All 
  $n$-Strings that Contain Exactly $m$ Nonnull Runs of Ones Under
  Different Constraints}
\label{sec:number-runs-oz-m}

In this section we address the problems of enumerating the number of
nonnull runs of ones and zeros over all $n$-strings that contain
exactly $m$ nonnull $(\ushort{k}\le\overline{k})$-runs of ones or
exactly $m$ nonnull $p$-parity runs of ones. We call these two
quantities $\beta_{\ushort{k}\le\overline{k}}(n,m)$ and $\beta_{[p]}(n,m)$,
respectively.  Special cases of both problems were previously studied
by Nyblom~\cite{nyblom12:_enumerating} and by Grimaldi and
Heubach~\cite{grimaldi05:_without_odd_runs}.

These two enumerations can be obtained in a similar way as the
equivalent enumerations in the previous section. The main difference
is that, rather than being respectively based on recurrences for
$w_{\ushort{k}\le\widebar{k}}(n,m)$ and $w_{[p]}(n,m)$ like in
Section~\ref{sec:number-runs-m}, the two recurrences are now based on
mutual recurrences for the same two quantities constrained to the
initial bit of the $n$-strings ---just like the ones used in
Section~\ref{sec:oz}, but for runs of ones only. The asymmetry of such
mutual recurrences implies that deriving ogfs, although certainly
feasible, is rather laborious
---cf.~\eqref{eq:s_n_m_klek_r_ogf}. %
Therefore we content ourselves with only providing recurrences for
$\beta_{\ushort{k}\le\overline{k}}(n,m)$ and $\beta_{[p]}(n,m)$.

\subsubsection[Number of Nonnull Runs of Ones and/or Zeros Over All 
n-Strings that Contain Exactly m Nonnull (ḵ≤k)-Runs of Ones]{Number of
  Nonnull Runs of Ones and/or Zeros Over All $n$-Strings that Contain
  Exactly $m$ Nonnull $(\ushort{k}\le \widebar{k})$-Runs of Ones}
\label{sec:h_n_m_klek}

As indicated we assume $\ushort{k}\ge 1$. We begin by finding mutual
recurrences for the number of $n$-strings that start with
$b\in\{0,1\}$ and contain exactly $m$
$(\ushort{k}\le\widebar{k})$-runs of ones, which we denote by
$w^b_{\ushort{k}\le\widebar{k}}(n,m)$. If an $n$-string starts with an
$i$-run of zeros then it contributes
$w^{1}_{\ushort{k}\le\widebar{k}}(n-i,m)$ to
$w^0_{\ushort{k}\le\widebar{k}}(n,m)$. On the other hand if an
$n$-string starts with an $i$-run of ones, if
$\ushort{k}\le i\le \widebar{k}$ then it contributes
$w^{0}_{\ushort{k}\le\widebar{k}}(n-i,m-1)$ to
$w^1_{\ushort{k}\le\widebar{k}}(n,m)$, but otherwise it contributes
$w^{0}_{\ushort{k}\le\widebar{k}}(n-i,m)$. Therefore
\begin{align}
  \label{eq:wb_n_m_klek_rec}
  w^b_{\ushort{k}\le\widebar{k}}(n,m)=&~\sum_{i=\ushort{k}}^{\widebar{k}}w^{\tilde{b}}_{\ushort{k}\le\widebar{k}}\big(n-i,m-\iv{b=1}\big)
                                      +\sum_{i=1}^{\ushort{k}-1}w^{\tilde{b}}_{\ushort{k}\le\widebar{k}}(n-i,m)
                                        +\sum_{i=\widebar{k}+1}^{n}w^{\tilde{b}}_{\ushort{k}\le\widebar{k}}(n-i,m),
\end{align}
where $b\in\{0,1\}$ and $\tilde{b}=\text{mod}(b+1,2)$.  When $n=1$ we
know by inspection that
\begin{equation}
  \label{eq:wb_n1_trivial}
  w^b_{\ushort{k}\le\widebar{k}}(1,m)=\iv{b=0}\iv{m=0}+\iv{b=1}\big(\iv{m=0}\iv{\ushort{k}>1}+\iv{m=1}\iv{\ushort{k}=1}\big).
\end{equation}
On the other hand, in the specific case of
$n=1$~\eqref{eq:wb_n_m_klek_rec} with becomes
\begin{equation}
  \label{eq:wb_n1_rec}
  w^b_{\ushort{k}\le\widebar{k}}(1,m)=w^{\tilde{b}}_{\ushort{k}\le\widebar{k}}\big(0,m-\iv{b=1}\iv{\ushort{k}=1}\big).
\end{equation}
We may verify that~\eqref{eq:wb_n1_rec}
equals~\eqref{eq:wb_n1_trivial} when the following initial values are used:
\begin{equation}
  \label{eq:wb_klek_init}
  w^b_{\ushort{k}\le\widebar{k}}(0,m)=\iv{m=0}.
\end{equation}
Of course,
$w_{\ushort{k}\le\widebar{k}}(n,m)=w^0_{\ushort{k}\le\widebar{k}}(n,m)+w^1_{\ushort{k}\le\widebar{k}}(n,m)$,
which allows us to recover the enumeration in
Section~\ref{sec:w_n_m_klek} for $\ushort{k}\ge 1$ in an alternative
way. However, it is not immediately obvious how to obtain a recurrence
for $w_{\ushort{k}\le\widebar{k}}(n,m)$ through this approach, due to
the asymmetry between the two mutual recurrences
in~\eqref{eq:wb_n_m_klek_rec} ---similarly to the enumeration in
Section~\ref{sec:s_n_m_klek}. Moreover, this alternative approach also
makes the derivation of the ogf~\eqref{eq:ogf_w_n_m_klek} less simple.

After these preliminaries we are ready to tackle the enumeration of
$\beta_{\ushort{k}\le\widebar{k}}(n,m)$. We begin by obtaining mutual
recurrences for the number of nonnull runs of ones and zeros in
$n$-strings that start with~$b$ and contain exactly $m$ nonnull
$({\ushort{k}\le\widebar{k}})$-runs (of ones), which we denote by
$\beta^b_{\ushort{k}\le\widebar{k}}(n,m)$. These two recurrences can
be reasoned as follows: the $n$-strings that start with an $i$-run of
zeros contribute $\beta^1_{\ushort{k}\le\widebar{k}}(n-i,m)$ plus
$w^1_{\ushort{k}\le\widebar{k}}(n,m)$ initial runs of zeros to
$\beta^0_{\ushort{k}\le\widebar{k}}(n,m)$. As for the $n$-strings that
start with an $i$-run of ones, if $\ushort{k}\le i\le \widebar{k}$
then they contribute $\beta^0_{\ushort{k}\le\widebar{k}}(n-i,m-1)$
plus $w^0_{\ushort{k}\le\widebar{k}}(n-i,m-1)$ initial runs of ones to
$\beta^1_{\ushort{k}\le\widebar{k}}(n,m)$. If they start with an
$i$-run of ones with $i<\ushort{k}$ or $i>\widebar{k}$ then they
contribute $\beta^0_{\ushort{k}\le\widebar{k}}(n-i,m)$ plus
$w^0_{\ushort{k}\le\widebar{k}}(n-i,m)$ initial runs of ones. This
yields
\begin{align}
  \label{eq:hb_klek_n_m_rec}
  \beta^b_{\ushort{k}\le\widebar{k}}(n,m)=&~\sum_{i=\ushort{k}}^{\widebar{k}}\Big(\beta^{\tilde{b}}_{\ushort{k}\le\widebar{k}}(n-i,m-\iv{b=1})+w^{\tilde{b}}_{\ushort{k}\le\widebar{k}}(n-i,m-\iv{b=1})\Big)\nonumber\\
                                      &+\sum_{i=1}^{\ushort{k}-1}\Big(\beta^{\tilde{b}}_{\ushort{k}\le\widebar{k}}(n-i,m)+w^{\tilde{b}}_{\ushort{k}\le\widebar{k}}(n-i,m)\Big)
  +\sum_{i=\widebar{k}+1}^{n}\Big(\beta^{\tilde{b}}_{\ushort{k}\le\widebar{k}}(n-i,m)+w^{\tilde{b}}_{\ushort{k}\le\widebar{k}}(n-i,m)\Big).
\end{align}
When $n=1$ we have by inspection that
$\beta^b_{\ushort{k}\le\widebar{k}}(1,m)=w^b_{\ushort{k}\le\widebar{k}}(1,m)$
---see~\eqref{eq:wb_n1_trivial}. Thus, the initialisation
of~\eqref{eq:hb_klek_n_m_rec} is
$\beta^b_{\ushort{k}\le\widebar{k}}(0,m)=0$. Finally, the enumeration that
we are interested in is
\begin{equation}
  \label{eq:h_n_m_klek}
  \beta_{\ushort{k}\le\widebar{k}}(n,m)=\beta^0_{\ushort{k}\le\widebar{k}}(n,m)+\beta^1_{\ushort{k}\le\widebar{k}}(n,m).
\end{equation}

\begin{remark}
  Nyblom~\cite[Thm.\ 3.2]{nyblom12:_enumerating} gave a recurrence
  based on the generalised Fibonacci
  numbers~\eqref{eq:nth_order_fibonacci_standard} to compute the total
  number of runs of ones and zeros over all the $w_{\ge k}(n,0)$
  $n$-strings that contain no runs of ones of length $k$ or longer for
  $k\ge 3$, which he denotes by~$\rho_k(n)$. Observe
  that~\eqref{eq:h_n_m_klek} allows us to obtain $\rho_k(n)=\beta_{k\le n}(n,0)$
  for all $1\le k\le n$.
\end{remark}

\subsubsection[Number of Nonnull Runs of Ones and/or Zeros Over All 
n-Strings that Contain Exactly m Nonnull p-Parity Runs of Ones]{Number
  of Nonnull Runs of Ones and/or Zeros Over All $n$-Strings that
  Contain Exactly $m$ Nonnull $p$-Parity Runs of Ones}
\label{sec:h_n_m_p}

As in the previous section, we find first mutual recurrences for the
number of $n$-strings that start with $b\in\{0,1\}$ and contain
exactly $m$ $p$-parity runs of ones, which we denote by
$w^b_{[p]}(n,m)$. If an $n$-string starts with an $i$-run of zeros
then it contributes $w^{1}_{[p]}(n-i,m)$ to $w^0_{[p]}(n,m)$. On the
other hand if an $n$-string starts with an $i$-run of ones, it
contributes $w^{0}_{[p]}(n-i,m-1)$ to $w^1_{[p]}(n,m)$ if
$p=\text{mod}(i,2)$; otherwise, it contributes
$w^{0}_{[p]}(n-i,m)$. Therefore
\begin{align}
  \label{eq:wb_n_m_p_rec}
  w^b_{[p]}(n,m)=&~\sum_{i=1}^{n}w^{\tilde{b}}_{[p]}\big(n-i,m-\iv{b=1}\iv{p=\text{mod}(i,2)}\big),
\end{align}
where $b\in\{0,1\}$ and $\tilde{b}=\text{mod}(b+1,2)$.  In the case $n=1$ we know by inspection
that
\begin{equation}
  \label{eq:wb_p_n1_trivial}
  w^b_{[p]}(1,m)=\iv{b=0}\iv{m=0}+\iv{b=1}\big(\iv{m=0}\iv{p=0}+\iv{m=1}\iv{p=1}\big),
\end{equation}
whereas in the specific case of $n=1$,
expression~\eqref{eq:wb_n_m_p_rec} yields
\begin{equation}
  \label{eq:wb_p_n1_rec}
  w^b_{[p]}(1,m)=w^{\tilde{b}}_{[p]}\big(0,m-\iv{b=1}\iv{p=1}\big).
\end{equation}
We wish~\eqref{eq:wb_p_n1_rec} to equal~\eqref{eq:wb_p_n1_trivial},
which happens the following initial values are used:
\begin{equation}
  \label{eq:wb_p_init}
  w^b_{[p]}(0,m)=\iv{m=0}.
\end{equation}
Thus the enumeration in Section~\ref{sec:wp_n_m} can alternatively be
obtained using~\eqref{eq:wb_n_m_p_rec} as
$w_{[p]}(n,m)=w^0_{[p]}(n,m)+w^1_{[p]}(n,m)$. But just like in the
previous section, the asymmetry between the mutual recurrences makes
it more difficult to get a recurrence for $w_{[p]}(n,m)$, or its ogf,
through this approach.

In any case, we only use the recurrences in~\eqref{eq:wb_n_m_p_rec} to
help us get $\beta_{[p]}(n,m)$.  We can do so by obtaining, in turn,
mutual recurrences for the number of nonnull runs of ones and zeros in
$n$-strings that start with $b$ and contain exactly $m$ nonnull $p$-parity
runs (of ones), which we denote by $\beta^b_{[p]}(n,m)$. These two mutual
recurrences can be reasoned in a similar way as in previous section:
the $n$-strings that start with an $i$-run of zeros contribute
$\beta^1_{[p]}(n-i,m)$ plus $w^1_{[p]}(n,m)$ initial runs of zeros to
$\beta^0_{[p]}(n,m)$. As for the $n$-strings that start with an $i$-run of
ones, if $p=\text{mod}(i,2)$ then they contribute $\beta^0_{[p]}(n-i,m-1)$
plus $w^0_{[p]}(n-i,m-1)$ initial runs of ones to $\beta^1_{[p]}(n,m)$;
otherwise they contribute $\beta^0_{[p]}(n-i,m)$ plus $w^0_{[p]}(n-i,m)$
initial runs of ones. This yields
\begin{align}
  \label{eq:hb_p_n_m_rec}
  \beta^b_{[p]}(n,m)=&~\sum_{i=1}^{n}\Big(\beta^{\tilde{b}}_{[p]}\big(n-i,m-b\,\iv{p=\text{mod}(i,2)}\big)+w^{\tilde{b}}_{[p]}\big(n-i,m-b\,\iv{p=\text{mod}(i,2)}\big)\Big)
\end{align}
When $n=1$ we have that $\beta^b_{[p]}(1,m)=w^b_{[p]}(1,m)$ by inspection
---see~\eqref{eq:wb_p_n1_trivial}. Thus, the initialisation
of~\eqref{eq:hb_p_n_m_rec} is $\beta^b_{[p]}(0,m)=0$. Finally, the
enumeration that we are interested in is
\begin{equation}
  \label{eq:h_n_m_p}
  \beta_{[p]}(n,m)=\beta^0_{[p]}(n,m)+\beta^1_{[p]}(n,m).
\end{equation}

\begin{remark}\label{rmk:grimaldi_r_2}
  Grimaldi and Heubach's enumeration mentioned in
  Remark~\ref{rmk:grimaldi_r}~\cite[Thm.\
  4]{grimaldi05:_without_odd_runs}, can be obtained as a special case
  of~\eqref{eq:h_n_m_p}: $t_n=\beta_{[1]}(n,0)$.
\end{remark}

\section{Number of Ones in Runs Under Different Constraints Over All
  $n$-Strings and Over  Restricted Subsets}
\label{sec:total-number-ones}

Our main goal in this section is enumerating how many ones are found
in runs under different constraints over all $n$-strings. The most
basic result, in Section~\ref{sec:total-number-ones-klek}, is the
enumeration of the ones in $(\ushort{k}\le\widebar{k})$-runs over all
$n$-strings, of which we consider several relevant special cases. We
also study in Section~\ref{sec:number-ones-m} the number of ones in
the $n$-strings enumerated in Sections~\ref{sec:w_n_m_klek}
and~\ref{sec:wp_n_m}.  Several of the enumerations considered here
were previously studied by Makri et
al.~\cite{makri12:_count_runs_ones_ones_runs},
Nyblom~\cite{nyblom12:_enumerating}, and Grimaldi and
Heubach~\cite{grimaldi05:_without_odd_runs}.

\subsection[Number of Ones in (ḵ≤k)-Runs Over All $n$-Strings]{Number
  of Ones in
  $(\ushort{k}\le\widebar{k})$-Runs Over All $n$-Strings}
\label{sec:total-number-ones-klek}
We denote by $t_{\ushort{k}\le\widebar{k}}(n)$ the number of ones in
$(\ushort{k}\le\widebar{k})$-runs over all $n$-strings. First of all,
it is not possible to establish a simple formula for
$t_{\ushort{k}\le\widebar{k}}(n)$ by directly exploiting the main
results in Section~\ref{sec:w_n_m_klek} ---i.e., using formulas
parallel to \eqref{eq:r_n_klek_from_w_n_m_klek} or
to~\eqref{eq:r_klek_through_ogf}. It is possible, though, to obtain
$t_{\ushort{k}\le\widebar{k}}(n)$ in a more roundabout way by using
the enumeration in Section~\ref{sec:w_n_m_k} as follows:
\begin{equation}\label{eq:tn_klek_from_w_n_m_k}
  t_{\ushort{k}\le\widebar{k}}(n)=\sum_{k=\ushort{k}}^{\widebar{k}} k\,\sum_{m=1}^{\big\lfloor
    \frac{n+1}{k+1}\big\rfloor} m\,   w_{k}(n,m).
\end{equation}
Nevertheless, if our intention is to evaluate this expression through
an explicit formula for $w_k(n,m)$ such as~\eqref{eq:w_gen}, then this
expression just adds another layer of misery with respect to the
already hard evaluation of~\eqref{eq:rnk_wknm}.  A kinder option is to
directly replace the second summation
in~\eqref{eq:tn_klek_from_w_n_m_k} with~$\rho_k(n)$. Denoting by
$t_k(n)$ the number of ones in $k$-runs over all $n$-strings, we
trivially have that
\begin{equation}
  \label{eq:tikei}
  t_k(n)=k \rho_k(n),
\end{equation}
and so we can write~\eqref{eq:tn_klek_from_w_n_m_k} as
\begin{equation}
  \label{eq:t_n_klek_from_r_n_k}
  t_{\ushort{k}\le\widebar{k}}(n)=\sum_{k=\ushort{k}}^{\widebar{k}} t_k(n).
\end{equation}
As from~\eqref{eq:tikei} and~\eqref{eq:r_n_k_explicit} we have that
\begin{equation}
  \label{eq:tkn}
  t_k(n)=k\,(n-k+3)\,2^{n-(k+2)}\iv{k<n}+n\,\iv{k=n},
\end{equation}
we can see that the evaluation of~\eqref{eq:t_n_klek_from_r_n_k} only
involves summations of the forms $\sum_k k a^k$ and $\sum_k k^2
a^k$. Thus, after some algebra, we can put the desired enumeration as
\begin{align}
  \label{eq:t_n_klek_general_khlenm1}
  t_{\ushort{k}\le\widebar{k}}(n)=~&\big(n(\ushort{k}+1)-\ushort{k}(\ushort{k}-1)\big)\,2^{n-(\ushort{k}+1)}\,\iv{\ushort{k}\le
                                     n}-\big(n\,(\widebar{k}+2)-\widebar{k}\,(\widebar{k}+1)\big)\,2^{n-(\widebar{k}+2)}\,\iv{\widebar{k}<n}
\end{align}
As usual, this expression is valid for
$0\le \ushort{k}\le \widebar{k}$ and $n\ge 1$ ---even though, of
course, null runs contain no ones and thus do not contribute to
$t_{\ushort{k}\le\widebar{k}}(n)$.

We will not show it here, but, like the main result in
Section~\ref{sec:total-number-runs-klek},
\eqref{eq:t_n_klek_general_khlenm1} can also be obtained by
establishing a recurrence for $t_{\ushort{k}\le\widebar{k}}(n)$ from
first principles (i.e, without resorting to any of the results in
previous sections), and then solving it directly by unrolling~it.

\subsubsection{Special Cases}
\label{sec:number-ones-klek-special-cases }
We denote the number of ones in $(\ge\! k)$-runs, the number of ones
in $(\le\! k)$-runs, and the total number of ones over all $n$-strings
by $t_{\ge k}(n)$, $t_{\le k}(n)$ and $t(n)$, respectively.  These are
special cases of the result in the previous section, since we can write
$t_{\ge k}(n)=t_{k\le n}(n)$, $t_{\le k}(n)=t_{1\le k}(n)$, and
$t(n)=t_{1\le n}(n)$. Thus, using~\eqref{eq:t_n_klek_general_khlenm1}
we have that
\begin{equation}
  \label{eq:tgtk_explicit}
  t_{\ge k}(n)=\big(n(k+1)-k(k-1)\big)\,2^{n-(k+1)}\,\iv{k\le n},%
\end{equation}
Because $t(n)=t_{\ge 1}(n)$, we have
\begin{equation}
  \label{eq:tn}
  t(n)=n\,2^{n-1}.
\end{equation}
Finally, 
\begin{equation}
  \label{eq:tltk_explicit}
  t_{\le k}(n)=t(n)-\big(n\,(k+2)-k\,(k+1)\big)\,2^{n-(k+2)}\,\iv{k<n}.
\end{equation}

\begin{remark}\label{rmk:onesandones}
  In the context of runs in binary strings, the expression for $t(n)$,
  i.e.,~\eqref{eq:tn}, was given by Nyblom~\cite[Lem.\
  3.1]{nyblom12:_enumerating}. 
  Of course, one can find~\eqref{eq:tn} without using runs at all. For
  example, $\sum_{j=1}^nj\, {n\choose j}=n\,2^{n-1}$. This quantity was
  also used by Gold~\cite{gold29:_note_frequen} in his computation of
  $\rho_k(n)$ ---see Remark~\ref{rmk:rk-r-alt}.

  Makri and Psillakis sketched a procedure in~\cite[Eq.\
  (25)]{makri11:_bernoul} to obtain~$t_{\ge k}(n)$ ---in their
  notation,~$R_{n,k}^{(s)}$. At that point they did not complete the
  calculation, probably due to their approach making it rather hard to
  obtain a closed-form expression. At the very end of the same paper
  they stated: ``\textit{A simple explicit form of $R_{n,k}^{(s)}$
    remains an open issue}''. However, soon afterwards, Makri et al.\
  were able to give the explicit closed-form
  expression~\eqref{eq:tgtk_explicit} for $t_{\ge k}(n)$~\cite[Eq.\
  (7), $a=s$]{makri12:_count_runs_ones_ones_runs}, by exploiting a
  previous probabilistic result of theirs about success runs
  statistics~\cite{makri07:_succes}. From~\eqref{eq:tgtk_explicit},
  they also deduced a recurrence for $t_{\ge k}(n)$~\cite[Eq.\ (9),
  $a=s$]{makri12:_count_runs_ones_ones_runs} ---see also the last
  paragraph in Section~\ref{sec:total-number-ones-klek}.
\end{remark}

\subsubsection{OEIS}
\label{sec:oeis-tklek-special-cases}
Below are the OEIS sequences that we have been able to find in
connection to the enumerations in this section. 
\begin{itemize}
\item Only four sequences emanating from $t_k(n)$ appear to be
    documented in the~OEIS:

  $t_1(n+1)$ is \seqnum{A045623} ---cf.\ $\rho_1(n+1)$.

  $t_2(n+1)$ is \seqnum{A087447} for $n\ge 2$. %
  
  $t_3(n+2)$ is \seqnum{A084860} for $n\ge 2$. %

  $t_4(n+3)$ is \seqnum{A001792} for $n\ge 2$. Incidentally, this is
  the same sequence as $\rho(n+1)=\rho_{\ge k}(n+k)$, but the connection
  between both is not immediately obvious. %
  
\item $t(n)$ is \seqnum{A001787}.

\item Sequences emanating from $t_{\ge k}(n)$: 

  $t_{\ge 1}(n)$ is \seqnum{A001787} ---cf.\ $t(n)$.
  
  $t_{\ge 2}(n)$ is \seqnum{A066373}.

  $t_{\ge 3}(n+1)$ is \seqnum{A128135} for $n\ge 2$.

  $t_{\ge 4}(n)$ is \seqnum{A386250}. %
  
  $t_{\ge 5}(n+3)$ is \seqnum{A053220} for $n\ge 2$.

\item Sequences emanating from $t_{\le k}(n)$:

  $t_{\le 1}(n)$ is \seqnum{A045623} ---cf.\ $t_1(n+1)$.

  $t_{\le 2}(n)$ is \seqnum{A386270}. %
  
\end{itemize}

\subsection{Number of Ones Over All $n$-Strings that Contain
  Exactly $m$ Nonnull Runs Under Different Constraints}
\label{sec:number-ones-m}
This section is the equivalent of
Section~\ref{sec:number-runs-m}. Here we deal with the problems of
enumerating how many ones are found over all $n$-strings that contain
exactly $m$ nonnull $(\ushort{k}\le\overline{k})$-runs or exactly $m$
nonnull $p$-parity runs. In keeping with our notation conventions, we
call these quantities $t_{\ushort{k}\le\overline{k}}(n,m)$ and
$t_{[p]}(n,m)$, respectively. Special cases of these problems were
previously studied by Nyblom~\cite{nyblom12:_enumerating} and by
Grimaldi and Heubach~\cite{grimaldi05:_without_odd_runs}.

Again, while relatively simple explicit expressions do not seem
possible in general, we give recurrence relations and generating
functions for these enumerations by relying on their counterparts in
Sections~\ref{sec:w_n_m_klek} and~\ref{sec:wp_n_m}, i.e.,
$w_{\ushort{k}\le\overline{k}}(n,m)$ and $w_{[p]}(n,m)$,
respectively. We must mention though that explicit expressions are
possible in some particular cases ---see Remark~\ref{rmk:grimaldi_t}
and Section~\ref{sec:oeis-tnm}.

\subsubsection[Number of Ones Over All n-Strings That Contain
Exactly m Nonnull (ḵ≤k)-Runs]{Number of Ones Over All  $n$-Strings That
  Contain Exactly $m$ Nonnull  $(\ushort{k}\le \widebar{k})$-Runs}
\label{sec:t_n_m_klek}

We obtain first a recurrence for $t_{\ushort{k}\le\widebar{k}}(n,m)$,
the number of ones over all $n$-strings that contain exactly $m$
nonnull $(\ushort{k}\le \widebar{k})$-runs. As indicated, we assume
$\ushort{k}\ge 1$. Following the customary strategy, the $n$-strings
that start with~$0$ contribute $t_{\ushort{k}\le\widebar{k}}(n-1,m)$
runs to $t_{\ushort{k}\le\widebar{k}}(n,m)$. As for the $n$-strings
that start with~$1$, if they start with a nonnull $i$-run with
$\ushort{k}\le i\le \widebar{k}$ then they contribute
$t_{\ushort{k}\le\widebar{k}}(n-(i+1),m-1)$ ones plus the
$i\,w_{\ushort{k}\le\widebar{k}}(n-(i+1),m-1)$ initial ones
themselves. If they start with a nonnull $i$-run with $i<\ushort{k}$
or $i>\widebar{k}$ then they contribute
$t_{\ushort{k}\le\widebar{k}}(n-(i+1),m)$ ones plus the
$i\,w_{\ushort{k}\le\widebar{k}}(n-(i+1),m)$ initial ones. This yields
the recurrence
\begin{align}
  \label{eq:t_klek_n_m_rec}
  t_{\ushort{k}\le\widebar{k}}(n,m)=t_{\ushort{k}\le\widebar{k}}(n-1,m)%
                                     &+ \sum_{i=\ushort{k}}^{\widebar{k}} \Big(
              t_{\ushort{k}\le\widebar{k}}(n-(i+1),m-1)+i\,w_{\ushort{k}\le\widebar{k}}(n-(i+1),m-1)\Big)\nonumber\\
            &+ \sum_{i=1}^{\ushort{k}-1} \Big(
              t_{\ushort{k}\le\widebar{k}}(n-(i+1),m)+i\, w_{\ushort{k}\le\widebar{k}}(n-(i+1),m)\Big)\nonumber\\
            &+ \sum_{i=\widebar{k}+1}^{n} \Big(
              t_{\ushort{k}\le\widebar{k}}(n-(i+1),m)+i\,w_{\ushort{k}\le\widebar{k}}(n-(i+1),m)\Big),
\end{align}
which is naturally very similar to \eqref{eq:r_klek_n_m_rec}. To
initialise the recurrence we use the case $n=1$, which leads to the
very same equations as in Section~\ref{sec:r_n_m_klek}, and thus, like
in~\eqref{eq:r_klek_n_m_init},
the initial values of~\eqref{eq:t_klek_n_m_rec} are
\begin{equation*}
  \label{eq:t_klek_n_m_init}
  t_{\ushort{k}\le\widebar{k}}(-1,m)=t_{\ushort{k}\le\widebar{k}}(0,m)=0.
\end{equation*}
In order to obtain the ogf
$T_{\ushort{k}\le\widebar{k}}(x,y)=\sum_{n,m}t_{\ushort{k}\le\widebar{k}}(n,m)\,x^ny^m$
we need a recurrence without an $n$-dependent summation. We can achieve this
by first obtaining the difference
$t_{\ushort{k}\le\widebar{k}}(n,m)-t_{\ushort{k}\le\widebar{k}}(n-1,m)$
using~\eqref{eq:t_klek_n_m_rec}, and then obtaining the same
difference again but using the recurrence resulting from the first
step. This is a straightforward but tedious computation and so,
instead of listing the intermediate steps, we directly give the final
recurrence ---which is of course equivalent
to~\eqref{eq:t_klek_n_m_rec}:
\begin{align}
  \label{eq:t_klek_n_m_rec_2}
  t_{\ushort{k}\le\widebar{k}}(n,m)=~&3\,t_{\ushort{k}\le\widebar{k}}(n-1,m)+t_{\ushort{k}\le\widebar{k}}(n-(\ushort{k}+1),m-1)
                                     -t_{\ushort{k}\le\widebar{k}}(n-(\widebar{k}+2),m-1)\nonumber\\&-t_{\ushort{k}\le\widebar{k}}(n-(\ushort{k}+1),m)
                                     +t_{\ushort{k}\le\widebar{k}}(n-(\widebar{k}+2),m)-2\,t_{\ushort{k}\le\widebar{k}}(n-2,m)\nonumber\\
                                     &-t_{\ushort{k}\le\widebar{k}}(n-(\ushort{k}+2),m-1)+t_{\ushort{k}\le\widebar{k}}(n-(\widebar{k}+3),m-1)\nonumber\\
                                     &+t_{\ushort{k}\le\widebar{k}}(n-(\ushort{k}+2),m)-t_{\ushort{k}\le\widebar{k}}(n-(\widebar{k}+3),m)\nonumber\\
                                     &-(\ushort{k}-1)\,w_{\ushort{k}\le\widebar{k}}(n-(\ushort{k}+2),m-1)-(\widebar{k}+1)\,w_{\ushort{k}\le\widebar{k}}(n-(\widebar{k}+2),m-1)\nonumber\\&+\widebar{k}\,w_{\ushort{k}\le\widebar{k}}(n-(\widebar{k}+3),m-1)+\ushort{k}\,w_{\ushort{k}\le\widebar{k}}(n-(\ushort{k}+1),m-1)\nonumber\\
                                     &+w_{\ushort{k}\le\widebar{k}}(n-2,m)+(\ushort{k}-1)\,w_{\ushort{k}\le\widebar{k}}(n-(\ushort{k}+2),m)%
                                     -\ushort{k}\,w_{\ushort{k}\le\widebar{k}}(n-(\ushort{k}+1),m)\nonumber\\&-\widebar{k}\,w_{\ushort{k}\le\widebar{k}}(n-(\widebar{k}+3),m)%
                                     +(\widebar{k}+1)\,w_{\ushort{k}\le\widebar{k}}(n-(\widebar{k}+2),m).
\end{align}
Now, multiplying~\eqref{eq:t_klek_n_m_rec_2} on both sides by $x^ny^m$
and summing over $n$ and $m$ we finally obtain the desired ogf in
terms of~\eqref{eq:ogf_w_n_m_klek}:
\begin{align}
  \label{eq:t_n_m_klek_ogf}
  T_{\ushort{k}\le\widebar{k}}(x,y)&=\frac{x^2-(1-y)\Big(\ushort{k}x^{\ushort{k}+1}-(\ushort{k}-1)x^{\ushort{k}+2}-(\widebar{k}+1)x^{\widebar{k}+2}+\widebar{k}x^{\widebar{k}+3}\Big)}{1-3x+2x^2+(1-y)\big(x^{\ushort{k}+1}-x^{\ushort{k}+2}-x^{\widebar{k}+2}+x^{\widebar{k}+3}\big)}W_{\ushort{k}\le\widebar{k}}(x,y).
\end{align}

\begin{remark}
  Nyblom~\cite[Thm.\ 3.1]{nyblom12:_enumerating} gave a recurrence to
  compute the number of \textit{zeros} in $n$-strings devoid of
  $(\ge k)$-runs, which he denotes by $\rho_k(n)$. We can also obtain
  this enumeration as a special case of our results in this section by
  using $t_{\ge k}(n,0)=t_{k\le n}(n,0)$, which allows us to write
  $\rho_k(n)=n\,w_{\ge k}(n,0)-t_{\ge k}(n,0)$. Also, this expression is
  valid for all $k\ge 1$ whereas the recurrence
  in~\cite{nyblom12:_enumerating} requires $k\ge 3$.
\end{remark}

\subsubsection{Number of Ones Over All  $n$-Strings That  Contain
  Exactly $m$ Nonnull $p$-Parity Runs}
\label{sec:t_n_m_p}
Finally, we find a recurrence relation for $t_{[p]}(n,m)$.  The
$n$-strings that start with~$0$ contribute $t_{[p]}(n-1,m)$ ones to
$t_{[p]}(n,m)$. As for the $n$-strings that start with $1$, if they
start with an odd $i$-run then they contribute
$t_{[p]}(n-(i+1),m-\iv{p=1})$ ones plus the
$i\, w_{[p]}(n-(i+1),m-\iv{p=1})$ initial ones themselves. If they
start with an even nonnull $i$-run then they contribute
$t_{[p]}(n-(i+1),m-\iv{p=0})$ runs plus the
$i\, w_{[p]}(n-(i+1),m-\iv{p=0})$ initial runs. This yields the
recurrence
\begin{align}
  \label{eq:t_n_m_p_rec}
  t_{[p]}(n,m)=~&t_{[p]}(n-1,m)+ \sum_{i=1}^{n}
                  t_{[p]}(n-(i+1),m-\iv{p=\text{mod}(i,2)})\nonumber\\
  &+ \sum_{i=1}^{n}i\, w_{[p]}(n-(i+1),m-\iv{p=\text{mod}(i,2)}).
\end{align}
To initialise~\eqref{eq:t_n_m_p_rec} we use special case $n=1$, which
leads to the very same equations as in Section~\ref{sec:r_n_m_p}
---i.e.,~\eqref{eq:r_n_m_p_trivial_n1} and~\eqref{eq:r_n_m_p_rec_n1}.  Thus
the initial values of~\eqref{eq:t_n_m_p_rec} are
\begin{equation*}
  \label{eq:t_n_m_p_init}
  t_{[p]}(-1,m)=t_{[p]}(0,m)=0.
\end{equation*}
In order to obtain the ogf
$T_{[p]}(x,y)=\sum_{n,m}t_{[p]}(n,m)\,x^ny^m$ we need a recurrence
without an $n$-dependent summation. We can achieve this by first obtaining
the difference $t_{[p]}(n,m)-t_{[p]}(n-2,m)$
using~\eqref{eq:t_n_m_p_rec}, and then obtaining again the same
difference using the recurrence resulting from the first step. We
directly give the final recurrence ---which is of course equivalent
to~\eqref{eq:t_n_m_p_rec}:
\begin{align}
    \label{eq:t_n_m_p_rec_2}
  t_{[p]}(n,m)=~&t_{[p]}(n-1,m)+2\, t_{[p]}(n-2,m)-2\, t_{[p]}(n-3,m)%
                +t_{[p]}(n-2,m-\iv{p=1})\nonumber\\&+t_{[p]}(n-3,m-\iv{p=0})%
                -t_{[p]}(n-4,m)+t_{[p]}(n-5,m)\nonumber\\
                &-t_{[p]}(n-4,m-\iv{p=1})-t_{[p]}(n-5,m-\iv{p=0})\nonumber\\
                &+w_{[p]}(n-4,m-\iv{p=1})+w_{[p]}(n-2,m-\iv{p=1})\nonumber\\
                &+2\,w_{[p]}(n-3,m-\iv{p=0}).
\end{align}
Lastly, multiplying~\eqref{eq:t_n_m_p_rec_2} on both sides by
$x^ny^m$ and summing over $n$ and $m$ we obtain the desired
ogf in terms of~\eqref{eq:wp_bivariate_ogf}:
\begin{equation}
  T_{[p]}(x,y) = \frac{2 x^3 y^{\iv{p = 0}} + (x^2+x^4) y^{\iv{p = 1}}}%
  {(x^5-x^3)y^{\iv{p = 0}}+(x^4-x^2)y^{\iv{p = 1}}+1-x-2x^2+2x^3+x^4-x^5}\, W_{[p]}(x,y).
\end{equation}

\begin{remark}\label{rmk:grimaldi_t}
  Grimaldi and Heubach~\cite[Thm.\ 3]{grimaldi05:_without_odd_runs}
  gave a recurrence, an explicit expression and a generating function
  for the number of zeros in $n$-strings devoid of odd runs of
  zeros---in their notation,~$z_n$. Notice that 
  $t_{[1]}(n,0)$ gives the number of ones in $n$-strings devoid of odd
  runs of ones, and thus, by symmetry, $z_n=t_{[1]}(n,0)$. Let us
  verify this fact by comparing the corresponding ogfs. In the
  specific case of $m=0$ and $p=1$,
  recurrence~\eqref{eq:t_n_m_p_rec_2} becomes
  \begin{align*}
    t_{[1]}(n,0)=~& t_{[1]}(n-1,0)+2\,t_{[1]}(n-2,0)-t_{[1]}(n-3,0)-t_{[1]}(n-4,0)
                  +2\, w_{[1]}\big(n-3,0\big).
  \end{align*}
  From here,  we have that
  \begin{align}
    [y^0]T_{[1]}(x,y)&=\frac{2x^3}{1-x-2x^2+x^3+x^4}\, [y^0]W_{[1]}(x,y)\label{eq:T1xy_a}\\
                     &=\frac{2x^2}{(1-x-x^2)^2},\label{eq:T1xy_b}
  \end{align}
  where we have used~\eqref{eq:ogf_wp_m} to get~\eqref{eq:T1xy_b}
  from~\eqref{eq:T1xy_a}.  As expected, the ogf~\eqref{eq:T1xy_b} is
  the same as $G_{z_n}(x)$ in~\cite[Thm.\  3]{grimaldi05:_without_odd_runs}. Grimaldi and Heubach also studied
  the number of ones in $n$-strings devoid of odd runs of zeros, which
  they denoted by~$w_n$. Using our own results, we can write
  $w_n=n\,w_{[1]}(n,0)-t_{[1]}(n,0)$.
\end{remark}

\subsubsection{OEIS}
\label{sec:oeis-tnm}
The only sequences resulting from the enumerations in
Section~\ref{sec:number-ones-m} that we have been able to find in the
OEIS are:

$t_{1}(n,0)$ is \seqnum{A259966} (Total binary weight of all
\seqnum{A005251}$(n)$ binary sequences of length $n$ not containing any
isolated $1$'s) ---cf.\ $w_{\le 1}(n-1,0)=w_1(n-1,0)$.

$t_{\ge 1}(n,1)$ is \seqnum{A000292} (Tetrahedral ---or triangular
pyramidal--- numbers).  %

$t_{\ge 2}(n,0)$ is \seqnum{A001629} ---cf.\ $\rho_{\ge 2}(n,0)$. %

$t_{\le 1}(n,0)$ is \seqnum{A259966} ---cf.\ $t_1(n,0)$.

\section{Conclusions}
\label{sec:conclusions}

We have studied fundamental enumeration problems concerning runs in
binary strings, where the runs conform to Mood's
criterion~\cite{mood40:_distr_theor_runs}. While many of these
problems had been solved over the years by different authors through a
variety of methods, we felt that there were gaps in the literature,
and that a uniform, systematic treatment would be beneficial to bring
together and connect the many contributions to the topic.  We believe
that the notation conventions that we have chosen, even if necessarily
unorthodox when it comes to probabilities (in order to accommodate
recurrences and pgfs), are effective when it comes to displaying the
relationships between enumerations concerning runs of ones, runs of
ones and/or zeros, and compositions, their probabilistic extensions,
and their associated generating functions. We consider that a relevant
contribution of this paper is the identification of the existence of
null runs of ones ($0$-runs). Probably due to their lack of practical
use, null runs seem to have escaped the attention of previous authors.
Nevertheless, we have seen that runs of ones of zero length are
essential to connect enumerations of runs of ones and zeros and
enumerations of runs of ones (Theorems~\ref{thm:chave_das_nozes}
and~\ref{thm:chave_das_nozes_2}). Null runs also naturally emerge in
problems involving the longest run of ones
(Sections~\ref{sec:number-ws_n_k_ltk} and~\ref{sec:ws_n_m_k_hamming}).
 
Hopefully, we have been able to show that approaching runs-related
problems by means of recurrences and generating functions is, for the
most part, more straightforward than doing so by means of direct
combinatorial analysis methods. The exception to this rule may be some
enumeration problems with a Hamming weight constraint. We are not
reinventing the wheel here, since, as discussed in the introduction,
this is an old and venerable approach ---even when it comes to
studying runs. However, this powerful technique has probably not been
exploited to its full advantage in this area of research.

It would be naive to think that we have exhausted all interesting
runs-related problems that one can address through the strategy
adopted in this paper. Extensions to probabilistic scenarios
---briefly discussed in Sections~\ref{sec:probability}
and~\ref{sec:note-prob-extens} and at some other points--- can be
investigated much further, following in the steps of previous authors who
went beyond the iid case to study Markov
dependencies~\cite{bateman48:_power,aki93:_markov,balakrishnan01:_runs_scans,fu03:_runs}. Many
purely enumerative problems not addressed here also merit
consideration, such as for example joint enumerations
---see~\cite[Sec.\
2.4.9]{goulden83:_combin_enumer}. 
Of course, explicit expressions based on binomial coefficients have
numerical limitations for large $n$. However, the generating functions
that we have provided can definitely be helpful in further studies
about the asymptotics the enumerations studied here. Austin and
Guy~\cite{austin78:_binar}, Sedgewick and
Flajolet~\cite{flajolet09:_analytic}, Bloom~\cite{bloom98:_singles},
Suman~\cite{suman94:_longest}, Prodinger~\cite{prodinger15}, and other
authors, previously followed this approach in some special
scenarios. Another interesting route in terms of asymptotics could be
using moment generating functions instead of probability generating
functions, as done by Wishart and Hirschfeld~\cite{wishart36:_theorem}
in one particular
problem. 
We should also remark that the mutual recurrences strategy on which
the results in Sections~\ref{sec:oz},~\ref{sec:note-prob-extens},
and~\ref{sec:oz_hamming} hinge ---which was first used by Wishart and
Hirschfeld in a probabilistic context~\cite{wishart36:_theorem}--- can
be readily extended to the enumeration of runs of different kinds in
$b$-ary strings with $b\ge 2$. Some of these $b$-ary generalisations
were already given by previous
authors~\cite{david62:_chance,suman94:_longest}, although not through
the procedure we just mentioned. Any of the research avenues sketched
above may prove fruitful: in Schilling's words, ``\textit{the variety
  of potential applications of runs theory is virtually
  boundless}''~\cite{schilling90:_longest}.

The problems addressed in this paper contain a very rich structure in
terms of the many distinct number sequences that they encompass
---many of them listed in Sloane's On-Line Encyclopedia of Integer
Sequences, but many others yet unexplored. The omissions may be as
relevant as the inclusions: observe that few sequences
corresponding to enumerations with $m\ge 2$ are in the OEIS. The many
connections unearthed thanks to the inestimable help of the OEIS may
open up new uses for the diverse formulas and results given
here.

To conclude, the recurrences for runs of ones that we have produced
compel us to say that, in the same sense that $\sqrt{-1}$ exists as a
mathematical object, a binary string with length $-1$ does also have
some kind of fleeting existence. Intriguingly, we are justified in
stating that the number of binary strings of length~$-1$ that are
devoid of runs of ones of length $k\ge 0$ is exactly one ---in the
same sense that we accept that there is a single binary string of
length~$0$ (empty binary string). We have also seen that, when binary
strings are drawn uniformly at random with probability $0<q<1$ of
drawing a $1$, the probability of a (or rather, ``the'') binary string
of length $-1$ without runs of ones of length~$k\ge 0$ must be
$1/(1-q)>1$. The majority of recurrences for runs of ones that we have
given would not be correct \textit{for all $n\ge 1$ and for all valid
  values of the parameters} without similar puzzling assertions being
true. Had we ignored or tried to skirt around these strange-looking
facts, then most of those recurrences would have been less general and
more piecemeal ---i.e., only valid for certain ranges of the
parameters and/or requiring ad-hoc initialisation. This fragmentary
character is in fact present in many recurrences concerning runs (or
compositions) previously given by other
authors~\cite{apostol88:_binary,chinn03:compositions,schilling90:_longest,bloom98:_singles,nyblom12:_enumerating,nej19:_binary}.
We are not the first ones to observe this phenomenon: the necessity of
considering a binary string of length~$-1$ in recurrences involving
runs of ones was implicit already in the work of of Austin and
Guy~\cite{austin78:_binar}.

\pdfbookmark[1]{References}{References}
\bibliographystyle{jis}

\bibliography{runs} %

\bigskip
\hrule
\bigskip

\noindent 2020 {\it Mathematics Subject Classification}: Primary 05A15 %

\noindent \emph{Keywords:} run, success run, recurrence,
ordinary generating function, probability generating function.

\end{document}